\documentclass[3p,computermodern,10pt]{elsarticle}
\usepackage{fullpage}

\usepackage{booktabs} 

\usepackage[tight-spacing=true]{siunitx}
\usepackage{siunitx} 
\usepackage{pgfplotstable} 
\usepackage{pifont}

\usepackage[ruled]{algorithm2e}
\usepackage{algpseudocode}
\usepackage[utf8]{inputenc}
\usepackage[english]{babel}
\usepackage{amsthm}

\newtheorem{remark}{Remark}[section] 
\usepackage{color}
\usepackage{amsmath}
\usepackage{amsfonts}
\usepackage{mathtools}
\usepackage[titletoc,title]{appendix}
\usepackage{float}
\usepackage{graphicx} 
\usepackage{epstopdf} 
\usepackage{caption}
\usepackage{subcaption}
\usepackage{pdfpages}
\usepackage{comment}

\usepackage{hyperref}
\hypersetup{
    colorlinks=true,
    linkcolor=blue,
    filecolor=magenta,      
    urlcolor=cyan,
}
\usepackage{tikz}
\usetikzlibrary{positioning, fit, arrows.meta, shapes}

\usetikzlibrary{
  arrows.meta, 
  fit, 
  positioning, 
}
\tikzset{
  neuron/.style={ 
    circle,draw,thick, 
    inner sep=0pt, 
    minimum size=3.5em, 
    node distance=1ex and 2em, 
  },
  group/.style={ 
    rectangle,draw,thick, 
    inner sep=0pt, 
  },
  io/.style={ 
    neuron, 
    fill=gray!15, 
  },
  conn/.style={ 
    -{Straight Barb[angle=60:2pt 3]}, 
    thick, 
  },
}

\newcommand{\LSPGG}{G--LSPG}

\newcommand{\LSPGSUPG}{SUPG--LSPG}
\newcommand{\LSPGGLS}{\SDGLS--LSPG}
\newcommand{\LSPGADJ}{\SDADJ--LSPG}
\newcommand{\LSPGSTGLS}{\STGLS--LSPG}
\newcommand{\LSPGSTADJ}{\STADJ--LSPG}
\newcommand{\STGLSLSPG}{\STGLS--LSPG}
\newcommand{\STADJLSPG}{\STADJ--LSPG}

\newcommand{\APGG}{G--APG}
\newcommand{\APGSUPG}{SUPG--APG}
\newcommand{\SUPGAPG}{SUPG--APG}

\newcommand{\APGSTADJ}{\STADJ--APG}
\newcommand{\APGSTGLS}{\STGLS--APG}

\newcommand{\tauApg}{\tau_{\mathrm{APG}}}
\newcommand{\APG}{APG}
\newcommand{\LSPG}{LSPG}
\newcommand{\SUPG}{SUPG}
\newcommand{\GLS}{\SDGLS}
\newcommand{\ADJ}{\SDADJ}
\newcommand{\STGLS}{GLS$_{\mathrm{ST}}$}
\newcommand{\STADJ}{ADJ$_{\mathrm{ST}}$}
\newcommand{\SDGLS}{GLS$_{\mathrm{DS}}$}
\newcommand{\SDADJ}{ADJ$_{\mathrm{DS}}$}

\newcommand{\Galerkin}{Gal.}

\newcommand{\bilinear}{\mathcal{B}}
\newcommand{\bilinearArg}[2]{\bilinear(#1,#2)}
\newcommand{\coerciveConstant}{C}
\newcommand{\coerciveNorm}{Y}

\newcommand{\normGalerkinArg}[1]{ \| #1 \|_\mathrm{G} }
\newcommand{\normSupgArg}[1]{ \| #1 \|_\mathrm{SUPG} }
\newcommand{\normGlsArg}[1]{ \| #1 \|_\mathrm{GLS} }
\newcommand{\normAdjArg}[1]{ \| #1 \|_\mathrm{ADJ} }

\newcommand{\bilinearGalerkin}{\mathcal{B}_\mathrm{G}}
\newcommand{\bilinearGalerkinArg}[2]{\bilinearGalerkin(#1,#2)}
\newcommand{\bilinearStabilized}{\mathcal{B}_\mathrm{S}}
\newcommand{\bilinearStabilizedArg}[2]{\bilinearStabilized(#1,#2)}

\newcommand{\cmark}{{\color{green}\ding{51}}}%
\newcommand{\xmark}{{\color{red}\ding{55}}}%
\newcommand{\pmark}{{\color{blue}\ding{60}}}%

\newcommand{\bz}{\boldsymbol 0}
\newcommand{\RPlusEqual}{\mathbb{R}_{\ge0}}
\newcommand{\RPlus}{\mathbb{R}_{>0}}
\newcommand{\energy}{\epsilon}

\newcommand{\tauSet}{\boldsymbol \tau}
\newcommand{\dtSet}{\boldsymbol \Delta \boldsymbol t}
\newcommand{\LTwo}{L^2(\cDomain)}
\newcommand{\LTwoK}{L^2(\cDomain_k)}

\newcommand{\errorLTwoBestFit}{\overline{e}_{L^2}}
\newcommand{\integratedErrorLTwoBestFit}{\overline{e}_{L^2}}
\newcommand{\projectorLTwoArg}[1]{\mathbb{P}_{L^2}(#1)}
\newcommand{\HOne}{\mathcal{H}^1(\cDomain)}

\newcommand{\errorHOneBestFit}{\overline{e}_{\mathcal{H}^1}}
\newcommand{\integratedErrorHOneBestFit}{\overline{e}_{\mathcal{H}^1}}
\newcommand{\projectorHOneArg}[1]{\mathbb{P}_{\mathcal{H}^1}(#1)}
\newcommand{\RR}[1]{\mathbb{R}^{#1}}
\newcommand{\RRS}[1]{\mathbb{S}^{#1}}

\newcommand{\EP}[1]{{\color{red}EP: #1}}
\newcommand{\MY}[1]{{\color{blue}MY: #1}}
\newcommand{\IT}[1]{{\color{purple}IT: #1}}
\newcommand{\TI}[1]{{\color{orange}TI: #1}}

\newcommand{\discreteROMs}{discrete ROMs}
\newcommand{\DiscreteROMs}{Discrete ROMs}
\newcommand{\continuousROMs}{continuous ROMs}
\newcommand{\ContinuousROMs}{Continuous ROMs}
\newcommand{\cDomain}{\Omega}
\newcommand{\cDomainClosed}{{\overline \Omega}}
\newcommand{\cDomainBoundary}{\Gamma}

\newcommand{\tDomain}{(0,T]}
\newcommand{\tDomainClosed}{[0,T]}
\newcommand{\femCoefficients}{FEM coefficients}
\newcommand{\cdcdrAcronym}{CDR}
\newcommand{\cdcrResidTimeDiscreteSpaceContinuous}{\mathsf{R}_{\mathrm{cdr}}}
\newcommand{\cdcrResidTimeDiscreteSpaceContinuousOrtho}{\mathsf{R}^{\parallel}_{\mathrm{G}}}
\newcommand{\cDomainDim}{d}
\newcommand{\residAPG}{\mathbf{r}_{\mathrm{APG}}}
\newcommand{\residAPGSUPG}{\mathbf{r}_{\mathrm{APG-S}}}

\newcommand{\PiFine}{\mathbb{A}'}

\newcommand{\energyCutoff}{\epsilon_{\mathrm{c}}}

\newcommand{\nSnap}{N_s}
\newcommand{\massRom}{\mass_{\mathrm{r}}}
\newcommand{\massRomArg}[1]{\mass_{\mathrm{r}_{#1}}}

\newcommand{\matFullRom}{\mathbf{\matFull}_{\mathrm{r}} }
\newcommand{\matFullRomArg}[1]{\mathbf{\matFull}_{{\mathrm{r}}_{#1}} }

\newcommand{\forcingVecRom}{\forcingVec_{\mathrm{r}} }
\newcommand{\forcingVecRomArg}[1]{  \forcingVec_{{\mathrm{r}}_{#1}} }

\newcommand{\residStabilizedRom}{\mathbf{r}_{\mathrm{S-r}}}
\newcommand{\stabilizationMatRom}{\mathbf{Q}_{\mathrm{r}}}
\newcommand{\stabilizationMatRomArg}[1]{ {\stabilizationMatRom}_{#1}}

\newcommand{\innerProductType}{\mathcal{X}}
\newcommand{\discreteInnerProductType}{\mathbf{P}}
\newcommand{\discreteInnerProductTypeTwo}{\mathbf{W}}

\newcommand{\snapshotMatrixDiscrete}{\mathbf{S}_{\genStateTimeDiscreteSpaceDiscrete}}
\newcommand{\snapshotMatrixContinuous}{\mathbf{S}_{\cStateTimeDiscreteSpaceDiscrete}}
\newcommand{\residualOp}{\mathcal{R}}
\newcommand{\stabilizationOpTest}{\mathcal{Q}}
\newcommand{\stabilizationOpTestArg}[1]{\stabilizationOpTest{#1}}
\newcommand{\stabilizationOpTestSUPG}{\mathcal{Q}_{\mathrm{DS-SUPG}}}
\newcommand{\stabilizationOpTestSTSUPG}{\mathcal{Q}_{\mathrm{ST-SUPG}}}

\newcommand{\stabilizationOpTestGLS}{\mathcal{Q}_{\mathrm{DS-GLS}}}
\newcommand{\stabilizationOpTestSTGLS}{\mathcal{Q}_{\mathrm{ST-GLS}}}

\newcommand{\stabilizationOpTestADJ}{\mathcal{Q}_{\mathrm{DS-ADJ}}}
\newcommand{\stabilizationOpTestSTADJ}{\mathcal{Q}_{\mathrm{ST-ADJ}}}

\newcommand{\stabilizationOpTrial}{\mathcal{L}}

\newcommand{\stabilizationMat}{\mathbf{Q}}
\newcommand{\stabilizationMatArg}[1]{\mathbf{Q}_{#1}}

\newcommand{\cSpace}{\mathcal{H}_0^1}
\newcommand{\resid}{\mathbf{r}}
\newcommand{\residGalerkinFEM}{\mathbf{r}_{\mathrm{G}}}
\newcommand{\residGalerkinFEMArg}[1]{[\mathbf{r}_{\mathrm{G}}]_{#1}}

\newcommand{\residGalerkinFEMOrtho}{\mathbf{r}^{\parallel}_{\mathrm{G}}}
\newcommand{\residGalerkinROM}{\mathbf{r}_{\mathrm{G-r}}}
\newcommand{\residGalerkinDiscreteRom}{\mathbf{r}_{\mathrm{G-DROM}}}

\newcommand{\residStabilizedFEM}{\mathbf{r}_{\mathrm{S}}}
\newcommand{\nel}{N_{\mathrm{el}}}
\newcommand{\cSpaceFEM}{\mathcal{V}_{\mathrm{h}}}
\newcommand{\cSpaceHRES}{[\cSpaceFEM]_{\mathrm{h-res}}}
\newcommand{\cSpaceFOM}{[\cSpaceFEM]_{\mathrm{fom}}}

\newcommand{\testSpaceFEM}{\mathcal{W}_{\mathrm{h}}}
\newcommand{\cBasis}{\mathsf{v}}
\newcommand{\cBasisVec}{\mathbf{\mathsf{v}}}
\newcommand{\cBasisDum}{\mathsf{v}}

\newcommand{\cState}{u}
\newcommand{\cStateStrong}{u_*}
\newcommand{\cStateStrongTimeDiscrete}{\mathsf{u}_*}
\newcommand{\cStateStrongTimeDiscreteArg}[1]{\cStateStrongTimeDiscrete^{#1}}

\newcommand{\cStateTimeDiscreteDum}{\mathsf{u}}
\newcommand{\cStateTimeDiscreteDumN}{\mathsf{w}}
\newcommand{\cStateTimeDiscreteDumNm}{\mathsf{z}}

\newcommand{\cStateTimeDiscrete}{\mathsf{u}}
\newcommand{\cStateTimeDiscreteArg}[1]{\cStateTimeDiscrete^{#1}}
\newcommand{\cStateTimeDiscreteSpaceDiscrete}{\mathsf{u}_{\mathrm{h}}}
\newcommand{\cStateTimeDiscreteSpaceDiscreteArg}[1]{\cStateTimeDiscreteSpaceDiscrete^{#1}}
\newcommand{\cStateFomTimeDiscreteSpaceDiscrete}{[\cStateTimeDiscreteSpaceDiscrete]_{\mathrm{fom}}}
\newcommand{\cStateFomTimeDiscreteSpaceDiscreteArg}[1]{\cStateFomTimeDiscreteSpaceDiscrete^{#1}}
\newcommand{\cStateHresTimeDiscreteSpaceDiscrete}{[\cStateTimeDiscreteSpaceDiscrete]_{\mathrm{h-res}}}
\newcommand{\cStateHresTimeDiscreteSpaceDiscreteArg}[1]{\cStateHresTimeDiscreteSpaceDiscrete^{#1}}

\newcommand{\genStateTimeDiscreteSpaceDiscrete}{\mathbf{a}_{\mathrm{h}}}
\newcommand{\genStateTimeDiscreteSpaceDiscreteDumN}{\mathbf{w}}
\newcommand{\genStateTimeDiscreteSpaceDiscreteDumNm}{\mathbf{z}}

\newcommand{\discreteROMGenStateTimeDiscreteSpaceDiscrete}{{\mathbf{a}}_{\mathrm{r}}}
\newcommand{\discreteROMGenStateTimeDiscreteSpaceDiscreteArg}[1]{\discreteROMGenStateTimeDiscreteSpaceDiscrete^{#1}}
\newcommand{\genStateTimeDiscreteSpaceDiscreteArg}[1]{\genStateTimeDiscreteSpaceDiscrete^{#1}}

\newcommand{\genStateFomTimeDiscreteSpaceDiscrete}{[\mathbf{a}_{\mathrm{h}}]_{\mathrm{fom}}}
\newcommand{\genStateFomTimeDiscreteSpaceDiscreteArg}[1]{\genStateFomTimeDiscreteSpaceDiscrete^{#1}}

\newcommand{\cSpaceRom}{\mathcal{V}_\mathrm{r}}
\newcommand{\cBasisRomcoefMat}{\mathbf{C}}
\newcommand{\cBasisRomMat}{\boldsymbol \phi}
\newcommand{\cBasisRom}{\phi}
\newcommand{\cBasisRomVec}{\boldsymbol \cBasisRom}

\newcommand{\genRomStateTimeDiscreteSpaceDiscrete}{\hat{\mathbf{x}}}
\newcommand{\genRomStateTimeDiscreteSpaceDiscreteDumN}{\mathbf{w}}
\newcommand{\genRomStateTimeDiscreteSpaceDiscreteDumNm}{\mathbf{z}}

\newcommand{\genRomStateTimeDiscreteSpaceDiscreteArg}[1]{\genRomStateTimeDiscreteSpaceDiscrete^{#1}}
\newcommand{\genDiscreteRomStateTimeDiscreteSpaceDiscreteDumN}{\mathbf{w}}
\newcommand{\genDiscreteRomStateTimeDiscreteSpaceDiscreteDumNm}{\mathbf{z}}
\newcommand{\genDiscreteRomStateTimeDiscreteSpaceDiscrete}{\hat{\mathbf{x}}}
\newcommand{\genDiscreteRomStateTimeDiscreteSpaceDiscreteArg}[1]{\genDiscreteRomStateTimeDiscreteSpaceDiscrete^{#1}}
\newcommand{\romStateTimeDiscreteSpaceDiscrete}{\mathsf{u}_{\mathrm{r}}}
\newcommand{\romStateTimeDiscreteSpaceDiscreteArg}[1]{\romStateTimeDiscreteSpaceDiscrete^{#1}}

\newcommand{\ddStateDum}{\mathbf{y}}
\newcommand{\ddgenStateDum}{\hat{\mathbf{y}}}

\newcommand{\defeq}{\vcentcolon=}

\newcommand{\reaction}{\sigma}

\newcommand{\wavespeed}{\mathbf{b}}

\newcommand{\forcing}{f}
\newcommand{\viscosity}{\epsilon}

\newcommand{\cipGenArbitrary}[3]{\cipGen{#2}{#3}_{#1}}
\newcommand{\cipGen}[2]{m\left(#1,#2\right)}
\newcommand{\cip}{m}
\newcommand{\cipGenElNoArg}{m_{\mathrm{el}}}
\newcommand{\cipGenEl}[2]{m_{\mathrm{el}}\left(#1,#2\right)}
\newcommand{\dipGen}[3]{m_d\left(#1,#2\right)_{#3}}
\newcommand{\dipGenNoArg}{{m_d}}
\newcommand{\dipGenArg}[3]{\left(#2,#3\right)_{#1}}

\newcommand{\fomdim}{N}

\newcommand{\romdim}{R}
\newcommand{\dBasisRomMat}{\mathbf{\Psi}}
\newcommand{\dBasisROMMat}{\mathbf{\Psi}}

\newcommand{\mass}{\mathbf{M}}
\newcommand{\massStabilized}{\mathbf{M}_{\mathrm{S}}}
\newcommand{\massStabilizedArg}[1]{\mathbf{M}_{\mathrm{S}_{#1}}}
\newcommand{\massStabilizedRom}{\mathbf{M}_{\mathrm{S-r}}}
\newcommand{\massStabilizedRomArg}[1]{ {\massStabilizedRom}_{#1}}

\newcommand{\massArg}[1]{\mass_{#1}}

\newcommand{\forcingVec}{\mathbf{f}}
\newcommand{\forcingVecStabilized}{\mathbf{f}_{\mathrm{S}}}
\newcommand{\forcingVecStabilizedArg}[1]{\mathbf{f}_{\mathrm{S}_{#1}}}

\newcommand{\forcingVecStabilizedRom}{\mathbf{f}_{\mathrm{S-r}}}
\newcommand{\forcingVecStabilizedRomArg}[1]{ { \forcingVecStabilizedRom }_{#1} }

\newcommand{\forcingVecArg}[1]{\forcingVec_{#1}}
\newcommand{\matFull}{\mathbf{B}}

\newcommand{\dSpaceROM}{\textit{V}_r}
\newcommand{\dSpaceROMFine}{\textit{V}_r^{\prime}}

\newcommand{\norm}[1]{||#1||_2}
\newcommand{\normArg}[2]{||#2||_{#1}}
\newcommand{\lspgWeight}{\mathbf{W}}

\newcommand{\cSnapCorr}{\mathbf{K}_{\cStateTimeDiscrete}}

\newcommand{\cSnapEigMat}{{\boldsymbol \Lambda}_{\cStateTimeDiscrete}}

\newcommand{\cSnapEigVecMat}{\mathbf{E}_{\cStateTimeDiscrete}}

\newcommand{\dSnapCorr}{\mathbf{K}_{\genStateTimeDiscreteSpaceDiscrete}}
\newcommand{\dSnapEig}{\Lambda}
\newcommand{\dSnapEigMat}{{ \boldsymbol \Lambda}_{\genStateTimeDiscreteSpaceDiscrete}}

\newcommand{\dSnapEigVecMat}{\mathbf{E}_{\genStateTimeDiscreteSpaceDiscrete}}
\newcommand{\dSnapEigMatTruncate}{\dSnapEigMat^R}
\newcommand{\dSnapEigVecMatTruncate}{\dSnapEigVecMat^R}

\newcommand{\nTimeSteps}{N_t}

\newcommand{\dynamicsMat}{\mathbf{G}}
\newcommand{\dynamicsMatLspg}{\mathbf{G}_{\mathrm{LSPG}}}
\newcommand{\dynamicsMatGalerkin}{\mathbf{G}_{\mathrm{G}}}

\newcommand{\diffusionMat}{\mathbf{D}}
\newcommand{\diffusionMatRom}{\diffusionMat_r}

\newcommand{\boundaryLayerFigureCaption}{Example 1, boundary layer. }
\newcommand{\advectingFrontFigureCaption}{Example 2, advecting front. }

\usepackage[normalem]{ulem}

\begin{document}
\begin{frontmatter}

\title{Residual-based stabilized reduced-order models of the transient convection-diffusion-reaction equation obtained through discrete and continuous projection}
\author[a]{Eric Parish}
\author[b]{Masayuki Yano}
\author[a]{Irina Tezaur}
\author[c]{Traian Iliescu}
\address[a]{Sandia National Laboratories,  Livermore, CA}
\address[b]{Institute for Aerospace Studies, University of Toronto, Toronto, ON, M3H 5T6, Canada}
\address[c]{Department of Mathematics, Virginia Tech, Blacksburg, VA 24061}
\begin{abstract}
Galerkin and Petrov--Galerkin projection-based reduced-order models (ROMs) of transient partial differential equations are typically obtained by performing a dimension reduction and projection process that is defined at either the spatially continuous or spatially discrete level. In both cases, it is common to add stabilization to the resulting ROM to increase the stability and accuracy of the method; the addition of stabilization is particularly common for advection-dominated systems 
when the ROM is under-resolved.
While these two approaches can be equivalent in certain settings, differing techniques have emerged in both contexts. This work outlines these two approaches within the setting of finite element method (FEM) discretizations (in which case a duality exists between the continuous and discrete levels) of the convection-diffusion-reaction equation, and compares residual-based stabilization techniques that have been developed in both contexts. In the spatially continuous case, we examine the Galerkin, streamline upwind Petrov--Galerkin (SUPG), Galerkin/least-squares (GLS), and adjoint (ADJ) stabilization methods. For the GLS and ADJ methods, we examine formulations constructed from both the ``discretize-then-stabilize" technique and the space--time technique. In the spatially discrete case, we examine the Galerkin, least-squares Petrov--Galerkin (LSPG), and adjoint Petrov--Galerkin (APG) methods. We summarize existing analyses for these methods, and provide numerical experiments, which demonstrate that residual-based stabilized methods developed via continuous and discrete processes yield substantial improvements over standard Galerkin methods when the underlying FEM model is under-resolved.
\end{abstract}

\end{frontmatter}

\section{Introduction}
The numerical solution of parameterized partial differential
equations (PDEs) plays a vital role in numerous fields, including
engineering design and optimization, financial analysis and climate
sciences. Despite advances in high-performance computing and
numerical methods, numerically solving ``full-order models" (FOMs)
comprising discretized PDEs remains prohibitively expensive for a
variety of systems due to the presence of a disparate range of
spatiotemporal scales. This challenge is exacerbated for many-query
problems, such as optimization and uncertainty quantification, in
which case many executions of the forward model are required. A
variety of techniques have thus been developed to generate
approximations to the PDE of interest at a reduced computational
cost
\cite{balanced_truncation_moore,balanced_truncation_roberts,krylov_rom,Hesthaven2016,willcox_benner_rev}.

Projection-based reduced-order models, like those of the Galerkin and Petrov--Galerkin type investigated in this paper, are one promising
technique for approximating solutions to PDEs at a reduced
computational cost. These methods operate by (1) approximating the
state variables in a low-dimensional ``trial" space, and either (2a)
executing a projection process to enforce the resulting residual to
be orthogonal to a ``test" space (analogously, this step can be viewed
as computing the Petrov--Galerkin approximation of the weak form of the PDE in the low-dimensional trial and test spaces) or (2b) executing a residual
minimization process that computes a residual-minimizing solution
within the trial space.  The result of this process is a
low-dimensional system which is referred to as the ``reduced-order
model" (ROM).

Two differing types of projection-based ROMs of PDEs have emerged
over the past several decades: projection-based ROMs applying a
residual orthogonalization/minimization process at the
spatially-\textit{continuous} level, and projection-based ROMs
applying a residual orthogonalization/minimization process at the
spatially-\textit{discrete} level. In this work, we refer to these approaches as
\textit{\continuousROMs}\ and \textit{\discreteROMs}, respectively.
We note that, in the transient case considered herein,
reduced-order models are typically formulated by reducing the
spatial dimension of the model and leveraging standard time-marching
schemes for temporal discretization. We restrict our attention to
this setting, but note that several pieces of work have examined the
construction of space--time
ROMs~\cite{choi_stlspg,URBAN2012203,Yano2014stBoussinesq,constantine_strom,benner_st,parish_wls},
in which case the same thematic similarities of discrete vs.\ continuous are present.

In \continuousROMs, the state variables are approximated at the
spatially-continuous level in a low-dimensional function space, and
generalized coordinates associated with the state representation are
then obtained via, e.g., a Galerkin projection process that enforces
the residual to be orthogonal to a test space in a continuous
$L^2(\Omega)$ inner product, where $\Omega$ denotes the physical
domain. This step can be analogously viewed as computing an approximate
solution to the weak form of the PDE by employing low-dimensional
trial and test spaces.
\ContinuousROMs\ are most often employed within the context of weighted residual methods (e.g., finite element methods, spectral methods), and examples of \continuousROMs\ can be found in~\cite{rb_1,rb_2,rb_3,Rozza2008,URBAN2012203,Yano2014stBoussinesq,rowley2004model,bergmann2009enablers,caiazzo2014numerical,iliescu2014variational,rovas_thesis,Wang_ROM_thesis,Balajewicz_rom0,basis_rotation} (and many other works). We do note that several pieces of work have examined extensions to finite volume methods~\cite{StHiMoLo17,StRo18,LoCaLuRo16}.   

 \DiscreteROMs, on the other hand, work directly with the spatially-discrete system emerging after discretization of the differential operators present in the
 PDE\footnote{For transient PDEs, discrete ROMs can be formulated at either the
 ordinary \emph{differential} equation (ODE) level (i.e., after spatial discretization) or the ordinary \emph{difference} equation (O$\Delta$E level) (i.e., after spatial and temporal discretization);
 Ref.~\cite{carlberg_lspg_v_galerkin} examines commutativity of the time-discretization step for Galerkin and least-squares Petrov Galerkin ROM formulations.}. \DiscreteROMs\ approximate the \textit{discretized} state variables within a low-dimensional Euclidean vector space, and the generalized coordinates associated with the state representation are then obtained via, e.g., a Galerkin projection process that enforces the \textit{discrete} residual to be orthogonal in a \textit{discrete} (e.g., Euclidean) inner product. \DiscreteROMs\ are most often employed within the context of finite volume and finite difference methods (although they can also be applied to 
 FEMs) and examples can be found in~\cite{legresley_1,legresley_2,legresley_3,bui_resmin_steady,bui_unsteady,bui_thesis,carlberg_gnat,carlberg_lspg_v_galerkin,carlberg_conservative_rom,choi_stlspg,l1,wentland_apg}.
  Discrete ROMs are typically less intrusive than their continuous counterparts and are more generic, but arguably comprise a less rigorous modeling approach than continuous ROMs, as they may neglect important information about the underlying PDE. 

Regardless of whether performed at the continuous or discrete level,
the choice of projection/residual minimization process dictates the
stability and accuracy of a ROM. Galerkin projection, where the
residual is restricted to be orthogonal to the trial space, is the
most popular projection method. This popularity likely stems from
the fact that Galerkin projection yields optimal results in a given
energy norm for symmetric coercive systems\footnote{Throughout this work, we define a coercive system as a system whose bilinear form $a(\cdot,\cdot)$ satisfies the condition $a(v,v) \geq \alpha \| v \|_\star^2$ for all $v$ in a function space $V$ endowed with the norm $\| \cdot \|_\star$.  The condition is referred to as ``strong coercivity'' in some literature.}.
It is well-known, however, that for the convection-dominated systems Galerkin projection often lacks robustness in the presence of sharp, under-resolved gradients.
As a result, a variety of stabilization techniques have been developed to increase the ROM stability and accuracy for both continuous and discrete ROMs.

While, in this paper, we focus
on ``residual-based''
stabilization techniques, we mention that a large body of work has
been dedicated to alternative approaches targeting this issue.  This work includes, but is not
limited to: stabilizing inner
products that guarantee a non-increasing energy \cite{KALASHNIKOVA2014569, kalash_convergence, BARONE20091932, rowley2004model, serre}
or non-decreasing entropy~\cite{kalash_aiaa, Chan:2020}; stabilizing
subspace rotations that account for truncated modes \textit{a priori}~\cite{basis_rotation, Balajewicz_rom0}; 
eigenvalue reassignment methods that calculate a stabilizing correction to a given linear \cite{kalash_eig_reassign} or nonlinear \cite{Rezaian:2021} ROM that is found to be unstable after it is constructed;
structure preserving
methods that guarantee that the ROM satisfy physical constraints~\cite{carlberg_conservative_rom,LALL2003304,carlberg2012spd,structurePreserveBeattie,chaturantabut2016structure,farhat2014dimensional, Gruber:2023};
spatial filtering-based stabilization methods~\cite{gunzburger2019evolve,iliescu2018regularized,kaneko2020towards,strazzullo2022consistency,wells2017evolve} that filter out unphysical high-frequency content, \textit{inf-sup} stabilization methods that enforce the \textit{inf-sup} condition in the incompressible Stokes and Navier--Stokes equations~\cite{bergmann2009enablers,caiazzo2014numerical,decaria2020artificial}; 
and closure modeling
approaches~\cite{ahmed2021closures,
wang2012proper,bergmann2009enablers,
Wang_ROM_thesis,Wang:269133,San2018} that add additional ``closure'' terms to the ROM so-as to account for the impact of truncated modes.

We focus our review herein on what are referred to as ``residual-based"
stabilization techniques. Within the finite element community,
various residual-based techniques have been proposed in an effort to
develop robust numerical methods for non-symmetric, non-coercive,
and under-resolved problems.
These methods, which include stabilized finite
elements~\cite{brooks_thesis,BROOKS1982199,HuFrHu89,
FRANCA1995299} (e.g., streamline upwind Petrov--Galerkin (SUPG), Galerkin/least-squares (GLS)) and variational
multiscale (VMS) methods~\cite{hughes1}\footnote{It is noted that
these methods are not mutually orthogonal, for example VMS methods
can recover several stabilized methods.}, are typically formulated
by adding terms involving a sum of element-wise
integrals to the Galerkin method. These terms typically comprise the
product of a test function with the residual of the governing
equations, and thus the stabilized formulations can be written as Petrov--Galerkin projections. 
 
These ``residual-based" methods have proven to be quite
successful, yielding, for example, robust solutions for the
convection--diffusion equation, incompressible Navier--Stokes
equations, and compressible Navier--Stokes
equations~\cite{BROOKS1982199,
HUGHES1984217,Johnson_1984_FEM_Hyperbolic,HuFrHu89,TEZDUYAR19911,
FRANCA1992253}.
We note that, in the context of VMS methods, a body of work additionally exists that examines the addition of phenomenologically-inspired terms to the weighted residual form (e.g., eddy viscosity methods); we do not consider these approaches here and restrict our attention to residual-based methods. The extension of residual-based approaches to ROMs obtained via continuous projection is straightforward as they operate in a similar variational setting. As a result, various works have examined the formulation of stabilized reduced-order models via classical finite element stabilization techniques~\cite{sotomayor_thesis,bergmann2009enablers,GIERE2015454,kragel2005streamline,PACCIARINI20141,rozza_stabilized,doi:10.1002/fld.3777,stabile_vms_rom,rovas_thesis}.

Another class of stabilized model reduction methods that is important to highlight is continuous minimum-residual methods (or least-squares
methods)~\cite{Maday_2002_RB_Noncoercive, rovas_thesis}. These methods compute solutions within the trial space that minimize the residual of the governing PDE in the least-squares sense. 
Minimum-residual methods can be interpreted as Petrov--Galerkin methods, where 
the (parameter-dependent) test space is defined to be the one that maximizes the inf-sup stability constant, and hence the
method is guaranteed to be stable. We also note a related double greedy algorithm~\cite{Dahmen_2014_Double_Greedy}, which constructs a fixed test space for a Petrov--Galerkin formulation that approximately maximizes the inf-sup constant for all parameter values in a greedy fashion. While minimum-residual methods are robust and display commonalities with discrete stabilization approaches, we do not consider them here.

The extension of classic stabilization techniques to \discreteROMs\ is less straightforward. This is a consequence of the fact that \discreteROMs\ start from the spatially discrete level, and, as such, do not operate in the same variational setting as their continuous counterparts. As a result, various stabilization techniques have been developed for discrete ROMs. One class of particularly popular stabilization techniques are discrete residual minimization approaches~\cite{legresley_1,legresley_2,legresley_3,bui_resmin_steady,carlberg_thesis,bui_thesis,bui_unsteady,carlberg_gnat,carlberg_lspg,carlberg_lspg_v_galerkin}. These approaches compute a solution within the trial space that minimizes the discrete 
FOM residual. In the time-varying case, which we consider here, this residual minimization process is typically formulated by sequentially minimizing the time-discrete residual arising at each time instance on a discrete time grid\footnote{We note that recent work has examined windowed least-squares minimization~\cite{parish_wls} and space--time residual minimization~\cite{choi_stlspg}. These approaches have better stability properties than LSPG, but here we restrict our focus to the standard LSPG approach, for simplicity.}. This formulation is commonly referred to as the least-squares Petrov--Galerkin (LSPG) approach~\cite{carlberg_gnat,carlberg_lspg_v_galerkin}, and can be written as a Petrov--Galerkin projection of the FOM O$\Delta$E (i.e., the FOM ODE after temporal discretization). In finite element language, LSPG most closely resembles a discrete least squares principle. We refer to Grimberg et al. for an overview of LSPG within the context of stabilization~\cite{GrFaYo20}.
The adjoint--Petrov Galerkin (APG) method~\cite{parish_apg} is an additional discrete model reduction approach that falls into the class of residual-based methods. In APG, the variational multiscale method is applied at the discrete level to decompose the Euclidean state-space into coarse and fine-scale components. The impact of the fine scales on the coarse scales is then accounted for by virtue of a residual-based stabilization term that is derived from the Mori--Zwanzig formalism~\cite{Chorin_predictwithmem}. APG differs from FEM stabilization techniques in that the residual is defined at the level of the FOM ODE (i.e., after spatial discretization).

\begin{figure}
\begin{center}
\includegraphics[trim={0cm 7.5cm 0cm 2cm},clip,width=0.95\textwidth]{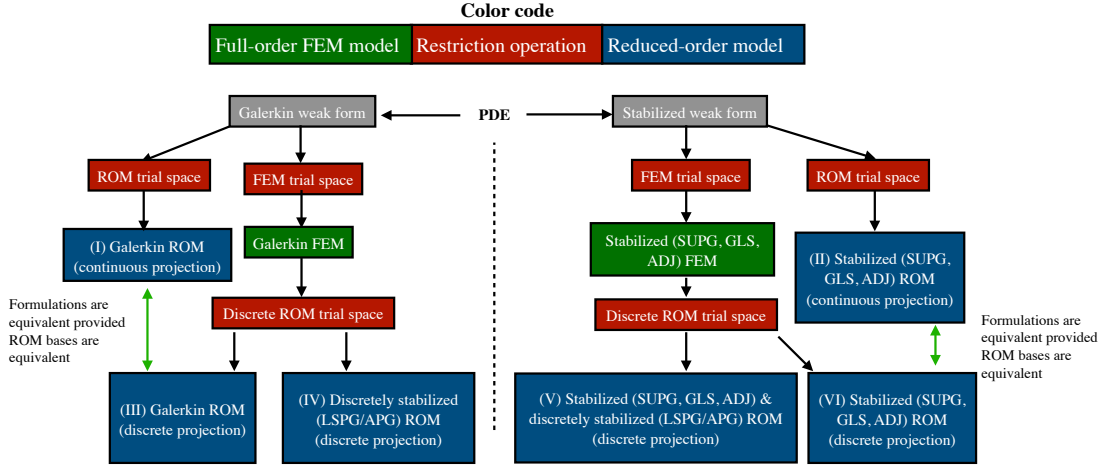}
\end{center}
\caption{Schematic of the various processes for constructing reduced-order models. Blocks in green comprise full-order models, blocks in red comprise restriction (projection) processes, and blocks in blue comprise reduced-order models.}
\label{fig:rom_types}
\end{figure}

Figure~\ref{fig:rom_types} provides a schematic of the various ROM
approaches just discussed within the context of finite element
discretizations. \ContinuousROMs\ (methods I and II in Figure~\ref{fig:rom_types}) rely on the definition of a weak
form and a ROM trial space, while \discreteROMs\ (methods III, IV, V, and VI in Figure~\ref{fig:rom_types}) rely on the
definition of a ``full-order" FEM system and a discrete ROM trial
space. In certain settings (e.g., finite
elements), it is straightforward to obtain the discrete form of the
ROM equations obtained via continuous projection. As a result,
provided proper selection of inner products at the discrete level,
it is well known that there is direct equivalence between continuous
and discrete ROMs (green arrows in Figure~\ref{fig:rom_types}); see, for instance,~\cite{rozza_monograph,volkwein2013proper}. This same duality does not exist within the context of finite volume methods, for instance. As a
result, the development and study of ROM methodologies has effectively forked
into bodies of work that start at the spatially continuous level,
and bodies of work that start at the spatially discrete level. This
fact is not well-documented in the community. To the best of the
authors' knowledge it is most clearly outlined in Refs.~\cite{KALASHNIKOVA2014569,kalashnikova2014reduced}.

The goal of this contribution is to help fill this gap and compare various stabilization techniques that have been developed in both settings. To this end, we outline the development of various continuous and discrete ROMs for the convection--diffusion-reaction equation. We outline these methods within the context of the 
FEM, in which case there is a duality between the spatially continuous and spatially discrete level. We establish precise and philosophical commonalities between the various approaches, and present several numerical examples assessing their performance. We also investigate the sensitivity of these methods with respect to the stabilization parameter and the time step.

The contributions of the present work are as follows:
\begin{enumerate}
\item We give the first side-by-side presentation of residual-based-stabilized ROMs developed through discrete and continuous projection.

\item We build on Ref.~\cite{KALASHNIKOVA2014569,kalashnikova2014reduced} to introduce a taxonomy for these various reduced-order modeling approaches. As will be seen in this manuscript, a ``Galerkin" approach may entail different ROMs to different communities.

\item We provide a summary of existing analyses for the various ROMs considered.

\item We present the first numerical comparison of discretely stabilized ROMs to continuously stabilized ROMs. This is particularly relevant in the context of the LSPG approach, which has existed in the community for over a decade.

\item We present the first study of the Galerkin/least-squares and adjoint stabilization methods applied to projection-based ROMs.

\item We provide a comprehensive study on the impact of the stabilization parameter, $\tau$, and time step, $\Delta t$, for all stabilized ROMs considered. This is the first such study that has been undertaken for the SUPG, GLS, and ADJ stabilization ROM methods.

\end{enumerate}

The layout of this manuscript is as follows.   In
Section~\ref{sec:math}, we outline the convection-diffusion-reaction
equation, the Galerkin 
FEM and stabilized FEM formulations for this equation.
In Sections~\ref{sec:croms} and~\ref{sec:droms}, we
outline ROMs constructed through continuous and discrete projection,
respectively. Alternate ROM stabilization approaches that do not
fall into either of these two categories are summarized in Section
\ref{sec:other}.
 Section~\ref{sec:analysis} summarizes the available theoretical support for the various methods.
Section~\ref{sec:experiments} presents numerical experiments, and
Section~\ref{sec:conclude} provides conclusions.

\section{Finite element discretizations for the convection-diffusion-reaction equation}\label{sec:math}
For concreteness, in this work, we consider the demonstrative example
of the convection-diffusion-reaction (\cdcdrAcronym) equation. We
emphasize that the concepts presented here generalize to other
systems, including nonlinear equations. The \cdcdrAcronym\ equation
is given by 
    \begin{equation}\label{eq:cdr}
       \begin{alignedat}{2}
 &\frac{\partial \cStateStrong}{\partial t} - \viscosity \Delta \cStateStrong + \wavespeed \cdot \nabla \cStateStrong + \reaction \cStateStrong = \forcing & \qquad & \text{in } \tDomain \times \Omega,  \\
&\cStateStrong(0,x) = \cState_0(x) , \qquad & & x \in \cDomain,  \\
&\cStateStrong(t,x) = 0  , \qquad & & x \in \cDomainBoundary , t \in \tDomain,
       \end{alignedat}
    \end{equation}
where $\cStateStrong:\tDomainClosed \times \cDomainClosed \rightarrow \RR{}$ is
the state implicitly defined as the solution to Eq.~\eqref{eq:cdr},
$\cDomain \subset \RR{\cDomainDim}$ is the physical domain,
$\cDomainBoundary$ is the domain boundary, $T \in \RPlus$ is the
final time, $\viscosity \in \RPlus$ is the viscosity, $\wavespeed
\in \RR{\cDomainDim}$ are the convection coefficients, $\reaction \in
\RPlusEqual$ is the reaction coefficient, $u_0 : \Omega \rightarrow
\RR{}$ is the initial condition, and $\forcing \in L^2(\cDomain)$ is
a forcing term. In what follows, we use the notation
$\cStateStrong(t) \equiv \cStateStrong(t,\cdot) : \cDomain
\rightarrow \RR{}$. We consider homogeneous boundary conditions, for
simplicity. In this setting, the system~\eqref{eq:cdr} is coercive.

We consider the standard weighted residual formulation of~\eqref{eq:cdr} in space, which reads as follows: find  
$\cState \in C^0((0,T]; L^2(\Omega)) \cap L^2((0,T]; \cSpace(\cDomain))$ such that $\forall t \in (0,T]$
\begin{equation}\label{eq:variational}
\cipGen{\cBasis}{\cState_t (t) } + \cipGen{ \viscosity \nabla \cBasis }{\nabla \cState(t)} + \cipGen{\cBasis}{\wavespeed \cdot \nabla \cState(t)} + \cipGen{\cBasis}{\reaction \cState(t)} = \cipGen{\cBasis}{\forcing}, \qquad \forall \cBasis \in \cSpace(\cDomain),
\end{equation}
and satisfies the initial condition $u(0) = u_0 \in L^2(\Omega)$, where $\cip : (v,w) \mapsto \int_{\Omega} v(x)  w(x)dx$ is the $L^2(\cDomain)$ inner product, and $\cSpace(\cDomain)$ is the standard Sobolev space of functions with square-integrable weak first derivatives that vanish on $\Gamma$.  The problem~\eqref{eq:variational} is well-posed; see, e.g.,~\cite{Quarteroni_2008_PDE}. For notational simplicity, we introduce the bilinear form
\begin{equation*}
\bilinearGalerkin : (\cBasisDum,\cStateTimeDiscreteDum) \mapsto  \cipGen{ \viscosity \nabla \cBasisDum }{\nabla \cStateTimeDiscreteDum} + \cipGen{\cBasisDum}{\wavespeed \cdot \nabla \cStateTimeDiscreteDum} + \cipGen{\cBasisDum}{\reaction \cStateTimeDiscreteDum}.
\end{equation*}


To transcribe~\eqref{eq:variational} into a discrete problem, we
need to introduce spatial and temporal discretizations. In this
work, we consider, for simplicity and without loss of generality, the implicit
Euler method for temporal discretization and a 
FEM 
for the spatial discretization; the concepts presented here can be extended to other time stepping schemes.
We introduce
without loss of generality a uniform
partition of the time domain $[0,T]$ into $\nTimeSteps+1$ time
instances $t^n = n \Delta t$, $n=0,\ldots,N_t$ with $\Delta t =
T/\nTimeSteps$. Application of the implicit Euler method yields the
series of strong form stationary PDEs for
$\cStateStrongTimeDiscreteArg{n}\left(\approx \cStateStrong \left(t^n\right) \right)$, $n=1,\ldots,\nTimeSteps$,
\begin{equation}\label{eq:cdr_td}
\begin{split}
&\frac{\cStateStrongTimeDiscreteArg{n} - \cStateStrongTimeDiscreteArg{n-1}}{\Delta t} - \viscosity \Delta \cStateStrongTimeDiscreteArg{n} + \wavespeed \cdot \nabla \cStateStrongTimeDiscreteArg{n} + \reaction \cStateStrongTimeDiscreteArg{n} = \forcing
\end{split}
\end{equation}
with $\cStateStrongTimeDiscreteArg{0} = \cState_0$ and $\cStateStrongTimeDiscreteArg{n} = 0$ on $\cDomainBoundary$, $n=1,\ldots,\nTimeSteps$. The weak form then yields the associated series of stationary problems: find $\cStateTimeDiscreteArg{n} (\approx \cStateStrongTimeDiscreteArg{n}) \in \cSpace(\cDomain)$, $n=1,\ldots,\nTimeSteps$, such that
\begin{equation}\label{eq:variational_timediscrete}
\cipGen{\cBasis}{\frac{\cStateTimeDiscreteArg{n} - \cStateTimeDiscreteArg{n-1}}{\Delta t} } + \bilinearGalerkinArg{ \cBasis}{\cStateTimeDiscreteArg{n}}  = \cipGen{\cBasis}{\forcing}, \qquad \forall \cBasis \in \cSpace(\cDomain),
\end{equation}
with initial condition $\cStateTimeDiscreteArg{0} = \cState_0$.

For spatial discretization, let $\cSpaceFEM \subset \cSpace(\cDomain)$ and $\testSpaceFEM \subset \cSpace(\cDomain)$ denote conforming trial and test spaces, respectively, obtained via a finite element discretization of $\cDomain$ into $\nel$ non-overlapping elements $\cDomain_k$, $k=1,\ldots,\nel$. The spatially discrete counterpart of~\eqref{eq:variational_timediscrete} reads:
 find $\cStateTimeDiscreteSpaceDiscreteArg{n} \in \cSpaceFEM$, $n=1,\ldots,\nTimeSteps$, such that
\begin{equation}\label{eq:variational_discrete}
\cipGen{\cBasis}{\frac{\cStateTimeDiscreteSpaceDiscreteArg{n} -
\cStateTimeDiscreteSpaceDiscreteArg{n-1}}{\Delta t} } +
\bilinearGalerkinArg{\cBasis}{\cStateTimeDiscreteSpaceDiscreteArg{n}}
 = \cipGen{\cBasis}{\forcing}, \qquad \forall \cBasis \in \testSpaceFEM,
\end{equation}
with (approximate) initial condition $\cStateTimeDiscreteSpaceDiscreteArg{n} = \cState_{0,h}$, where $\cState_{0,h}$ is, e.g., the $L^2(\Omega)$ projection of $\cState_0$ onto $\cSpaceFEM$.

\subsection{Galerkin approach}
The standard Galerkin approach is obtained by setting $\testSpaceFEM = \cSpaceFEM$ in~\eqref{eq:variational_discrete}.
We introduce the basis $\{\cBasis_i\}_{i=1}^\fomdim$ for $\cSpaceFEM$, which yields the following FOM basis vector $\forall x \in \cDomain$:
\begin{equation}\label{eq:fem_basis}
\cBasisVec(x) \equiv \begin{bmatrix} \cBasis_1 (x)& \cdots & \cBasis_{\fomdim} (x)\end{bmatrix}. 
\end{equation}
The time-discrete state at time-instance $t^n$ is described with these basis functions as
$\cStateTimeDiscreteSpaceDiscreteArg{n}(x) = \cBasisVec(x) \genStateTimeDiscreteSpaceDiscreteArg{n}$,
where $\genStateTimeDiscreteSpaceDiscreteArg{n} \in \RR{\fomdim}$, $n=0,\ldots,\nTimeSteps$. We refer to $\genStateTimeDiscreteSpaceDiscreteArg{n}$ as the \femCoefficients.
The Galerkin method yields the O$\Delta$E system to be solved for $\genStateTimeDiscreteSpaceDiscreteArg{n}$, $n=1,\ldots,\nTimeSteps,$
\begin{equation}\label{eq:g_fom_odeltae}
\residGalerkinFEM(\genStateTimeDiscreteSpaceDiscreteArg{n};\genStateTimeDiscreteSpaceDiscreteArg{n-1}) = \mathbf{0},
\end{equation}
where
$$\residGalerkinFEM: (\genStateTimeDiscreteSpaceDiscreteDumN;\genStateTimeDiscreteSpaceDiscreteDumNm) \mapsto \mass \left[\frac{ \genStateTimeDiscreteSpaceDiscreteDumN - \genStateTimeDiscreteSpaceDiscreteDumNm }{\Delta t}\right] + \matFull \genStateTimeDiscreteSpaceDiscreteDumN  - \forcingVec.$$
In the above, $\massArg{ij} \equiv \cipGen{\cBasis_i  }{\cBasis_j}  \in \RRS{\fomdim} $ is the FEM mass matrix, $\matFull_{ij}
\equiv 
\bilinearGalerkinArg{ \cBasis_i}{\cBasis_j} \in  \RR{\fomdim \times \fomdim}$ is a dynamics matrix resulting from the bilinear form, and
$\forcingVecArg{i}  \equiv \cipGen{ \cBasis_i}{f}$, $\forcingVec \in
\RR{\fomdim}$ is the discrete forcing; we denote the space of $N\times N$ symmetric positive
definite matrices by $\RRS{N}$. We refer to~\eqref{eq:g_fom_odeltae}
as the Galerkin FOM O$\Delta$E.

\begin{remark}
Obtaining the discrete problem~\eqref{eq:g_fom_odeltae} requires
evaluating the inner products in the
system~\eqref{eq:variational_discrete}. In general, evaluating these
inner products requires introducing a discrete quadrature rule, as done in \cite{BARONE20091932}. We
note that for linear problems (and problems displaying polynomial
nonlinearities) with piecewise polynomial forcing operators, the
inner products can be evaluated exactly with an appropriate
quadrature rule, e.g., Gaussian quadrature.
\end{remark}

\subsection{Residual-based stabilization}\label{sec:resid_stab}
\begin{figure}
\begin{center}
\begin{subfigure}[t]{0.49\textwidth}
\includegraphics[ trim={1cm 1cm 1cm 1cm},clip, width=1.\linewidth]{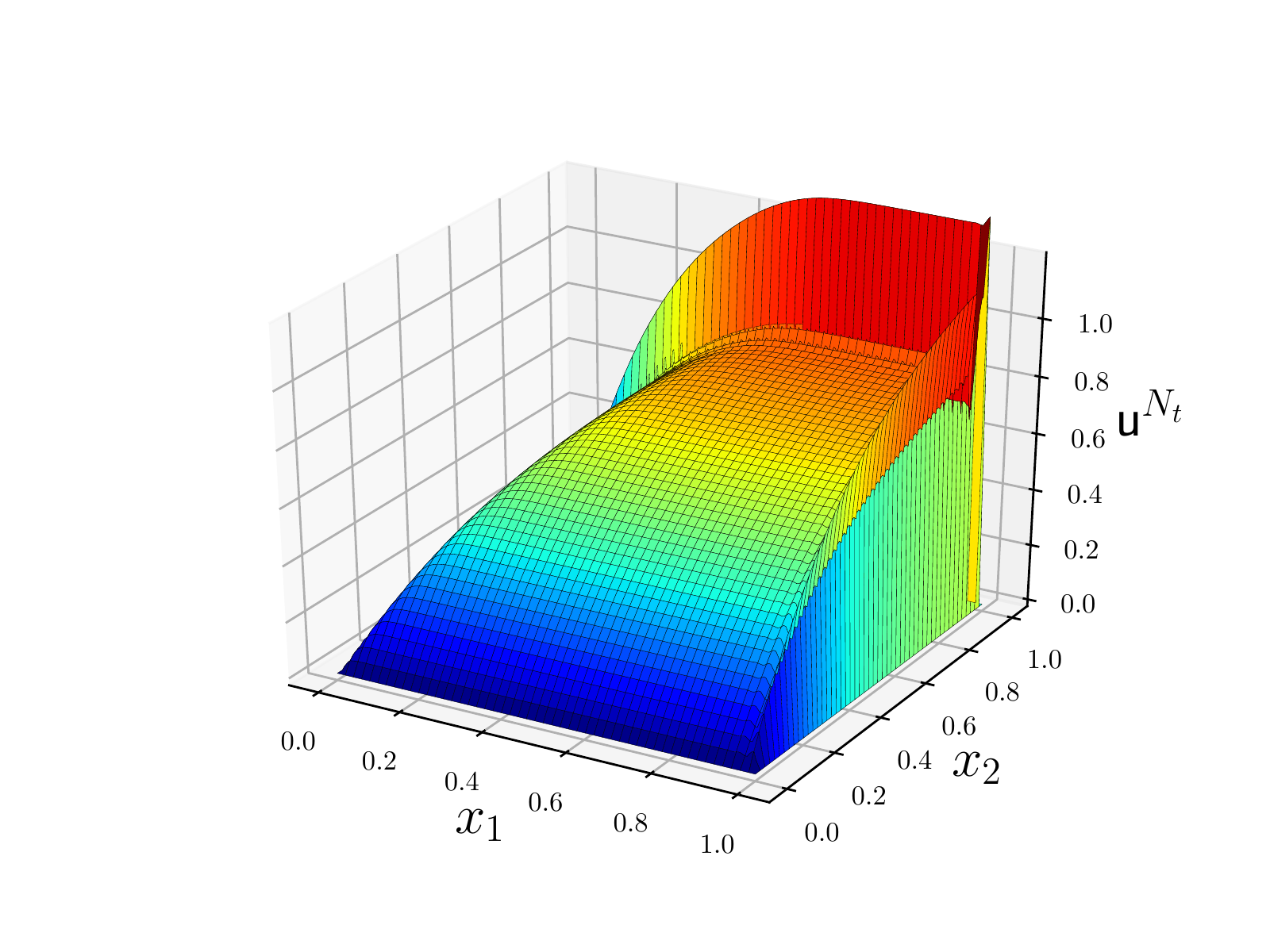}
\caption{Galerkin FEM}
\label{fig:galerkin_dem}
\end{subfigure}
\begin{subfigure}[t]{0.49\textwidth}
\includegraphics[trim={1cm 1cm 1cm 1cm},clip, width=1.\linewidth]{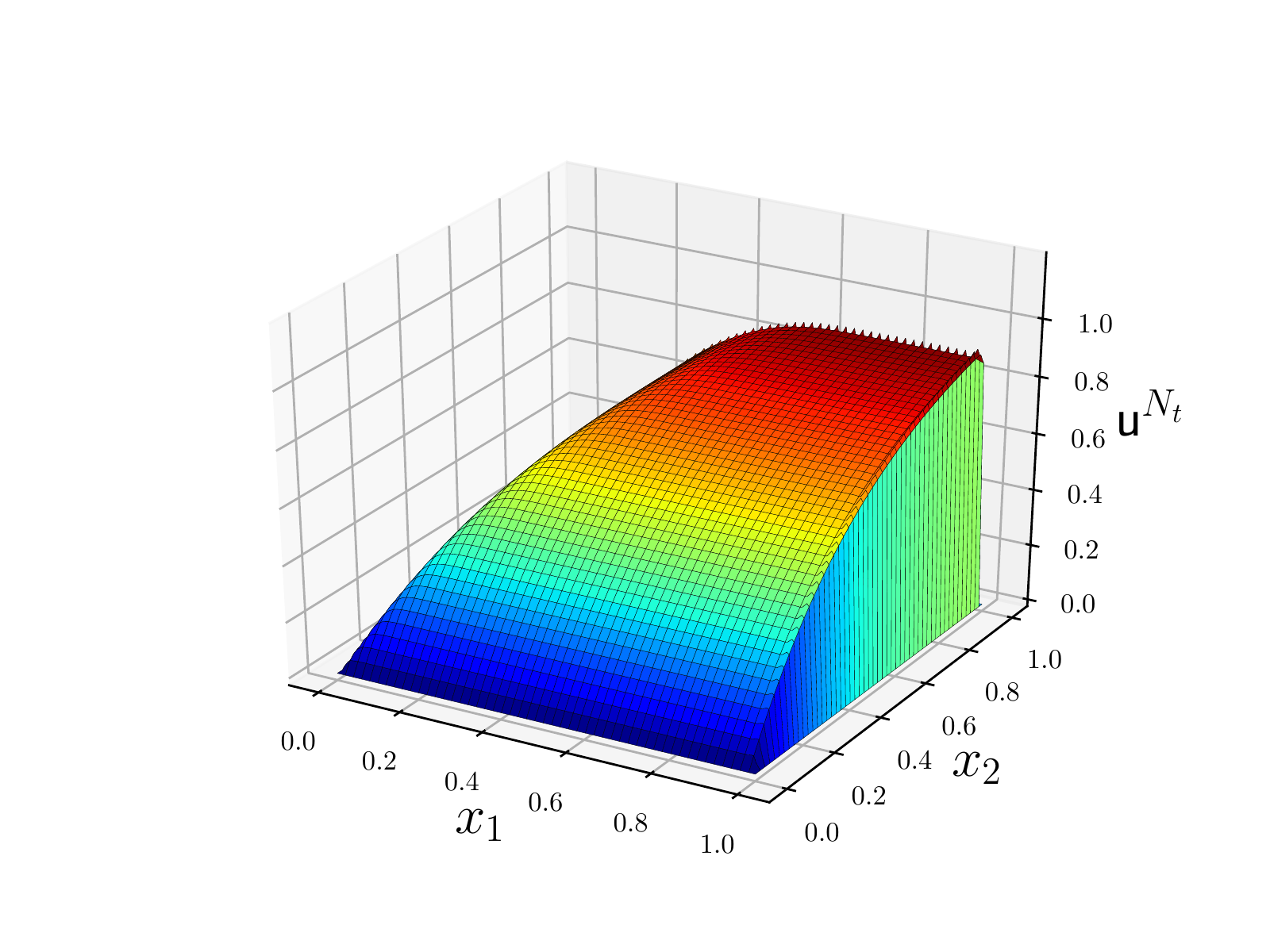}
\caption{SUPG FEM}
\label{fig:supg_dem}
\end{subfigure}
\caption{Finite element solutions to the \cdcdrAcronym\ equation, with the setup described in Section~\ref{sec:example1}, at $t=5$.}
\label{fig:fem_dem}
\end{center}
\end{figure}

The Galerkin approach is known to perform well for symmetric positive definite systems, in which case the Galerkin approach comprises a minimization principle in a system-specific energy norm.
In the presence of sharp, under-resolved gradients, however, it is well-known that the Galerkin approach can lack robustness. In the present
context, this poor performance is most pronounced for large
grid Peclet numbers (i.e., Pe$_{\text{h}}:=\norm{ \wavespeed} h/\viscosity \gg 1$,
where $h$ is a measure of the element size and $\norm{ \cdot}$ is the Euclidian norm). We note that large grid Peclet numbers occur, e.g., for coarse meshes (i.e., in the under-resolved regime), small diffusion coefficients, or a combination of both.
In this regime, the
skew-symmetric convection operator dominates the symmetric diffusion
operator. Figure~\ref{fig:fem_dem} demonstrates this by showing
finite element solutions to the \cdcdrAcronym\ equation obtained
using the Galerkin FEM, as well as a stabilized FEM.  The Galerkin approach is seen to yield large oscillations
near the boundary of the computational domain, while the stabilized
approach suppresses these oscillations and yields accurate
solutions\footnote{We note that one negative consequence of stabilized methods is that often the convergence rates for the methods are often lower than the unstabilized FEM.}.

To improve the performance of the numerical method in such regimes,
it is common to introduce stabilization to smooth the numerical solution. Various stabilization techniques exist, including flux limiters,
artificial viscosity, etc. In the finite element community,
\textit{residual-based stabilization} is a popular
stabilization technique. Residual-based FEMs, which include the
likes of the SUPG~\cite{BROOKS1982199,brooks_supg0,brooks_thesis}, GLS~\cite{HuFrHu89}, and
adjoint~\cite{FRANCA1995299,FRANCA1992253,hughes1} (ADJ) (also known as
unusual or subgrid-scale) stabilization methods, are typically
formulated by adding terms involving a sum of
element-wise integrals to the Galerkin method. These terms usually
comprise the product of a test function with the residual of the
governing equations. These approaches have been successful in
providing robust methodologies for a variety of systems, including
the CDR, incompressible Navier--Stokes,
and compressible Navier--Stokes equations~\cite{BROOKS1982199,
HUGHES1984217,Johnson_1984_FEM_Hyperbolic,HuFrHu89,TEZDUYAR19911,
FRANCA1992253,roos2008robust}.

For transient systems, like the \cdcdrAcronym\ equation described in this work, stabilized methods are typically employed in one of two ways:
\begin{itemize}
\item \textit{Space--time discretizations:} Space--time finite elements are employed in both time and space. The temporal dimension is then viewed as an additional spatial dimension, 
and standard stabilization approaches can be applied. 
This is the approach that was first employed for stabilized and variational multiscale methods of unsteady problems~\cite{SHAKIB199135,HUGHES1996217}.

\item \textit{Discretize-then-stabilize:} The PDE is first discretized in time, and then a stabilized method is applied to the time-discrete, spatially-continuous system. This approach is quite popular as it can be more computationally efficient than space--time discretizations and is compatible with numerous time marching schemes~\cite{CODINA20072413}.
\end{itemize}
In this work, we explore both approaches. 

A stabilized form of~\eqref{eq:variational_timediscrete} can be written generally as:
find $\cStateTimeDiscreteArg{n} \in \cSpaceFEM$, $n=1,\ldots,\nTimeSteps$, such that
\begin{equation}\label{eq:variational_timediscrete_stable}
m\left(\cBasisDum,\frac{\cStateTimeDiscreteArg{n} -
\cStateTimeDiscreteArg{n-1}}{\Delta t} \right) + 
\bilinearGalerkinArg{\cBasis}{\cStateTimeDiscreteArg{n}}
  + 
 m_{\text{el}}\left(\stabilizationOpTestArg{\cBasis} , \tau  \residualOp\left( \cStateTimeDiscreteArg{n} ,  \cStateTimeDiscreteArg{n-1} \right) \right) = \cipGen{\cBasis}{\forcing}, \qquad \forall \cBasis \in \cSpaceFEM,
\end{equation}
where $\cipGenElNoArg: (u,v) \mapsto \sum_{K=1}^{\nel}
\int_{\Omega_k} u(x)  v(x)dx$ denotes the sum of element-wise
$L^2(\cDomain)$ inner products, $\tau : \Omega \rightarrow \RR{}$
is a grid-dependent stabilization parameter,
$\residualOp : \left( \cStateTimeDiscreteArg{n} ,  \cStateTimeDiscreteArg{n-1} \right) \rightarrow 
\frac{\cStateTimeDiscreteArg{n} - \cStateTimeDiscreteArg{n-1}}{\Delta t} +  \stabilizationOpTrial \cStateTimeDiscreteArg{n} - \forcing$ 
is the strong form residual operator, 
 $\stabilizationOpTrial:
\cStateTimeDiscreteArg{n} \mapsto  -\viscosity \Delta
\cStateTimeDiscreteArg{n} + \wavespeed \cdot \nabla
\cStateTimeDiscreteArg{n} + \reaction \cStateTimeDiscreteArg{n}$,
and $\stabilizationOpTest$ is a linear stabilization operator that
is scheme dependent. 
Note that $\left[ \stabilizationOpTrial
\cStateTimeDiscreteArg{n} + \frac{ \cStateTimeDiscreteArg{n} - 
\cStateTimeDiscreteArg{n-1}}{\Delta t} - \forcing \right]$ yields
the strong form of the time-discrete residual. 
For notational simplicity, we denote the bilinear form associated with the stabilized formulations as
\begin{equation*}
\bilinearStabilized : (\cBasisDum,\cStateTimeDiscreteDum) \mapsto \bilinearGalerkinArg{\cBasisDum}{\cStateTimeDiscreteDum} + 
 m_{\text{el}}\left(\stabilizationOpTestArg{\cBasisDum} , \tau \left( \stabilizationOpTrial \cStateTimeDiscreteDum + \frac{\cStateTimeDiscreteDum }{\Delta t} \right) \right).
\end{equation*}

\begin{remark}
We note that the stabilized
form~\eqref{eq:variational_timediscrete_stable} is consistent with
respect to the continuous equations \eqref{eq:cdr} semi-discretized
in time using a implicit Euler scheme. That is, if one substitutes in the
exact solution, $\cStateTimeDiscreteArg{n} \leftarrow
\cStateStrongTimeDiscreteArg{n}$, the additional stabilization terms
in~\eqref{eq:variational_timediscrete_stable} vanish. We discuss
consistency in more detail in Section~\ref{sec:consistency}.
\end{remark}
The most popular types of stabilization methods are
SUPG~\cite{brooks_supg0,brooks_thesis}, GLS~\cite{HuFrHu89}, and ADJ \cite{HARARI20041491}. If a discretize-then-stabilize approach is taken, the
operator $\stabilizationOpTest$ takes the form
\begin{align}\label{eq:stabilizationOpsOne}
 \stabilizationOpTestSUPG &: \cBasisDum \mapsto \frac{1}{2}\left( \stabilizationOpTrial \cBasisDum - \stabilizationOpTrial^* \cBasisDum \right) \defeq \wavespeed \cdot \nabla
\cBasisDum,\\
 \stabilizationOpTestGLS &: \cBasisDum \mapsto  \frac{\cBasisDum}{\Delta t} + \stabilizationOpTrial \cBasisDum \defeq \frac{
\cBasisDum}{\Delta t}  -\viscosity \Delta
\cBasisDum + \wavespeed \cdot \nabla
\cBasisDum + \reaction \cBasisDum,\\
 \stabilizationOpTestADJ &: \cBasisDum \mapsto   -\frac{\cBasisDum}{\Delta t}   -\stabilizationOpTrial^* \cBasisDum \defeq -\frac{
\cBasisDum}{\Delta t}  +\viscosity \Delta
\cBasisDum + \wavespeed \cdot \nabla
\cBasisDum - \reaction \cBasisDum,
\end{align}
where the subscript ``DS" denotes ``discretize-then-stablize",  and $\stabilizationOpTrial^* : \cBasisDum \mapsto -\viscosity \Delta
\cBasisDum - \wavespeed \cdot \nabla \cBasisDum + \reaction \cBasisDum$ denotes the adjoint of
$\stabilizationOpTrial$.
For the space--time approach, we note that the implicit Euler method is equivalent to a space--time method with $p=0$ discontinuous Galerkin (DG) finite elements, and in this setting the space--time stabilization approach results in the operator $\stabilizationOpTest$ taking the form
\begin{align}\label{eq:stabilizationOpsTwo}
 \stabilizationOpTestSTSUPG &\equiv \stabilizationOpTestSUPG,\\
 \stabilizationOpTestSTGLS &: \cBasisDum \mapsto  \stabilizationOpTrial \cBasisDum \defeq   -\viscosity \Delta
\cBasisDum + \wavespeed \cdot \nabla
\cBasisDum + \reaction \cBasisDum,\\
 \stabilizationOpTestSTADJ &: \cBasisDum \mapsto   -\stabilizationOpTrial^* \cBasisDum \defeq \viscosity \Delta
\cBasisDum + \wavespeed \cdot \nabla
\cBasisDum - \reaction \cBasisDum, \label{eq:stabilizationOpsFinal}
\end{align}
where the subscript ``ST" denotes space--time.
\begin{remark}
We emphasize that, in this work, for space--time stabilization methods with the implicit Euler method, we employ the full residual operator $\residualOp$ in the right slot of the stabilization term as in \eqref{eq:variational_timediscrete_stable}. In Refs.~\cite{CODINA1998185,HuFrHu89}, the authors do not include the $\left( \cStateTimeDiscreteArg{n} - \cStateTimeDiscreteArg{n-1} \right) / \Delta t$ term for $p=0$ DG. We include this term for consistency at the time-discrete level. In our numerical experiments, we observed this term to make very little difference in practice.
\end{remark} 
Employing a finite element discretization
in space and leveraging the basis vector
\eqref{eq:fem_basis} yields the stabilized O$\Delta$E system to be solved for
$\genStateTimeDiscreteSpaceDiscreteArg{n}$,
$n=1,\ldots,\nTimeSteps$,
\begin{equation}\label{eq:stab_fom_odeltae}
\residStabilizedFEM(\genStateTimeDiscreteSpaceDiscreteArg{n};\genStateTimeDiscreteSpaceDiscreteArg{n-1}) = \mathbf{0}.
\end{equation}
The discrete residual of the stabilized discretization is given by
$$\residStabilizedFEM : (\genStateTimeDiscreteSpaceDiscreteDumN;\genStateTimeDiscreteSpaceDiscreteDumNm)  \mapsto \residGalerkinFEM(\genStateTimeDiscreteSpaceDiscreteDumN;\genStateTimeDiscreteSpaceDiscreteDumNm) + \stabilizationMat \genStateTimeDiscreteSpaceDiscreteDumN - \forcingVecStabilized - \massStabilized \frac{\genStateTimeDiscreteSpaceDiscreteDumNm}{\Delta t},$$
with $ \stabilizationMatArg{ij} = \cipGenEl{ \stabilizationOpTest \cBasis_i}{\tau \left( \stabilizationOpTrial \cBasis_j + \frac{\cBasis_j}{\Delta t} \right)}$, $\forcingVecStabilizedArg{i} = \cipGenEl{ \stabilizationOpTest \cBasis_i}{\tau \forcing},$ and $\massStabilizedArg{ij} = \cipGenEl{ \stabilizationOpTest \cBasis_i}{\tau \cBasis_j}.$

\subsection{Selection of the stabilization parameter, $\tau$}
The stabilized form~\eqref{eq:variational_timediscrete_stable}
requires specification of the stabilization parameter $\tau$. The
\textit{a priori} selection of suitable stabilization parameters has
been a topic of much research; see, for
example,~\cite{HSU2010828,HARARI20041491,CODINA20072413} and
references therein. Traditionally, the stabilization constants are
obtained through asymptotic scaling arguments~\cite{HSU2010828}, and
depend on, e.g., the grid size,  the diffusion coefficient and, for transient problems, the time step.

Relevant to the current work is the fact that classical definitions
of the stabilization parameters for transient problems are subject
to several issues. First and foremost, in addition to depending on
the spatial grid resolution, classical definitions of $\tau$ depend
on the time step. These definitions become poorly behaved in
both the low time step \textit{and} steady-state regimes. In
all numerical experiments considered in the work we present
results for numerous values of $\tau$ and as such do not restrict
ourselves to a particular definition. We refer to~\cite{CODINA1998185,HSU2010828,roos2008robust} for examples of
definitions for $\tau$.

\subsection{Sensitivity to the time step, $\Delta t$}
In addition to depending on the stabilization parameter $\tau$, it
is well-known that stabilized formulations depend on the time step
$\Delta t$. Stability analyses have demonstrated, for example, that
stabilized formulations may become unstable at low CFL
numbers~\cite{BOCHEV20042301}. This sensitivity to the time step can
be understood intuitively for GLS and ADJ stabilization,  where
changing the time step size changes the nature of the stabilization
operator $\stabilizationOpTest$. Thus, changing the time step
modifies both the error incurred due to temporal discretization
\textit{and} the properties of the stabilized scheme. As will be
seen later in the manuscript, the least-squares Petrov--Galerkin
approach suffers from similar issues. In our numerical examples, we present results for a variety of time steps to quantify this
dependence.

\section{Continuous projection reduced-order models}\label{sec:croms}
We now develop ROMs of the \cdcdrAcronym\ system via
continuous projection. Continuous ROMs generate approximate
solutions $\romStateTimeDiscreteSpaceDiscreteArg{n} (\approx \cStateTimeDiscreteArg{n})$ within a low-dimensional
spatial trial space $\romStateTimeDiscreteSpaceDiscreteArg{n} \in \cSpaceRom  \subset \cSpaceFEM
\subset \cSpace(\cDomain)$, and have been studied in a number of
references including~\cite{Balajewicz_rom0, basis_rotation,
BARONE20091932, GIERE2015454, iliescu2014variational,
KALASHNIKOVA2014569, kalashnikova2014reduced, rowley2004model,
wang2012proper,
Wang_ROM_thesis}. 
 Various techniques exist for constructing this trial space, and here
we consider proper orthogonal decomposition
(POD)~\cite{berkooz_turbulence_pod}. To construct the trial space
through POD, we assume access to an ensemble of snapshots at time
instances $t^n$, $n=0,\ldots,\nTimeSteps$\footnote{In practice, snapshots are often collected at only a subset of the time steps.  Additionally, snapshots can be collected for a variety of parameter values, in the case of a parametrized PDE.}. We collect these snapshots
into the matrix
$$\snapshotMatrixContinuous = \begin{bmatrix}
\cStateTimeDiscreteArg{0} & \cdots & \cStateTimeDiscreteArg{\nTimeSteps}
\end{bmatrix}.
$$
The POD method seeks to
find an $\innerProductType$-orthonormal basis of rank $\romdim \ll \fomdim$ (where $\fomdim$ is the size of the FOM from which the reduced basis is built) that minimizes the projection error
\begin{equation}\label{eq:pod_min_problem}
\underset{\{\cBasisRom_i\}_{i=1}^R, \phi_i \in \cSpaceFEM}{\text{minimize}}\;  \sum_{n=0}^{\nTimeSteps} \left|
\left| \cStateTimeDiscreteArg{n} - \sum_{j=1}^{\romdim}
\cipGenArbitrary{\innerProductType}{\cStateTimeDiscreteArg{n}}{\cBasisRom_j}\cBasisRom_j
\right|\right|^2_{\innerProductType},
\end{equation}
where $\phi_j : \Omega \rightarrow \RR{}$, $j=1,\ldots,\romdim$, are ROM basis functions and $\innerProductType$ denotes inner product type (e.g, $\HOne$, $L^2(\cDomain)$, weighted $L^2(\cDomain)$~\cite{BARONE20091932}). 
The minimization problem~\eqref{eq:pod_min_problem} can be solved via the eigenvalue problem
$$\cSnapCorr \cSnapEigVecMat = \cSnapEigVecMat  \cSnapEigMat,$$
where $[\cSnapCorr]_{ij}  = \cipGenArbitrary{\innerProductType}{[\snapshotMatrixContinuous]_i}{[\snapshotMatrixContinuous]_j}\in \RRS{\nSnap \times \nSnap}$, and $\cSnapEigVecMat$ and $\cSnapEigMat$ are the matrices associated with the eigenvectors and eigenvalues, respectively. Assuming the snapshot matrix is full rank, it can be shown that the minimizer of the problem~\eqref{eq:pod_min_problem} is
\begin{equation}\label{eq:pod_min_solution_continuous}
\cBasisRomVec = \snapshotMatrixContinuous \cSnapEigVecMat  \sqrt{\cSnapEigMat^{-1}}.
\end{equation}
For each $x \in \cDomain$, we evaluate these basis functions at $x$ and assemble the ROM basis vector  
$\cBasisRomMat(x) \equiv \begin{bmatrix} \cBasisRom_1(x) & \cdots & \cBasisRom_{\romdim}(x) \end{bmatrix}$ with
$\cBasisRomMat(x) \in \RR{1 \times \romdim}$.  We then set $\cSpaceRom \equiv \text{span}\{ \phi_1, \dots, \phi_R \}$.
We additionally note that, as $\cSpaceRom \subset \cSpaceFEM$, it directly follows that the ROM basis vectors can be described with a linear combination of the FOM basis vectors, i.e.,
$\cBasisRomMat(x) = \cBasisVec(x) \cBasisRomcoefMat$, where $\cBasisRomcoefMat \in \RR{\fomdim \times \romdim}$ is a coefficient matrix.

\subsection{Galerkin reduced-order models}
The Galerkin ROM achieved after time discretization is: find $\romStateTimeDiscreteSpaceDiscreteArg{n} \in \cSpaceRom$, $n=1,\ldots,\nTimeSteps$, such that
\begin{equation}\label{eq:variational_discrete_grom}
\cipGen{\cBasisRom}{\frac{\romStateTimeDiscreteSpaceDiscreteArg{n} - \romStateTimeDiscreteSpaceDiscreteArg{n-1}}{\Delta t} } + \cipGen{ \viscosity \nabla \cBasisRom }{\nabla \romStateTimeDiscreteSpaceDiscreteArg{n}} + \cipGen{\cBasisRom}{\wavespeed \cdot \nabla \romStateTimeDiscreteSpaceDiscreteArg{n}} + \cipGen{\cBasisRom}{\reaction
 \romStateTimeDiscreteSpaceDiscreteArg{n}} = \cipGen{\cBasisRom}{\forcing}, \qquad \forall \cBasisRom \in \cSpaceRom.
\end{equation}
Leveraging the ROM basis vector $\cBasisRomVec$, the  Galerkin ROM can be cast as the sequence of O$\Delta$Es
to be solved for $\genRomStateTimeDiscreteSpaceDiscreteArg{n}$, $n=1,\ldots,\nTimeSteps$,
\begin{equation}\label{eq:g_rom_odeltae}
\residGalerkinROM(\genRomStateTimeDiscreteSpaceDiscreteArg{n};\genRomStateTimeDiscreteSpaceDiscreteArg{n-1}) = \mathbf{0},
\end{equation}
where $\genRomStateTimeDiscreteSpaceDiscreteArg{n} \in \RR{\romdim}$ are the ROM ``generalized coordinates" such that the approximate state is defined as $\romStateTimeDiscreteSpaceDiscreteArg{n} = \cBasisRomMat \genRomStateTimeDiscreteSpaceDiscreteArg{n}$, and the residual operator is given by 
$$\residGalerkinROM: (\genRomStateTimeDiscreteSpaceDiscreteDumN;\genRomStateTimeDiscreteSpaceDiscreteDumNm) \mapsto \massRom \left[ \frac{ \genRomStateTimeDiscreteSpaceDiscreteDumN - \genRomStateTimeDiscreteSpaceDiscreteDumNm }{\Delta t} \right] + \matFullRom \genRomStateTimeDiscreteSpaceDiscreteDumN  - \forcingVecRom,$$
with $\matFullRomArg{ij}  \equiv  \bilinearGalerkinArg{ \cBasisRom_i}{\cBasisRom_j} \in \RR{\romdim \times \romdim} $, $\massRomArg{ij} \equiv \cipGen{\cBasisRom_i  }{\cBasisRom_j} \in \RRS{\romdim}$ and $\forcingVecRomArg{i}  \equiv \cipGen{ \cBasisRom_i}{\forcing } \in \RR{\romdim}$. The Galerkin ROM formulation~\eqref{eq:variational_discrete_grom} corresponds to method I in Figure~\ref{fig:rom_types}.


\subsection{Residual-based stabilized reduced-order models}\label{sec:resid_stab_rom}
Analogously to the FEM case, the Galerkin ROM can lack robustness in the presence of sharp, under-resolved gradients. In the present context, such a situation may arise when the diffusion constant is small (relatively to ${\bf b}$) and not enough ROM basis vectors are employed to capture the behavior of the solution of interest.
A stabilized ROM formulation can be written generally as:
find $\romStateTimeDiscreteSpaceDiscreteArg{n} \in \cSpaceRom$,
$n=1,\ldots,\nTimeSteps$, such that
\begin{multline}\label{eq:variational_timediscrete_stable_rom}
\cipGen{\cBasisRom}{\frac{\romStateTimeDiscreteSpaceDiscreteArg{n} - \romStateTimeDiscreteSpaceDiscreteArg{n-1}}{\Delta t} } + \cipGen{ \viscosity \nabla \cBasisRom}{\nabla \romStateTimeDiscreteSpaceDiscreteArg{n}} + \cipGen{\cBasisRom}{\wavespeed \cdot \nabla \romStateTimeDiscreteSpaceDiscreteArg{n}} + \cipGen{\cBasisRom}{\reaction
 \romStateTimeDiscreteSpaceDiscreteArg{n}}  +
 \\  m_{\text{el}}\left(\mathcal{Q}\phi,  \tau \residualOp \left( \romStateTimeDiscreteSpaceDiscreteArg{n} , \romStateTimeDiscreteSpaceDiscreteArg{n-1}\right) \right) = \cipGen{\cBasisRom}{\forcing}, \qquad \forall \cBasisRom \in \cSpaceRom
\end{multline}
where $\mathcal{Q}$ can be any of the forms given in Eqns.~\eqref{eq:stabilizationOpsOne}--\eqref{eq:stabilizationOpsFinal}. We consider the SUPG, \SDGLS, \SDADJ, \STGLS, and \STADJ\ ROMs as defined by the operators in Eqns.~\eqref{eq:stabilizationOpsOne}--\eqref{eq:stabilizationOpsFinal}.
Leveraging the ROM basis vectors yields the stabilized O$\Delta$E system
to be solved for $\genRomStateTimeDiscreteSpaceDiscreteArg{n}$, $n=1,\ldots,\nTimeSteps$,
\begin{equation}\label{eq:stab_rom_odeltae}
\residStabilizedRom(\genRomStateTimeDiscreteSpaceDiscreteArg{n};\genRomStateTimeDiscreteSpaceDiscreteArg{n-1}) = \mathbf{0}.
\end{equation}
The discrete residual of the stabilized discretization is given by
$$\residStabilizedRom : (\genRomStateTimeDiscreteSpaceDiscreteDumN;\genRomStateTimeDiscreteSpaceDiscreteDumNm)  \mapsto \residGalerkinROM(\genRomStateTimeDiscreteSpaceDiscreteDumN;\genRomStateTimeDiscreteSpaceDiscreteDumNm) + \stabilizationMatRom \genRomStateTimeDiscreteSpaceDiscreteDumN - \massStabilizedRom \frac{\genRomStateTimeDiscreteSpaceDiscreteDumNm}{\Delta t} - \forcingVecStabilizedRom ,$$
with $ \stabilizationMatRomArg{ij} = \cipGenEl{ \stabilizationOpTest \cBasisRom_i}{\tau \left( \stabilizationOpTrial \cBasisRom_j + \frac{\cBasisRom_j}{\Delta t} \right)} \in \RR{\romdim \times \romdim}$,
$ \forcingVecStabilizedRomArg{i} =  \cipGenEl{ \stabilizationOpTest \cBasisRom_i}{\tau \forcing } \in \RR{\romdim}$,
and  $\massStabilizedRomArg{ij}= \cipGenEl{ \stabilizationOpTest \cBasisRom_i}{\tau \cBasisRom_j} \in \RR{\romdim \times \romdim}$. The stabilized ROM formulations~\eqref{eq:variational_timediscrete_stable_rom} corresponds to method II in Figure~\ref{fig:rom_types}.

In the literature, SUPG ROMs have been considered in~\cite{bergmann2009enablers,GIERE2015454,JoMoNo22,PACCIARINI20141,mclaughlin2016stabilized,kragel2005streamline,zoccolan2023streamline,zoccolan2023stabilized}, amongst others. To the best of our knowledge, ROMs based on the other stabilized formulations discussed above have not been studied in the literature.

\subsection{Selection of the stabilization parameter, $\tau$}
    \label{sec:tau}
Like in the standard FEM, the stabilized ROM form~\eqref{eq:variational_timediscrete_stable_rom} requires
specification of the stabilization parameter $\tau$. \textit{A priori} selection of this parameter in the context of
 ROMs is not as well explored as in the standard FEM case. It remains unclear, for instance, if the stabilization parameter should take
 on a different value for the ROM as opposed to the full-order FEM discretization.  While selecting a different value of $\tau$ for the
 ROM as the one used in the FOM would lead to a lack of consistency between the ROM and the FOM, it may nonetheless improve the
 stability and accuracy of the ROM. To the best of the authors knowledge,~\cite{GIERE2015454}, which explores the selection of the stabilization parameter in the context of the SUPG approach, 
 is the only present study that examines the selection of $\tau$ within the ROM context. Specifically,~\cite{GIERE2015454} employs the same strategy as that used in standard FEM: a SUPG-ROM error bound is first proved, and $\tau$ is chosen to minimize this bound.  
 However, since the ROM space is a subspace of the FEM space, two types of inverse inequalities are used to prove the SUPG-ROM bound: a standard FEM inverse inequality, and a ROM inverse inequality~\cite{kunisch1}.
 These two inverse inequalities yield two SUPG-ROM error bounds, which in turn yield two $\tau$ scalings: a standard FEM scaling in which $\tau$ depends on the FEM mesh size, and a new ROM scaling in which $\tau$ depends on the ROM parameters (e.g., the ROM dimension, the POD basis functions, and the corresponding eigenvalues).
 We note that other approaches leveraging residual-based stabilization for
 ROMs (see, e.g.,~\cite{PACCIARINI20141}) use standard definitions of $\tau$ inherited from the FEM community in which case $\tau$ depends on the FEM discretization (as opposed to the
 resolution of the ROM basis). We additionally note that~\cite{parish_apg} numerically explored the selection of $\tau$ within the context of the APG method, where a relationship between the optimal value of $\tau$ and the spectral radius of the Jacobian of the right-hand side operator (which varies with ROM dimension) was observed.


\section{Discrete projection reduced-order models}\label{sec:droms}
ROMs developed through continuous projection operate in a weighted residual setting defined at the
spatially-continuous level. ROMs developed through discrete projection, however, perform model
reduction at either the level of the FOM ODE or FOM O$\Delta$E; Ref.~\cite{carlberg_lspg_v_galerkin} shows that these
two approaches are equivalent for the Galerkin method. In this work, we restrict our discussion to discrete projection-based ROMs developed at
the O$\Delta$E level.

Discrete projection ROMs approximate the degrees of freedom
associated with the spatial discretization in a low-dimensional
(vector) trial space,
$\discreteROMGenStateTimeDiscreteSpaceDiscreteArg{n} (\approx
\genStateTimeDiscreteSpaceDiscreteArg{n}) \in \dSpaceROM \subseteq
\RR{\fomdim}$, $n=0,\ldots,\nTimeSteps$, where $\dSpaceROM$ is the discrete ROM trial subspace. 
We again employ POD to construct this space. Towards this end, on the vector space
$\RR{\fomdim}$ we first define the
$\discreteInnerProductType$-weighted inner product $\dipGenNoArg(\cdot,\cdot)_{\discreteInnerProductType}:
(\mathbf{U},\mathbf{V}) \mapsto \mathbf{U}^T
\discreteInnerProductType \mathbf{V}$, where
$\discreteInnerProductType  \in \RRS{\fomdim}$ is a
symmetric-positive definite weighting matrix. The associated
$\discreteInnerProductType$-weighted norm is $\| \mathbf{x}
\|_{\discreteInnerProductType}^2 = \mathbf{x}^T
\discreteInnerProductType \mathbf{x}$. Next, we assume access to an
ensemble of snapshots of the \femCoefficients\ at time instances
$t^n$, $n=0,\ldots,\nTimeSteps$. We then seek a
$\discreteInnerProductType$-orthonormal basis of rank $\romdim$ that
minimizes the projection error
\begin{equation}\label{eq:discrete_pod_min_problem}
\underset{\dBasisRomMat \in \RR{\fomdim \times \romdim}, \; \dBasisRomMat^T \discreteInnerProductType \dBasisRomMat = \mathbf{I}}{ \text{minimize} }\; \sum_{n=0}^{\nTimeSteps} \| \genStateTimeDiscreteSpaceDiscreteArg{n} - \dBasisRomMat\dBasisRomMat^T \discreteInnerProductType \genStateTimeDiscreteSpaceDiscreteArg{n}   \|_{\discreteInnerProductType} ^2,
\end{equation}
The solution to the minimization problem~\eqref{eq:discrete_pod_min_problem} can be obtained via an eigenvalue problem or via the generalized singular value decomposition; we present the former here. We denote the snapshots of FEM coefficients as
$$\snapshotMatrixDiscrete \equiv \begin{bmatrix}
\genStateTimeDiscreteSpaceDiscreteArg{0} & \cdots & \genStateTimeDiscreteSpaceDiscreteArg{\nTimeSteps}
\end{bmatrix} \in \RR{N \times N_t + 1}
.$$
We note that $\snapshotMatrixContinuous = \cBasisVec \snapshotMatrixDiscrete$.
Defining the time correlation matrix as $\dSnapCorr = \dipGen{ \snapshotMatrixDiscrete}{\snapshotMatrixDiscrete}{\discreteInnerProductType},$
we can leverage the eigenvalue problem
$$\dSnapCorr \dSnapEigVecMat = \dSnapEigVecMat  \dSnapEigMat$$
to obtain the POD bases. The solution to the minimization problem~\eqref{eq:discrete_pod_min_problem} can be shown to be
$$\dBasisRomMat =  \snapshotMatrixDiscrete  \dSnapEigVecMatTruncate \sqrt{[\dSnapEigMatTruncate]^{-1}},$$
where $\dSnapEigVecMatTruncate$ and $\dSnapEigMatTruncate$ comprise the first $\romdim$ columns of $\dSnapEigVecMat$ and the first $\romdim$ columns and rows of $\dSnapEigMat$, respectively.
\begin{remark}
Setting $\discreteInnerProductType_{ij} \leftarrow \cipGenArbitrary{\innerProductType}{\cBasisVec_i}{\cBasisVec_j}$, we can express the correlation matrix as
 $$ [\dSnapCorr]_{ij} =  \cipGenArbitrary{\innerProductType}{[ \snapshotMatrixContinuous]_i }{[\snapshotMatrixContinuous ]_j},$$
which recovers the correlation matrix used in continuous projection ROMs. Further, as $\snapshotMatrixContinuous = \cBasisVec \snapshotMatrixDiscrete$, we can express Eq.~\eqref{eq:pod_min_solution_continuous} as $\cBasisRomVec(x) = \cBasisVec(x) \snapshotMatrixDiscrete \cSnapEigVecMat  \sqrt{\cSnapEigMat^{-1}}$ and we see that
$$\cBasisRomVec(x) = \cBasisVec(x) \dBasisRomMat.$$
We emphasize that this result is well-documented in the community, see, e.g.,~\cite{volkwein2013proper,rozza_monograph}.

\end{remark}

\subsection{Galerkin reduced-order model}
The Galerkin ROM developed through discrete projection is obtained by (i) making the substitution $\genStateTimeDiscreteSpaceDiscreteArg{n} \leftarrow \dBasisRomMat \genDiscreteRomStateTimeDiscreteSpaceDiscreteArg{n}$, $n=0,\ldots,\nTimeSteps$, and (ii) restricting the residual of the FOM O$\Delta$E to be $\discreteInnerProductTypeTwo$-orthogonal to the vector trial space $\dSpaceROM$. Here $\genDiscreteRomStateTimeDiscreteSpaceDiscreteArg{n} \in \RR{\romdim}$, $n=0,\ldots,\nTimeSteps$ are the ROM generalized coordinates and $\discreteInnerProductTypeTwo \in \mathbb{S}^{\fomdim}$ is a weighting matrix inducing the inner product $\dipGenNoArg(\cdot,\cdot)_{\discreteInnerProductTypeTwo}:
(\mathbf{U},\mathbf{V}) \mapsto \mathbf{U}^T
\discreteInnerProductTypeTwo \mathbf{V}$ with the associated
$\discreteInnerProductTypeTwo$-weighted norm is $\| \mathbf{x}
\|_{\discreteInnerProductTypeTwo}^2 = \mathbf{x}^T
\discreteInnerProductTypeTwo \mathbf{x}$; $\discreteInnerProductTypeTwo$ may or may not be equivalent to $\discreteInnerProductType$. It is critical to note that the Galerkin ROM developed via discrete projection can be developed for any FOM O$\Delta$E; e.g., the FOM O$\Delta$E could associate with the Galerkin FOM O$\Delta$E~\eqref{eq:g_fom_odeltae}, \textit{or} it could associate with the stabilized FOM O$\Delta$E~\eqref{eq:stab_fom_odeltae}.

We denote the residual of a generic FOM O$\Delta$E\footnote{We assume the FOM O$\Delta$E to depend only on the state at the current time instance and previous time instance, as would be the case with an implicit Euler temporal discretization.} as

$$\resid: (\genDiscreteRomStateTimeDiscreteSpaceDiscreteDumN;\genDiscreteRomStateTimeDiscreteSpaceDiscreteDumNm) \mapsto
\resid(\genDiscreteRomStateTimeDiscreteSpaceDiscreteDumN;\genDiscreteRomStateTimeDiscreteSpaceDiscreteDumNm). $$
Examples of this residual are $\resid = \residGalerkinFEM$ for association with the Galerkin FOM O$\Delta$E~\eqref{eq:g_fom_odeltae} and $\resid = \residStabilizedFEM$ for association with the stabilized FOM O$\Delta$E~\eqref{eq:stab_fom_odeltae}.
The Galerkin ROM obtained via discrete projection yields the O$\Delta$E system to be solved for $\genDiscreteRomStateTimeDiscreteSpaceDiscreteArg{n}$, $n=1,\ldots,\nTimeSteps$,
\begin{equation}\label{eq:drom_odeltae}
\residGalerkinDiscreteRom(\genDiscreteRomStateTimeDiscreteSpaceDiscreteArg{n};\genDiscreteRomStateTimeDiscreteSpaceDiscreteArg{n-1}) = \mathbf{0},
\end{equation}
where the residual of the discretely projected Galerkin ROM is given by
\begin{equation*}
\residGalerkinDiscreteRom :(\genDiscreteRomStateTimeDiscreteSpaceDiscreteDumN;\genDiscreteRomStateTimeDiscreteSpaceDiscreteDumNm) \mapsto
\dipGen{\dBasisRomMat}{\resid(\dBasisRomMat \genDiscreteRomStateTimeDiscreteSpaceDiscreteDumN;\dBasisRomMat \genDiscreteRomStateTimeDiscreteSpaceDiscreteDumNm)}{\discreteInnerProductTypeTwo}.
\end{equation*}
The discrete Galerkin ROM formulation~\eqref{eq:drom_odeltae} corresponds to methods III and VI in Figure~\ref{fig:rom_types}, depending on the underlying FEM model.
\begin{remark}\label{remark:discrete_g_remark}
Setting $\discreteInnerProductType_{ij} \leftarrow \cipGenArbitrary{\innerProductType}{\cBasis_i}{\cBasis_j}$ in
optimization problem~\eqref{eq:discrete_pod_min_problem},
and $\discreteInnerProductTypeTwo \leftarrow \mathbf{I}$, $\resid \leftarrow \residGalerkinFEM$ in problem~\eqref{eq:drom_odeltae}, the Galerkin ROM obtained via discrete projection~\eqref{eq:drom_odeltae} recovers the Galerkin ROM obtained via continuous projection~\eqref{eq:g_rom_odeltae}.
\end{remark}
\begin{remark}\label{remark:discrete_stab_remark}
Analogously to Remark~\ref{remark:discrete_g_remark}, setting $\discreteInnerProductType_{ij} \leftarrow \cipGenArbitrary{\innerProductType}{\cBasis_i}{\cBasis_j}$ in
optimization problem~\eqref{eq:discrete_pod_min_problem}, and $\discreteInnerProductTypeTwo \leftarrow \mathbf{I}$, $\resid \leftarrow \residStabilizedFEM$ in problem~\eqref{eq:drom_odeltae}, the Galerkin ROM obtained via discrete projection~\eqref{eq:drom_odeltae} recovers the stabilized ROM obtained via continuous projection~\eqref{eq:stab_rom_odeltae}.
\end{remark}

\subsection{Least-squares Petrov--Galerkin reduced-order model}\label{sec:LSPG}
Similar to the continuous Galerkin ROM, the discrete Galerkin ROM has been observed to yield inaccurate or unstable solutions in a variety of settings and thus various stabilization approaches have been developed for discrete ROMs. The LSPG approach comprises one particularly popular stabilization approach for discrete ROMs. LSPG operates by computing a sequence of solutions $\genDiscreteRomStateTimeDiscreteSpaceDiscreteArg{n}$, $n=1,\ldots,\nTimeSteps$, that satisfy the minimization problem
\begin{equation}\label{eq:LSPG_minproblem}
\genDiscreteRomStateTimeDiscreteSpaceDiscreteArg{n} = \underset{\ddgenStateDum  \in \RR\romdim}{\text{arg\,min}} \; \normArg{\lspgWeight}{ \resid(\dBasisROMMat \ddgenStateDum;\dBasisROMMat \genDiscreteRomStateTimeDiscreteSpaceDiscreteArg{n-1}) }^2,
\end{equation}
where $\resid$ is again the residual of the FOM O$\Delta$E.
The optimization problem~\eqref{eq:LSPG_minproblem} can be solved via the first-order optimality conditions, which yield the sequence of algebraic equations for $\genDiscreteRomStateTimeDiscreteSpaceDiscreteArg{n}, n=1,\ldots,\nTimeSteps$,
$$\dipGen{ \frac{\partial \resid}{\partial \ddStateDum}(\dBasisROMMat \genDiscreteRomStateTimeDiscreteSpaceDiscreteArg{n}) \dBasisROMMat }{\resid(\dBasisROMMat \genDiscreteRomStateTimeDiscreteSpaceDiscreteArg{n};\dBasisROMMat \genDiscreteRomStateTimeDiscreteSpaceDiscreteArg{n-1})}{\discreteInnerProductTypeTwo} = \mathbf{0}, $$
where $ \frac{\partial \resid}{\partial \ddStateDum}$ is the Jacobian of the residual $\resid(\cdot,\cdot)$ with respect to the first argument.
In the case $\resid \leftarrow \residGalerkinFEM$, the optimality conditions become
$$\dipGen{\frac{\mass \dBasisROMMat}{\Delta t} + \matFull \dBasisROMMat }{\residGalerkinFEM(\dBasisROMMat \genDiscreteRomStateTimeDiscreteSpaceDiscreteArg{n};\dBasisROMMat \genDiscreteRomStateTimeDiscreteSpaceDiscreteArg{n-1})}{\discreteInnerProductTypeTwo} = \mathbf{0}.$$
We see that LSPG takes the form of a Petrov--Galerkin ROM and hence we classify it as a residual-based method.
The LSPG ROM formulation~\eqref{eq:LSPG_minproblem} corresponds to methods IV and V in Figure~\ref{fig:rom_types}, depending on the underlying FEM.

\begin{remark} (LSPG can correspond to a continuous minimization principal.) 
Setting $\resid \leftarrow \residGalerkinFEM$ and $\discreteInnerProductTypeTwo \leftarrow \mass^{-1}$ in optimization problem~\eqref{eq:LSPG_minproblem}, LSPG corresponds to the continuous minimization principle for $\romStateTimeDiscreteSpaceDiscreteArg{n}$, $n=1,\ldots,\nTimeSteps$,
\begin{equation}\label{eq:LSPG_minproblem_continuous}
\romStateTimeDiscreteSpaceDiscreteArg{n} =
\underset{\cStateTimeDiscreteDum \in \cSpaceRom  }{\text{arg\,min}}
\; \int_{\cDomain} \left(
\cdcrResidTimeDiscreteSpaceContinuousOrtho(\cStateTimeDiscreteDum;\romStateTimeDiscreteSpaceDiscreteArg{n-1})
\right)^2 dx,
\end{equation}
where
$$
\cdcrResidTimeDiscreteSpaceContinuousOrtho : (\cStateTimeDiscreteDumN ;\cStateTimeDiscreteDumNm) \mapsto
\cBasisVec \mass^{-1} \cipGen{\cBasisVec}{ \cdcrResidTimeDiscreteSpaceContinuous(\cStateTimeDiscreteDumN ;\cStateTimeDiscreteDumNm)}
$$
and 
$$\cdcrResidTimeDiscreteSpaceContinuous : (\cStateTimeDiscreteDumN ;\cStateTimeDiscreteDumNm) \mapsto
\frac{\cStateTimeDiscreteDumN - \cStateTimeDiscreteDumNm}{\Delta t} - \nu \nabla^2 \cStateTimeDiscreteDumN + \wavespeed
\cdot \nabla \cStateTimeDiscreteDumN + \reaction
 \cStateTimeDiscreteDumN - \forcing.
$$
LSPG computes the solution $\romStateTimeDiscreteSpaceDiscreteArg{n}$ within the ROM trial space $\cSpaceRom$ that minimizes the $\LTwo$-norm of the time-discrete, spatially continuous residual projected onto the finite element trial space $\cSpaceFEM$. The full derivation for this equivalence is presented in Appendix~\ref{appendix:lspg_proof}.
\end{remark}

\begin{remark}
For the case $\resid \leftarrow \residStabilizedFEM$, it is not clear if LSPG corresponds to an underlying residual minimization principle defined at the continuous level.
\end{remark}

\subsubsection{Selection of the time step, $\Delta t$}
While LSPG does not contain a stabilization parameter, its performance depends on the time step and time integration scheme~\cite{carlberg_lspg_v_galerkin}. This is due to the fact that changing the time step (i) modifies the error incurred due to temporal discretization and (ii) modifies the LSPG minimization problem (i.e., the time-discrete residual changes). 
As a result, LSPG yields best results at an intermediary time step~\cite{carlberg_lspg_v_galerkin}. LSPG lacks robustness for too small a time step (in the limit $\Delta t \rightarrow 0$ LSPG recovers the Galerkin approach~\cite{carlberg_lspg_v_galerkin}) and too large a time step. Minimal work has examined the \textit{a priori} selection of an appropriate time step.

\subsection{Adjoint Petrov--Galerkin reduced-order model}
The final residual-based stabilization technique considered in this work is the APG method~\cite{parish_apg}. APG is a VMS-based approach for constructing \discreteROMs, and is derived from a time-continuous ODE setting. 
APG is derived via a multiscale decomposition of $\mathbb{R}^N$ into a coarse-scale, resolved trial space $\dSpaceROM$ and a fine-scale, unresolved trial space, $\dSpaceROMFine$ such that $\mathbb{R}^N = \dSpaceROM \oplus \dSpaceROMFine$. The impact of fine scales on the coarse-scale dynamics is then accounted for by virtue of the Mori--Zwanzig formalism and the variational multiscale method~\cite{parish_apg}. Setting $\discreteInnerProductType \leftarrow \mass$ in~\eqref{eq:discrete_pod_min_problem}, associating with the Galerkin FOM O$\Delta$E~\eqref{eq:g_fom_odeltae}, and using the implicit Euler method for time-discretization, APG yields the sequence of O$\Delta$E's to be solved for $\genDiscreteRomStateTimeDiscreteSpaceDiscreteArg{1}, \ldots, \genDiscreteRomStateTimeDiscreteSpaceDiscreteArg{\nTimeSteps}$,
\begin{equation}\label{eq:APG_problem}
\residAPG(\genDiscreteRomStateTimeDiscreteSpaceDiscreteArg{n}; \genDiscreteRomStateTimeDiscreteSpaceDiscreteArg{n-1}) =
\mathbf{0}.
\end{equation}
The APG residual is given by
$$
\residAPG: (\genDiscreteRomStateTimeDiscreteSpaceDiscreteDumN;
\genDiscreteRomStateTimeDiscreteSpaceDiscreteDumNm) \mapsto
\dipGenArg{\mathbf{I}}{ \left( \mathbf{I} - \tauApg [\PiFine]^T \big[
\matFull \big]^T \right) \dBasisROMMat
}{\residGalerkinFEM(\dBasisROMMat
\genDiscreteRomStateTimeDiscreteSpaceDiscreteDumN;\dBasisROMMat
\genDiscreteRomStateTimeDiscreteSpaceDiscreteDumNm)} = \mathbf{0},$$
where $\PiFine = \mass^{-1} - \dBasisROMMat \dBasisROMMat^T$,
$\PiFine : \RR{\fomdim} \rightarrow \dSpaceROMFine$, and $\tauApg \in \RPlus$ is a stabilization parameter. The full derivation for 
APG is provided in Appendix~\ref{appendix:derive_apg}.


\begin{remark}
The APG approach displays conceptual similarities with adjoint stabilization and the (quasi-static) orthogonal subscales (OSS) approach from the variational multiscale method~\cite{codina_oss,BAIGES2015173,REYES2020112844}; see Ref.~\cite{parish_apg} for details. There is no clear \textbf{direct} equivalence between these approaches, however. This is a result of APG being formulated at the discrete level, while adjoint stabilization and orthogonal subscales are formulated at the continuous level.
\end{remark}

Like LSPG, APG could also associate with a stabilized FOM. Some of the stabilized FEM formulations considered in this work, however, are developed at the time-discrete level (e.g., \ADJ\ and \GLS). As APG is derived from a time-continuous setting, it is not straightforward to construct an APG ROM of all stabilized formulations, and we only consider the \SUPG, \STGLS, and \STADJ\ FEM models.  We note that this is due to the fact that the test functions in these formulations do not contain terms of the form $\frac{v}{\Delta t}$. The APG ROM associating with one of these stabilized FEM models is obtained by setting $\discreteInnerProductType \leftarrow \mass$ in~\eqref{eq:discrete_pod_min_problem}, associating with the stabilized FEM O$\Delta$E~\eqref{eq:g_fom_odeltae} with $\stabilizationOpTest \leftarrow \stabilizationOpTestSUPG,\stabilizationOpTestSTGLS$, or $\stabilizationOpTestSTADJ$, and using the implicit Euler method for time-discretization. This process yields the sequence of O$\Delta$E's to be solved for $\genDiscreteRomStateTimeDiscreteSpaceDiscreteArg{1}, \ldots, \genDiscreteRomStateTimeDiscreteSpaceDiscreteArg{\nTimeSteps}$,
\begin{equation}\label{eq:APG_problem_stabilized}
\residAPGSUPG(\genDiscreteRomStateTimeDiscreteSpaceDiscreteArg{n}; \genDiscreteRomStateTimeDiscreteSpaceDiscreteArg{n-1}) =
\mathbf{0},
\end{equation}
where
$$
\residAPGSUPG: (\genDiscreteRomStateTimeDiscreteSpaceDiscreteDumN;
\genDiscreteRomStateTimeDiscreteSpaceDiscreteDumNm) \mapsto
\dipGenArg{\mathbf{I}}{ \left( \mathbf{I} - \tauApg [\PiFine]^T \left[
\matFull + \stabilizationMat \right]^T \right) \dBasisROMMat
}{\residStabilizedFEM(\dBasisROMMat
\genDiscreteRomStateTimeDiscreteSpaceDiscreteDumN;\dBasisROMMat
\genDiscreteRomStateTimeDiscreteSpaceDiscreteDumNm)} = \mathbf{0}.$$
The APG ROM formulation~\eqref{eq:APG_problem} corresponds to method IV Figure~\ref{fig:rom_types}, while the APG ROM formulation~\eqref{eq:APG_problem_stabilized} corresponds to method V in Figure~\ref{fig:rom_types}.

%

\subsection{Summary of remarks for discrete ROMs}
A summary of the remarks provided in this section is as follows:
\begin{itemize}

\item \textbf{The discrete POD basis recovers the continuous POD basis}.  The POD bases obtained through discrete projection recover the POD bases obtained through continuous projection under the conditions $\discreteInnerProductType_{ij} \leftarrow \cipGenArbitrary{\innerProductType}{\cBasisVec_i}{\cBasisVec_j}$ in optimization problem~\eqref{eq:discrete_pod_min_problem}.

\item \textbf{The discrete Galerkin ROM recovers the continuous Galerkin ROM.} The discrete Galerkin ROM recovers the continuous Galerkin ROM under the conditions $\discreteInnerProductTypeTwo \leftarrow \mathbf{I}$ in problem~\eqref{eq:drom_odeltae}, $\discreteInnerProductType_{ij} \leftarrow \cipGenArbitrary{\innerProductType}{\cBasisVec_i}{\cBasisVec_j}$ in
optimization problem~\eqref{eq:discrete_pod_min_problem}, and $\resid \leftarrow \residGalerkinFEM$ in problem~\eqref{eq:drom_odeltae}.

\item \textbf{The discrete Galerkin ROM recovers stabilized ROMs.}
The discrete Galerkin ROM recovers the stabilized continuous ROM under the conditions $\discreteInnerProductTypeTwo \leftarrow \mathbf{I}$ in problem~\eqref{eq:drom_odeltae}, $\discreteInnerProductType_{ij} \leftarrow \cipGenArbitrary{\innerProductType}{\cBasisVec_i}{\cBasisVec_j}$ in
optimization problem~\eqref{eq:discrete_pod_min_problem}, and $\resid \leftarrow \residStabilizedFEM$ in problem~\eqref{eq:drom_odeltae}.

\item \textbf{LSPG mimics a continuous $L^2(\cDomain)$ minimization principle.} LSPG mimics a continuous $L^2(\cDomain)$ minimization principle under the conditions $\resid \leftarrow \residGalerkinFEM$ and $\discreteInnerProductTypeTwo \leftarrow \mass^{-1}$ in optimization problem~\eqref{eq:LSPG_minproblem}.

\item \textbf{APG displays similarities to adjoint stabilization.} Similar to APG, ADJ can also be derived from the variational multiscale method. For transient systems ADJ results in a set of equations that are conceptually similar to APG, but without the appearance of an orthogonal projector. We note that FEM approaches for orthogonal subscales do exist, e.g.,~\cite{codina_oss,BAIGES2015173,REYES2020112844}, and APG also displays similarities with these approaches. 
\end{itemize}

\section{Other stabilized ROMs} \label{sec:other}


Although the main focus of this paper is on residual-based stabilized ROMs for the CDR equation, in this section we briefly outline some of the stabilized ROMs not covered herein.
This work includes, but is not limited to, ROM stabilizations that are not residual-based, and ROM stabilizations for equations different from the convection-diffusion-reaction equation (e.g., the incompressible Navier-Stokes equations and, especially, the compressible Euler equations). We outline several such techniques here:

\begin{itemize}

    \item {\it Closure models}
    that add additional ``closure" terms to the ROM in order to account for the impact of truncated modes.
    For classical numerical discretizations (e.g., finite element, finite volume, or spectral methods) of turbulent incompressible or compressible flows, there is an extensive literature on closure models, especially in large eddy simulation (LES)~\cite{sagaut2006large}.
    Closure models for ROMs (see~\cite{ahmed2021closures} for a survey) have also been developed, using ideas from different fields, e.g., image processing~\cite{xie2017approximate}, data-driven modeling~\cite{xie2018data,hijazi2019data}, machine learning~\cite{san2018machine}, information
      theory~\cite{majda2018model} the Mori-Zwanzig formalism from statistical mechanics~\cite{lin2019data}, or dynamical systems~\cite{chekroun2019variational}.
    We emphasize, however, that (just as in LES) arguably the most popular type of ROM closure models are the eddy viscosity models~\cite{HLB96,rebollo2017certified,wang2012proper}, which add a dissipative term to the standard ROM.
    These eddy viscosity ROM closures are generally built by invoking physical arguments, i.e., the concept of energy cascade, which states that in three-dimensional (3D) turbulent flows energy is transferred from large scales to small scales where the energy is dissipated~\cite{CSB03,sagaut2006large}.

    The eddy viscosity ROM closure models are similar {\it in spirit} to residual-based ROM stabilization methods since they both increase the numerical stability of the ROM.
    There are, however, notable differences.
    Probably the most important difference is that the two strategies target different aspects of the ROM inaccuracy in under-resolved simulations: The ROM closures target the {\it cause} of the problem, i.e., they model the term representing the effect of the discarded modes on the ROM dynamics.
    In contrast, stabilized ROMs target the {\it symptoms}: since under-resolved ROMs generally yield spurious numerical oscillations, adding numerical stabilization can often cure the problem. 
    Another difference between the two strategies is that 
    ROM eddy viscosity closure models are generally built for nonlinear equations, e.g., the Navier-Stokes equations or the quasi-geostrophic equations that model the large scale ocean circulation~\cite{san2015stabilized} (see, however,~\cite{koc2019commutation} for a notable exception regarding the commutation error).
    Residual-based stabilized ROMs, on the other hand, are built for  both linear (as in this paper) and nonlinear equations.
    Furthermore, the eddy viscosity ROM closure models are based on physical arguments, and, as a result, eddy viscosity ROMs are often designed to be deployed on a specific equation set or discretization. Residual-based stabilized ROMs constitute a more generic modeling strategy, broadly speaking, but may be less effective in domain-specific applications. 

    \item Stabilizing {\it inner products}, which are used to construct more stable ROMs.
    One of the earliest examples of stabilizing inner products is the $\HOne$ inner product that is used in~\cite{iollo2000stability} instead of the standard $L^2$ inner product.
    Other examples of stabilizing inner products are present in the literature for ROMs applied to multistate 
    systems where a classic vector $L^2$ inner product does not result in a physically meaningful energy principle. For example, energy-based projections are proposed for compressible flows in Refs.~\cite{BARONE20091932,iollo2000stability,kalashnikova2014reduced,rowley2004model, serre}, which use different inner products and flow variables to construct stabilized ROMs.  Recently, Ref.~\cite{kaptanoglu2020physics} has proposed similar ideas within the context of magnetohydrodynamics.  We note that the LSPG ROM-based preconditioning approach developed in \cite{LindsayEtAl2022} can be interpreted as a modification to the underlying inner product that has the effect of scaling different solution components to ensure that they are all of roughly the same magnitude.  

    \item {\it Inf-sup} stabilizations, that aim at enforcing the \textit{inf-sup} (or the Ladyzhenskaya-Babuska-Brezzi (LBB)) condition in incompressible Stokes and Navier-Stokes equations.
    We emphasize that the \textit{inf-sup} condition is used to ensure the well-posedness of saddle-point problems, such as the incompressible Stokes and Navier-Stokes equations.
    Thus, the \textit{inf-sup} stabilizations are different from the stabilization of convection-dominated systems, such as those we consider in this paper.

    In standard (e.g., FEM) numerical discretizations 
    of the Stokes and the Navier-Stokes equations, it is well-known that not enforcing the
    \textit{inf-sup} condition can yield spurious numerical oscillations in the pressure field.
    There are two main approaches to tackle this issue:
    (i) choose finite elements that do satisfy the \textit{inf-sup} condition, or
    (ii) choose finite elements that do not satisfy the \textit{inf-sup} condition and add pressure stabilization.

    In the ROM community, the first \textit{inf-sup} stabilizations have been proposed in pioneering work by the reduced basis methods (RBM) group~\cite{Rozza_2007_RB_Stokes}, 
    who developed ROMs that satisfy the
    \textit{inf-sup} condition (which is a significantly more difficult task than for finite elements, since the velocity and pressure ROM bases are problem dependent).
    Recognizing that enforcing the \textit{inf-sup} condition at a ROM level can be prohibitively expensive~\cite{ballarin2015supremizer}, more efficient stabilized ROMs that do not satisfy the
    \textit{inf-sup} condition were devised by using, e.g., the penalty method~\cite{bergmann2009enablers,caiazzo2014numerical}, 
    artificial compressibility~\cite{decaria2020artificial}, or local projection stabilization~\cite{rubino2020numerical}.
    

    \item {\it Structure preserving} methods that guarantee that the ROM satisfies the same physical constraints as those satisfied by the underlying equations.
    As for classical numerical discretizations, preserving these physical constraints generally yield more stable ROMs.
    For example, for the incompressible Navier-Stokes equations, ROMs in which the nonlinear terms preserve the kinetic energy are more stable than standard ROMs~\cite{kondrashov2015data,loiseau2018constrained} (see also~\cite{mohebujjaman2019physically} for ROM closure modeling and~\cite{kaptanoglu2020physics} for work in magnetohydrodynamics).
    Furthermore, ROMs that preserve Lagrangian structure were developed in~\cite{carlberg2015preserving,LALL2003304}, and ROMs that preserve Hamiltonian structure were constructed in~\cite{afkham2017structure,gong2017structure,peng2016symplectic, Sockwell:2019}.  The recent work by Gruber \textit{et al.} is the first to construct ROMs in which the more general metriplectic structure is preserved \cite{Gruber:2023}.

    \item Stabilizing {\it basis modification} 
    methods designed to remedy the so-called ``mode truncation instability", that is, to account for truncated modes \textit{a priori}
    ~\cite{amsallem_stab, basis_rotation, Balajewicz_rom0}.  In \cite{amsallem_stab}, 
    Amsallem and Farhat develop a non-intrusive method for stabilizing linear time-invariant (LTI) ROMs through the minimal modification of the left ROM basis.  The new reduced-order basis is obtained by formulating and solving a small-scale convex constrained optimization problem in which the constraint 
    imposes asymptotic stability of the modified ROM.   In \cite{basis_rotation, Balajewicz_rom0}, Balajewicz \textit{et al.} demonstrate that a ROM for (nonlinear) fluid flow can be stabilized through a \textit{stabilizing rotation} of the projection subspace.  
Specifically, the projection subspace is ``rotated" into a more dissipative regime by modifying the eigenvalue distribution of the linear operator.
Mathematically, the approach
is formulated as a trace minimization on the Stiefel manifold.  Like the approach in \cite{amsallem_stab},
the methods in \cite{Balajewicz_rom0, basis_rotation}
are non-intrusive and do not require any empirical eddy viscosity closure modeling terms. 

    \item {\it Spatial filtering}-based stabilization~\cite{gunzburger2019evolve,iliescu2018regularized,kaneko2020towards,wells2017evolve}, in which explicit filtering performed either in the physical space or in the ROM space is used to regularize/smooth different terms in the underlying equations, e.g., the convective term in the incompressible Navier-Stokes equations.
    Due to their simplicity, modularity, and effectiveness, spatially-filtered regularized models have been extensively studied in standard CFD (e.g., with finite element discretizations, surveyed in~\cite{layton2012approximate}), only a few regularized ROMs have been proposed in both deterministic~\cite{kaneko2020towards,sabetghadam2012alpha,wells2017evolve} and stochastic~\cite{gunzburger2019evolve,iliescu2018regularized} settings.



\item  \textit{Eigenvalue reassignment methods}, which calculate a stabilizing correction to an unstable ROM after the ROM has been constructed.
The correction is computed offline by solving a constrained optimization problem.
The approach was originally developed in Kalashnikova \textit{et al.} \cite{kalash_eig_reassign} in the context of 
LTI systems, for which it is natural to impose a constraint on the Lyapunov stability of the ROM system by requiring that the eigenvalues of the ROM matrix defining the problem have negative real parts.  The approach was subsequently extended to the nonlinear compressible flow equations by Rezaian and Wei in \cite{Rezaian}.  Here, appropriate constraints on the system energy, namely that it is non-increasing, were developed and applied.  
Eigenvalue reassignment methods are non-intrusive by construction, as they operate on a ROM \textit{a posteriori} (i.e., after the ROM has been constructed), and can be effective regardless
of 
the nature of the instability.  The methods can also be used to assimilate data into a given ROM, again after the model has been constructed. 
\\



\end{itemize}

\section{Brief survey of numerical analysis of residual-based ROM stabilizations}\label{sec:analysis}

In this section, we summarize the  numerical analysis results that are currently available for the residual-based ROM  stabilizations presented above. Specifically, we discuss the consistency, stability, and error bounds for these methods.
We emphasize that this is just a brief summary of the existing results and reflects only our own view on the topic.
Furthermore, we note that these definitions are not necessarily agreed upon.

\subsection{Consistency} \label{sec:consistency}
We start by considering consistency. For ROMs, two
types of consistency can be considered, and for concreteness we use
the following terminology:
\begin{itemize}
\item \textbf{Type 1: (Time-discrete) PDE consistency.} The ROM weak form holds when evaluated at the PDE solution, $\romStateTimeDiscreteSpaceDiscreteArg{n} \leftarrow \cStateStrongTimeDiscreteArg{n}$, assuming $\cStateStrongTimeDiscreteArg{n} \in \mathcal{H}^2(\Omega)$. We note that 
Type 1 consistency is only relevant for \continuousROMs, as discrete ROMs have no notion of the underlying PDE. We also note that Type 1 consistency is the consistency concept used for classical numerical methods (e.g., FEM).
 
\item \textbf{Type 2: FOM consistency:} The ROM weak form holds when evaluated at the FOM solution from which it is constructed. For continuous ROMs, this condition states that the weak form holds under the substitution $\romStateTimeDiscreteSpaceDiscreteArg{n}  \leftarrow \cStateTimeDiscreteSpaceDiscreteArg{n}$. Analogously for discrete ROMs, the ``discrete weak form" holds under the substitution $ \discreteROMGenStateTimeDiscreteSpaceDiscreteArg{n}  \leftarrow \genStateTimeDiscreteSpaceDiscreteArg{n}$. 

\end{itemize}

\begin{remark}[Model Consistency]
{\it Model consistency} is the setting in which the same stabilization method is used in the FOM and ROM.
We note that, when the same parameters are used in the FOM and ROM (i.e., we have parameter FOM-ROM consistency~\cite{strazzullo2022consistency}), model consistency is a special class of Type 2 consistency.
In~\cite{PACCIARINI20141} (see also~\cite{GIERE2015454}), the authors have argued both numerically and theoretically (in particular, see Section 3.3 and Proposition 3.1 in~\cite{PACCIARINI20141}) that using the same type of stabilization (i.e., SUPG) in the FOM and ROM yields more accurate ROM results.
More recently, model consistency for the evolve-filter-relax ROM~\cite{strazzullo2022consistency}  
(which is a spatial filtering-based stabilization, such as those described in Section~\ref{sec:other}) was shown to increase the ROM accuracy.

\end{remark}

\subsubsection{Continuous residual-based ROM stabilizations}


\noindent \textit{Continuous Galerkin ROMs}: 
The continuous Galerkin ROM is Type 1 and Type 2 consistent. Type 1 consistency follows directly from setting $\cBasis = \cBasisRom$ in the weak form~\eqref{eq:variational_timediscrete}, where we have leveraged $\cSpaceRom \subset \cSpace(\cDomain)$. Analogously, Type 2 consistency is shown from setting $\cBasis = \cBasisRom$ in the weak form~\eqref{eq:variational_discrete}.  
We note that, assuming a consistent FOM, Type 2 consistency automatically implies Type 1 consistency. This allows for \textit{a priori} convergence
analyses of continuous Galerkin ROMs with respect to the FOM solution as well as the solution to the governing continuous PDEs, as discussed in Refs. \cite{kalash_convergence,rozza2008reduced,Hesthaven2016}.
We emphasize  that, in order to maintain Type 2 consistency for a continuous Galerkin ROM,
it is important to employ the same spatial and temporal discretization method in building
the ROM as the one employed in building the FOM. \\

\noindent \textit{Continuous stabilized ROMs}:
Like continuous Galerkin ROMs, the \SUPG, \GLS, and \ADJ\ stabilized ROMs (developed through both the discretize-then-stabilize and space--time formulations)
are Type 1 and Type 2 consistent.  These stabilized formulations display Type 1 consistency as 
the stabilization term vanishes when evaluated about the PDE solution; the term vanishes 
as the residual evaluates to zero for the PDE solution  (assuming the solution is sufficiently regular). Type 2 consistency follows directly 
from setting $\cBasis = \cBasisRom$ in the stabilized weak form~\eqref{eq:variational_timediscrete_stable}, where we have leveraged $\cSpaceRom \subset \cSpaceFEM$.

\subsubsection{Discrete residual-based ROM stabilizations}
\noindent \textit{Discrete Galerkin ROMs}: The discrete Galerkin ROM
displays Type 2 consistency with the FEM model from which it is
constructed. We show this by first setting $
\genStateTimeDiscreteSpaceDiscreteArg{n}$, $n=1,\ldots,\nTimeSteps$
to be the FOM solution obtained from the Galerkin
FOM~\eqref{eq:g_fom_odeltae}. Making the substitution
$\discreteROMGenStateTimeDiscreteSpaceDiscreteArg{n}  \leftarrow
\genStateTimeDiscreteSpaceDiscreteArg{n}$, $n=0,\ldots,\nTimeSteps$,
it is straightforward to see that Eq.~\eqref{eq:drom_odeltae} is
satisfied under the conditions $\resid \leftarrow \residGalerkinFEM$
as $\residGalerkinFEM( \genStateTimeDiscreteSpaceDiscreteArg{n} ,
\genStateTimeDiscreteSpaceDiscreteArg{n-1}) = \mathbf{0}$,
$n=1,\ldots,\nTimeSteps$. Analogously, let
$\genStateTimeDiscreteSpaceDiscreteArg{n}$, $n=1,\ldots,\nTimeSteps$
be the FOM solution obtained from a stabilized
FOM~\eqref{eq:stab_fom_odeltae}. Making the substitution
$\discreteROMGenStateTimeDiscreteSpaceDiscreteArg{n}  \leftarrow
\genStateTimeDiscreteSpaceDiscreteArg{n}$, $n=0,\ldots,\nTimeSteps$,
it is again straightforward to see that~\eqref{eq:drom_odeltae}
is satisfied under the conditions $\resid \leftarrow
\residStabilizedFEM$ as $\residStabilizedFEM(
\genStateTimeDiscreteSpaceDiscreteArg{n} ,
\genStateTimeDiscreteSpaceDiscreteArg{n-1}) = \mathbf{0}$,
$n=1,\ldots,\nTimeSteps$.  
\vspace{0.2 in}\\
\noindent \textit{Discrete Stabilized ROMs}: Like the discrete Galerkin ROM, discrete stabilized ROMs display Type
2 consistency with the FEM model from which they are constructed. As both LSPG and APG 
can be written as a Petrov--Galerkin method, Type 2 consistency
follows from the same arguments as the discrete Galerkin ROM. \\


\begin{remark}
The residual-based stabilizations examined here all display Type 2 consistency\footnote{Formally, these methods are Type 2 consistent only if the same stabilization parameters and time steps are used in the FOM and ROM.}, and all continuous ROMs display Type 1 consistency. We emphasize that, while all methods considered here are consistent within the setting described above, this does not hold for all stabilized methods. Stabilization approaches based on, for example, eddy viscosity approaches typically \textit{do not} display Type 1 consistency. 
\end{remark}
%
\subsection{Stability}
Stability properties of the ROM depend on the type of projection, and are a driving factor in the ROM development. 
Here, we highlight the stabilization properties of the various ROMs considered. For concreteness, we restrict our discussion to stability within the context of the CDR equation. 

\subsubsection{Continuous ROMs}
For continuous ROMs, we define a stable formulation as one whose spatial bilinear form is \textit{strongly coercive}. For continuous ROMs, coercivity is defined by 
\begin{equation}\label{eq:stability}
\bilinearArg{\cStateTimeDiscreteDum}{\cStateTimeDiscreteDum} \ge \coerciveConstant \| \cStateTimeDiscreteDum \|_{\coerciveNorm}^2
\end{equation}
for some $\coerciveConstant \in \RPlus$. In the above, $\bilinear$ denotes a bilinear form, and $\| \cdot \|_{\coerciveNorm}^2$ denotes a norm associated with the formulation.  

\begin{remark}
If the spatial bilinear form is strongly coercive, then the fully discrete bilinear form associated with an implicit Euler discretization in time is also strongly coercive. 
\end{remark}

Coercivity of the spatial bilinear form guarantees boundedness of the solution, i.e., for finite $n$ and some $\beta > 0$,
\begin{equation}\label{eq:stability-definition}
\|\cStateTimeDiscreteArg{n} \|_{\coerciveNorm} \le \beta \|\forcing^n  \|_{\coerciveNorm},
\end{equation}
where $\forcing^n$ is the data at the $n$th time step. 

Before proceeding, we make the important point that a stable ROM does not necessarily imply an accurate ROM, and, often times, terminologies between instabilities and inaccuracies are mixed. As an example, for the CDR equation, the constants $C$ and $\beta$ in~\eqref{eq:stability} and~\eqref{eq:stability-definition} depend on the  diffusion parameter, $\epsilon$. It can be shown that the continuous Galerkin ROM has $\beta \to \infty$ as $\epsilon \to 0$.
Thus, for small $\epsilon$ values, the stability constant $\beta$ can be very large. As a result, although the standard Galerkin ROM may be formally stable (in the sense of~\eqref{eq:stability-definition}), it can be extremely inaccurate. We emphasize that this is not just a theoretical issue.
In practical ROM computations of convection-dominated systems (i.e., when $\epsilon$ is very small), the standard Galerkin ROM approximation --- while mathematically stable --- can indeed display large, spurious numerical oscillations (just as in the FEM setting~\cite{roos2008robust}). Although these oscillations are often referred to as instabilities, we emphasize here that they are large, but bounded (by the large stability constant $\beta$).

With this in mind, we now outline stability properties of the various formulations in the sense of the definition~\eqref{eq:stability}.
\begin{itemize}
\item \textbf{Continuous Galerkin ROM.} Coercivity of the continuous Galerkin FEM model has been demonstrated in numerous contexts (see, for example, Ref.~\cite{St05}), and it is straightforward to show that the continuous Galerkin FEM model is stable in the sense 
\begin{equation}\label{eq:Galerkin_coercivity}
 \bilinearGalerkinArg{\cStateTimeDiscreteDum}{\cStateTimeDiscreteDum} \ge \coerciveConstant \normGalerkinArg{ \cStateTimeDiscreteDum},
\end{equation}
where $\normGalerkinArg{\cStateTimeDiscreteDum}^2 = \reaction \| \cStateTimeDiscreteDum \|_{\LTwo}^2 + \viscosity \|  \nabla \cStateTimeDiscreteDum \|_{\LTwo}^2$. We emphasize that this coercivity property guarantees boundedness of the solution. For instance, under a suitable time-step restriction, the Galerkin method with a $\theta$-scheme for temporal discretization can be equipped with the stability bound 
\begin{equation*}
    \| \cStateTimeDiscreteSpaceDiscreteArg{n} \|_{L^2(\Omega)} \leq     \| u_{0,h} \|_{L^2(\Omega)} + C \sqrt{\frac{t^n}{\epsilon}} \max_{t \in [0,T]} \| f(t) \|_{L^2(\Omega)}
\end{equation*}
for a constant $C$ that depends only on $\Omega$.  See, e.g., Proposition 12.2.1 in \cite{Quarteroni_2008_PDE}. 
 In the limit of $\epsilon \rightarrow 0$, the stability statement~\eqref{eq:Galerkin_coercivity} loses control over the gradient. Hence, the Galerkin method is \textit{stable}, but not \textit{robust} in the limit of $\epsilon \rightarrow 0$. The stability of the continuous Galerkin ROM follows along the same lines as that for the continuous Galerkin FEM model~\cite[Theorem 5]{kunisch1}.



\item \textbf{Continuous SUPG ROM.} Coercivity of the SUPG FEM model has additionally been demonstrated in various contexts (see again Ref.~\cite[pg. 494]{St05}, 
or Refs.~\cite{JoNo11,GIERE2015454}). Coercivity of the SUPG FEM model depends on inverse estimates, and it is fairly straightforward to show that for some $\tau \in [0,\tau^*_{\text{SUPG}}]$, where $\tau^*_{\text{SUPG}}$ is a grid and parameter dependent upper threshold on $\tau$, the continuous SUPG FEMmodel is stable in the sense 
\begin{equation}\label{eq:SUPG_coercivity}
 \bilinearStabilizedArg{\cStateTimeDiscreteDum}{\cStateTimeDiscreteDum} \ge \coerciveConstant \normSupgArg{ \cStateTimeDiscreteDum},
\end{equation}
where 
$\normSupgArg{\cStateTimeDiscreteDum}^2 = \reaction \| \cStateTimeDiscreteDum \|_{\LTwo}^2 +  \viscosity\| \nabla \cStateTimeDiscreteDum \|_{\LTwo}^2 + \tau \sum_{k=1}^{\nel} \| \wavespeed \cdot \nabla \cStateTimeDiscreteDum \|_{\LTwoK}^2$.  We note that, in the limit that $\epsilon \rightarrow 0$, the stability statement~\eqref{eq:SUPG_coercivity} maintains control over the gradient of the state in the streamline direction. Thus, the SUPG method is \textit{stable} and \textit{robust} in the limit of $\epsilon \rightarrow 0$. As a result, we do not expect the accuracy of the method to deteriorate for small $\viscosity$. We additionally note that $\normSupgArg{\cStateTimeDiscreteDum} \ge \normGalerkinArg{\cStateTimeDiscreteDum}$, so that SUPG is more dissipative than Galerkin. The stability of the continuous SUPG ROM follows along the same lines as that for the continuous SUPG FEM model (see, e.g.,~\cite[Lemma 3.3]{GIERE2015454}). 

\item \textbf{Continuous \SDGLS\ ROM.} Stability of GLS is given in Ref.~\cite{HuFrHu89} in the steady case. Coercivity is straightforward to demonstrate as GLS adds a symmetric non-negative term to the bilinear form. As \SDGLS\ is equivalent to the steady case but with a modified source term, the analysis in Ref.~\cite{HuFrHu89} is directly applicable and results in the stability statement 
\begin{equation}\label{eq:GLS_coercivity}
 \bilinearStabilizedArg{\cStateTimeDiscreteDum}{\cStateTimeDiscreteDum} \ge \coerciveConstant \normGlsArg{ \cStateTimeDiscreteDum},
\end{equation}
where 
$\normGlsArg{\cStateTimeDiscreteDum}^2 = \reaction \|  \cStateTimeDiscreteDum \|_{\LTwo}^2 +  \viscosity \| \nabla \cStateTimeDiscreteDum \|_{\LTwo}^2 + \tau \sum_{k=1}^{\nel}\| \wavespeed \cdot \nabla \cStateTimeDiscreteDum + \left( \reaction + \frac{1}{\Delta t} \right) \cStateTimeDiscreteDum  - \epsilon \Delta \cStateTimeDiscreteDum \|_{\LTwoK}^2$.  We note that \GLS\ is stable for non-negative values of $\tau$. Like SUPG, the stability statement~\eqref{eq:GLS_coercivity} maintains control over the gradient of the state in the streamwise direction in the limit $\epsilon \rightarrow 0$. Thus, GLS is \textit{stable} and \textit{robust} in the limit of $\epsilon \rightarrow 0$. We additionally note that $\normGlsArg{\cStateTimeDiscreteDum}^2 $ depends on the time step $\Delta t$. For very small time steps, the stability statement~\eqref{eq:GLS_coercivity} will be dominated by the $\frac{1}{\Delta t}$ term and we are thus not robust in this limit. The stability of the continuous \SDGLS\ ROM follows along the same lines as that for the continuous \SDGLS\ FEM model.

\item \textbf{Continuous \SDADJ\ ROM.} Coercivity of the ADJ FEM model has been demonstrated for the steady convection diffusion reaction equation~\cite{FrVa20}. As \SDADJ\ is equivalent to the steady case but with a modified forcing term, the stability statement presented in Ref.~\cite{FrVa20} applies. The stability statement is given as: for $0 \le \tau \le \tau^*_{\text{ADJ}}$, 
\begin{equation}\label{eq:GLS_coercivity}
 \bilinearStabilizedArg{\cStateTimeDiscreteDum}{\cStateTimeDiscreteDum} \ge \coerciveConstant \normAdjArg{ \cStateTimeDiscreteDum},
\end{equation}
where $\normAdjArg{\cStateTimeDiscreteDum}^2 = \sum_{k=1}^{\nel} \left( \left( \reaction + \frac{1}{\Delta t} \right) \alpha_k \|  \cStateTimeDiscreteDum \|_{\LTwoK}^2 +  \viscosity \alpha_k \| \nabla \cStateTimeDiscreteDum \|_{\LTwoK}^2 + \tau  \| \wavespeed \cdot \nabla \cStateTimeDiscreteDum  \|_{\LTwoK}^2 \right)$ with $\alpha_k$ being a constant that depends on the mesh, parameters, and inverse estimates. We again observe more robust behavior in the limit of $\viscosity \rightarrow 0$ as well as a dependence on the time step $\Delta t$. We again expect poor behavior in the limit of $\Delta t \rightarrow 0$ as coercivity is dominated by the $\frac{1}{\Delta t}$ term. The stability of the continuous \SDADJ\ ROM follows along the same lines as that for the continuous \SDADJ\ FEM model.

\item \textbf{Continuous \STGLS\ ROM.} Coercivity of the \STGLS\ FEM  model was demonstrated in one of the original references on GLS~\cite{HuFrHu89} by virtue of the formulation adding a symmetric term. We note that, here, we include the $\left( \cStateTimeDiscreteArg{n} - \cStateTimeDiscreteArg{n-1} \right) / \Delta t$ term in the definition of our residual to retain consistency for the $p=0$ DG trial space, and as a result the analysis in~\cite{HuFrHu89} does not directly extend to the current case.\footnote{We found that this term makes little difference in practice.} We also note that the space--time formulation was advocated in the original reference~\cite{HuFrHu89}.

\item \textbf{Continuous \STADJ\ ROM.} Coercivity of the \STADJ\ FEM model applied to the unsteady convection-diffusion-reaction equation has not been demonstrated to the best of 
our knowledge. 

%
\end{itemize}

\subsubsection{Discrete ROMs}
We now consider stability of the various discrete ROMs discussed above. For the following analysis, we introduce the following notation for a generic discrete ROM as
$$ \mathbf{E} \frac{ \genDiscreteRomStateTimeDiscreteSpaceDiscreteArg{n}  - \genDiscreteRomStateTimeDiscreteSpaceDiscreteArg{n-1} }{\Delta t}  
+ \mathbf{G} \genDiscreteRomStateTimeDiscreteSpaceDiscreteArg{n} = \mathbf{f},$$
where $\mathbf{E}\in \RR{\romdim \times \romdim}$ is a ``mass matrix" (e.g., $\mathbf{E} = \dBasisROMMat^T \mathbf{M} \dBasisROMMat$ for the discrete Galerkin ROM of the Galerkin FEM), $\mathbf{G} \in \RR{\romdim \times \romdim}$ is a ``dynamics" matrix, $\mathbf{f} \in \RR{\romdim}$ is a forcing vector, and $\hat{\mathbf{x}}^n \in \RR{\romdim}$ are the reduced coordinates. 
We define a stable discrete ROM as one whose ``mass" matrix $\mathbf{E}$ is symmetric positive definite and whose dynamics matrix $\dynamicsMat$ is positive definite, i.e.,
\begin{equation}\label{eq:discrete_stability}
\mathbf{v}^T \dynamicsMat \mathbf{v} > 0, \;  \; \forall \mathbf{v} 
\in   \RR{\romdim} \setminus \{ {\bf 0} \}.
\end{equation} 
We emphasize that this is analogous to coercivity in finite dimensional spaces since all norms are equivalent.

To proceed, it is first helpful to note that the continuous Galerkin FEM model~\eqref{eq:g_fom_odeltae} results in the system 
$$ \mass \frac{\genStateTimeDiscreteSpaceDiscreteArg{n}  - \genStateTimeDiscreteSpaceDiscreteArg{n-1} }{\Delta t}  + \left(\mathbf{A} + \viscosity \mathbf{D} + \reaction \mathbf{M} \right) \genStateTimeDiscreteSpaceDiscreteArg{n} = \mathbf{f},$$ 
where $\mathbf{A} = m\left( \cBasis_i , \mathbf{b} \cdot \nabla \cBasis_j\right)$ is the convection matrix, $\mathbf{D} = m\left( \nabla \cBasis_i, \nabla \cBasis_j \right)$ the symmetric positive definite diffusion matrix, and $\mathbf{M} = m \left( \cBasis_i,\cBasis_j\right)$ the symmetric positive definite mass matrix. We note that $\mathbf{v}^T \mathbf{A} \mathbf{v} = 0$, $\mathbf{v}^T \mathbf{D} \mathbf{v} > 0$, and $\mathbf{v}^T \mathbf{M} \mathbf{v} > 0$ $\forall \mathbf{v} 
\in \RR{\fomdim} \setminus \{ {\bf 0} \}$. For notational simplicity, we define $\matFull = \mathbf{A} + \epsilon \mathbf{D} + \sigma \mass$. 
\begin{itemize}
\item \textbf{Discrete Galerkin ROM with $\mathbf{W}=\mathbf{I}$}. As the discrete Galerkin ROM of the continuous Galerkin FEM model with $\mathbf{W}=\mathbf{I}$ is equivalent to the continuous Galerkin ROM, stability is implied. It is further straightforward to show that $\mathbf{v}^T \mathbf{B}\mathbf{v} > 0$ $\forall \mathbf{v} \in \RR{\fomdim} \setminus \{ {\bf 0} \}$ at the discrete level directly. The Galerkin discrete ROM constructed from the Galerkin continuous ROM results in the dynamics matrix
$$\dynamicsMatGalerkin = \dBasisROMMat^T \left(\mathbf{A} + \viscosity \mathbf{D} + \reaction \mathbf{M} \right) \dBasisROMMat.$$ 
It is straightforward to see
$$\mathbf{v}^T \dynamicsMatGalerkin \mathbf{v} =  \viscosity \| \mathbf{v} \|_{\diffusionMatRom}^2 + \reaction \| \mathbf{v} \|_{\massRom}^2 > 0,$$
where $\diffusionMatRom = \dBasisROMMat^T \mathbf{D} \dBasisROMMat$.
We note the above is simply the discrete equivalent to the inequality~\eqref{eq:Galerkin_coercivity}. We additionally note that $ \viscosity \| \mathbf{v} \|_{\diffusionMatRom}^2 + \reaction \| \mathbf{v} \|_{\massRom}^2 $ is the discrete statement of a weighted $\HOne$ norm 
that approaches a $\sigma$-weighted discrete $\LTwo$ norm as $\viscosity \rightarrow 0$.

It is less straightforward to demonstrate that $\dynamicsMat$ is positive definite for the discrete Galerkin ROM constructed from a stabilized FEM model. As the discrete Galerkin ROM dynamics recover the continuous Galerkin ROM dynamics under the condition $\mathbf{W} \leftarrow \mathbf{I}$, however, these ROMs can be expected to obey the same stability properties as their continuous counterparts.

\item \textbf{Discrete LSPG ROMs}. To the best of our knowledge, no result exists in the literature demonstrating stability of a discrete LSPG ROM constructed from a Galerkin FEM FOM in the sense of~\eqref{eq:discrete_stability}. Stability analyses for LSPG have been carried out in other contexts; for instance, Ref.~\cite{HuWeDu22} demonstrates that, for LTI systems, LSPG with orthonormal bases results in an asymptotically stable ROM if the underlying FOM is asymptotically stable.\footnote{No mass matrix was considered.}

In the present context, one can show that a discrete LSPG ROM of the continuous Galerkin FEM model constructed in the inner product $\discreteInnerProductTypeTwo \leftarrow \mass^{-1}$ results in a mass matrix $\mathbf{E} = \dBasisROMMat^T \mass \dBasisRomMat$ and a dynamics matrix
$$\dynamicsMatLspg =  \dBasisROMMat^T \left[ \matFull + \matFull^T \right] \dBasisROMMat + \Delta t \dBasisROMMat^T \matFull^T \mass^{-1} \matFull  \dBasisROMMat.$$
As the inequality $\mathbf{v}^T \matFull \mathbf{v} > 0$ $\forall \mathbf{v} \in \RR{\fomdim} \setminus \{ {\bf 0} \}$, stability is implied.

From our analysis of the LSPG ROM applied to the Galerkin FEM model, it is straightforward to see that if the FOM has a dynamics matrix that is positive definite, then the resulting LSPG ROM constructed in the inner product $\discreteInnerProductTypeTwo = \mathbf{I}$ will additionally have a positive definite dynamics matrix. Thus, LSPG with $\discreteInnerProductTypeTwo = \mathbf{I}$ will be stable if deployed on a stabilized FEM model that is also stable.
 \item \textbf{APG ROMs}. The APG ROM is derived from a formulation of the Mori--Zwanzig formalism, and a stability analysis has been undertaken for a model displaying a structural equivalence to APG in Ref.~\cite{HaSt07}. This analysis demonstrates that the so-called $t$-model (which is equivalent to APG for $\tau$ set to $t$) will be be dissipative when applied to a system that is energy conserving; this result directly implies a positive definite dynamics matrix. However, no result exists in the literature demonstrating stability of the APG ROM for systems that dissipate energy. 

\end{itemize}

\subsection{Error bounds}
We now summarize existing numerical analyses that attempt to bound the ROM error for the CDR system. 
We note that \textit{a priori} error bounds for POD-based methods are typically limited to reproductive cases (unless some assumption is made on the solution manifold), while \textit{a posteriori} error bounds are typically valid in both the reproductive and predictive regime.

\subsubsection{Continuous ROMs}
\begin{itemize}
\item  \textbf{Continuous Galerkin ROMs:} 
As a continuous Galerkin ROM arises from the Galerkin approximation of the CDR equation in a POD (or RBM) subspace, its \textit{a priori} error bound can be derived by leveraging the FEM error analysis for parabolic PDEs~\cite{thomee2006galerkin} together with the approximability properties of the POD~\cite{kunisch1} or RBM~\cite{Hesthaven2016,quarteroni2015reduced} space. For example, error bounds for the POD-Galerkin ROM constructed using continuous projection were derived for parabolic linear systems 
and certain nonlinear systems by Kunisch and Volkwein in Ref.~\cite{kunisch1}. In a follow-up paper, the authors 
derived error bounds for equations pertaining to fluid dynamics~\cite{kunisch2}, e.g., the two-dimensional incompressible Navier-Stokes equations. New error bounds were proved by Singler~\cite{singler}, who derived exact expressions for the POD data approximation errors considering four different POD projections and two different Hilbert space error norms.
Error bounds for the RBM-Galerkin ROM constructed using continuous projection were derived for parabolic problems in~\cite{Grepl_2005_RB_Parabolic, Haasdonk_2008_RB_FV_Evolution,Grepl_2007_RB_EIM}, the convergence of POD-Greedy algorithm was analyzed in \cite{Haasdonk_2013_POD_Greedy_Convergence}, and sharper error bounds using space-time formulations were obtained in~\cite{Urban_2014_Space_Time_RB,Yano_2014_Space_Time_RB_Boussinesq}.
We emphasize that these Galerkin ROM 
error bounds for parabolic PDEs grow with the inverse of the coercivity constant, which scales with $\epsilon$ for the CDR system. Hence, 
there is no guarantee that the continuous Galerkin ROMs will provide an approximation comparable to the best-fit approximation in the limit of $\epsilon \to 0$. 

It is worth noting that error bounds and convergence analyses exist for ROMs built using continuous Galerkin projection for PDEs other than the CDR equation, e.g., hyperbolic equations. 
In \cite{kalash_convergence}, for example, Kalashnikova and
Barone derived \textit{a priori} error estimates for an
energy-stability-preserving ROM formulation developed in
\cite{BARONE20091932} for linearized compressible flow.  These error
bounds 
were derived by adapting techniques traditionally used in the
numerical analysis of spectral approximations to PDEs \cite{funaro}
and employed a carefully constructed stable penalty-like
implementation of the relevant boundary conditions in the ROM.  

\item \textbf{Continuous SUPG ROMs:} 
Error analysis of the SUPG ROM was undertaken recently in Ref.~\cite{JoMoNo22}, where it was demonstrated that the SUPG ROM could be equipped with robust error estimates that do not deteriorate as 
$\viscosity \rightarrow 0$. These estimates bound the error between the SUPG ROM solution and a corresponding SUPG FEM solution. To obtain $\Delta t$ independent error bounds, the authors employed POD snapshots that included the time-difference quotients~\cite{kunisch1,koc2021optimal}; without these coefficients the resulting error bounds depend on $\Delta t$. In numerical experiments, however, it was observed that including the time-difference quotients did not lead to improved results, and thus the authors believe an important open question is the derivation of $\Delta t$-independent bounds for the case where the time-difference quotients are not included in the bases.  
Lastly, we note that Ref.~\cite{JoMoNo22} supports previous analysis of the SUPG ROM in~\cite{GIERE2015454}.  

\item \textbf{Other stabilized ROMs:} To the best of our knowledge, no error analysis exists for the other stabilized ROMs considered here. We do note that error analysis exists for the corresponding FEM formulations. 
We also note that error analysis exists for stabilized ROMs that are not residual-based, such as those outlined in Section~\ref{sec:other}. To our knowledge, the first numerical analysis of ROM closures was performed in~\cite{borggaard2011artificial}, where error bounds for the time discretization of the Smagorinsky eddy viscosity POD-Galerkin ROM were proven. 
Error bounds for the space and time discretizations of the Smagorinsky RBM-Galerkin ROM were later proven in~\cite{rebollo2017certified}.
The first numerical analysis of variational multiscale eddy viscosity ROM closures was performed in~\cite{iliescu2014variational} for the incompressible Navier-Stokes equations and in~\cite{iliescu2013variational} for the CDR equation, where  stability and convergence were proven.
Related work was performed in~\cite{roop2013proper,eroglu2017modular,rubino2020numerical,azaiez2021cure}. 


To our knowledege, the only numerical analysis for spatial filtering-based ROM stabilizations was performed in~\cite{xie2018numerical}, where error bounds for the Leray regularzied ROM were proved.

\bigskip


\end{itemize}

\subsubsection{Discrete ROMs}
\begin{itemize}
\item \textbf{Discrete Galerkin ROM}: Due to the equivalence between a discrete Galerkin ROM and its continuous counterpart, the error bounds derived for the
continuous Galerkin ROM are applicable to the discrete Galerkin ROM constructed on top of their corresponding continuous FEM system. Various authors have
derived error bounds for the discrete Galerkin ROM in a more generic context. In
Ref.~\cite{doi:10.1137/S0036142901389049}, error bounds are derived for the
discrete Galerkin ROM within the context of a linear and nonlinear dynamical
system $\dot{\boldsymbol x} = \boldsymbol f(\boldsymbol x)$.
Analogously, Ref.~\cite{carlberg_lspg_v_galerkin} derives error
bounds for the discrete Galerkin ROM for O$\Delta $Es arriving from
linear multistep and Runge--Kutta time discretizations of
$\dot{\boldsymbol x} = \boldsymbol f(\boldsymbol x)$. In the general nonlinear case, these error bounds 
depend on difficult-to-compute Lipshitz constants, grow exponentially in time, and lack sharpness. No error analysis of the discrete Galerkin ROM 
specialized to the CDR system exists to the best of our knowledge. 

\item \textbf{LSPG ROMs}: \textit{A priori} and \textit{a posteriori} upper error bounds for LSPG ROMs applied to generic time-discrete nonlinear dynamical systems are derived in Ref.~\cite{carlberg_lspg_v_galerkin}. These error bounds rely on the assumption of Lipshitz continuity of the nonlinear right-hand side velocity operator, which can be related to coercivity in the linear setting. The bounds demonstrate that the upper error bound of the LSPG ROM grows exponentially with the number of time steps, and that this upper error bound can be bounded by the maximum residual over a given time step. Further,~\cite{carlberg_lspg_v_galerkin} shows that LSPG can be equipped with an \textit{a posteriori} upper error bound lower than the Galerkin ROM. However, bounds presented in Ref.~\cite{carlberg_lspg_v_galerkin} are derived for the case where there is no mass matrix, and sharpness of the bounds is not addressed (thus, LSPG having a lower ``upper bound" than Galerkin is not a robust statement of accuracy). Again, no error analysis of LSPG ROMs specialized to the CDR system exists to the best of our knowledge.
 \\

\item \textbf{APG ROMs}: Ref.~\cite{parish_apg} derives
\textit{a priori} error bounds for the APG ROM for nonlinear
dynamical systems and linear time-invariant dynamical systems.
In the linear case, it is shown that, for sufficiently small $\tau$,
the upper bound on the error in the APG ROM is lower than in the
Galerkin ROM. Similar to LSPG, however, the bounds are presented for the case where there is no mass matrix and sharpness is additionally not addressed (again, APG having a lower ``upper bound" than Galerkin is not a robust statement of accuracy).
Once again, no error analysis of APG ROMs
specialized to the CDR system exists to the best of our knowledge.
\end{itemize}

\subsection{Selection of the stabilization parameter, $\tau$}
    \label{sec:tau-analysis}
   
Lastly, we comment on analyses for selecting the ROM stabilization parameter. In standard FEM, numerical analysis arguments are generally used to determine the scaling of the stabilization parameter, $\tau$, in residual-based stabilizations  (see, e.g., the survey in~\cite{roos2008robust}).
The general approach used to determine the $\tau$ scaling is to (i) prove error bounds for the stabilized method, and (ii) choose a $\tau$ scaling with respect to the discretization parameters (e.g., the mesh size $h$ and the time step $\Delta t$) that ensures an optimal error bound.  
    
For residual-based ROM stabilizations, one heuristic approach for choosing the stabilization parameter, $\tau$, is to use the same value as that used in the standard FEM (see, e.g., equation (13) in~\cite{PACCIARINI20141}, equation (20) in~\cite{kragel2005streamline}, and equation (10) in~\cite{mclaughlin2016stabilized}).
We note that this approach is purely heuristic and 
does not use the numerical analysis arguments generally employed for standard FEM residual-based stabilized methods~\cite{roos2008robust}.

A fundamentally different approach, which utilizes numerical analysis arguments to determine the $\tau$ scaling, was proposed in~\cite{GIERE2015454} for the SUPG-ROM.
As explained in Section 3.5.1 in~\cite{GIERE2015454}, since the ROM space is a subspace of the FEM space, two types of inverse estimates can be used to prove optimal error estimates for the SUPG-ROM: (i) a FEM inverse estimate, which yields the standard FEM 
scaling in which $\tau$ depends on the FEM mesh size, and (ii) a ROM inverse estimate~\cite{kunisch1}, which yields a new ROM 
scaling in which $\tau$ depends on the ROM parameters (e.g., the ROM dimension, the POD basis functions, and the corresponding eigenvalues).
The preliminary numerical investigation in~\cite{GIERE2015454} suggests that the FEM $\tau$ scaling yields more accurate results for large $R$ values, but the ROM $\tau$ scaling is competitive for low $R$ values.
Further theoretical and numerical investigation is needed in order to determine optimal $\tau$ scalings for residual-based stabilized ROMs.



\section{Numerical experiments}\label{sec:experiments}
\subsection{Overview}
We now present several studies to numerically assess the various ROM
formulations for the CDR equation \eqref{eq:cdr}. We first provide
specifics on the setup of the numerical experiments.

\subsubsection{Investigated ROMs and implementation details}
Table~\ref{tab:crom_summary} details the \continuousROMs\ we
investigate, while Tables~\ref{tab:lspg_summary} and~\ref{tab:apg_summary} detail the
LSPG and APG ROMs investigated. For continuous ROMs, we investigate the Galerkin, SUPG, ADJ, and GLS ROMs, as detailed in Section~\ref{sec:resid_stab_rom}. For ADJ and GLS, we investigate formulations developed both through the ``discretize-then-stabilize" approach (DS) and the ``space--time" approach (ST), as discussed in Section~\ref{sec:resid_stab}. For discrete ROMs, we examine (1) LSPG ROMs based on the Galerkin FEM, SUPG FEM, \ADJ\  FEM, \STADJ\ FEM, \GLS\ FEM, and \STGLS\ FEM,  and (2) APG ROMs based on the Galerkin FEM, SUPG FEM, \STGLS\ FEM, and \STADJ\ FEM. In what follows, we will abbreviate these discrete ROM formulations as ``FEM type--discrete ROM type", e.g., ``SUPG--LSPG" denotes an LSPG ROM of the SUPG FEM.  
All in all, we consider 16 ROM formulations. All ROMs employ the implicit Euler method for
temporal discretization. The numerical experiments are carried out
in the {\tt FEniCS} package~\cite{AlnaesBlechta2015a,LoggMardalEtAl2012a,LoggWells2010a,LoggWellsEtAl2012a}.
We note that stabilization will be carried out both at the FEM level and the ROM level; this will be detailed in subsequent sections. Lastly, we additionally note that all experiments will focus only on reproductive ROMs.

\begin{table}[H]
\begin{centering}
\begin{tabular}{l c c c c c c c c}
\hline
& \Galerkin & \SUPG &  \GLS & \ADJ & \STGLS & \STADJ\\
\hline
ROM O$\Delta$E & Eq.~\eqref{eq:g_rom_odeltae} & Eq.~\eqref{eq:stab_rom_odeltae} & Eq.~\eqref{eq:stab_rom_odeltae} & Eq.~\eqref{eq:stab_rom_odeltae}  & Eq.~\eqref{eq:stab_rom_odeltae} & Eq.~\eqref{eq:stab_rom_odeltae}\\
Conditions & N/A &  $\stabilizationOpTest = \stabilizationOpTestSUPG$ & $\stabilizationOpTest= \stabilizationOpTestGLS$ & $\stabilizationOpTest = \stabilizationOpTestADJ$  & $\stabilizationOpTest= \stabilizationOpTestSTGLS$ & $\stabilizationOpTest = \stabilizationOpTestSTADJ$ & \\ \hline
\end{tabular}
\caption{Summary of \continuousROMs\ investigated 
}
\label{tab:crom_summary}
\end{centering}
\end{table}

\begin{table}[H]
\begin{centering}
\begin{tabular}{l c c c c c c }
\hline
&  \LSPGG &  \LSPGSUPG  & \LSPGGLS & \LSPGADJ & \LSPGSTGLS & \LSPGSTADJ \\
\hline
ROM O$\Delta$E & Eq.~\eqref{eq:LSPG_minproblem} & Eq.~\eqref{eq:LSPG_minproblem} & Eq.~\eqref{eq:LSPG_minproblem} & Eq.~\eqref{eq:LSPG_minproblem} & Eq.~\eqref{eq:LSPG_minproblem}  & Eq.~\eqref{eq:LSPG_minproblem} \\
FOM O$\Delta$E & Eq.~\eqref{eq:g_fom_odeltae} & Eq.~\eqref{eq:stab_fom_odeltae} & Eq.~\eqref{eq:stab_fom_odeltae} & Eq.~\eqref{eq:stab_fom_odeltae} & Eq.~\eqref{eq:stab_fom_odeltae} & Eq.~\eqref{eq:stab_fom_odeltae} \\
Conditions &  
N/A & 
 $\stabilizationOpTest = \stabilizationOpTestSUPG$ &
 $\stabilizationOpTest = \stabilizationOpTestGLS$ &
 $\stabilizationOpTest = \stabilizationOpTestADJ$ &
 $\stabilizationOpTest = \stabilizationOpTestSTGLS$ &
 $\stabilizationOpTest = \stabilizationOpTestSTADJ$ \\
Inner product &  $\discreteInnerProductTypeTwo = \mass^{-1}$ & $ \discreteInnerProductTypeTwo = \mass^{-1}$ &
$\discreteInnerProductTypeTwo = \mass^{-1}$ &
$\discreteInnerProductTypeTwo = \mass^{-1}$ &
$\discreteInnerProductTypeTwo = \mass^{-1}$ &
$\discreteInnerProductTypeTwo = \mass^{-1} $
 \\ \hline
\end{tabular}
\caption{Summary of LSPG ROMs investigated}
\label{tab:lspg_summary}
\end{centering}
\end{table}

\begin{table}[H]
\begin{centering}
\begin{tabular}{l c c c c  }
\hline
&  \APGG &  \APGSUPG  &  \APGSTGLS & \APGSTADJ \\
\hline
ROM O$\Delta$E & Eq.~\eqref{eq:APG_problem} &  Eq.~\eqref{eq:APG_problem_stabilized} & Eq.~\eqref{eq:APG_problem_stabilized}  & Eq.~\eqref{eq:APG_problem_stabilized} \\
FOM O$\Delta$E & Eq.~\eqref{eq:g_fom_odeltae} &  Eq.~\eqref{eq:stab_fom_odeltae} & Eq.~\eqref{eq:stab_fom_odeltae} & Eq.~\eqref{eq:stab_fom_odeltae} \\
Conditions &  
N/A & 
 $\stabilizationOpTest = \stabilizationOpTestSUPG$ &
 $\stabilizationOpTest = \stabilizationOpTestSTGLS$ &
 $\stabilizationOpTest = \stabilizationOpTestSTADJ$ \\
Inner product &  $\discreteInnerProductTypeTwo = \mathbf{I}$ &
$ \discreteInnerProductTypeTwo = \mathbf{I}$ &
$\discreteInnerProductTypeTwo = \mathbf{I}$ &
$\discreteInnerProductTypeTwo = \mathbf{I} $
 \\ \hline
\end{tabular}
\caption{Summary of APG ROMs investigated}
\label{tab:apg_summary}
\end{centering}
\end{table}

\subsubsection{Metrics}\label{sec:metrics}
We use as metrics the (discrete) time-integrated relative $\LTwo$  error and the (discrete) time-integrated relative $\HOne$ error between the ROM solution and best-fit solution (i.e., the error between the ROM solution and the FOM solution projected onto the trial space). The time-integrated $\LTwo$ relative error is defined as
\begin{equation}\label{eq:til2error}
\errorLTwoBestFit =  \frac{ \sum_{n=1}^{\nTimeSteps}  \| \romStateTimeDiscreteSpaceDiscreteArg{n} - \projectorLTwoArg{\cSpaceRom} \cStateStrongTimeDiscreteArg{n} \|_{\LTwo}^2 }{
  \sum_{n=1}^{\nTimeSteps} \| \projectorLTwoArg{\cSpaceRom} \cStateStrongTimeDiscreteArg{n} \|_{\LTwo}^2 },
\end{equation}
where $\projectorLTwoArg{\cSpaceRom}$ is the orthogonal $\LTwo$ projector onto $\cSpaceRom$. Analogously, the relative $\HOne$ error is defined as
\begin{equation}\label{eq:tih1error}
\errorHOneBestFit =  \frac{ \sum_{n=1}^{\nTimeSteps} | \romStateTimeDiscreteSpaceDiscreteArg{n} - \projectorHOneArg{\cSpaceRom} \cStateStrongTimeDiscreteArg{n} |_{\HOne}^2}{  \sum_{n=1}^{\nTimeSteps} |\projectorHOneArg{\cSpaceRom} \cStateStrongTimeDiscreteArg{n} |_{\HOne}^2 },
\end{equation}
where $\| \cdot \|_{\HOne}$ is the $\HOne$ semi-norm and $\projectorHOneArg{\cSpaceRom}$ is the orthogonal projector onto $\cSpaceRom$ in the $\HOne$ semi-norm. We note that in our studies, we execute the ROMs for varying time step sizes\footnote{ This is done because the time step size impacts the stabilization.}, and as a result the ROM is executed on a different time grid than the high resolution solution. We design our studies such that the ratio of the ROM time step to the high resolution time step is always a positive integer. The summations in Eqs.~\eqref{eq:til2error} and~\eqref{eq:tih1error} are then performed on the coarser ROM time grid.
%

\subsubsection{Construction of ROM trial space}\label{sec:basis_construction}

In this work, we are interested in examining cases where the
underlying FEM requires stabilization (this is the relevant case for
real-world applications), and as such examine the scenario where the FEM trial space is not fully resolved such that our FEM requires stabilization to be accurate. The natural approach to generate the ROM trial space in this setting 
is to (i) solve a FOM, which comprises a
stabilized FEM model, and (ii) leverage the solution data to
construct the ROM trial space. The ROM is then executed with the
same stabilized form used to generate the snapshots; this is the
so-called ``offline--online stabilization strategy" 
outlined in
Ref.~\cite{rozza_stabilized} (see also~\cite{strazzullo2022consistency}), which comprises a form of model consistency.
This procedure, unfortunately, makes
it difficult to compare the performance of different stabilization
techniques; i.e., if we generated the FOM solution with SUPG, then
comparing SUPG ROM solutions to LSPG ROM solutions becomes unfair.
As a result, we are not using this approach in our numerical investigation.

To circumvent this issue, we generate the reduced-order trial spaces by projecting ``truth data" onto a FOM trial space. 
As the problems we analyze in this section do not have analytical solutions, we generate the truth data via a high fidelity Galerkin FEM model that uses a high resolution trial space $\cSpaceHRES$. We generate bases by projecting this truth data onto the (lower-dimensional) FOM trial space $\cSpaceFOM$ and, e.g., performing an SVD. We emphasize that the high-fidelity Galerkin FEM model and the high resolution trial space are used only for generating reference solutions and snapshots. At a high level, one may think of this approach as the setting where the truth data comes from experimental measurements. We emphasize that a consequence of this approach is that none of our ROMs will be consistent with the truth data,  meaning that adding ROM basis vectors will not necessarily result in a more accurate solution. Algorithm~\ref{alg:gen_rom_trialspace} outlines the algorithm used to generate the ROM trial spaces. 

\begin{algorithm}
\caption{Algorithm for generating high fidelity solutions and the ROM trial space}
\label{alg:gen_rom_trialspace}
\textbf{Input:} Energy cutoff criterion, $\energyCutoff$
 \\
\textbf{Output:} ROM bases, $\dBasisROMMat$\\
\textbf{Steps:}

\begin{enumerate}
\item Solve the Galerkin O$\Delta$E~\eqref{eq:g_fom_odeltae} with $\cSpaceFEM \leftarrow \cSpaceHRES$ for $n=1,\ldots,\nTimeSteps$ to generate solutions $[\cStateTimeDiscreteSpaceDiscreteArg{n}]_{\mathrm{h-res}}$, $n=1,\ldots,\nTimeSteps$. 

\item Perform the restriction $\cStateFomTimeDiscreteSpaceDiscreteArg{n} = \projectorLTwoArg{\cSpaceFOM}\cStateHresTimeDiscreteSpaceDiscreteArg{n}$, $n=1,\ldots,\nTimeSteps$.

\item Collect the \femCoefficients\ $\genStateFomTimeDiscreteSpaceDiscreteArg{n}$ associated with $[\cStateTimeDiscreteSpaceDiscreteArg{n}]_{\text{FOM}}$, $n=0,\ldots,\nTimeSteps$ into the snapshot matrix
$$\snapshotMatrixDiscrete  = \begin{bmatrix}
\genStateFomTimeDiscreteSpaceDiscreteArg{1}& \cdots & \genStateFomTimeDiscreteSpaceDiscreteArg{\nTimeSteps}
\end{bmatrix}
.$$
\item Execute Algorithm~\ref{alg:pod} with inputs $\snapshotMatrixDiscrete$, $\energyCutoff$, $\discreteInnerProductType = \mass$ to obtain the ROM basis $\dBasisROMMat$ and trial space $\cSpaceRom = \text{Range}(\dBasisROMMat)$.
\end{enumerate}

\end{algorithm}

\subsubsection{Selection of stabilization parameters and time step}\label{sec:timestep_information}
In addition to depending on the choice of inner product, all
stabilized methods considered depend on the stabilization parameter
$\tau$, the time step $\Delta t$, or both. As discussed earlier, the \textit{a priori}
selection of $\tau$ and $\Delta t$ is an area that is receiving
attention in both the FEM and ROM
communities~\cite{HSU2010828,GIERE2015454}, and is still an
outstanding issue for ROMs in particular. Here we perform a
grid sweep to explore this sensitivity and to select ``optimal" time step sizes and stabilization
parameters. The grid sweep is obtained by executing ROM solves for
$(\tau,\Delta t) \in \tauSet \times \dtSet$, where $\tauSet = \dtSet =
\{10^{-4}, \; 2.5 \times 10^{-4} , \; 5 \times 10^{-4}, \; 10^{-3},
\; 2 \times 10^{-3}, \; 3 \times 10^{-3}, \; 4 \times 10^{-3}, \;
5\times 10^{-3},  \; 6\times10^{-3},  \; 7\times10^{-3}, \; 8 \times
10^{-3}, \; 9 \times 10^{-3}, \; 10^{-2}, \; 1.5\times10^{-2}, \;
2\times10^{-2}, \; 2.5\times10^{-2}, \; 3\times10^{-2}, \;
4\times10^{-2}, \;5\times 10^{-2} , 6\times 10^{-2},8 \times
10^{-2},10^{-1},2\times 10^{-1},3\times10^{-1}, 4 \times 10^{-1},
5\times10^{-1} \}$. We note that the \APGSUPG, \APGSTGLS, and \APGSTADJ\ ROMs depend upon both
the APG stabilization parameter and the FEM stabilization
parameter. As will be seen, these methods are well-behaved in the low time-step limit 
and as such we execute these methods for
$\tau_{\text{APG}},\tau  \in \tauSet \times
\tauSet$ with a fixed time step $\Delta t = 10^{-3}$ equivalent to the FOM. As will be detailed in the following section, we note that all ROMs and FOMs will be performed on the same spatial grid. 

subsection{Example 1: boundary layer}\label{sec:example1}

\begin{figure}
\begin{center}
\includegraphics[width=0.45\linewidth]{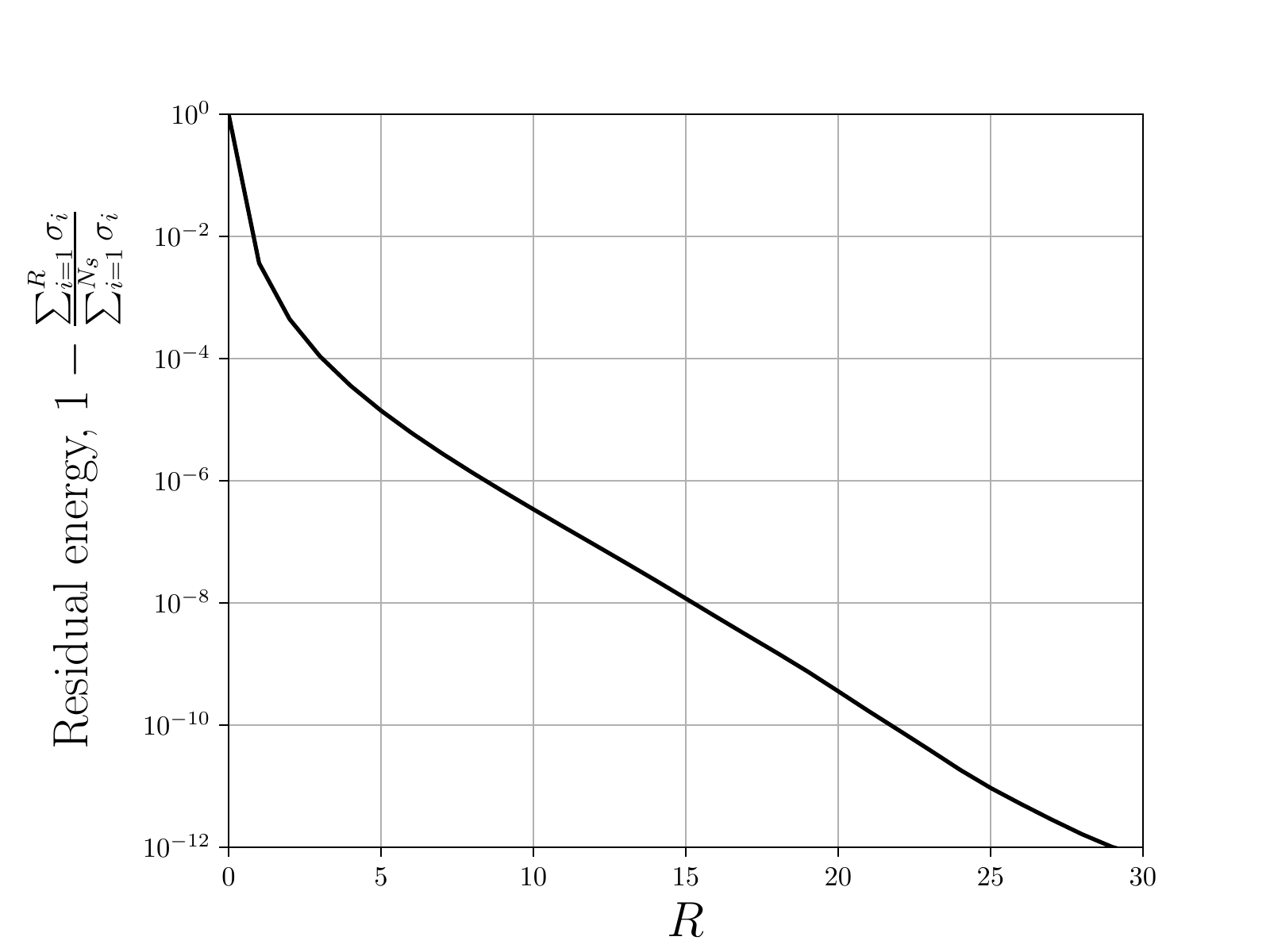}
\caption{\boundaryLayerFigureCaption Residual statistical energy as a function of basis dimension.}
\label{fig:rom_case1_svd}
\end{center}
\end{figure}

The first numerical experiment we consider is a transient version of the setup 
used by Codina in~\cite{CODINA1998185}. 
We solve Eq.~\eqref{eq:cdr} with a final time $T=5$ and a physical domain $x \in \Omega = (0,1) \times (0,1)$. We take the parameters to be a slightly modified version of those used in by Codina in~\cite{CODINA1998185}, where we set $\viscosity = 10^{-3}$, $\forcing = 1$, $\reaction = 1$, and $\wavespeed = \frac{1}{2} \begin{bmatrix} \cos( \pi / 3) & \sin(\pi / 3) \end{bmatrix}^T$. The high-resolution trial space $\cSpaceHRES$ is obtained via a uniform triangulation of $\Omega$ into $N_{\text{el}} = 2 \times 128^2$ elements equipped with a $\mathcal{C}^0(\cDomain)$ continuous discretization with polynomials of order $p=2$. The grid Peclet number is $\text{Pe}_{\text{h}} =  1.953125$, where we used $h=(128 p)^{-1}$. Analogously, the full-order model trial space $\cSpaceFOM$ is obtained via a uniform triangulation of $\Omega$ based on 32 nodes in each direction into $N_{\text{el}} = 2048$ elements equipped with a $\mathcal{C}^0(\cDomain)$ continuous discretization with polynomials of order $p=2$. The grid Peclet number is $\text{Pe}_{\text{h}} =  7.8125$, where we used $h=(32 p)^{-1}$. In both cases, the triangulations are obtained via a partition of $\Omega$ into uniform square cells. The triangles are then cut from bottom-left to top-right of each cell. The Galerkin method equipped with the high resolution trial space is mesh-converged and accurate. The Galerkin method equipped with the FOM trial space yields inaccurate solutions; this error will be quantified later in this section and we consider this case as this is representative of practical problems. Lastly, the ROM trial space is obtained by executing Algorithm~\ref{alg:gen_rom_trialspace}. Figure~\ref{fig:rom_case1_svd} presents the residual statistical energy as a function of basis dimension, where it is seen that the first five basis vectors capture over $99.999\%$ of the statistical energy (in $\LTwo$).

\subsubsection{Full-order model results}\label{sec:example1-fom}
\begin{figure}
\begin{center}
\begin{subfigure}[t]{0.32\textwidth}
\includegraphics[trim={1cm 1cm 1cm 1cm},clip, width=1.\linewidth]{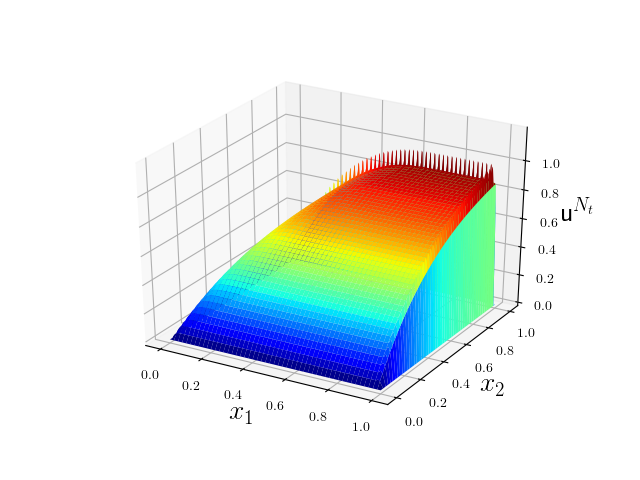}
\caption{H-Res}
\end{subfigure}
\begin{subfigure}[t]{0.32\textwidth}
\includegraphics[trim={1cm 1cm 1cm 1cm},clip, width=1.\linewidth]{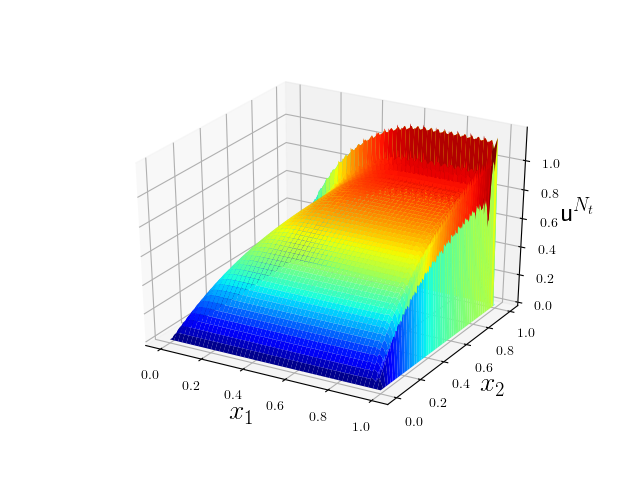}
\caption{Galerkin}
\end{subfigure}
\begin{subfigure}[t]{0.32\textwidth}
\includegraphics[trim={1cm 1cm 1cm 1cm},clip, width=1.\linewidth]{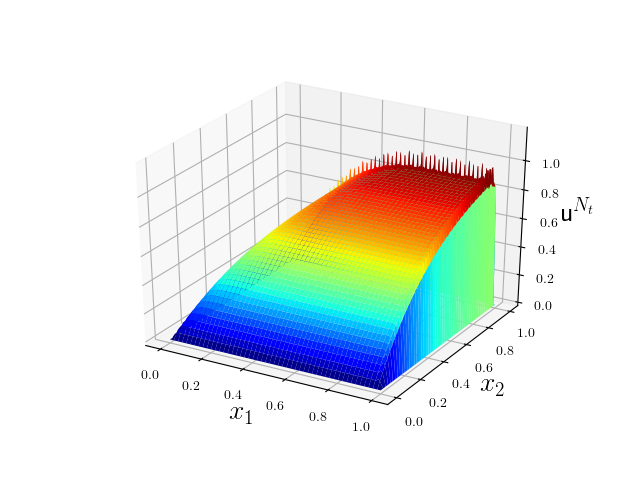}
\caption{\SUPG}
\end{subfigure}
\begin{subfigure}[t]{0.32\textwidth}
\includegraphics[trim={1cm 1cm 1cm 1cm},clip, width=1.\linewidth]{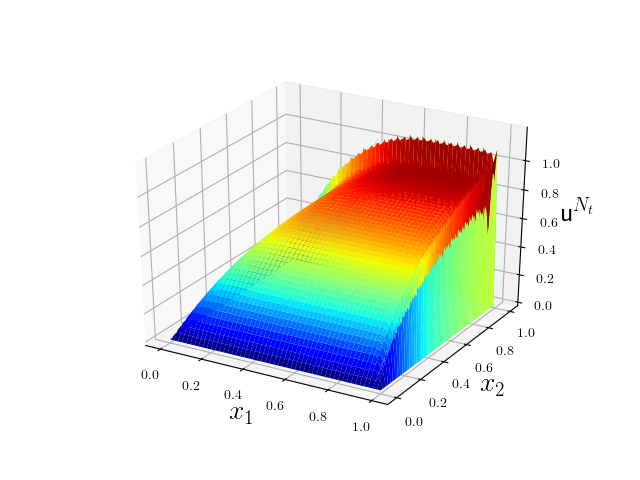}
\caption{\GLS}
\end{subfigure}
\begin{subfigure}[t]{0.32\textwidth}
\includegraphics[trim={1cm 1cm 1cm 1cm},clip,  width=1.\linewidth]{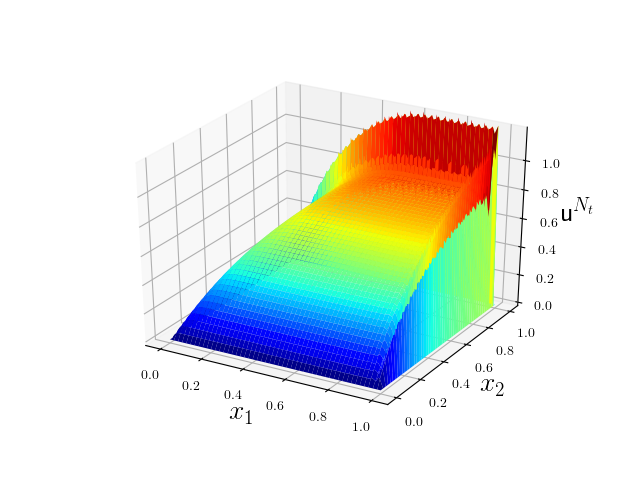}
\caption{\ADJ}
\end{subfigure}
\begin{subfigure}[t]{0.32\textwidth}
\includegraphics[trim={1cm 1cm 1cm 1cm},clip,  width=1.\linewidth]{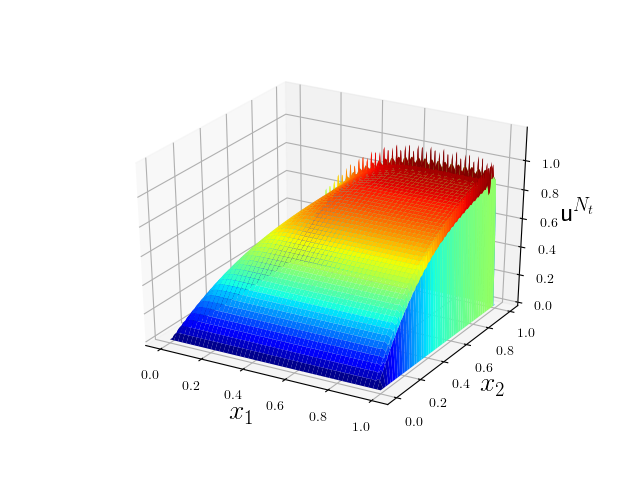}
\caption{\STGLS}
\end{subfigure}
\begin{subfigure}[t]{0.32\textwidth}
\includegraphics[trim={1cm 1cm 1cm 1cm},clip,  width=1.\linewidth]{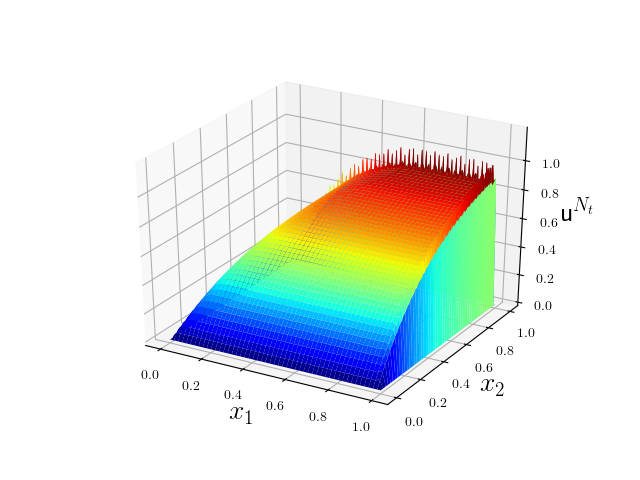}
\caption{\STADJ}
\end{subfigure}

\caption{\boundaryLayerFigureCaption FOM solutions to the \cdcdrAcronym\ equation at $t=5$.} 
\label{fig:ex1_fom_figs}
\end{center}
\end{figure}

We first present results of the various FOMs considered and re-emphasize that (1) these FOMs are executed on a coarse trial space such that the FOM requires stabilization to be accurate and (2) the data from these FOMs are \textit{not} used to construct the ROM trial subspace; instead, we employ high-fidelity data as described in Section~\ref{sec:basis_construction}. The FOM solutions are presented to quantify the underlying FEM error of a given method on this coarse mesh. 
We re-emphasize that we examine the case where the FOM requires stabilization as this is representative of practical problems.
Figure~\ref{fig:ex1_fom_figs} presents the various FOM solutions at the final time, $t=5$, while Table~\ref{tab:ex1_fom_errors} tabulates the solution errors and stabilization parameters employed in the simulations; these parameters were selected by executing the FOMs on the $\tauSet,\dtSet$ grid described above and extracting solutions with the lowest $\LTwo$ error. We observe the following:
\begin{itemize}
\item The Galerkin, \GLS, and \ADJ\ FEM FOMs are the worst performing methods, and all result in solutions with large oscillations at the boundary.
\item The \SUPG, \STGLS, and \STADJ\ FEM FOMs provide the best solutions. These methods result in solution errors that are approximately an order of magnitude better than the Galerkin FOM and provide qualitatively accurate solutions. 
\item The ``space--time" \STGLS\ and \STADJ\ FEM FOMs perform much better than the ``discretize-then-stabilize" \GLS\ and \ADJ\ FEM FOMs.
\end{itemize} 

\begin{table}
  \begin{center}
    \caption{\boundaryLayerFigureCaption Integrated (relative) $\LTwo$ and $\HOne$ errors of various FOMs presented in Figure~\ref{fig:ex1_fom_figs}, along with the stabilization parameters at which the FOMs were executed.}
\label{tab:ex1_fom_errors}
    \begin{filecontents*}{figs/example1/results_fem_foms.csv}
Codes, G, SUPG, GLS, ADJ, STGLS, STADJ
0.000000000000000000e+00,1.834824722163699373e-03,2.105654040677873006e-04,1.692612450636087496e-03,2.094449067099121795e-03,4.342305215844730081e-04,2.375536028119398986e-04
1.000000000000000000e+00,1.559585502722802023e-03,1.616356341809318696e-04,1.482233412636241700e-03,1.605757476403169275e-03,3.522400109615543388e-04,2.040959449879668959e-04
2.000000000000000000e+00,nan,1.000000000000000021e-02,1.000000000000000048e-04,1.000000000000000048e-04,1.000000000000000021e-02,1.000000000000000021e-02
3.000000000000000000e+00,1.000000000000000021e-03,1.000000000000000021e-03,1.000000000000000021e-03,1.000000000000000021e-03,1.000000000000000021e-03,1.000000000000000021e-03
\end{filecontents*}
\pgfplotstabletypeset[
      multicolumn names, 
      col sep=comma, 
      string replace*={0.000000000000000000e+00}{$\integratedErrorLTwoBestFit$},
      string replace*={1.000000000000000000e+00}{$\integratedErrorHOneBestFit$},
      string replace*={2.000000000000000000e+00}{$\tau$},
      string replace*={3.000000000000000000e+00}{$\Delta t$},
      string replace*={errorH}{$\integratedErrorHOneBestFit$},
      string replace*={nan}{N/A},
      display columns/0/.style={
		column name=,string type},  
      display columns/1/.style={
		column name=Galerkin, 
		column type={S[table-parse-only]},string type},  
     display columns/2/.style={
		column name={\begin{centering} \SUPG \end{centering}},
		column type={S[table-parse-only]},string type},
      display columns/3/.style={
		column name=\GLS,
		column type={S[table-parse-only]},string type},
     display columns/4/.style={
		column name=\ADJ,
		column type={S[table-parse-only]},string type},
     display columns/5/.style={
		column name=\STGLS,
		column type={S[table-parse-only]},string type},
     display columns/6/.style={
		column name=\STADJ,
		column type={S[table-parse-only]},string type},
    every head row/.style={
		before row={\toprule}, 
		after row={
			\midrule} 
			},
		every last row/.style={after row=\bottomrule}, 
    ]{figs/example1/results_fem_foms.csv} 
  \end{center}
\end{table}

\subsubsection{Results as a function of basis dimension}
\begin{figure}
\begin{center}
\begin{subfigure}[t]{0.49\textwidth}
\includegraphics[trim={1.5cm 8.8cm 0cm 0cm},clip, width=1.\linewidth]{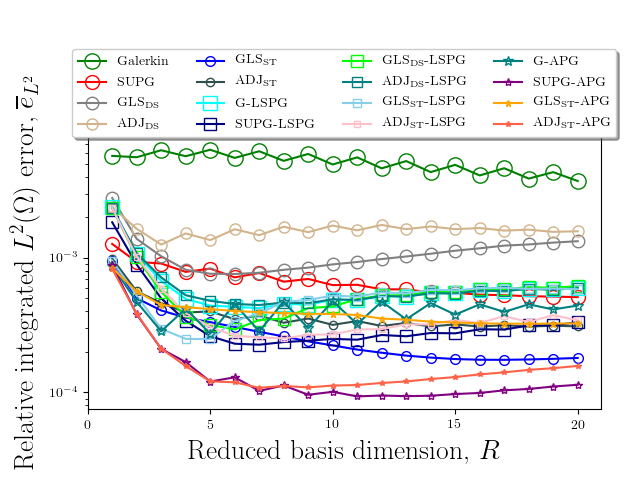}
\end{subfigure}

\begin{subfigure}[t]{0.49\textwidth}
\includegraphics[width=1.\linewidth]{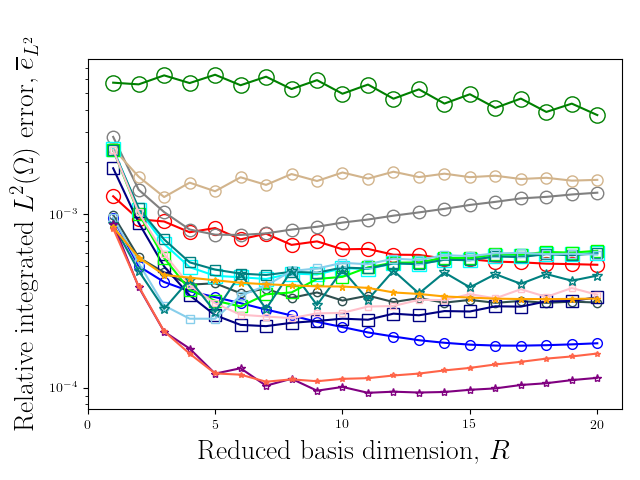}
\caption{$L^2(\cDomain)$ error.}
\end{subfigure}
\begin{subfigure}[t]{0.49\textwidth}
\includegraphics[width=1.\linewidth]{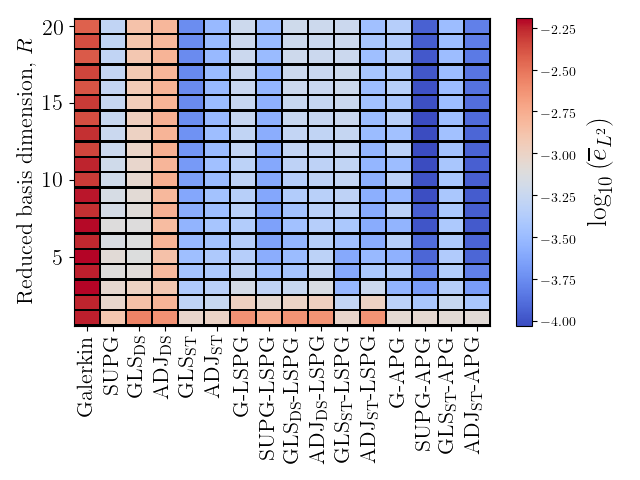}
\caption{$L^2(\cDomain)$ error.}
\end{subfigure}

\begin{subfigure}[t]{0.49\textwidth}
\includegraphics[width=1.\linewidth]{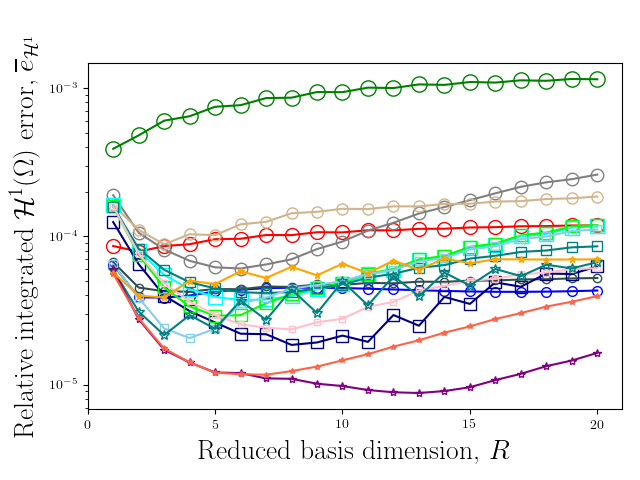}
\caption{$\HOne$ error.}
\end{subfigure}
\begin{subfigure}[t]{0.49\textwidth}
\includegraphics[width=1.\linewidth]{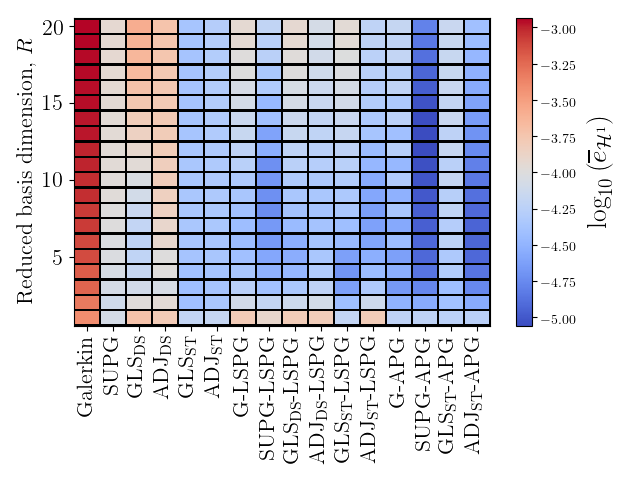}
\caption{$\HOne$ error.}
\end{subfigure}
\caption{\boundaryLayerFigureCaption $\LTwo$ (top) and $\HOne$ (bottom) error as a function of ROM basis dimension for the various ROMs evaluated. We note that the left and right figures show the same quantities, but with different visualization techniques. Results are shown for optimal values of $t,\tau$ as discussed in Section~\ref{sec:timestep_information}.}
\label{fig:bl_romconverge}
\end{center}
\end{figure}

\begin{figure}
\begin{center}
\begin{subfigure}[t]{0.49\textwidth}
\includegraphics[trim={1.5cm 8.8cm 0cm 0cm},clip, width=1.\linewidth]{figs/example1/ex1ConvergenceLegend.png}
\end{subfigure}

\begin{subfigure}[t]{0.49\textwidth}
\includegraphics[width=1.\linewidth]{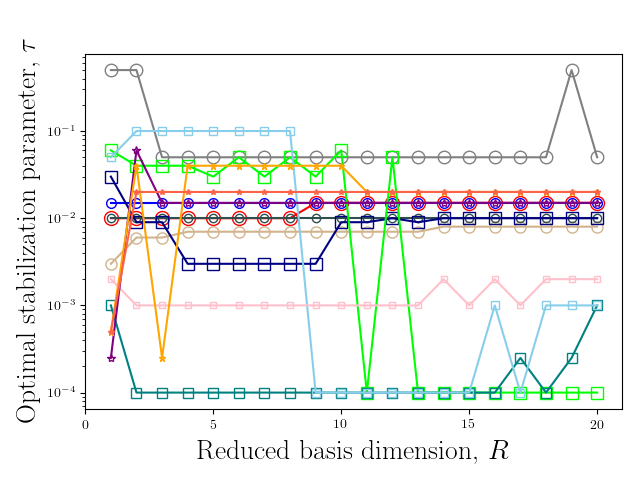}
\caption{Optimal parameter, $\tau$.}
\end{subfigure}
\begin{subfigure}[t]{0.49\textwidth}
\includegraphics[width=1.\linewidth]{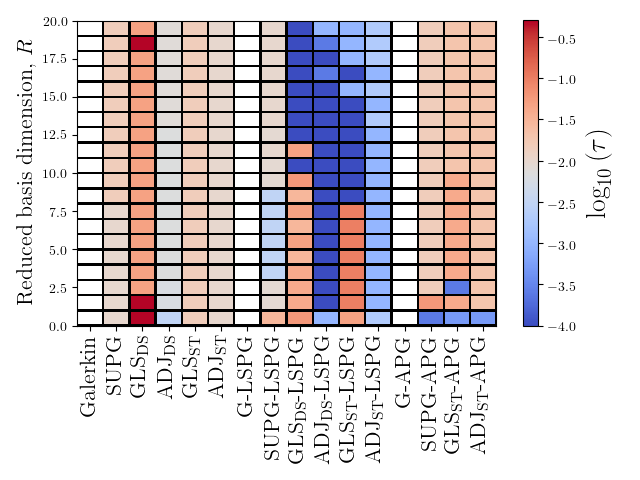}
\caption{Optimal parameter, $\tau$.}
\end{subfigure}

\begin{subfigure}[t]{0.49\textwidth}
\includegraphics[width=1.\linewidth]{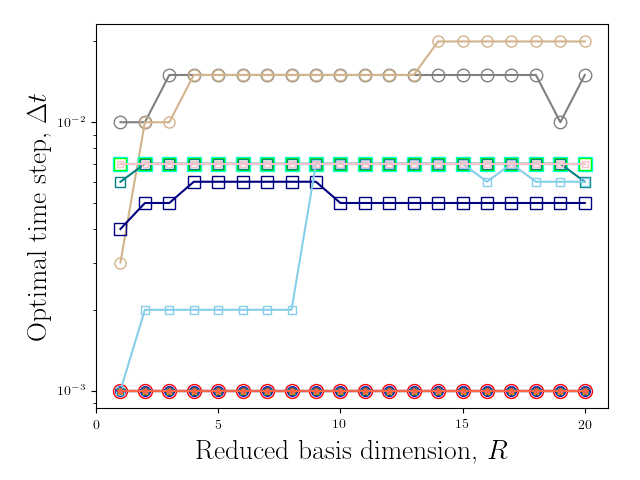}
\caption{Optimal time step, $\Delta t$.}
\end{subfigure}
\begin{subfigure}[t]{0.49\textwidth}
\includegraphics[width=1.\linewidth]{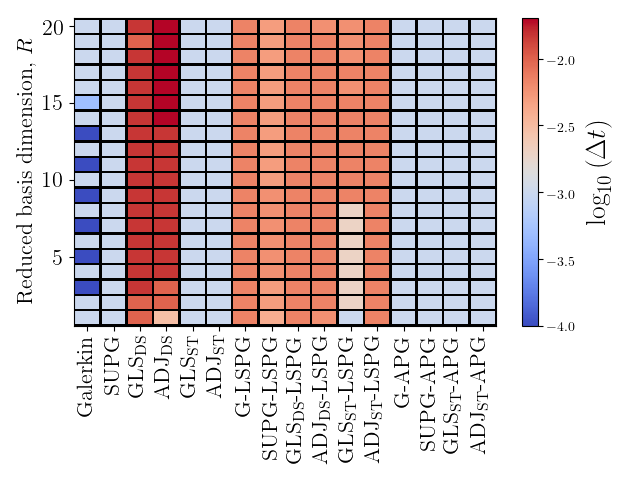}
\caption{Optimal time step, $\Delta t$.}
\end{subfigure}
\caption{\boundaryLayerFigureCaption Optimal stabilization parameter (top) and time step (bottom) as a function of ROM basis dimension for the various ROMs evaluated. We note that the left and right figures show the same quantities, but with different visualization techniques. Results are shown for optimal values of $t,\tau$ as discussed in Section~\ref{sec:timestep_information}}.
\label{fig:bl_romconverge_params}
\end{center}
\end{figure}

We next examine the performance of the various ROMs as the dimension of the ROM basis is varied for $1 \le \romdim \le 20$. For each basis dimension, we present results for optimal ($\tau$, $\Delta t$) as measured by the $\LTwo$ error. Figure~\ref{fig:bl_romconverge} shows the convergence of the $\LTwo$-error and $\HOne$-error as as a function of RB size, while Figure~\ref{fig:bl_romconverge_params} shows the corresponding optimal stabilization parameters and time steps. We make the following observations about the \textit{accuracy} of the various ROMs:
\begin{itemize}
\item No method results in a monotonic decrease in error in both the $\LTwo$ and $\HOne$ norms. This is, in part, a result of the fact that the ROMs are not consistent with the high resolution reference solution (see the discussion in Section~\ref{sec:basis_construction}). 
\item The Galerkin ROM performs poorly for all basis dimensions.  
\item \APGSUPG, \APGSTADJ, \STGLS, and \LSPGSUPG\ are the best performing ROMs.
\item \STGLS\ and \STADJ\ are consistently more accurate than \SUPG.
\item The ``discretize-then-stabilize" ROMs outperform the Galerkin ROM, but are consistently worse than their space--time counterparts. 
\item LSPG ROMs outperform their continuous counterparts in all cases. 
\item APG ROMs outperform their continuous counterparts in all cases except for \APGSTGLS. 
\item It is interesting to note that the \APGG\ and \LSPGG\ ROMs perform significantly better than the standard Galerkin ROM for large number of bases, even though these two methods formally converge to the Galerkin ROM in the limit of a full basis.  
\item A decrease in error in the $\LTwo$ norm does not always correspond to a decrease in error in the $\HOne$ norm, and vice versa.
\item Comparing Figure~\ref{fig:bl_romconverge} to Table~\ref{tab:ex1_fom_errors}, it is interesting to observe that, while the same trends are observed, some ROMs are more accurate than their corresponding FEM models. This is again likely a result of the inconsistency between the ROMs and the high-resolution reference solution.
\end{itemize}

Examining Figure~\ref{fig:bl_romconverge_params}, we make the following observations about the behavior of the stabilization parameters of the various ROMs:

\begin{itemize}
\item The optimal stabilization parameters for the \SUPG, \STADJ, \STGLS, and APG ROMs are more or less constant for all reduced basis dimensions (with the exception of a few APG solutions at very small basis dimensions). 
\item The \GLS, \ADJ, and LSPG-based ROMs are optimal for time step sizes larger than the FOM. In particular, the optimal time step for almost all LSPG ROMs occurs at an intermediate time step. This is well documented in the literature~\cite{carlberg_lspg_v_galerkin}.
\item It is difficult to decipher any pattern in the optimal stabilization parameters for the LSPG-based ROMs. We expect that this is, in part, due to the complex interplay between the dependence of the time-step \text{and} stabilization parameters on the LSPG ROM performance.
\end{itemize}

\begin{figure}
\begin{center}
\begin{subfigure}[t]{0.28\textwidth}
\includegraphics[trim={1cm 1cm 1cm 1cm},clip, width=1.\linewidth]{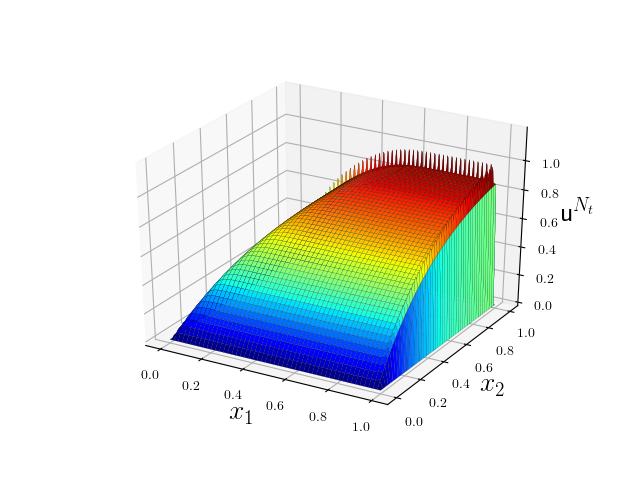}
\caption{H-Res ($L^2(\cDomain)$ best fit)}
\label{fig:rom_fig1b}
\end{subfigure}
\begin{subfigure}[t]{0.28\textwidth}
\includegraphics[trim={1cm 1cm 1cm 1cm},clip, width=1.\linewidth]{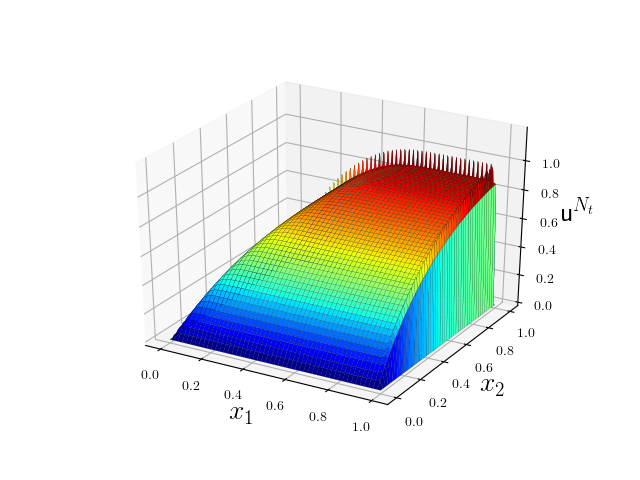}
\caption{H-Res ($H_1(\cDomain)$ best fit)}
\label{fig:rom_fig1c}
\end{subfigure}
\begin{subfigure}[t]{0.28\textwidth}
\includegraphics[trim={1cm 1cm 1cm 1cm},clip,  width=1.\linewidth]{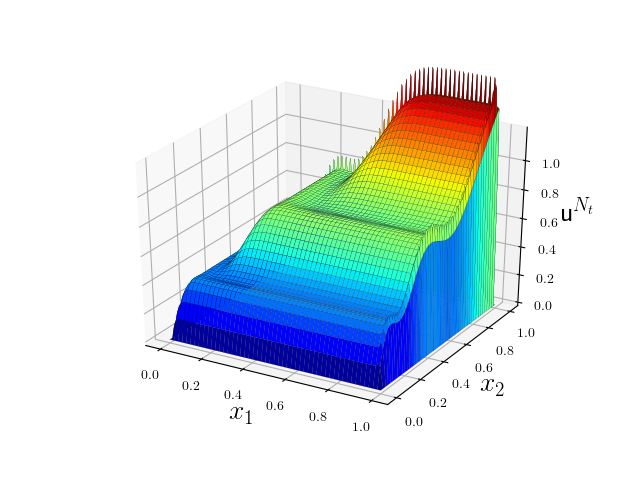}
\caption{Galerkin}
\label{fig:rom_fig1d}
\end{subfigure}
\begin{subfigure}[t]{0.28\textwidth}
\includegraphics[trim={1cm 1cm 1cm 1cm},clip,  width=1.\linewidth]{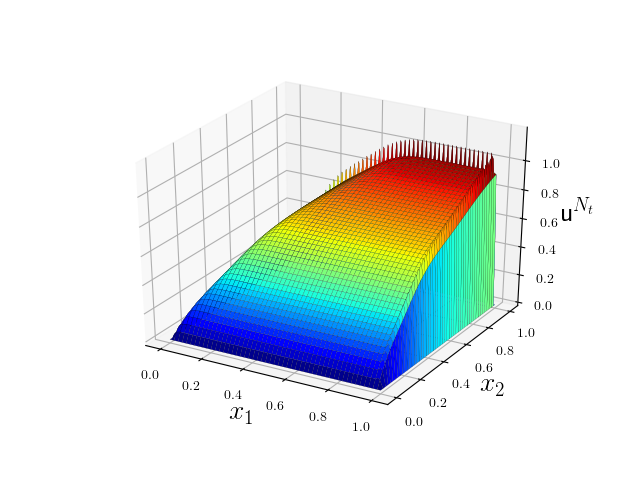}
\caption{\SUPG}
\label{fig:rom_fig1e}
\end{subfigure}
\begin{subfigure}[t]{0.28\textwidth}
\includegraphics[trim={1cm 1cm 1cm 1cm},clip,  width=1.\linewidth]{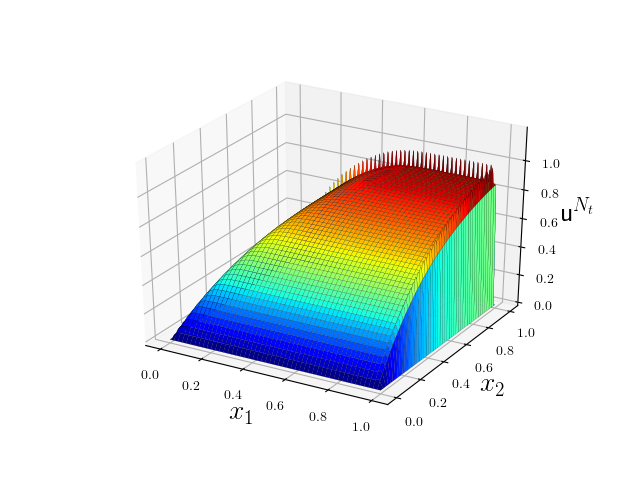}
\caption{\GLS}
\label{fig:rom_fig1f}
\end{subfigure}
\begin{subfigure}[t]{0.28\textwidth}
\includegraphics[trim={1cm 1cm 1cm 1cm},clip,  width=1.\linewidth]{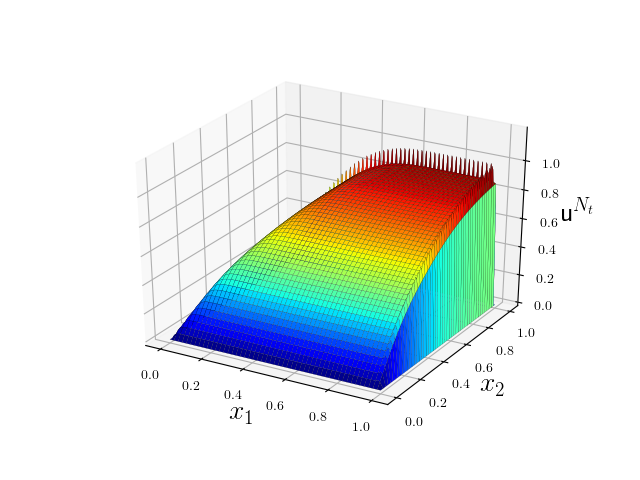}
\caption{\ADJ}
\label{fig:rom_fig1g}
\end{subfigure}
\begin{subfigure}[t]{0.28\textwidth}
\includegraphics[trim={1cm 1cm 1cm 1cm},clip,  width=1.\linewidth]{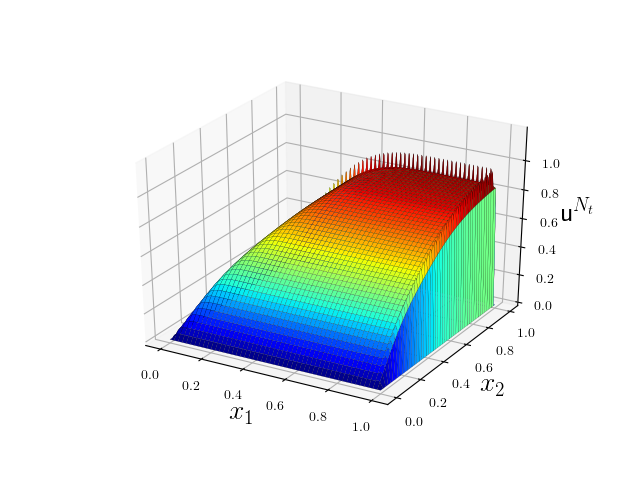}
\caption{\STGLS}
\label{fig:rom_fig1f}
\end{subfigure}
\begin{subfigure}[t]{0.28\textwidth}
\includegraphics[trim={1cm 1cm 1cm 1cm},clip,  width=1.\linewidth]{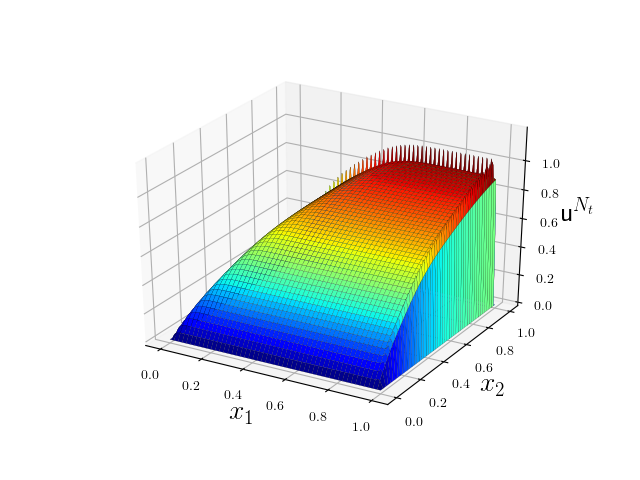}
\caption{\STADJ}
\label{fig:rom_fig1g}
\end{subfigure}
\begin{subfigure}[t]{0.28\textwidth}
\includegraphics[ trim={1cm 1cm 1cm 1cm},clip, width=1.\linewidth]{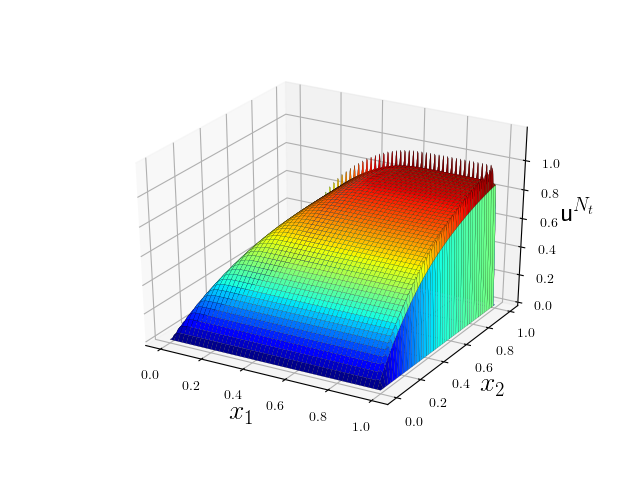}
\caption{\LSPGG}
\label{fig:rom_fig1h}
\end{subfigure}
\begin{subfigure}[t]{0.28\textwidth}
\includegraphics[ trim={1cm 1cm 1cm 1cm},clip, width=1.\linewidth]{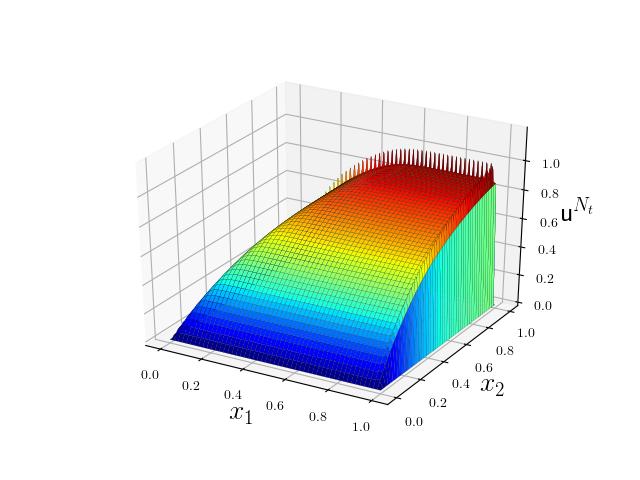}
\caption{\LSPGSUPG}
\label{fig:rom_fig1i}
\end{subfigure}
\begin{subfigure}[t]{0.28\textwidth}
\includegraphics[ trim={1cm 1cm 1cm 1cm},clip, width=1.\linewidth]{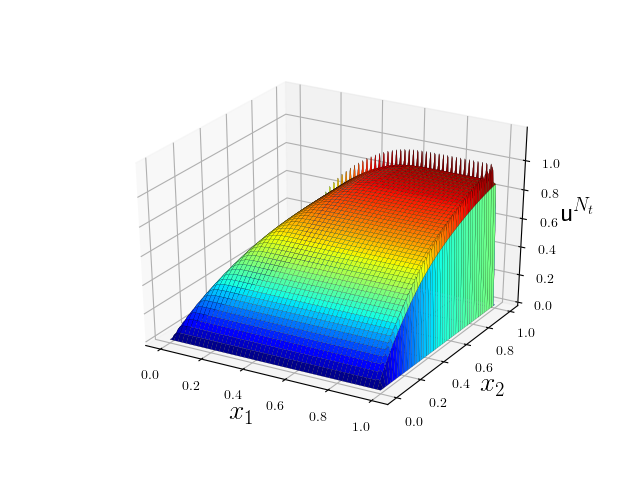}
\caption{\LSPGADJ}
\label{fig:rom_fig1i}
\end{subfigure}
\begin{subfigure}[t]{0.28\textwidth}
\includegraphics[ trim={1cm 1cm 1cm 1cm},clip, width=1.\linewidth]{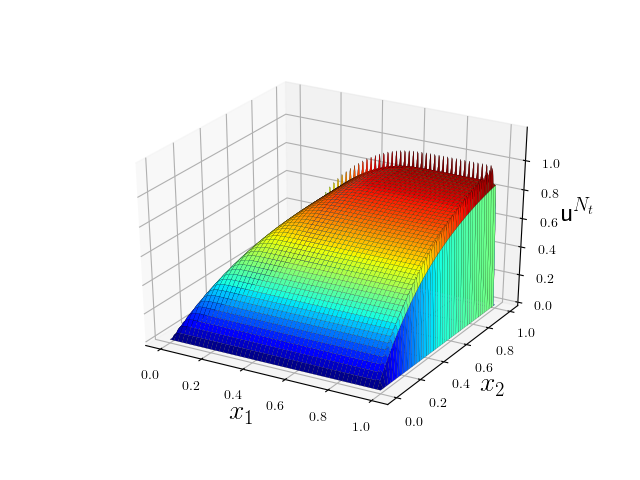}
\caption{\LSPGADJ}
\label{fig:rom_fig1i}
\end{subfigure}
\begin{subfigure}[t]{0.28\textwidth}
\includegraphics[ trim={1cm 1cm 1cm 1cm},clip, width=1.\linewidth]{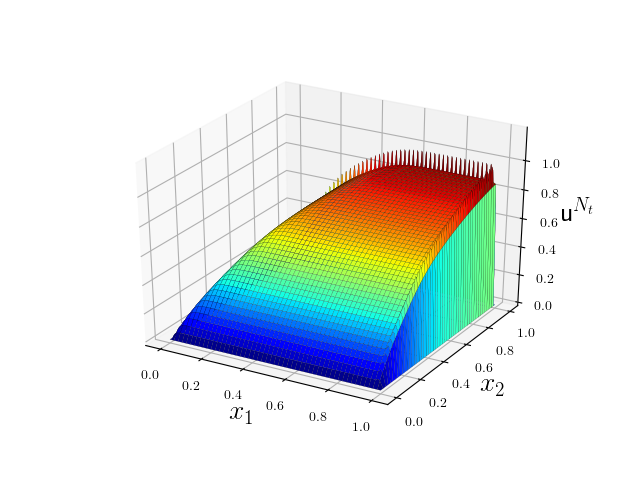}
\caption{\LSPGSTGLS}
\label{fig:rom_fig1i}
\end{subfigure}
\begin{subfigure}[t]{0.28\textwidth}
\includegraphics[ trim={1cm 1cm 1cm 1cm},clip, width=1.\linewidth]{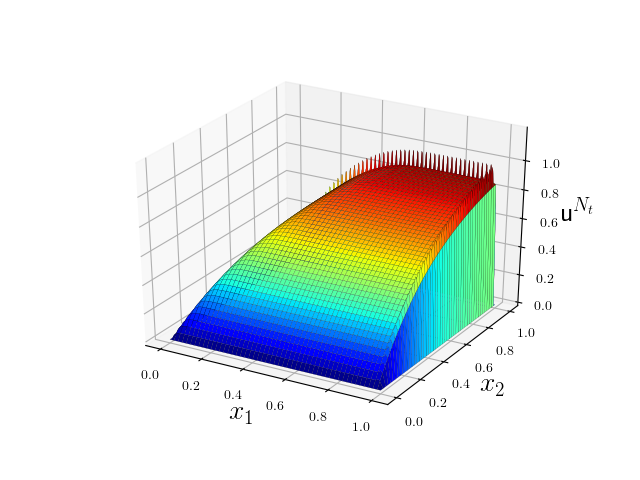}
\caption{\LSPGSTADJ}
\label{fig:rom_fig1i}
\end{subfigure}
\begin{subfigure}[t]{0.28\textwidth}
\includegraphics[ trim={1cm 1cm 1cm 1cm},clip, width=1.\linewidth]{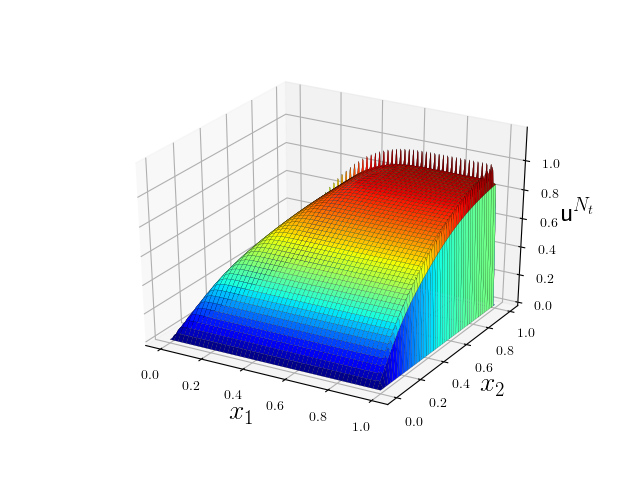}
\caption{\APGG}
\label{fig:rom_fig1i}
\end{subfigure}
\begin{subfigure}[t]{0.28\textwidth}
\includegraphics[ trim={1cm 1cm 1cm 1cm},clip, width=1.\linewidth]{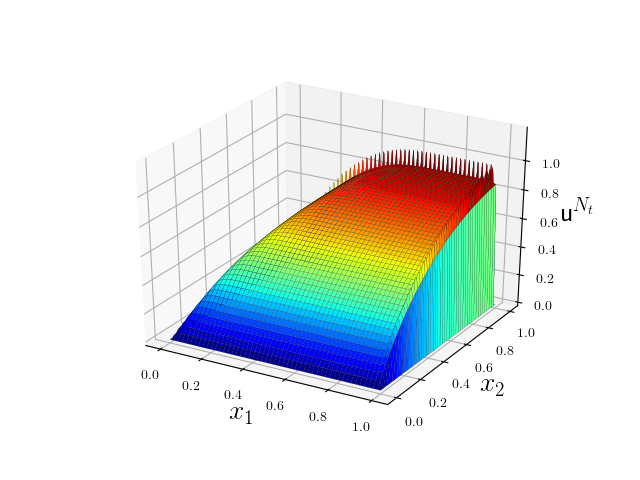}
\caption{\APGSUPG}
\label{fig:rom_fig1i}
\end{subfigure}
\begin{subfigure}[t]{0.28\textwidth}
\includegraphics[ trim={1cm 1cm 1cm 1cm},clip, width=1.\linewidth]{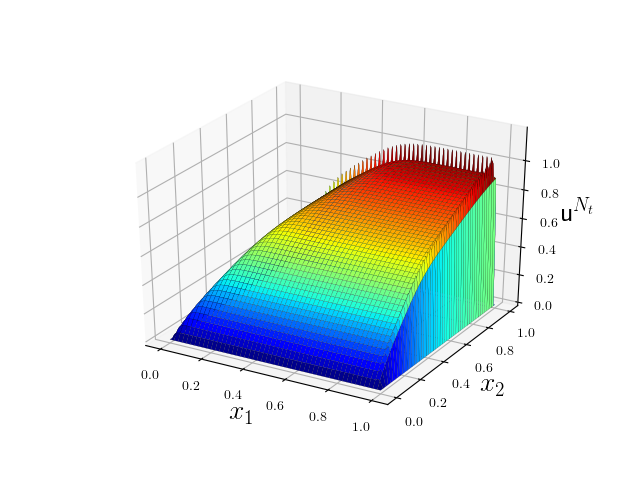}
\caption{\APGSTGLS}
\label{fig:rom_fig1i}
\end{subfigure}
\begin{subfigure}[t]{0.28\textwidth}
\includegraphics[ trim={1cm 1cm 1cm 1cm},clip, width=1.\linewidth]{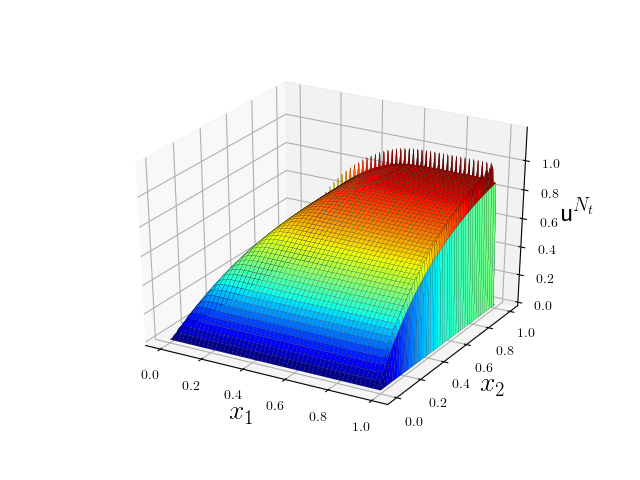}
\caption{\APGSTADJ}
\label{fig:rom_fig1i}
\end{subfigure}
\caption{\boundaryLayerFigureCaption ROM solutions to the \cdcdrAcronym\ equation at $t=5$. Results are shown for solutions obtained with $\romdim = 5$ and with the optimal time step and stabilization parameter as measured by the $\LTwo$-error.} 
\label{fig:rom_fig1}
\end{center}
\end{figure}

Next, Figure~\ref{fig:rom_fig1} presents solution profiles for the various ROMs at a reduced basis dimension of $\romdim=5$, which corresponds to an energy criterion of $\energyCutoff=0.99999$, and at optimal values of $\tau$ and $\Delta t$ (and $\tauApg$) for APG ROMs) for the final time instance, $t=5$. We observe that the projected truth solution displays a small oscillation at the boundary. This oscillation is a result of the FOM trial space $\cSpaceFOM$ being unable to fully resolve the boundary layer. Next, we see that the Galerkin ROM results in inaccurate solutions with large-scale oscillations. All stabilized ROMs are qualitatively accurate with minimal variation between their solutions. 

\subsubsection{Sensitivity to time step and stabilization parameters}
\begin{figure}
\begin{center}
\begin{subfigure}[t]{0.32\textwidth}
\includegraphics[width=1.\linewidth]{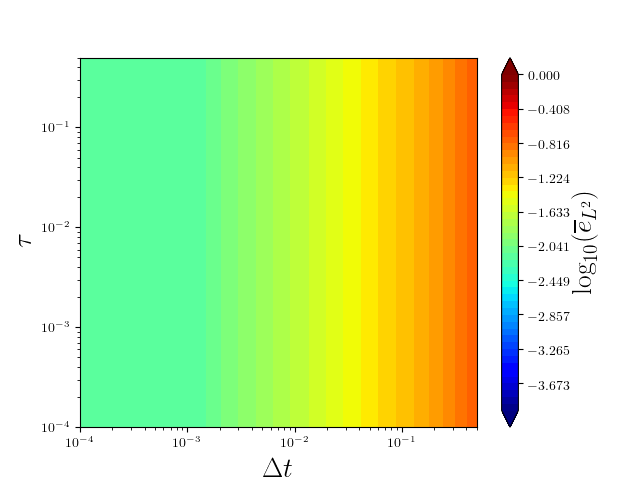}
\caption{Galerkin}
\label{fig:rom_bl_taudtvary_g}
\end{subfigure}
\begin{subfigure}[t]{0.32\textwidth}
\includegraphics[width=1.\linewidth]{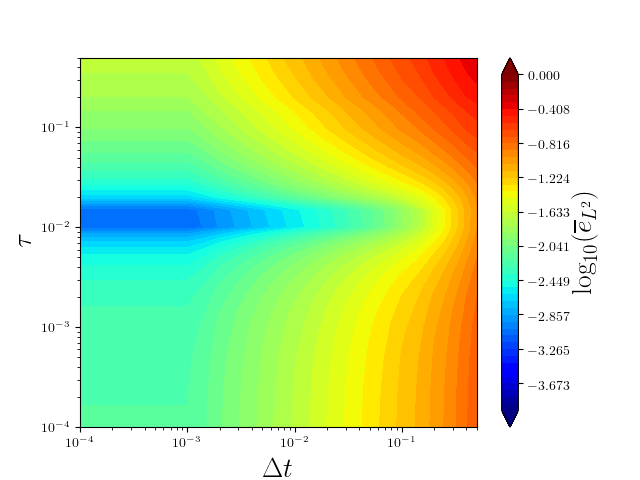}
\caption{\SUPG}
\label{fig:rom_bl_taudtvary_supg}
\end{subfigure}
\begin{subfigure}[t]{0.32\textwidth}
\includegraphics[width=1.\linewidth]{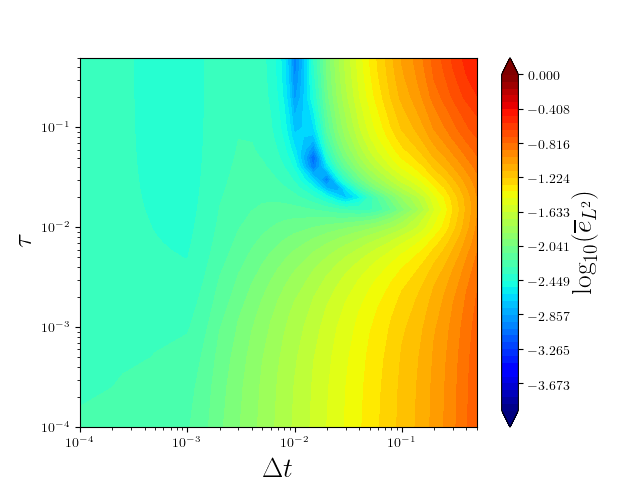}
\caption{\GLS}
\label{fig:rom_bl_taudtvary_gls}
\end{subfigure}
\begin{subfigure}[t]{0.32\textwidth}
\includegraphics[width=1.\linewidth]{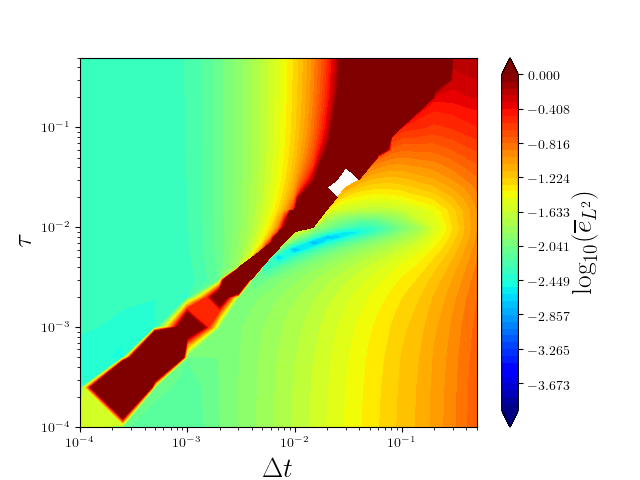}
\caption{\ADJ}
\label{fig:rom_bl_taudtvary_adj}
\end{subfigure}
\begin{subfigure}[t]{0.32\textwidth}
\includegraphics[width=1.\linewidth]{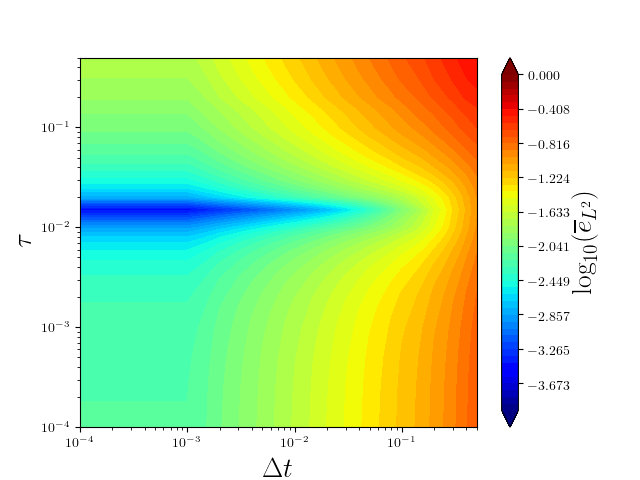}
\caption{\STGLS}
\label{fig:rom_bl_taudtvary_stgls}
\end{subfigure}
\begin{subfigure}[t]{0.32\textwidth}
\includegraphics[width=1.\linewidth]{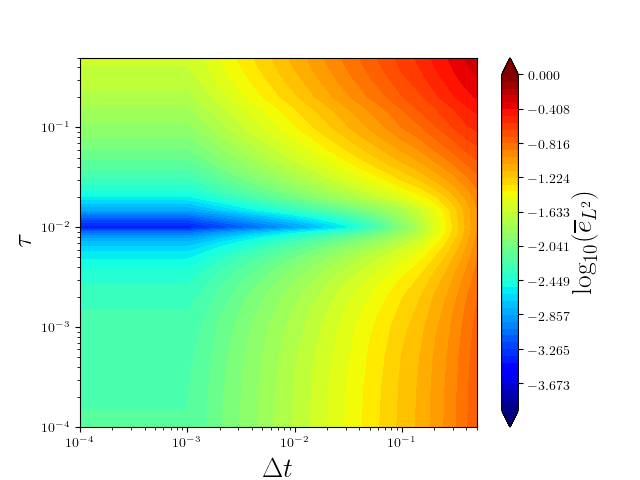}
\caption{\STADJ}
\label{fig:rom_bl_taudtvary_stadj}
\end{subfigure}
\begin{subfigure}[t]{0.32\textwidth}
\includegraphics[width=1.\linewidth]{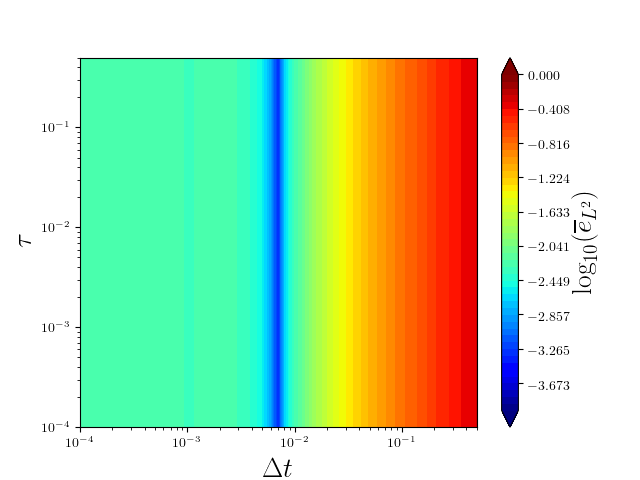}
\caption{\LSPGG}
\label{fig:rom_bl_taudtvary_lspg}
\end{subfigure}
\begin{subfigure}[t]{0.32\textwidth}
\includegraphics[width=1.\linewidth]{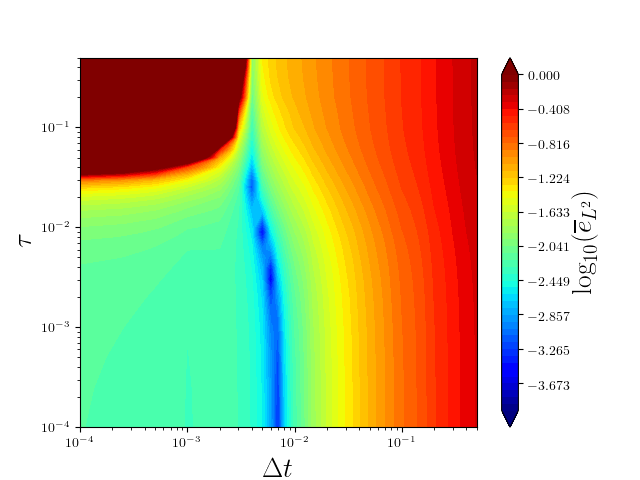}
\caption{\LSPGSUPG}
\end{subfigure}
\begin{subfigure}[t]{0.32\textwidth}
\includegraphics[width=1.\linewidth]{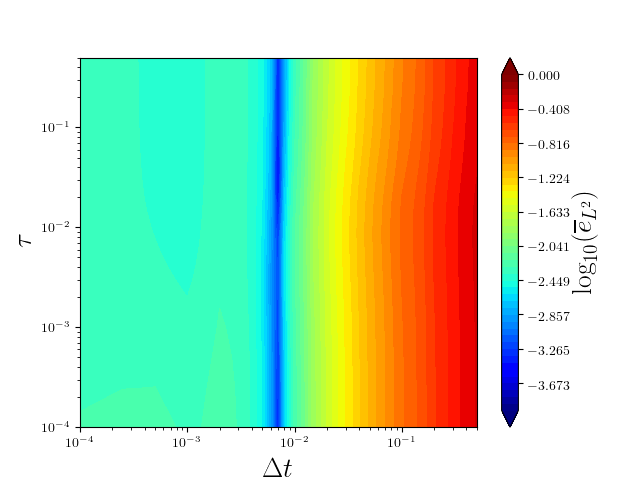}
\caption{\LSPGGLS}
\end{subfigure}
\begin{subfigure}[t]{0.32\textwidth}
\includegraphics[width=1.\linewidth]{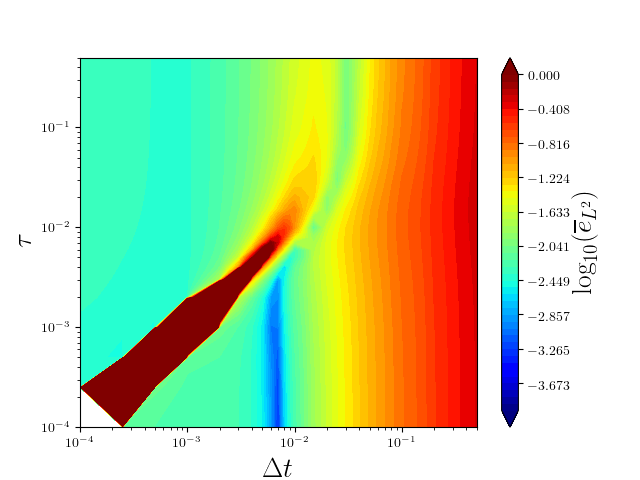}
\caption{\LSPGADJ}
\end{subfigure}
\begin{subfigure}[t]{0.32\textwidth}
\includegraphics[width=1.\linewidth]{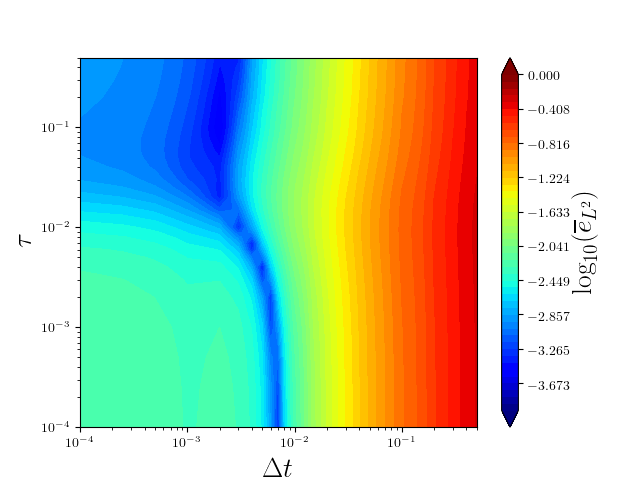}
\caption{\STGLSLSPG}
\end{subfigure}
\begin{subfigure}[t]{0.32\textwidth}
\includegraphics[width=1.\linewidth]{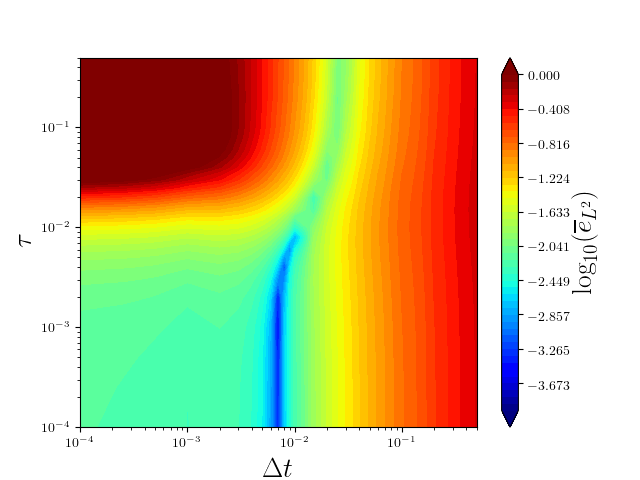}
\caption{\STADJLSPG}
\end{subfigure}
\caption{\boundaryLayerFigureCaption Time integrated $\LTwo$ best-fit error as a function of time step and stabilization parameter for the various ROMs evaluated. Note that Galerkin and LSPG display no dependence on the stabilization parameter. White regions indicate regions where the solution diverged to NaN.}
\label{fig:rom_bl_taudtvary}
\end{center}
\end{figure}

\begin{figure}
\begin{center}
\begin{subfigure}[t]{0.32\textwidth}
\includegraphics[width=1.\linewidth]{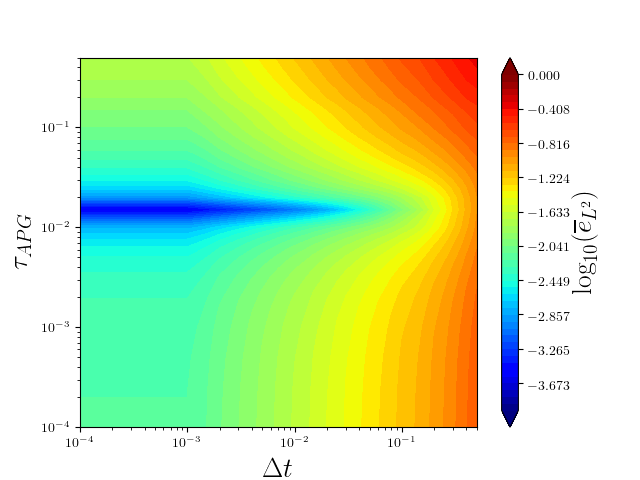}
\caption{\APGG}
\label{fig:rom_bl_taudtvary_apgg}
\end{subfigure}
\begin{subfigure}[t]{0.32\textwidth}
\includegraphics[width=1.\linewidth]{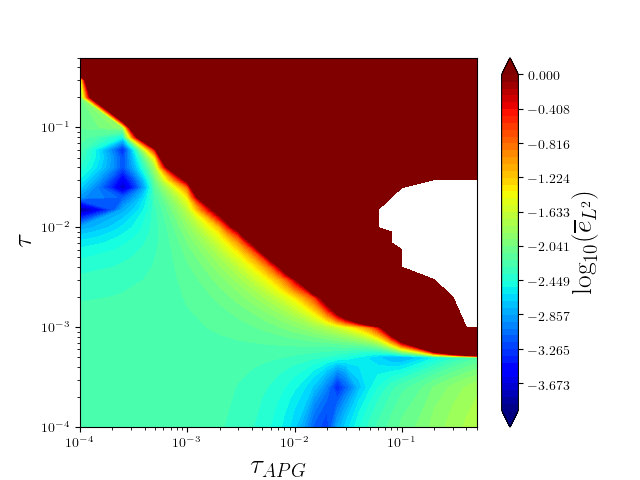}
\caption{\APGSUPG}
\label{fig:rom_bl_taudtvary_apgsupg}
\end{subfigure}
\begin{subfigure}[t]{0.32\textwidth}
\includegraphics[width=1.\linewidth]{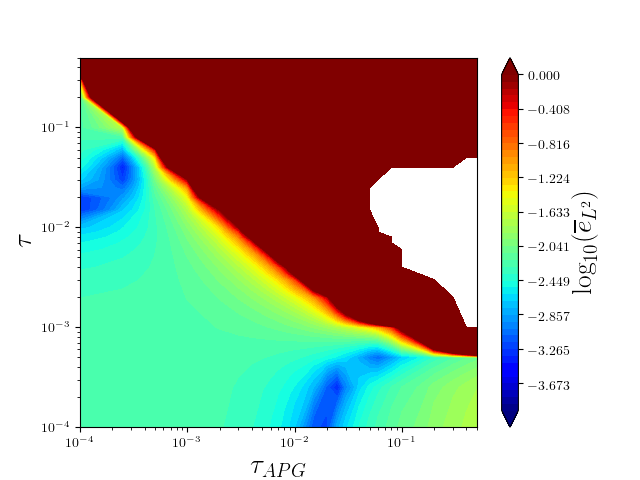}
\caption{\APGSTGLS}
\end{subfigure}
\begin{subfigure}[t]{0.32\textwidth}
\includegraphics[width=1.\linewidth]{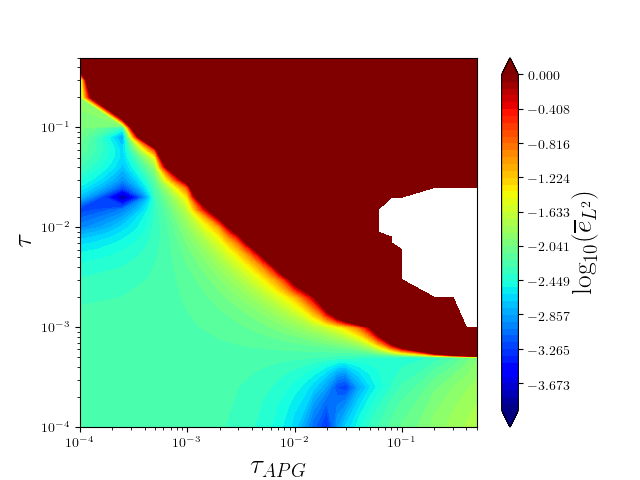}
\caption{\APGSTADJ}
\label{fig:rom_bl_taudtvary_apgstadj}
\end{subfigure}
\caption{\boundaryLayerFigureCaption Time integrated $\LTwo$ error as a function of time step and stabilization parameter for the APG ROMs evaluated. White regions indicate regions where the solution diverged to NaN.}
\label{fig:rom_bl_taudtvary_apg}
\end{center}
\end{figure}

\begin{figure}
\begin{center}
\begin{subfigure}[t]{0.32\textwidth}
\includegraphics[width=1.\linewidth]{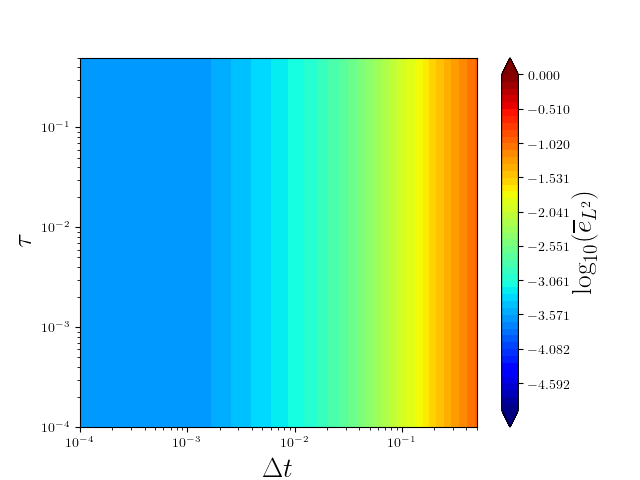}
\caption{Galerkin FEM}
\label{fig:fom_bl_taudtvary_g}
\end{subfigure}
\begin{subfigure}[t]{0.32\textwidth}
\includegraphics[width=1.\linewidth]{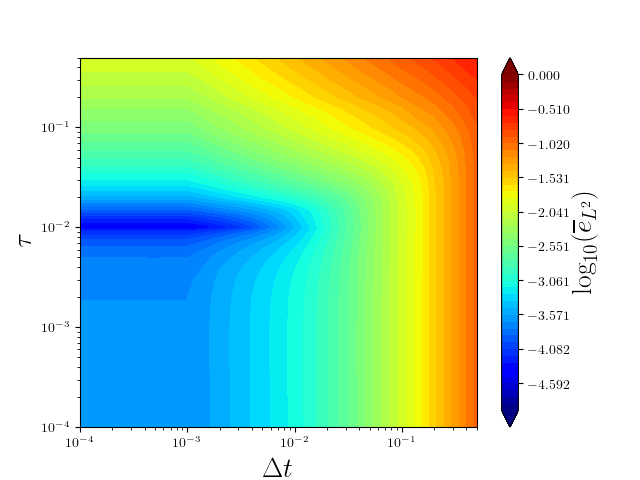}
\caption{\SUPG\ FEM}
\end{subfigure}
\begin{subfigure}[t]{0.32\textwidth}
\includegraphics[width=1.\linewidth]{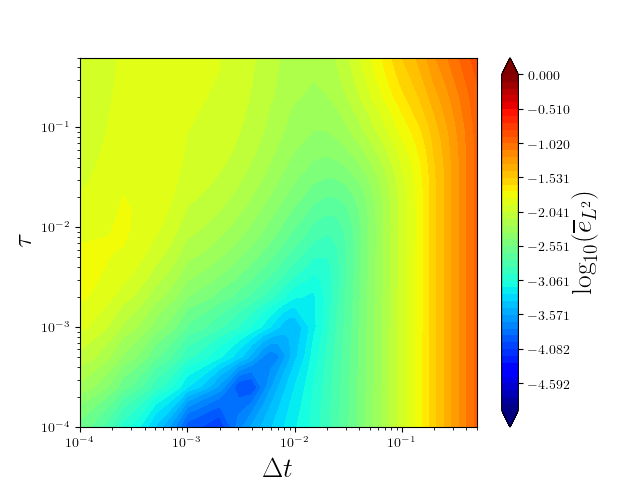}
\caption{\GLS\ FEM}
\end{subfigure}
\begin{subfigure}[t]{0.32\textwidth}
\includegraphics[width=1.\linewidth]{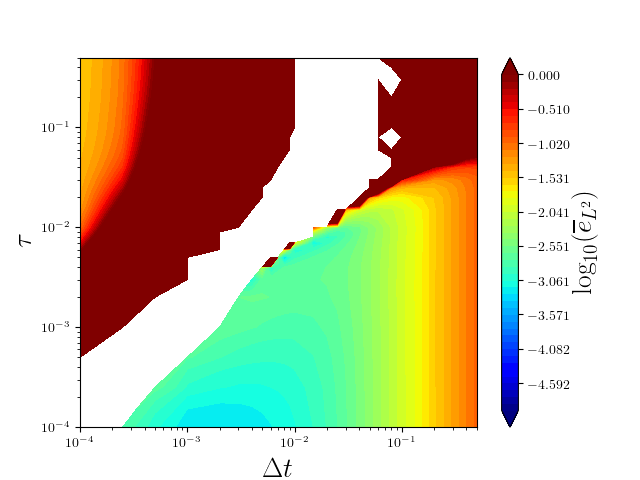}
\caption{\ADJ\ FEM}
\end{subfigure}
\begin{subfigure}[t]{0.32\textwidth}
\includegraphics[width=1.\linewidth]{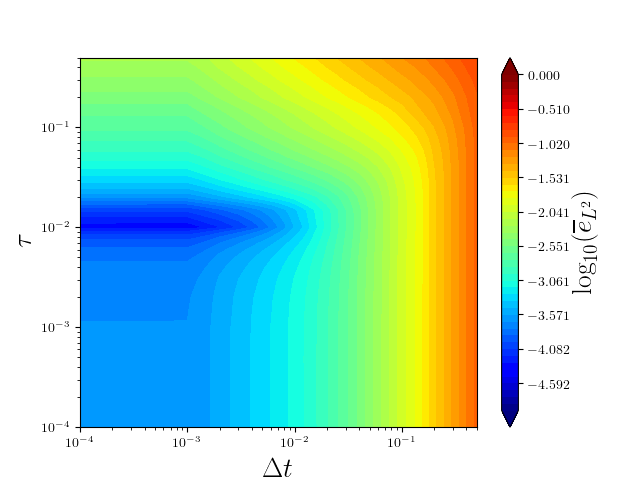}
\caption{\STGLS\ FEM}
\label{fig:fom_bl_taudtvary_stgls}
\end{subfigure}
\begin{subfigure}[t]{0.32\textwidth}
\includegraphics[width=1.\linewidth]{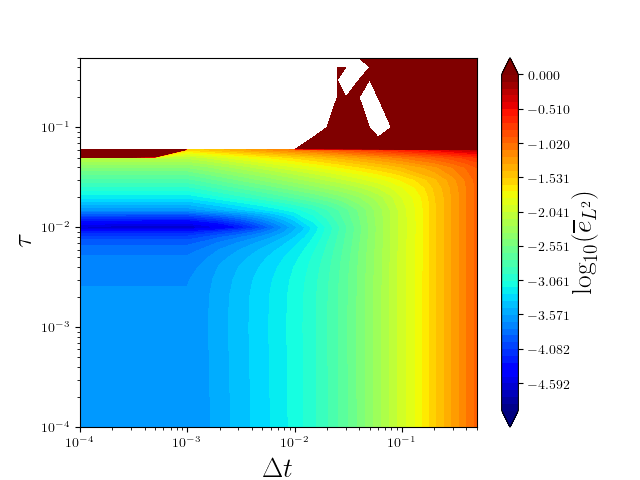}
\caption{\STADJ\ FEM}
\label{fig:fom_bl_taudtvary_stadj}
\end{subfigure}
\caption{\boundaryLayerFigureCaption Time integrated $\LTwo$ error as a function of time step and stabilization parameter for various FEM models. Results are shown for full-order FEM solutions executed on the FOM trial space. White regions indicate regions where the solution diverged to NaN.}
\label{fig:fom_bl_taudtvary}
\end{center}
\end{figure}

The performance of the stabilized methods can depend on both the stabilization parameter $\tau$ and the time step $\Delta t$ (and, for APG, the APG stabilization parameter $\tauApg$). To quantify this sensitivity, Figure~\ref{fig:rom_bl_taudtvary} presents results for the continuous and LSPG ROM solutions obtained on the parameter grid $(\tau, \Delta t ) \in \tauSet \times \dtSet$. 
Figure~\ref{fig:rom_bl_taudtvary_apg} presents results for the \APGG\ ROM solution obtained on the parameter grid $(\tauApg, \Delta t ) \in \tauSet \times \dtSet$ and the remaining APG ROM solutions (which depend on three parameters, $\Delta t, \tau,$ and $\tauApg$) obtained on the parameter grid $(\tau,\tauApg) \in \tauSet \times \tauSet$ with a fixed time step $\Delta t = 10^{-3}$. All ROM results are shown for a reduced basis dimension $\romdim = 5$. 
As a reference, Figure~\ref{fig:fom_bl_taudtvary} shows the same results, but for full-order finite element simulations executed on the FOM trial space. We observe the following:
\begin{itemize}
\item In the limit that $\tau \rightarrow 0$ (or $\Delta t \rightarrow 0$ for \LSPG), all ROMs converge to the standard G-ROM with the exception of \ADJ. This ROM displays poor behavior in the low time-step limit when $\tau \approx \Delta t$.
\item The \SUPG\ (Figure~\ref{fig:rom_bl_taudtvary_supg}), \STGLS\ (Figure~\ref{fig:rom_bl_taudtvary_stgls}), \STADJ\ (Figure~\ref{fig:rom_bl_taudtvary_stadj}), and \APGG\ (Figure~\ref{fig:rom_bl_taudtvary_apgg}) ROMs again all display a similar dependence on the time step and stabilization parameter. Optimal results are obtained for an intermediate value of $\tau$, and the solutions all converge in the limit of $\Delta t \rightarrow 0$.
\item All LSPG ROMs (Figure~\ref{fig:rom_bl_taudtvary_lspg}-\ref{fig:rom_bl_taudtvary_stadj}) yield optimal results at an intermediate time step, and are thus not robust in the low time-step limit. LSPG's optimality at an intermediate time step is well documented in the community~\cite{carlberg_lspg_v_galerkin}, and it is well known that LSPG converges to Galerkin in the low time-step limit. 
\item Errors in the \SDGLS (Figure~\ref{fig:rom_bl_taudtvary_gls}) and \SDADJ\ (Figure~\ref{fig:rom_bl_taudtvary_adj}) ROMs start to increase once the time step becomes small enough, and thus these ROMs are not robust in the low time-step limit. 
\item The Galerkin (Figure~\ref{fig:rom_bl_taudtvary_g}), \SUPG\ (Figure~\ref{fig:rom_bl_taudtvary_supg}), \STGLS (Figure~\ref{fig:rom_bl_taudtvary_stgls}), and \STADJ\ (Figure~\ref{fig:rom_bl_taudtvary_stadj}) ROMs display a similar dependence to the stabilization parameter and time step as their corresponding FOMs (Figure~\ref{fig:fom_bl_taudtvary_g}-\ref{fig:rom_bl_taudtvary_stadj}). The behaviors of the \SDGLS\ and \SDADJ\ ROMs display some qualitative similarities with their corresponding FEM solutions, but in general are different.
\item For APG ROMs built on top of a stabilized FEM model (Figures~\ref{fig:rom_bl_taudtvary_apgsupg}-\ref{fig:rom_bl_taudtvary_apgstadj}), optimal results are obtained for either an intermediate value of $\tauApg$ and low value of $\tau$, or vice versa. It is interesting to note that the solutions are almost symmetric with respect to these two parameters. In addition, we see regions of instability for high values of $\tauApg$. 
\end{itemize} 

\subsection{Example 2: Advecting front}\label{sec:example2}
The second numerical experiment we consider examines the \cdcdrAcronym\ equation in a setting that yields an advecting front. ROMs of this problem require more basis vectors to accurately characterize the system and it is easier to examine the regime where the ROM itself is under-resolved.

\subsubsection{Description of problem setup, full-order model, and generation of trial spaces}\label{ex:example2_description}
\begin{figure}
\begin{center}
\includegraphics[width=0.45\linewidth]{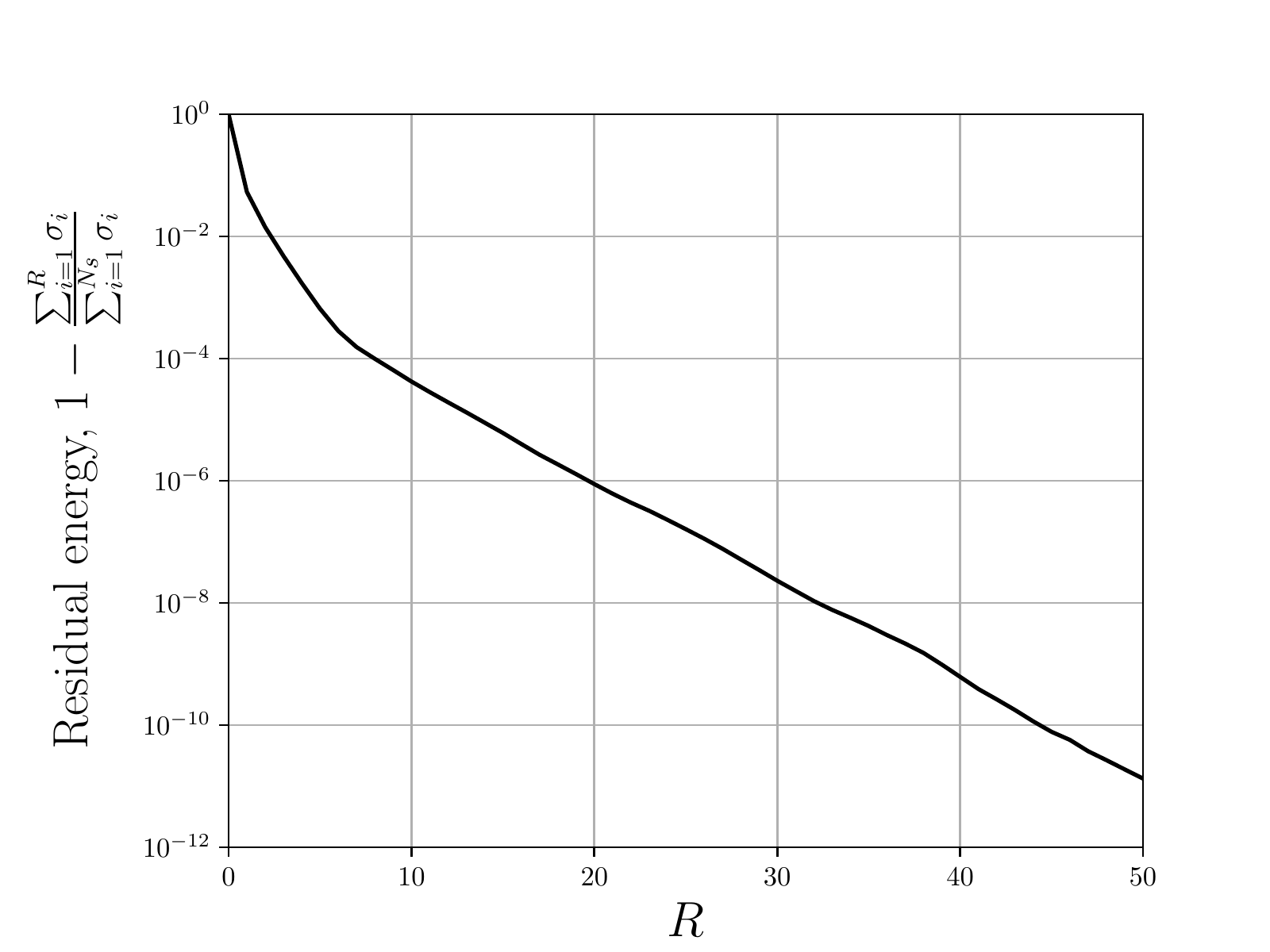}
\caption{\advectingFrontFigureCaption Residual statistical energy as a function of basis dimension.}
\label{fig:rom_case2_svd}
\end{center}
\end{figure}

We now describe the problem setup. We solve Eq.~\eqref{eq:cdr} with a final time $T=2$ 
and a spatial domain $\Omega = (0,1) \times (0,1)$. We take $\viscosity = 10^{-4}$, $\reaction = 1$, and $\wavespeed = \frac{1}{2} \begin{bmatrix} \cos( \pi / 3) & \sin(\pi / 3) \end{bmatrix}^T$. The Peclet number is $\text{Pe}:=\norm{ \wavespeed} /\viscosity = 5000$. The forcing is set as
$$\forcing = \begin{cases} 1 & 0 \le x \le 0.5 \;\; \text{and}\;\; 0 \le y \le 0.25 \\ 0 & x > 0.5 \; \; \text{and} \;\; y > 0.25. \end{cases}$$
The high-resolution trial space $\cSpaceHRES$ is obtained via a uniform triangulation of $\Omega$ into $N_{\text{el}} = 2 \times 256^2$ elements equipped with a $\mathcal{C}^0(\cDomain)$ continuous discretization with polynomials of order $p=2$. The grid Peclet number is $\text{Pe}_{\text{g}} = 9.77$, where we used $h=(256 p)^{-1}$. Analogously, the FOM trial space $\cSpaceFOM$ is obtained via a uniform triangulation of $\Omega$ into $N_{\text{el}} = 2 \times 32^2$ elements equipped with a $\mathcal{C}^0(\cDomain)$ continuous discretization with polynomials of order $p=2$. The grid Peclet number is $\text{Pe}_{\text{g}} =  78.125$, where we used $h=(32 p)^{-1}$. The triangulations are obtained in the same manner as in the previous experiment. The ROM trial space is obtained by executing Algorithm~\ref{alg:gen_rom_trialspace}; Figure~\ref{fig:rom_case2_svd} presents the residual statistical energy as a function of basis dimension. As compared to the previous numerical example, here it is seen that more basis vectors are required to characterize the system.

\subsubsection{Full-order model results}
\begin{figure}
\begin{center}
\begin{subfigure}[t]{0.32\textwidth}
\includegraphics[trim={1cm 1cm 1cm 1cm},clip, width=1.\linewidth]{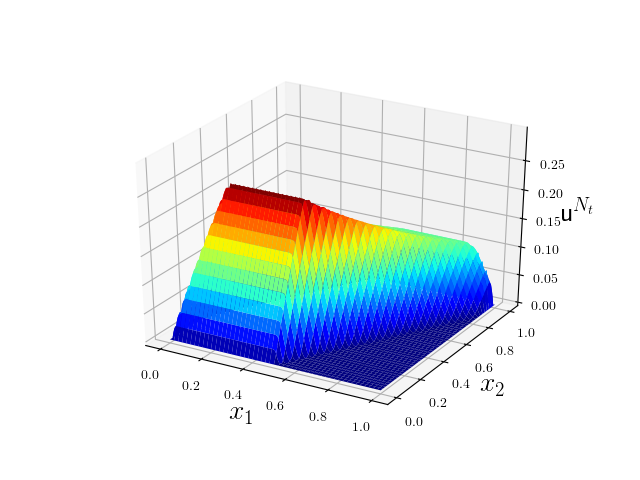}
\caption{H-Res}
\end{subfigure}
\begin{subfigure}[t]{0.32\textwidth}
\includegraphics[trim={1cm 1cm 1cm 1cm},clip, width=1.\linewidth]{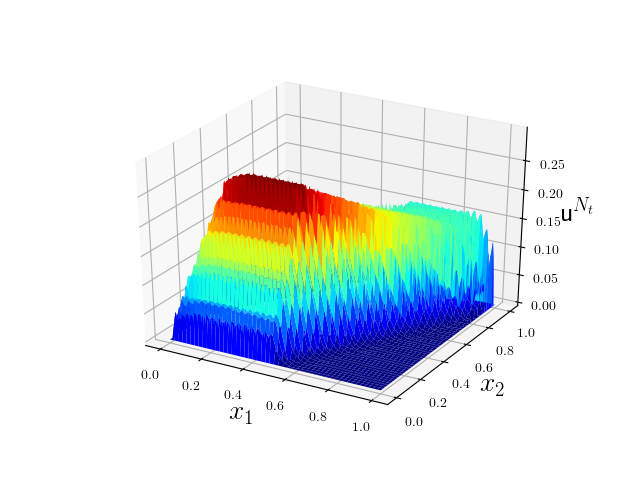}
\caption{Galerkin}
\end{subfigure}
\begin{subfigure}[t]{0.32\textwidth}
\includegraphics[trim={1cm 1cm 1cm 1cm},clip, width=1.\linewidth]{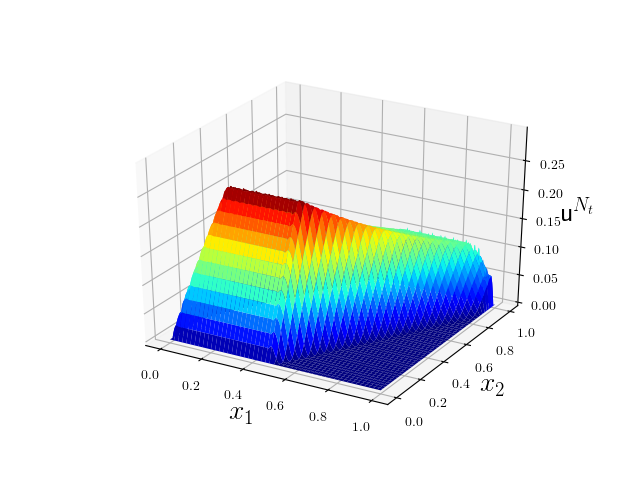}
\caption{\SUPG}
\end{subfigure}
\begin{subfigure}[t]{0.32\textwidth}
\includegraphics[trim={1cm 1cm 1cm 1cm},clip, width=1.\linewidth]{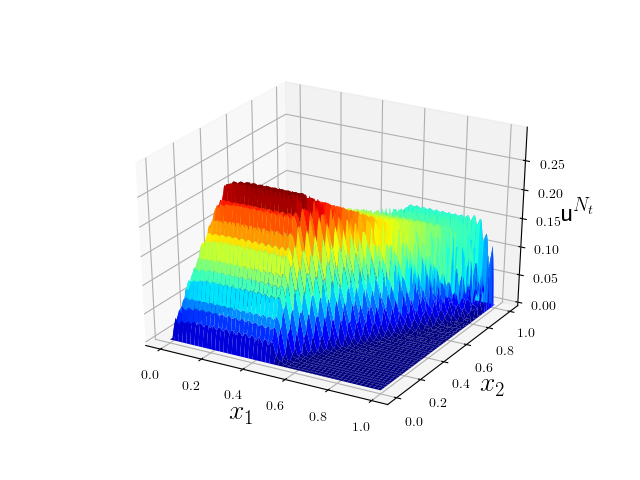}
\caption{\GLS}
\end{subfigure}
\begin{subfigure}[t]{0.32\textwidth}
\includegraphics[trim={1cm 1cm 1cm 1cm},clip,  width=1.\linewidth]{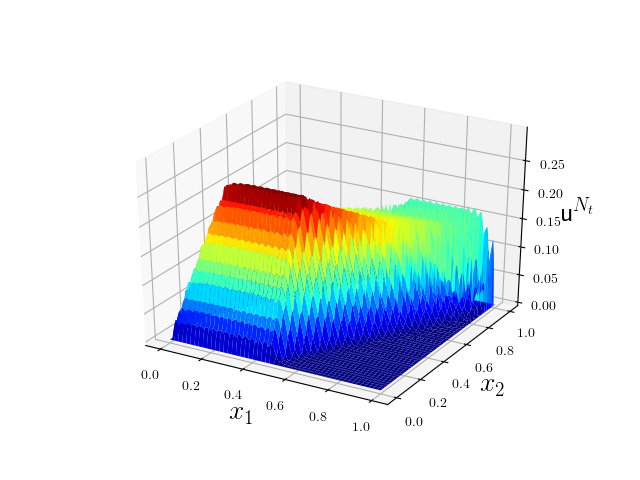}
\caption{\ADJ}
\end{subfigure}
\begin{subfigure}[t]{0.32\textwidth}
\includegraphics[trim={1cm 1cm 1cm 1cm},clip,  width=1.\linewidth]{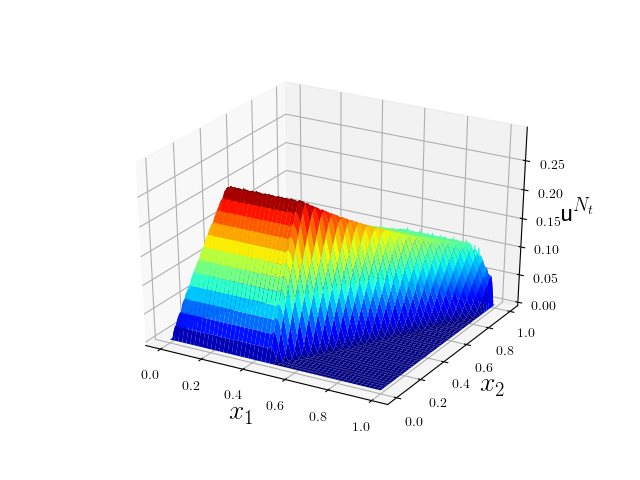}
\caption{\STGLS}
\end{subfigure}
\begin{subfigure}[t]{0.32\textwidth}
\includegraphics[trim={1cm 1cm 1cm 1cm},clip,  width=1.\linewidth]{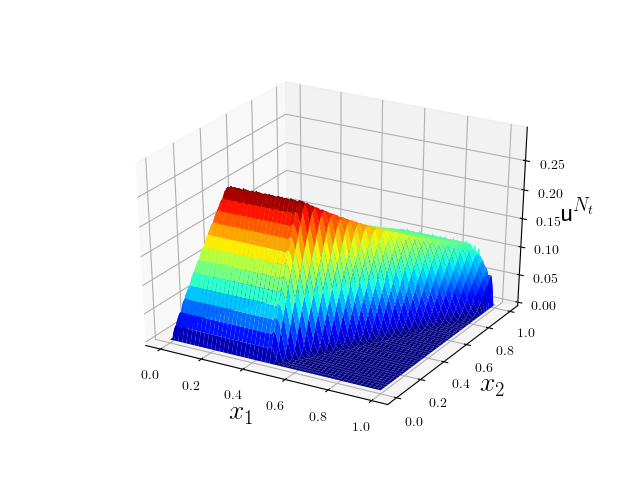}
\caption{\STADJ}
\end{subfigure}

\caption{
\advectingFrontFigureCaption FOM solutions to the \cdcdrAcronym\ equation at $t=5$.} 
\label{fig:ex2_fom_figs}
\end{center}
\end{figure}

We again first present results of the various FOMs considered (again for optimal stabilization parameters and time steps). Figure~\ref{fig:ex2_fom_figs} presents the various FOM solutions at the final time, $t=2$, while Table~\ref{tab:ex2_fom_errors} tabulates the solution errors and stabilization parameters employed in the simulations. We observe the following:
\begin{itemize}
\item The Galerkin, \ADJ, and \GLS\ FEM FOMs are the worst performing methods and all result in oscillatory solutions. 
\item \SUPG, \STGLS, and \STADJ\ all provide solutions of a similar qualitative and quantitative quality, and are all able to suppress the oscillations seen in the Galerkin FEM solution.  
\item The ``space--time" \STGLS\ and \STADJ\ FEMs again perform much better than the ``discretize-then-stabilize" \GLS\ and \ADJ\ FEM FOMs.
\end{itemize} 
These results are similar to those obtained for example 1 (Section~\ref{sec:example1-fom}).

\begin{table}
  \begin{center}
    \caption{\advectingFrontFigureCaption Integrated (relative) $\LTwo$ and $\HOne$ errors of various FOMs presented in Figure~\ref{fig:ex2_fom_figs}, along with the stabilization parameters at which the FOMs were executed.}
\label{tab:ex2_fom_errors}
    \begin{filecontents*}{figs/example2/results_fem_foms.csv}
Codes, G, SUPG, GLS, ADJ, STGLS, STADJ
0.000000000000000000e+00,2.759173028166726771e-02,6.620296672187269728e-03,2.239904833384482805e-02,2.192068632024630689e-02,6.657746643360291233e-03,6.587166969012498048e-03
1.000000000000000000e+00,9.300839286897473468e-02,5.371882630527806748e-03,7.530943471626841412e-02,7.319738029699919546e-02,5.806980897480060536e-03,5.020729023394916030e-03
2.000000000000000000e+00,nan,1.000000000000000021e-02,1.000000000000000048e-04,1.000000000000000048e-04,1.000000000000000021e-02,1.000000000000000021e-02
3.000000000000000000e+00,1.000000000000000021e-03,1.000000000000000021e-03,1.000000000000000021e-03,1.000000000000000021e-03,1.000000000000000021e-03,1.000000000000000021e-03
\end{filecontents*}
\pgfplotstabletypeset[
      multicolumn names, 
      col sep=comma, 
      string replace*={0.000000000000000000e+00}{$\integratedErrorLTwoBestFit$},
      string replace*={1.000000000000000000e+00}{$\integratedErrorHOneBestFit$},
      string replace*={2.000000000000000000e+00}{$\tau$},
      string replace*={3.000000000000000000e+00}{$\Delta t$},
      string replace*={errorH}{$\integratedErrorHOneBestFit$},
      string replace*={nan}{N/A},
      display columns/0/.style={
		column name=,string type},  
      display columns/1/.style={
		column name=Galerkin, 
		column type={S[table-parse-only]},string type},  
     display columns/2/.style={
		column name={\begin{centering} \SUPG \end{centering}},
		column type={S[table-parse-only]},string type},
      display columns/3/.style={
		column name=\GLS,
		column type={S[table-parse-only]},string type},
     display columns/4/.style={
		column name=\ADJ,
		column type={S[table-parse-only]},string type},
     display columns/5/.style={
		column name=\STGLS,
		column type={S[table-parse-only]},string type},
     display columns/6/.style={
		column name=\STADJ,
		column type={S[table-parse-only]},string type},
    every head row/.style={
		before row={\toprule}, 
		after row={
			\midrule} 
			},
		every last row/.style={after row=\bottomrule}, 
    ]{figs/example2/results_fem_foms.csv} 
  \end{center}
\end{table}

\subsubsection{Results as a function of RB dimension}
\begin{figure}
\begin{center}
\begin{subfigure}[t]{0.49\textwidth}
\includegraphics[trim={1.5cm 8.8cm 0cm 0cm},clip, width=1.\linewidth]{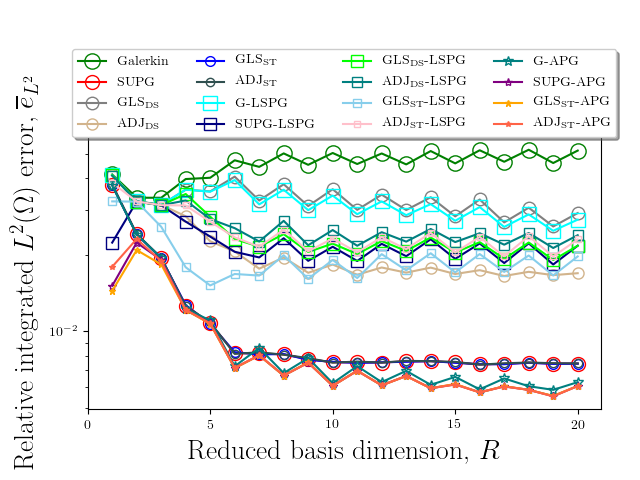}
\end{subfigure}

\begin{subfigure}[t]{0.49\textwidth}
\includegraphics[width=1.\linewidth]{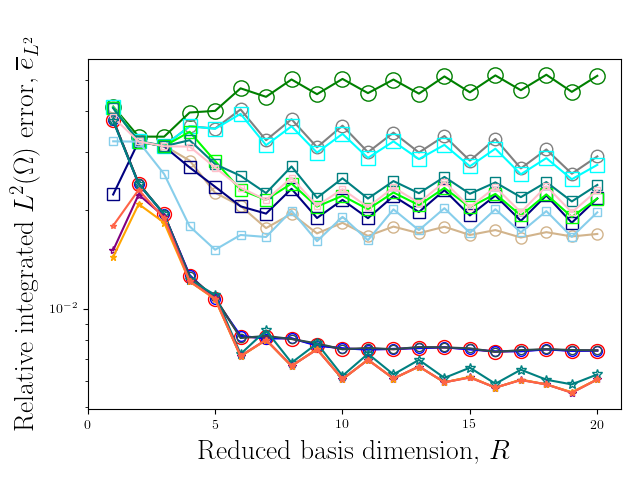}
\caption{$L^2(\cDomain)$ error.}
\end{subfigure}
\begin{subfigure}[t]{0.49\textwidth}
\includegraphics[width=1.\linewidth]{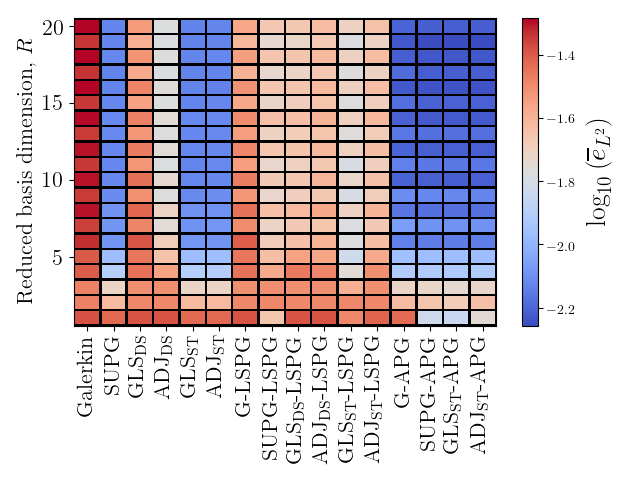}
\caption{$L^2(\cDomain)$ error.}
\end{subfigure}

\begin{subfigure}[t]{0.49\textwidth}
\includegraphics[width=1.\linewidth]{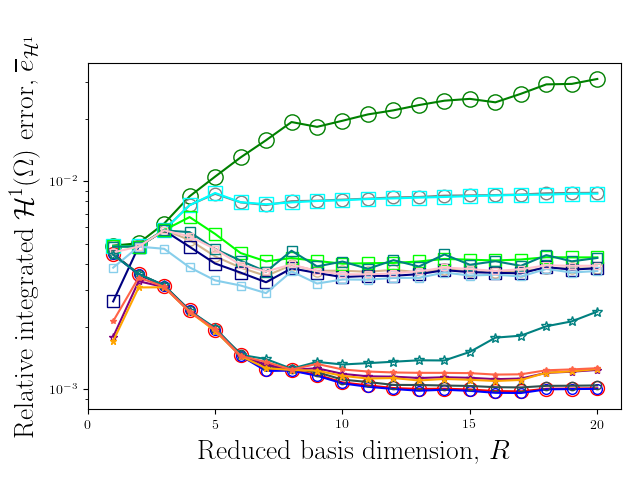}
\caption{$\HOne$ error.}
\end{subfigure}
\begin{subfigure}[t]{0.49\textwidth}
\includegraphics[width=1.\linewidth]{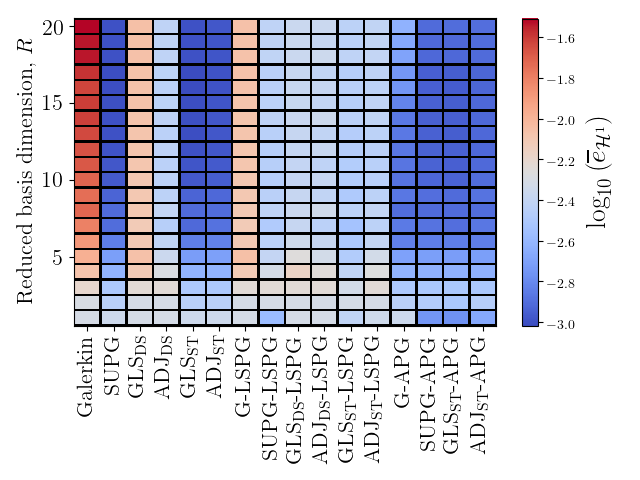}
\caption{$\HOne$ error.}
\end{subfigure}

\caption{\advectingFrontFigureCaption $\LTwo$ (top) and $\HOne$ (bottom) error as a function of RB dimension for the various ROMs evaluated. We note that the left and right figures show the same quantities, but with different visualization techniques. Results are shown for optimal values of $t,\tau$ as discussed in Section~\ref{sec:timestep_information}.}
\label{fig:piecewise_forcing_romconverge}
\end{center}
\end{figure}

\begin{figure}
\begin{center}

\begin{subfigure}[t]{0.49\textwidth}
\includegraphics[trim={1.5cm 8.8cm 0cm 0cm},clip, width=1.\linewidth]{figs/example2/ex1ConvergenceLegend.png}
\end{subfigure}

\begin{subfigure}[t]{0.49\textwidth}
\includegraphics[width=1.\linewidth]{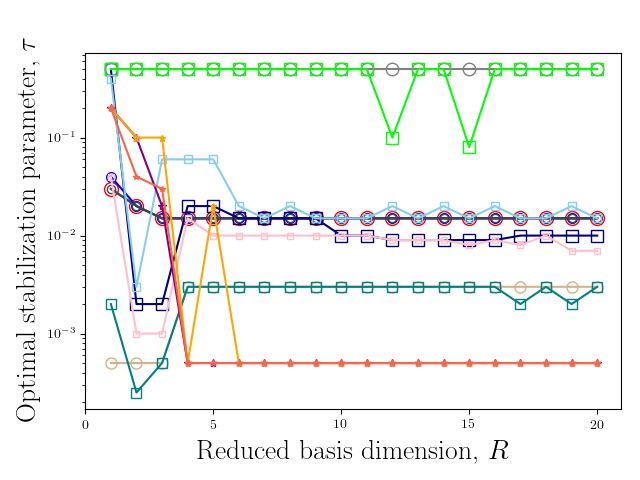}
\caption{Optimal parameter, $\tau$.}
\label{fig:piecewise_forcing_romconverge_c}
\end{subfigure}
\begin{subfigure}[t]{0.49\textwidth}
\includegraphics[width=1.\linewidth]{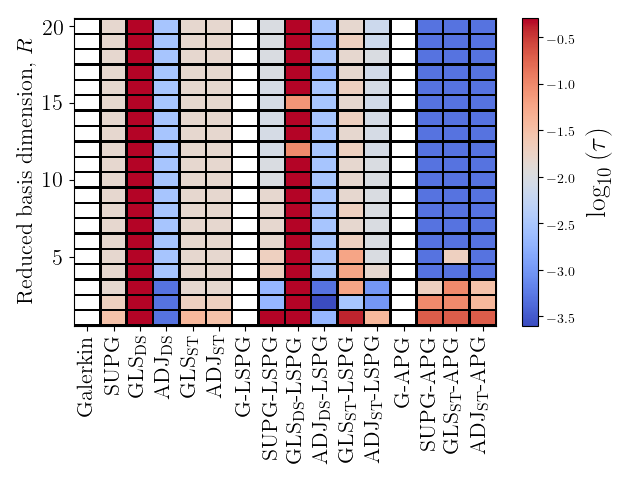}
\caption{Optimal parameter, $\tau$.}
\end{subfigure}

\begin{subfigure}[t]{0.49\textwidth}
\includegraphics[width=1.\linewidth]{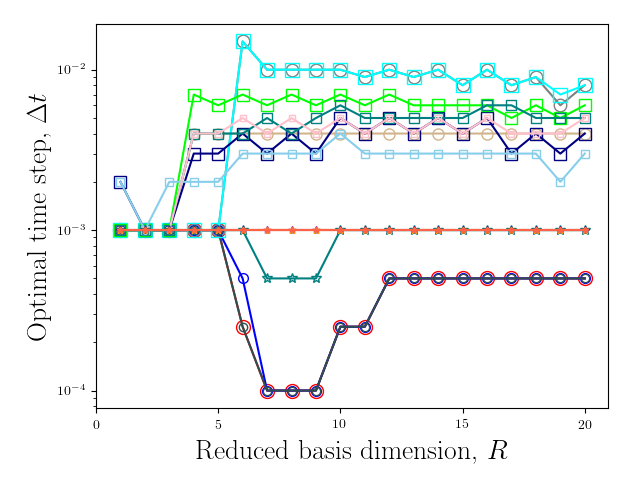}
\caption{Optimal time step, $\Delta t$.}
\label{fig:piecewise_forcing_romconverge_d}
\end{subfigure}
\begin{subfigure}[t]{0.49\textwidth}
\includegraphics[width=1.\linewidth]{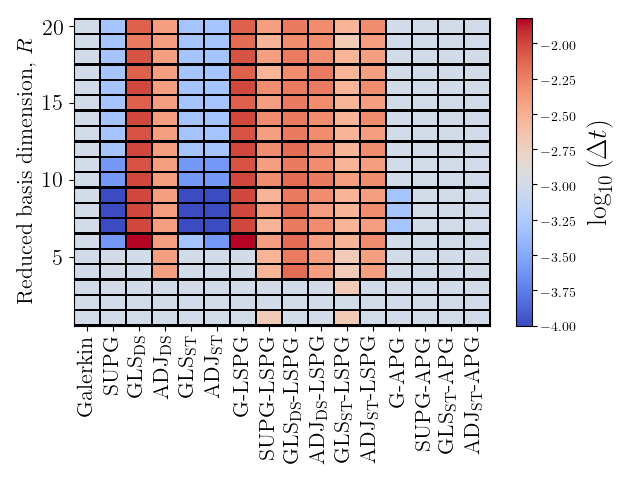}
\caption{Optimal time step, $\Delta t$.}
\label{fig:piecewise_forcing_romconverge_optimaldt_pcolor}
\end{subfigure}
\caption{\advectingFrontFigureCaption Optimal stabilization parameter (top) and time step (bottom) as a function of RB dimension. We note that the left and right figures show the same quantities, but with different visualization techniques. 
Results are shown for optimal values of $t,\tau$ as discussed in Section~\ref{sec:timestep_information}.} 
\label{fig:piecewise_forcing_romconverge_params}
\end{center}
\end{figure}

Figure~\ref{fig:piecewise_forcing_romconverge} shows the $\LTwo$ error and $\HOne$ error for each of the ROMs considered as a function of RB size.
Results are presented for the value of $\tau$ and $\Delta t$ that led to the lowest $\LTwo$ error;
Figure~\ref{fig:piecewise_forcing_romconverge_params} shows these  
optimal stabilization parameters and
time steps for each basis dimension. 

Examining Figure~\ref{fig:piecewise_forcing_romconverge}, we make the following observations about the \textit{accuracy} of the various ROMs: 
\begin{itemize}
\item The APG-based, \SUPG, \STADJ, and \STGLS\ ROMs are the best overall performing methods. 
\item The ``discretize-then-stabilize" ROMs outperform the Galerkin ROM, but are consistently worse than their space--time stabilized counterparts for all ROM dimensions. In general, even when equipped with LSPG, the discretize-then-stabilize methods perform quite poorly. 
\item When applied to the standard Galerkin and \SDGLS\ FEM models, LSPG leads to slightly improved solutions. LSPG does not lead to improved solutions for \SUPG, \STGLS, or \STADJ. 
\item It is again interesting to note the improved performance of \APGG\ over the standard Galerkin method at high reduced basis dimensions, given that \APGG\ will converge to the Galerkin method in the limit of a full reduced basis. The improved performance is also seen for \LSPGG, but to a lesser extent.
\item Once again, a decrease in error in the $\LTwo$ norm does not always correspond to a decrease in error in the $\HOne$ norm, and vice versa.
\item 
Comparing Figure~\ref{fig:piecewise_forcing_romconverge} to Table~\ref{tab:ex2_fom_errors}, we again observe that some ROMs are more accurate than their corresponding FEM.
\end{itemize}

Examining Figure~\ref{fig:piecewise_forcing_romconverge_params}, where we show the optimal stabilization parameters and time steps associated with Figure~\ref{fig:piecewise_forcing_romconverge}, we make the following observations about the behavior of the stabilization parameters of the various ROMs: 
\begin{itemize}
\item The optimal stabilization parameters for the \SUPG, \STADJ, \STGLS, and APG ROMs increase for $R < 4$. This result is consistent with analyses performed for SUPG in~\cite{GIERE2015454} and APG in~\cite{parish_apg}, which suggested that the optimal stabilization parameter decreases with increasing ROM dimension.
\item The optimal stabilization parameter for \SDGLS\ is quite high for all ROM dimensions. 
\item The optimal time step for all LSPG ROMs occurs at an intermediate time step larger than the FEM FOM.
\item Interestingly, the optimal time step for \LSPGG\ is the same as the optimal time step for \GLS. 
\item The reader may observe that the optimal time step for \SUPG, \STGLS, and \SDADJ\ decreases for moderately high ROM dimensions; we note that this improvement is very minor as can will be seen in Figure~\ref{fig:rom_piecewise_forcing_taudtvary}.
\end{itemize}


\begin{figure}
\begin{center}
\begin{subfigure}[t]{0.28\textwidth}
\includegraphics[trim={1cm 1cm 1cm 1cm},clip, width=1.\linewidth]{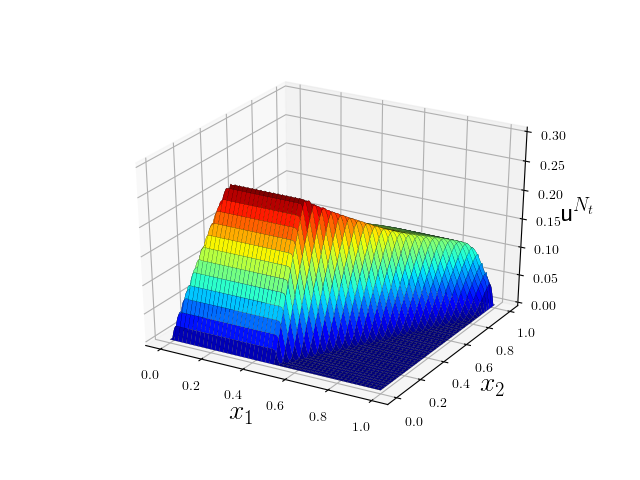}
\caption{H-Res ($L^2(\cDomain)$ best fit)}
\label{fig:rom_piecwise_forcing_b}
\end{subfigure}
\begin{subfigure}[t]{0.28\textwidth}
\includegraphics[trim={1cm 1cm 1cm 1cm},clip, width=1.\linewidth]{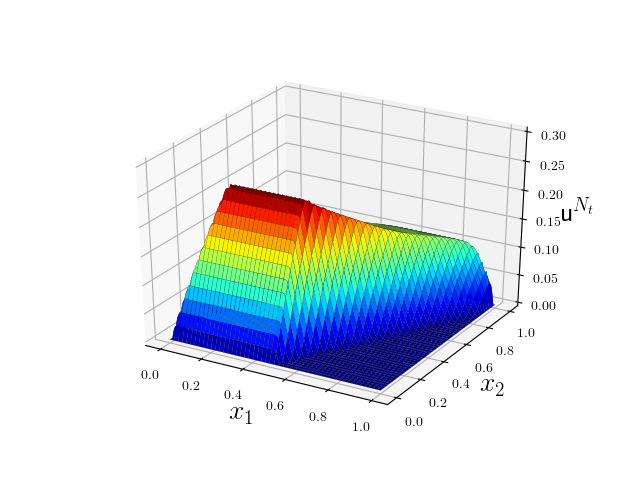}
\caption{H-Res ($H_1(\cDomain)$ best fit)}
\label{fig:rom_piecewise_forcing_c}
\end{subfigure}
\begin{subfigure}[t]{0.28\textwidth}
\includegraphics[trim={1cm 1cm 1cm 1cm},clip,  width=1.\linewidth]{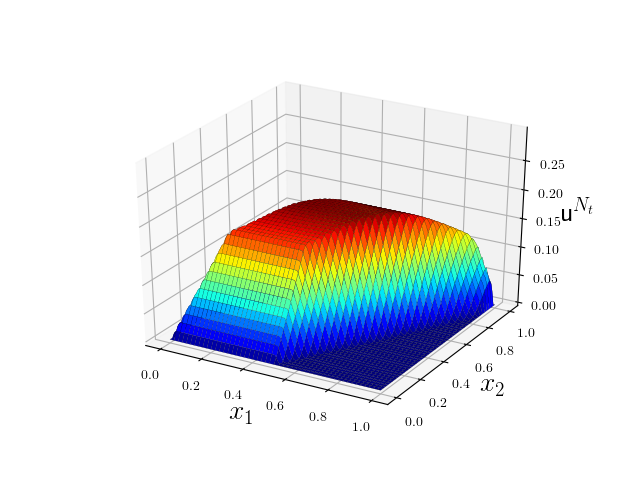}
\caption{Galerkin}
\label{fig:rom_piecewise_forcing_d}
\end{subfigure}
\begin{subfigure}[t]{0.28\textwidth}
\includegraphics[trim={1cm 1cm 1cm 1cm},clip,  width=1.\linewidth]{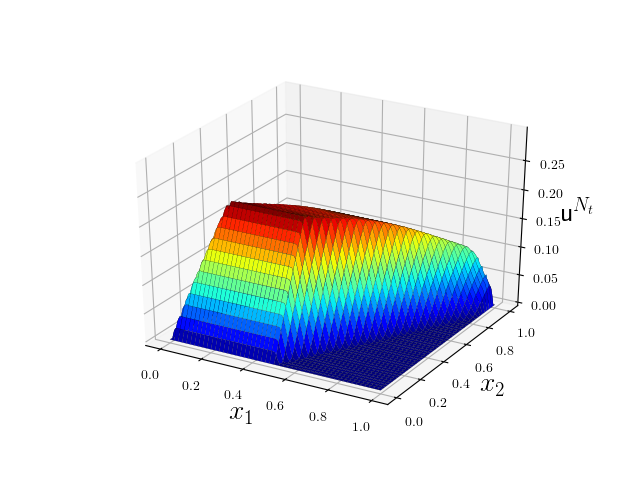}
\caption{\SUPG}
\label{fig:rom_piecewise_forcing_e}
\end{subfigure}
\begin{subfigure}[t]{0.28\textwidth}
\includegraphics[trim={1cm 1cm 1cm 1cm},clip,  width=1.\linewidth]{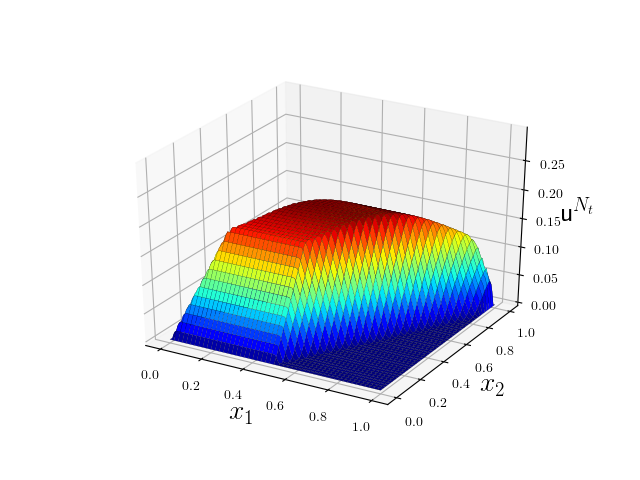}
\caption{\GLS}
\label{fig:rom_piecewise_forcing_f}
\end{subfigure}
\begin{subfigure}[t]{0.28\textwidth}
\includegraphics[trim={1cm 1cm 1cm 1cm},clip,  width=1.\linewidth]{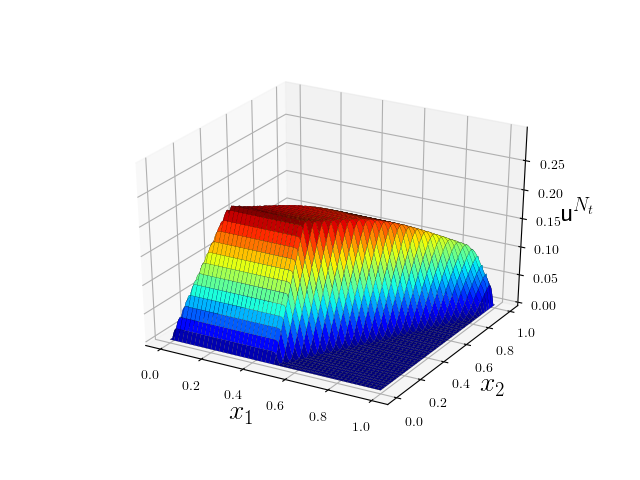}
\caption{\ADJ}
\label{fig:rom_piecewise_forcing_g}
\end{subfigure}
\begin{subfigure}[t]{0.28\textwidth}
\includegraphics[trim={1cm 1cm 1cm 1cm},clip,  width=1.\linewidth]{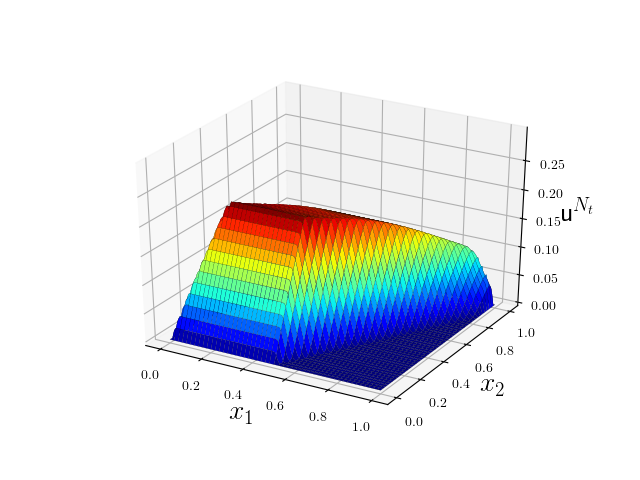}
\caption{\STGLS}
\label{fig:rom_piecewise_forcing_f}
\end{subfigure}
\begin{subfigure}[t]{0.28\textwidth}
\includegraphics[trim={1cm 1cm 1cm 1cm},clip,  width=1.\linewidth]{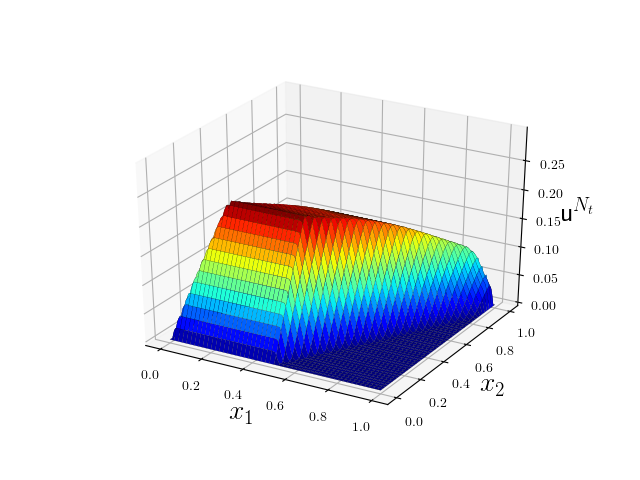}
\caption{\STADJ}
\label{fig:rom_piecewise_forcing_g}
\end{subfigure}
\begin{subfigure}[t]{0.28\textwidth}
\includegraphics[ trim={1cm 1cm 1cm 1cm},clip, width=1.\linewidth]{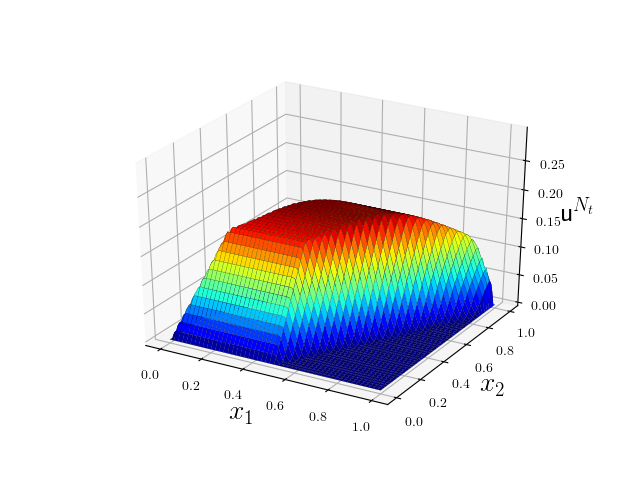}
\caption{\LSPGG}
\label{fig:rom_piecewise_forcing_h}
\end{subfigure}
\begin{subfigure}[t]{0.28\textwidth}
\includegraphics[ trim={1cm 1cm 1cm 1cm},clip, width=1.\linewidth]{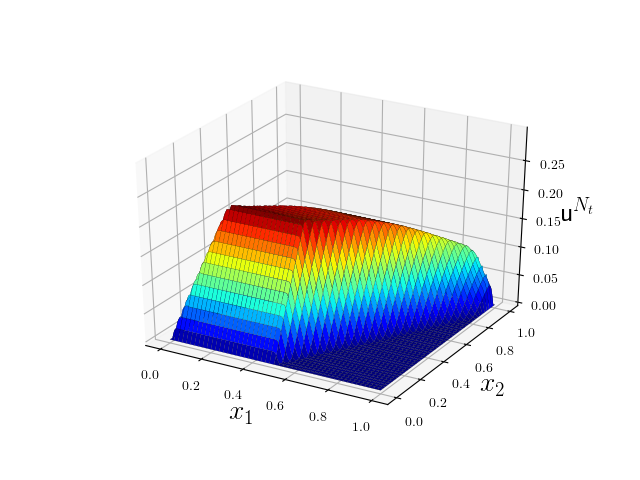}
\caption{\LSPGSUPG}
\label{fig:rom_piecewise_forcing_h}
\end{subfigure}
\begin{subfigure}[t]{0.28\textwidth}
\includegraphics[ trim={1cm 1cm 1cm 1cm},clip, width=1.\linewidth]{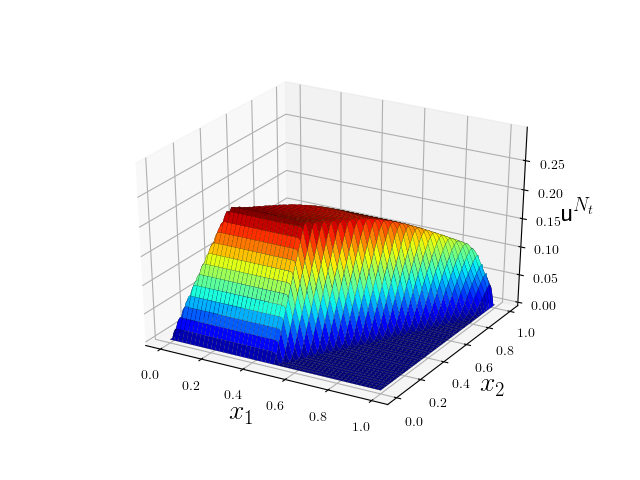}
\caption{\LSPGGLS}
\label{fig:rom_piecewise_forcing_h}
\end{subfigure}
\begin{subfigure}[t]{0.28\textwidth}
\includegraphics[ trim={1cm 1cm 1cm 1cm},clip, width=1.\linewidth]{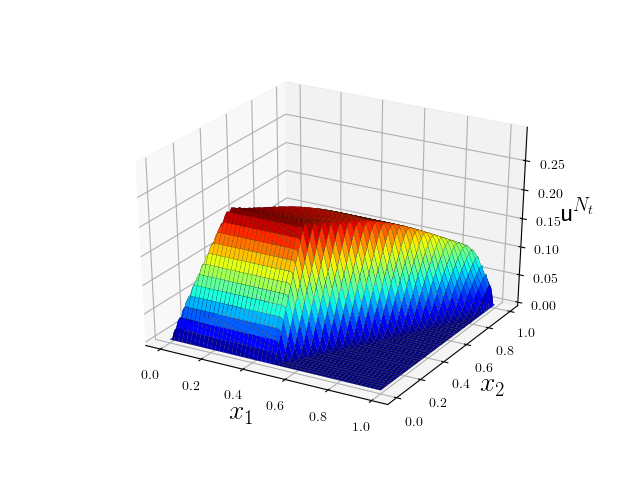}
\caption{\LSPGADJ}
\label{fig:rom_piecewise_forcing_h}
\end{subfigure}
\begin{subfigure}[t]{0.28\textwidth}
\includegraphics[ trim={1cm 1cm 1cm 1cm},clip, width=1.\linewidth]{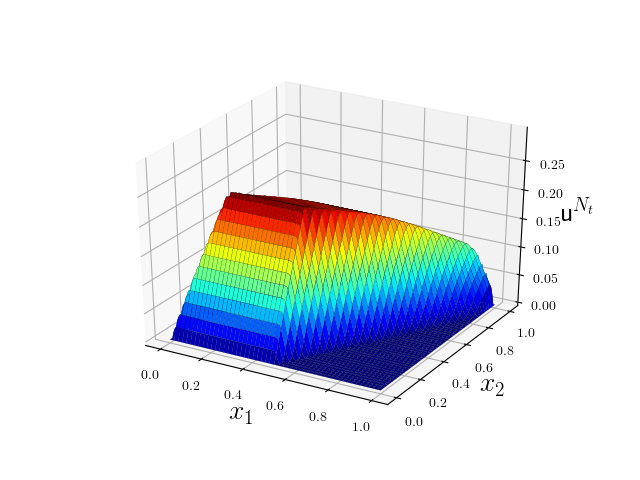}
\caption{\LSPGSTGLS}
\label{fig:rom_piecewise_forcing_h}
\end{subfigure}
\begin{subfigure}[t]{0.28\textwidth}
\includegraphics[ trim={1cm 1cm 1cm 1cm},clip, width=1.\linewidth]{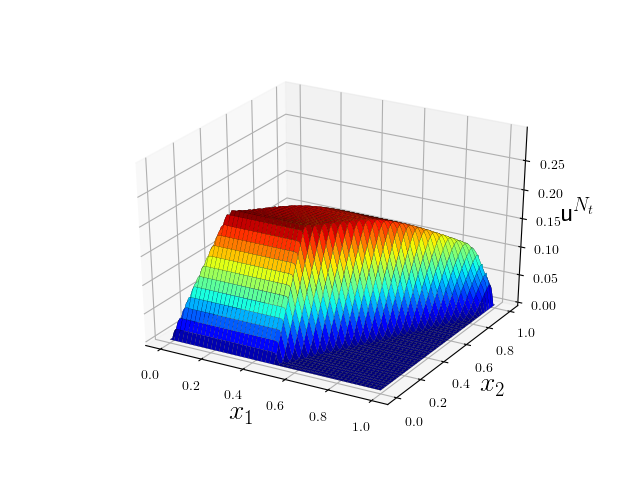}
\caption{\LSPGSTADJ}
\label{fig:rom_piecewise_forcing_h}
\end{subfigure}
\begin{subfigure}[t]{0.28\textwidth}
\includegraphics[ trim={1cm 1cm 1cm 1cm},clip, width=1.\linewidth]{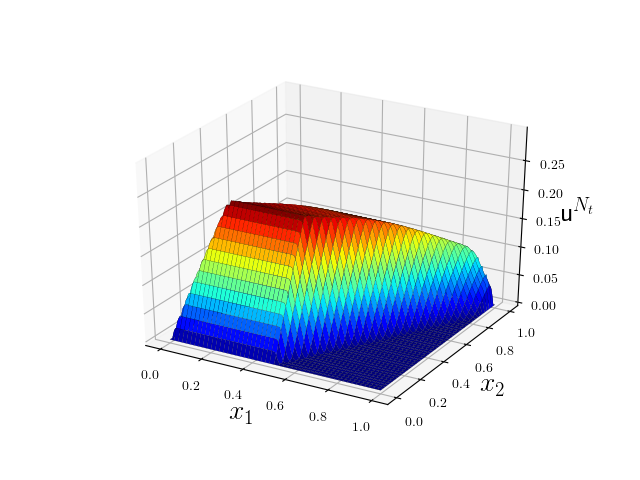}
\caption{\APGG}
\label{fig:rom_piecewise_forcing_i}
\end{subfigure}
\begin{subfigure}[t]{0.28\textwidth}
\includegraphics[ trim={1cm 1cm 1cm 1cm},clip, width=1.\linewidth]{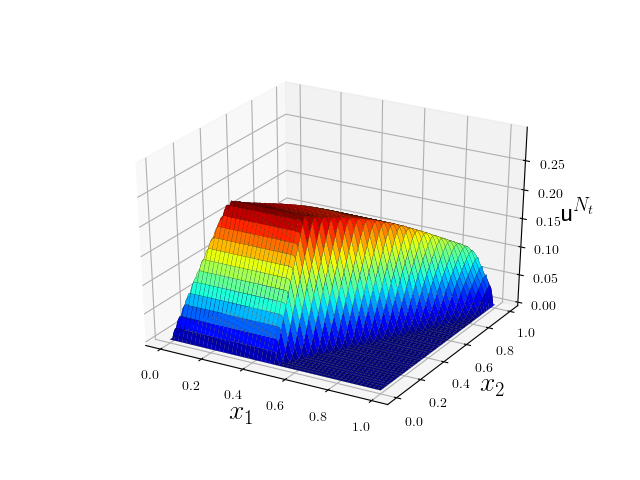}
\caption{\APGSUPG}
\label{fig:rom_piecewise_forcing_i}
\end{subfigure}
\begin{subfigure}[t]{0.28\textwidth}
\includegraphics[ trim={1cm 1cm 1cm 1cm},clip, width=1.\linewidth]{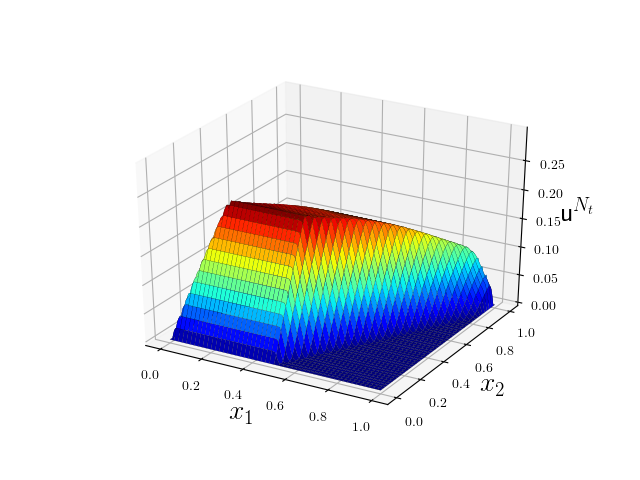}
\caption{\APGSTGLS}
\label{fig:rom_piecewise_forcing_i}
\end{subfigure}
\begin{subfigure}[t]{0.28\textwidth}
\includegraphics[ trim={1cm 1cm 1cm 1cm},clip, width=1.\linewidth]{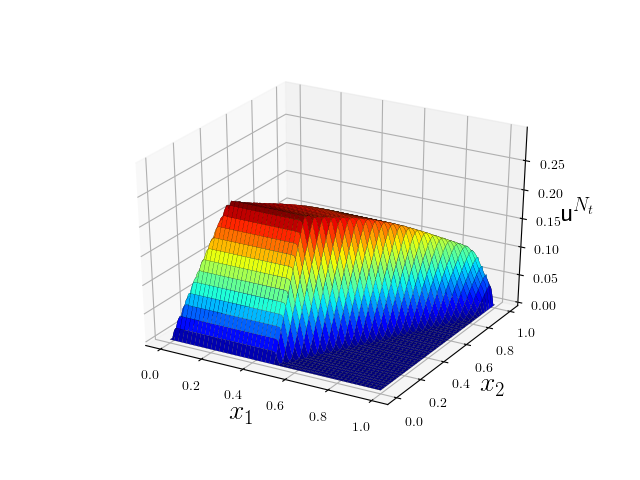}
\caption{\APGSTADJ}
\label{fig:rom_piecewise_forcing_i}
\end{subfigure}
\caption{\advectingFrontFigureCaption ROM solutions to the \cdcdrAcronym\ equation at $t=2.0$. Results are shown for $\romdim=5$.}
\label{fig:rom_piecewise_forcing}
\end{center}
\end{figure}

Next, Figure~\ref{fig:rom_piecewise_forcing} presents physical space solution profiles for the various ROMs for a reduced basis dimension of $R = 5$ at the final time instance, $t=2.0$. 
We observe that all methods yield qualitatively accurate solutions with the exception of Galerkin, \LSPGG, and \SDGLS; these three methods under-predict the magnitude of the 
solution in the lower-left quadrant of the domain ($x_1,x_2 \le 0.4$). It is interesting to note that, although Figure~\ref{fig:piecewise_forcing_romconverge} showed that \SDADJ\ and the non-Galerkin LSPG ROMs clearly perform less well than the 
other formulations, Figure~\ref{fig:rom_piecewise_forcing} shows that their physical space solutions still show a significant improvement over the standard Galerkin ROM. 

\subsubsection{Sensitivity to time step and stabilization parameters}
\begin{figure}
\begin{center}
\begin{subfigure}[t]{0.32\textwidth}
\includegraphics[width=1.\linewidth]{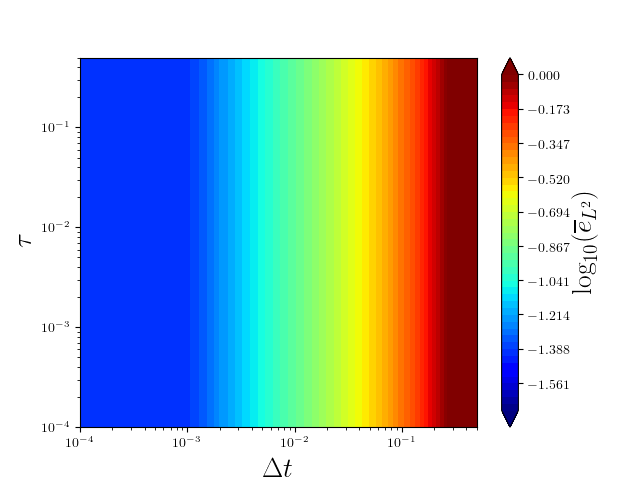}
\caption{Galerkin}
\label{fig:rom_piecewise_forcing_taudtvary_g}
\end{subfigure}
\begin{subfigure}[t]{0.32\textwidth}
\includegraphics[width=1.\linewidth]{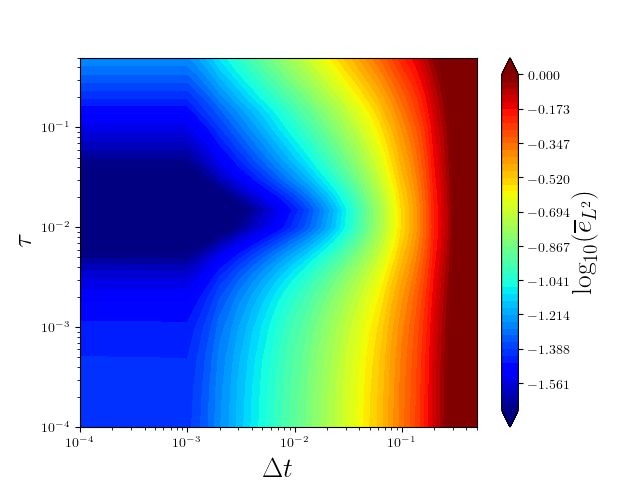}
\caption{\SUPG}
\label{fig:rom_piecewise_forcing_taudtvary_supg}
\end{subfigure}
\begin{subfigure}[t]{0.32\textwidth}
\includegraphics[width=1.\linewidth]{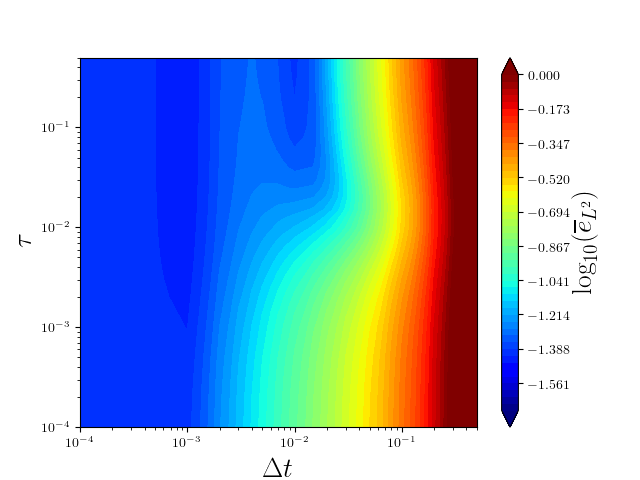}
\caption{\GLS}
\label{fig:rom_piecewise_forcing_taudtvary_gls}
\end{subfigure}
\begin{subfigure}[t]{0.32\textwidth}
\includegraphics[width=1.\linewidth]{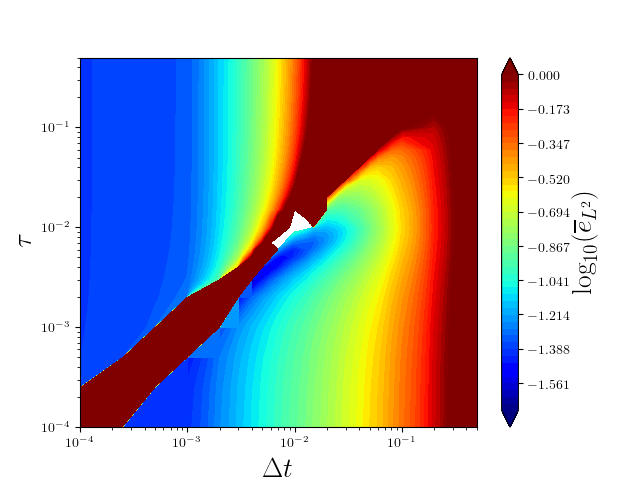}
\caption{\ADJ}
\label{fig:rom_fig1_taudtvary_adj}
\end{subfigure}
\begin{subfigure}[t]{0.32\textwidth}
\includegraphics[width=1.\linewidth]{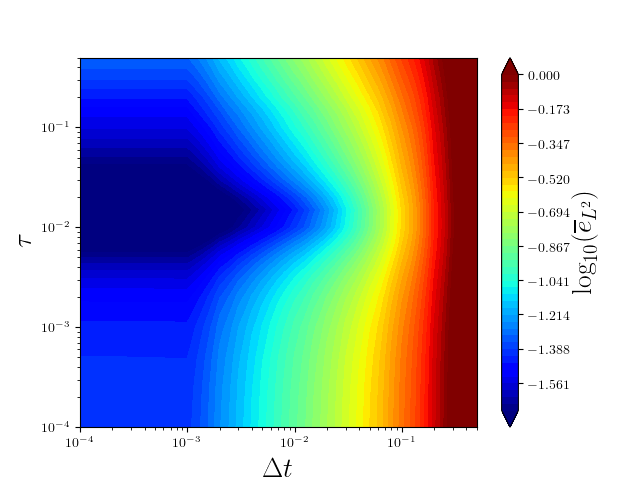}
\caption{\STGLS}
\label{fig:rom_piecewise_forcing_taudtvary_stgls}
\end{subfigure}
\begin{subfigure}[t]{0.32\textwidth}
\includegraphics[width=1.\linewidth]{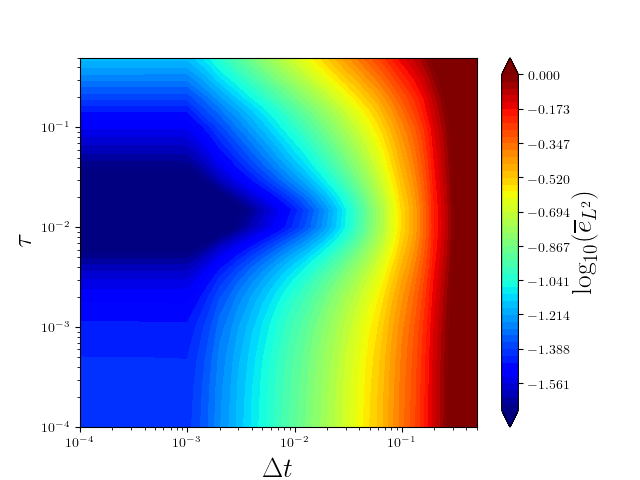}
\caption{\STADJ}
\label{fig:rom_piecewise_forcing_taudtvary_stadj}
\end{subfigure}
\begin{subfigure}[t]{0.32\textwidth}
\includegraphics[width=1.\linewidth]{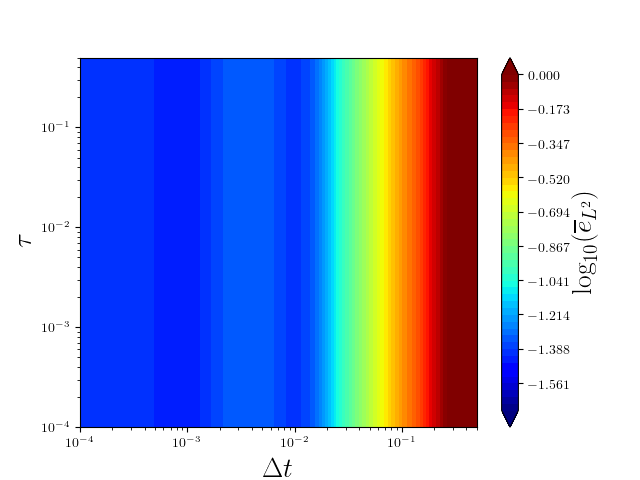}
\caption{\LSPGG}
\label{fig:rom_piecewise_forcing_taudtvary_lspg}
\end{subfigure}
\begin{subfigure}[t]{0.32\textwidth}
\includegraphics[width=1.\linewidth]{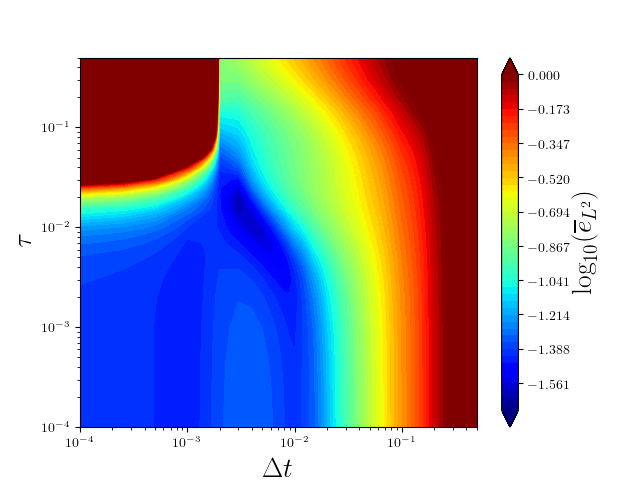}
\caption{\LSPGSUPG}
\end{subfigure}
\begin{subfigure}[t]{0.32\textwidth}
\includegraphics[width=1.\linewidth]{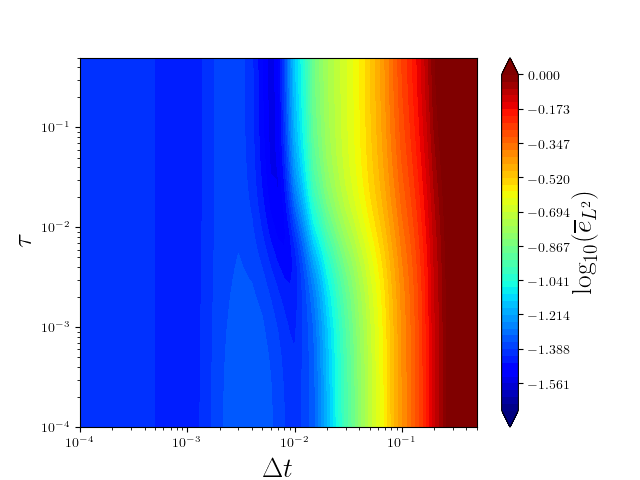}
\caption{\LSPGGLS}
\end{subfigure}
\begin{subfigure}[t]{0.32\textwidth}
\includegraphics[width=1.\linewidth]{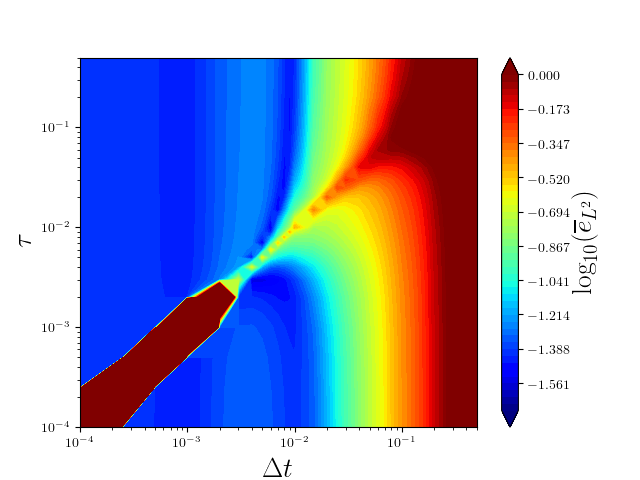}
\caption{\LSPGADJ}
\end{subfigure}
\begin{subfigure}[t]{0.32\textwidth}
\includegraphics[width=1.\linewidth]{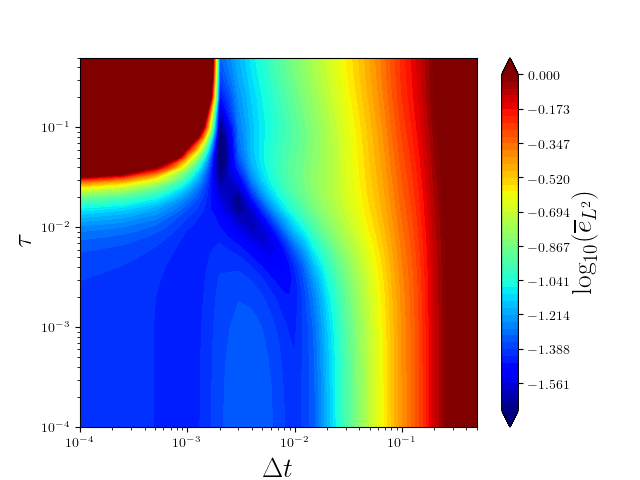}
\caption{\STGLSLSPG}
\end{subfigure}
\begin{subfigure}[t]{0.32\textwidth}
\includegraphics[width=1.\linewidth]{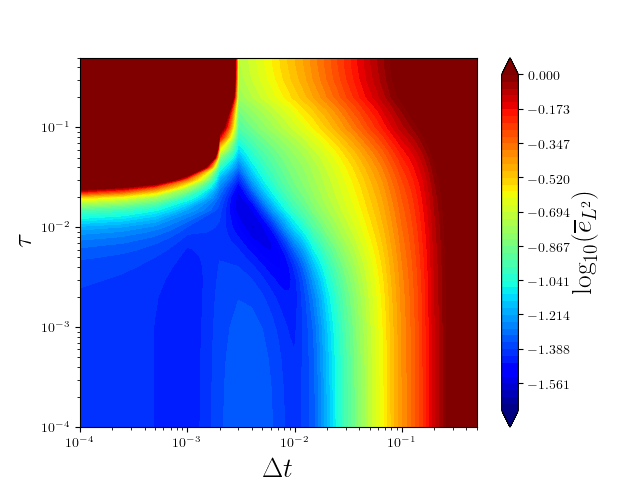}
\caption{\STADJLSPG}
\end{subfigure}
\caption{\advectingFrontFigureCaption Time integrated $\LTwo$ error as a function of time step and stabilization parameter for the various ROMs evaluated. Note that Galerkin and LSPG display no dependence on the stabilization parameter. White regions indicate solutions that diverged to NaN.}
\label{fig:rom_piecewise_forcing_taudtvary}
\end{center}
\end{figure}

\begin{figure}
\begin{center}
\begin{subfigure}[t]{0.32\textwidth}
\includegraphics[width=1.\linewidth]{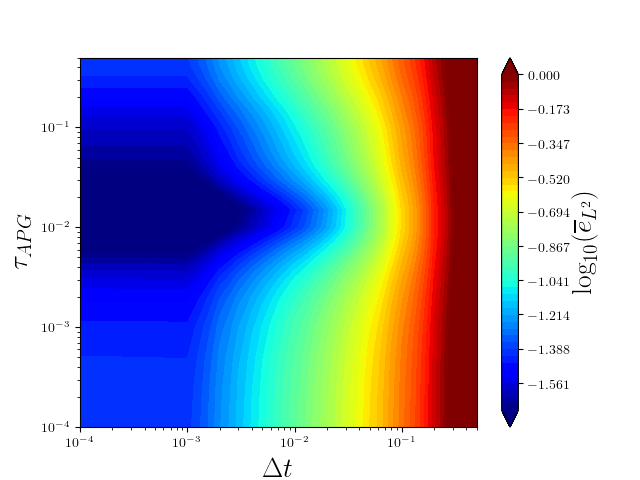}
\caption{\APGG}
\label{fig:rom_piecewise_forcing_taudtvary_apgg}
\end{subfigure}
\begin{subfigure}[t]{0.32\textwidth}
\includegraphics[width=1.\linewidth]{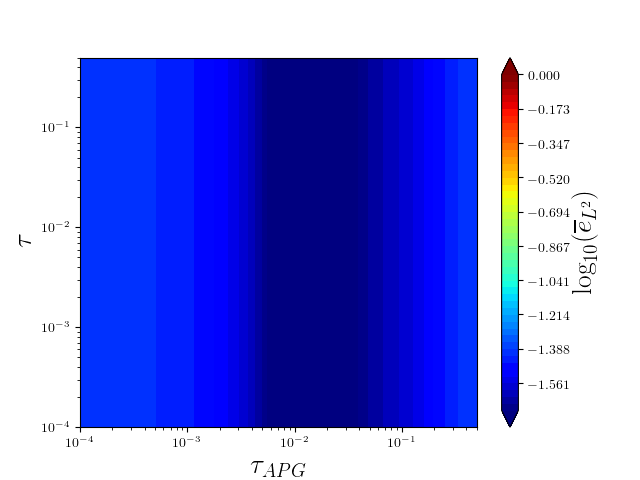}
\caption{\APGG}
\end{subfigure}
\begin{subfigure}[t]{0.32\textwidth}
\includegraphics[width=1.\linewidth]{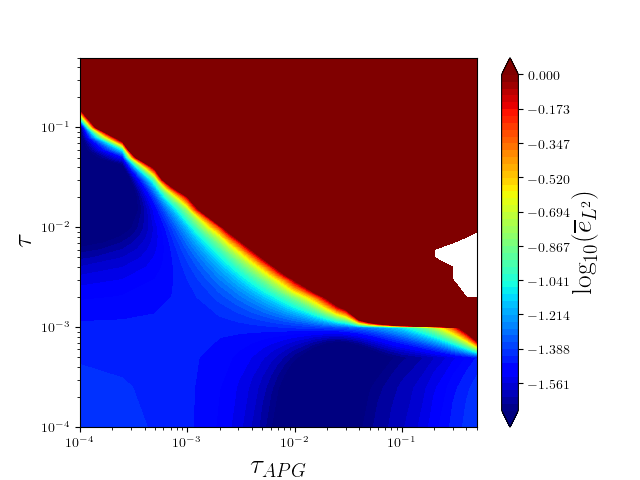}
\caption{\APGSUPG}
\label{fig:rom_piecewise_forcing_taudtvary_apgsupg}
\end{subfigure}
\begin{subfigure}[t]{0.32\textwidth}
\includegraphics[width=1.\linewidth]{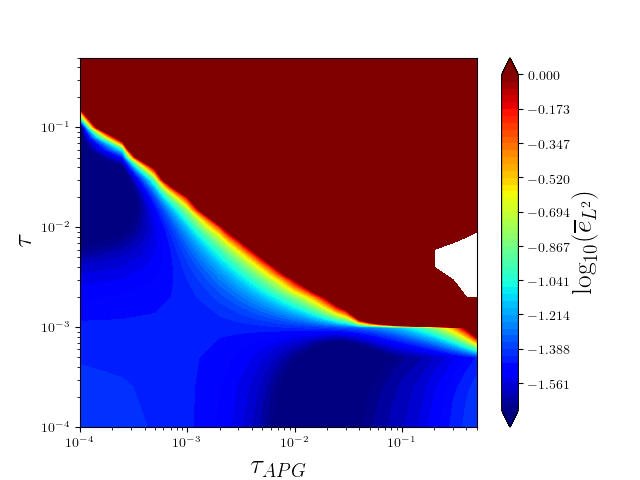}
\caption{\STGLS}
\end{subfigure}
\begin{subfigure}[t]{0.32\textwidth}
\includegraphics[width=1.\linewidth]{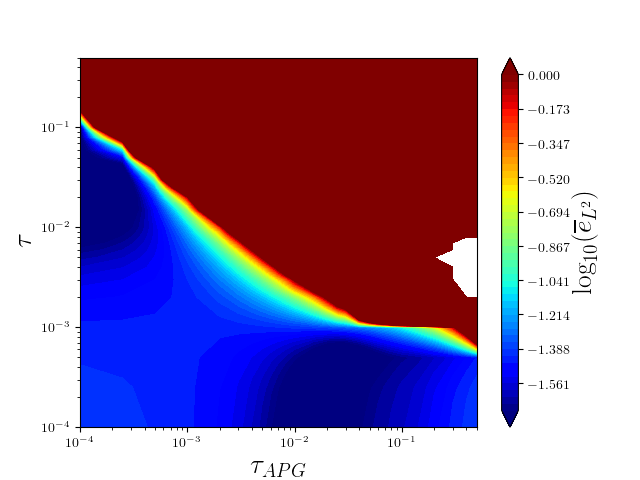}
\caption{\STADJ}
\label{fig:rom_piecewise_forcing_taudtvary_apgstadj}
\end{subfigure}
\caption{\advectingFrontFigureCaption Time integrated $\LTwo$ error as a function of time step and stabilization parameter for the various ROMs evaluated. Note that Galerkin and LSPG display no dependence on the stabilization parameter. White regions indicate errors higher than the color limit.}
\label{fig:rom_piecewise_forcing_taudtvary_apg}
\end{center}
\end{figure}

\begin{figure}
\begin{center}
\begin{subfigure}[t]{0.32\textwidth}
\includegraphics[width=1.\linewidth]{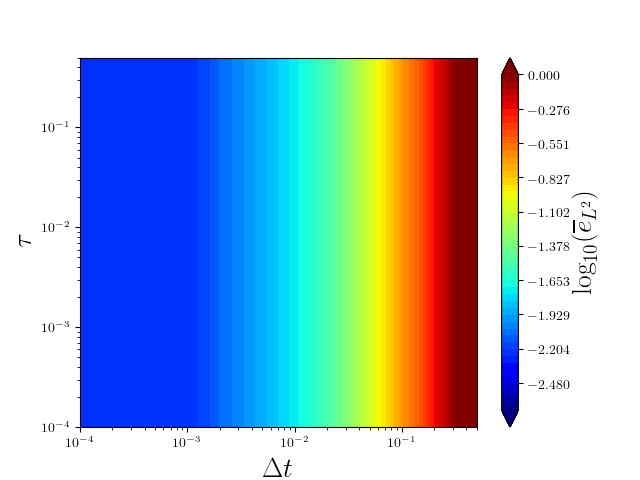}
\caption{Galerkin-FEM}
\label{fig:fom_piecewise_forcing_taudtvary_g}
\end{subfigure}
\begin{subfigure}[t]{0.32\textwidth}
\includegraphics[width=1.\linewidth]{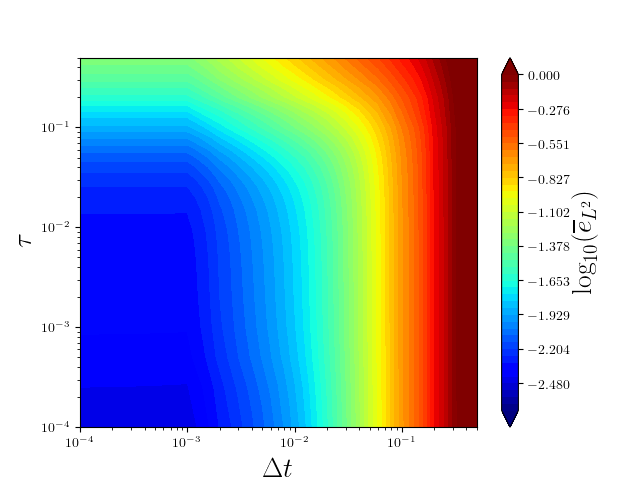}
\caption{\SUPG-FEM}
\end{subfigure}
\begin{subfigure}[t]{0.32\textwidth}
\includegraphics[width=1.\linewidth]{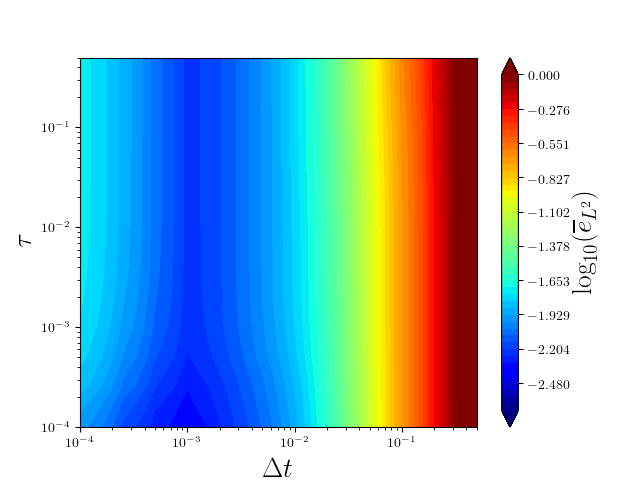}
\caption{\GLS-FEM}
\end{subfigure}
\begin{subfigure}[t]{0.32\textwidth}
\includegraphics[width=1.\linewidth]{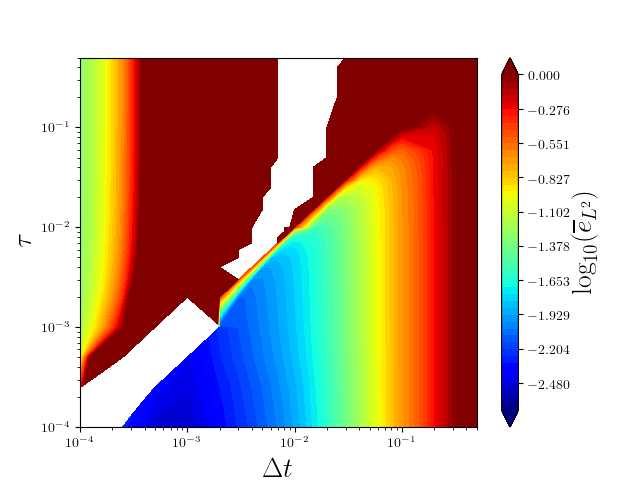}
\caption{\ADJ-FEM}
\end{subfigure}
\begin{subfigure}[t]{0.32\textwidth}
\includegraphics[width=1.\linewidth]{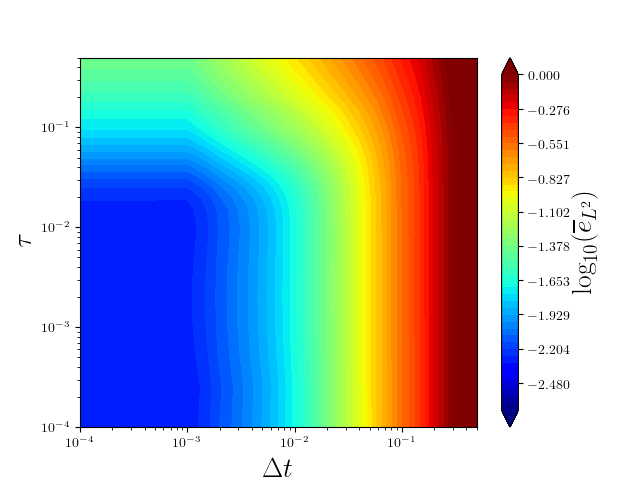}
\caption{\STGLS-FEM}
\end{subfigure}
\begin{subfigure}[t]{0.32\textwidth}
\includegraphics[width=1.\linewidth]{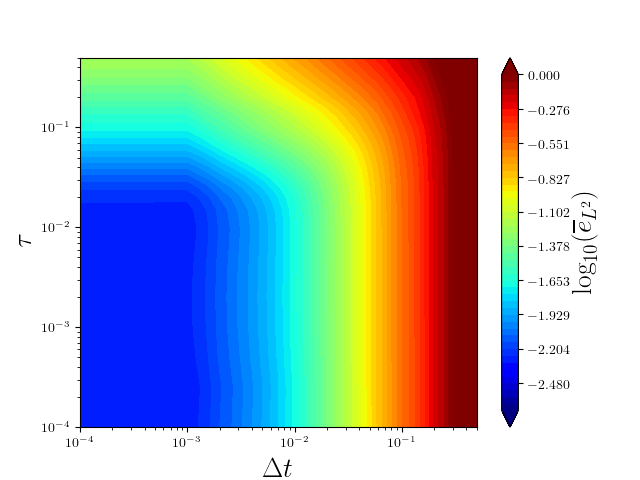}
\caption{\STADJ-FEM}
\label{fig:fom_piecewise_forcing_taudtvary_stadj}
\end{subfigure}

\caption{\advectingFrontFigureCaption Time integrated $\LTwo$ error as a function of time step and stabilization parameter. Results are shown for full-order FEM solutions executed on the FOM trial space. White regions indicate errors higher than the color limit.}
\label{fig:fom_piecewise_forcing_taudtvary}
\end{center}
\end{figure}

We now quantify the sensitivity of the various methods to their stabilization parameters and the time step.
Figure~\ref{fig:rom_piecewise_forcing_taudtvary} presents results for the continuous and LSPG ROM solutions obtained at $\romdim=5$ on the parameter grid $(\tau, \Delta t ) \in \tauSet \times \dtSet$, while 
Figure~\ref{fig:rom_piecewise_forcing_taudtvary_apg} presents results for the \APGG\ ROM solutions obtained at $\romdim=5$ on the parameter grid $(\tauApg, \Delta t ) \in \tauSet \times \dtSet$ and the remaining APG ROM solutions (which depend on three parameters, $\Delta t, \tau,$ and $\tauApg$) obtained at $\romdim=5$ on the parameter grid $(\tau,\tauApg) \in \tauSet \times \tauSet$ with a fixed time step $\Delta t = 10^{-3}$.  As a reference, Figure~\ref{fig:fom_piecewise_forcing_taudtvary} shows the same sensitivities but for full-order FEM simulations executed on the FOM trial space. We make the following observations:

\begin{itemize}
\item The \SUPG\ (Figure~\ref{fig:rom_piecewise_forcing_taudtvary_supg}), \STGLS\ (Figure~\ref{fig:rom_piecewise_forcing_taudtvary_stgls}), \STADJ\ (Figure~\ref{fig:rom_piecewise_forcing_taudtvary_stadj}), and \APGG\ (Figure~\ref{fig:rom_piecewise_forcing_taudtvary_apgg}) ROMs again all display a dependence on the time step and stabilization parameter that is similar to the first example. Optimal results are obtained for an intermediate value of $\tau$, and the solutions all converge in the limit of $\Delta t \rightarrow 0$.
\item \LSPGG\ (Figure~\ref{fig:rom_piecewise_forcing_taudtvary_lspg}), which contains no dependence on $\tau$, yields optimal results at an intermediate time step. LSPG's optimality at an intermediate time step is well documented in the community~\cite{carlberg_lspg_v_galerkin,parish_apg}. 
\item For sufficiently large $\tau$, the dependence of \GLS\ (Figure~\ref{fig:rom_piecewise_forcing_taudtvary_gls}) on the time step is similar to that of LSPG (Figure~\ref{fig:rom_piecewise_forcing_taudtvary_lspg}).
\item Neither \SDGLS\ nor \SDADJ\ ROMs have optimal results at the minimum time step, suggesting that these methods are not well behaved in the low time-step limit. This result reinforces those presented in Figure~\ref{fig:piecewise_forcing_romconverge_optimaldt_pcolor}.  
\item All continuous ROMs (Figures~\ref{fig:rom_piecewise_forcing_taudtvary_g}-\ref{fig:rom_piecewise_forcing_taudtvary_stadj}) display a similar dependence to the stabilization parameter and time step as their corresponding FOMs (Figures~\ref{fig:fom_piecewise_forcing_taudtvary_g}-\ref{fig:fom_piecewise_forcing_taudtvary_stadj}).  
\item Lastly, for APG ROMs built on top of a stabilized FEM (Figures~\ref{fig:rom_piecewise_forcing_taudtvary_apgsupg}-\ref{fig:rom_piecewise_forcing_taudtvary_apgstadj}), optimal results are obtained for either an intermediate value of $\tauApg$ and low value of $\tau$, or vice versa. It is interesting to note that the solutions are almost symmetric with respect to these two parameters. 
\end{itemize}

\subsection{Summary of numerical experiments and empirical findings}
Sections~\ref{sec:example1} and~\ref{sec:example2} presented results for stabilized ROMs applied to the CDR system for two different configurations. Across both cases, 
we observed that the ``space--time" stabilized continuous ROM formulations were superior to their ``discretize-then-stabilize" counterparts: the space--time stabilization formulations had lower errors, were well-behaved in the low time-step limit (the discretize-then-stabilized methods were not), and had a smoother, more intuitive dependence on the time step and stabilization parameter. 

ROMs constructed from LSPG projection  had mixed results in terms of accuracy. In the first numerical experiment, ROM solutions computed via LSPG projection tended to be slightly better than their non-stabilized (e.g., discrete Galerkin) ROM counterparts. In the second example, however, the performance was mixed. While \LSPGG\ and \LSPGGLS\ led to better solutions than Galerkin and \GLS, respectively, \LSPGSUPG, \LSPGADJ, \LSPGSTGLS, and \LSPGSTADJ\ all performed worse than \SUPG, \ADJ, \STGLS, and \STADJ, respectively. Further, all LSPG ROMs were optimal for intermediate time steps and, when built on top of a stabilized FEM solution, displayed a complex sensitivity to both the stabilization parameter and time step. The \textit{a priori} selection of an optimal time step appears difficult, and is likely problem dependent.   

Constructing ROMs from APG stabilization tended to have a positive result in terms of accuracy. In the first example, APG ROMs outperform their non-stabilized (i.e., discrete Galerkin) counterpart for all cases with the exception of \STGLS. In the second example, \APG-based ROMs were better than their non-stabilized counterparts for all FEM formulations. Further, unlike LSPG, APG was well behaved in the low time step limit. However, building an \APG\ ROM on top of a stabilized FEM solution led to complex behavior for the stabilization parameter: high values of $\tau$ for the stabilized FEM models required low values of $\tauApg$ for the ROM, and vice versa. The \textit{a priori} selection of these parameters again appears difficult.  

\begin{figure}
\begin{center}
\begin{subfigure}[t]{0.48\textwidth}
\includegraphics[trim={0cm 0cm 0cm 0cm},clip,  width=1.\linewidth]{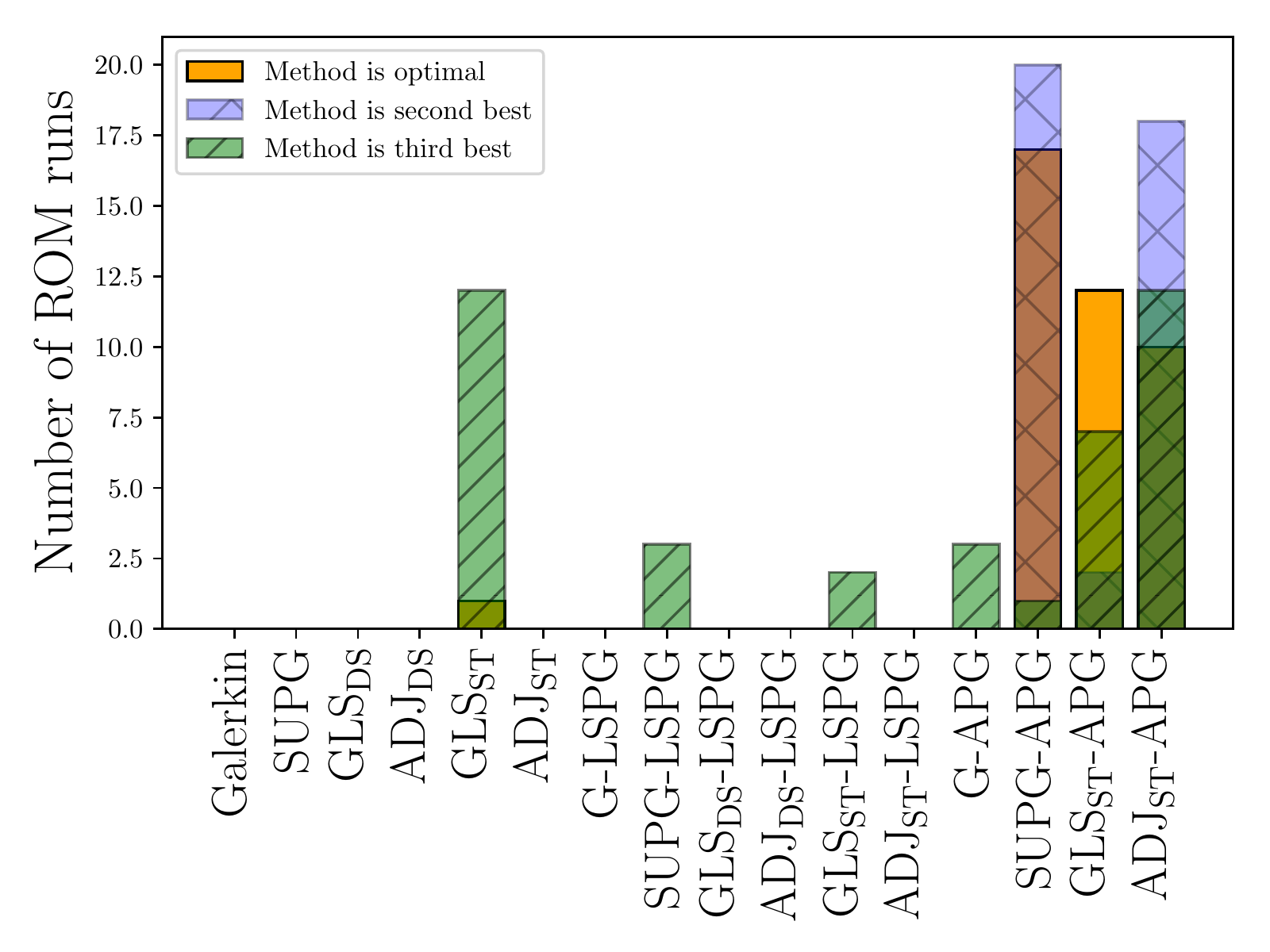}
\caption{$\LTwo$-norm}
\end{subfigure}
\begin{subfigure}[t]{0.48\textwidth}
\includegraphics[trim={0cm 0cm 0cm 0cm},clip,  width=1.\linewidth]{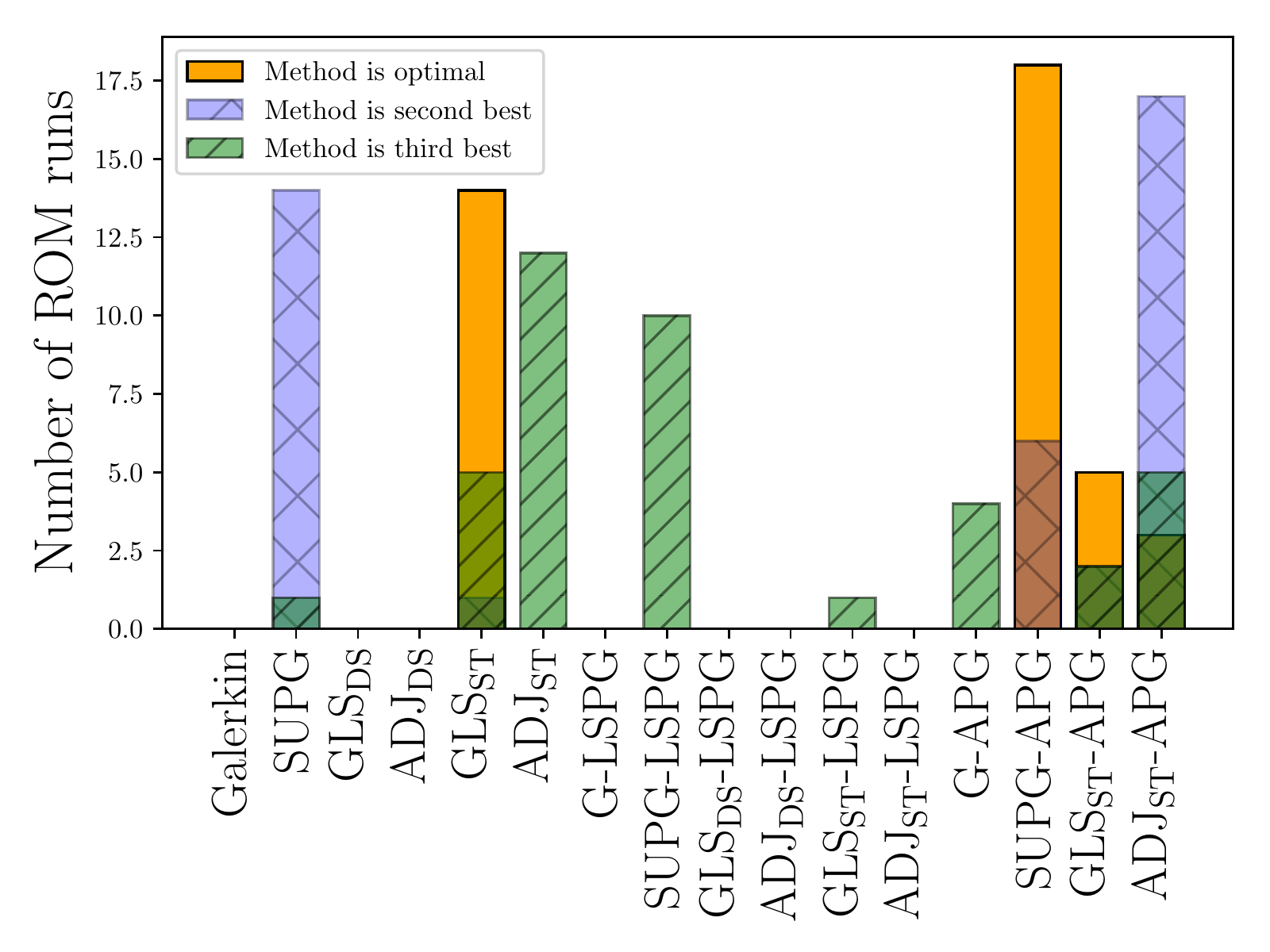}
\caption{$\HOne$-norm}
\end{subfigure}
\caption{Summary of numerical experiments. Number of times a given ROM formulation led to lowest errors in the $\LTwo$-norm (left) and $\HOne$-norm (right). Results are compiled across both numerical examples for all basis dimensions.}
\label{fig:rom_optimal_summary}
\end{center} 
\end{figure}

Figure~\ref{fig:rom_optimal_summary} summarizes the performance of the various ROMs by tabulating the number of times a given ROM formulation led to the lowest errors in the $\LTwo$ and $\HOne$-norm. Results are compiled for ROMs of basis dimensions $R=1,\ldots,20$ across both numerical experiments. We observe that \SUPGAPG\ was the best performing ROM for both error metrics. The next best performing ROMs were \APGSTADJ\ and \APGSTGLS, followed closely by \STGLS, \SUPG, and \STADJ. It is interesting to note that the two different error measures lead to slightly different measure of optimality (SUPG is never a top-performing method in $\LTwo$ but is consistently a top-performing method in $\HOne$). 

\begin{figure}
\begin{center}
\begin{subfigure}[t]{0.48\textwidth}
\includegraphics[trim={0cm 0cm 0cm 0cm},clip,  width=1.\linewidth]{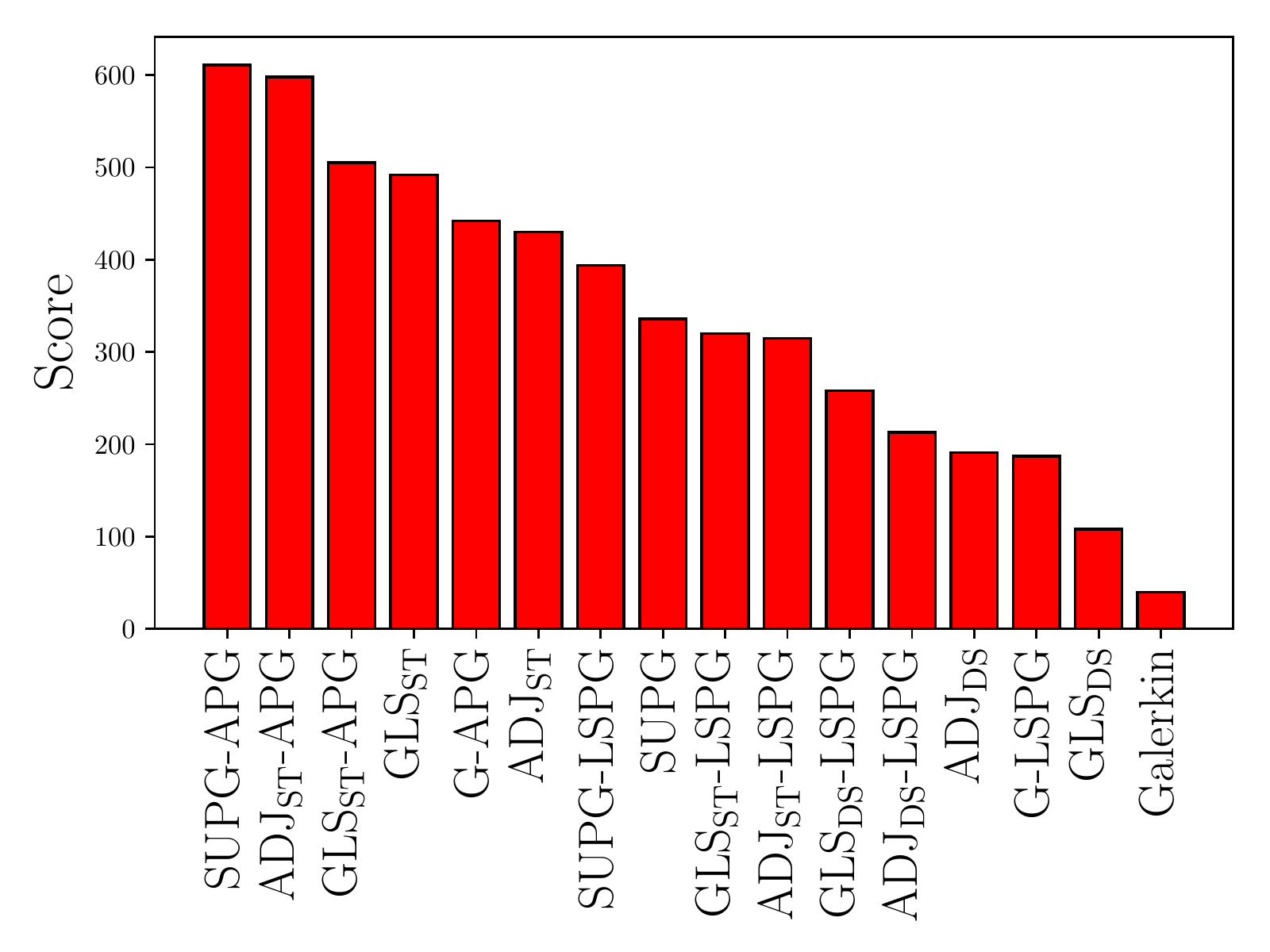}
\caption{Score based on the $\LTwo$-norm}
\end{subfigure}
\begin{subfigure}[t]{0.48\textwidth}
\includegraphics[trim={0cm 0cm 0cm 0cm},clip,  width=1.\linewidth]{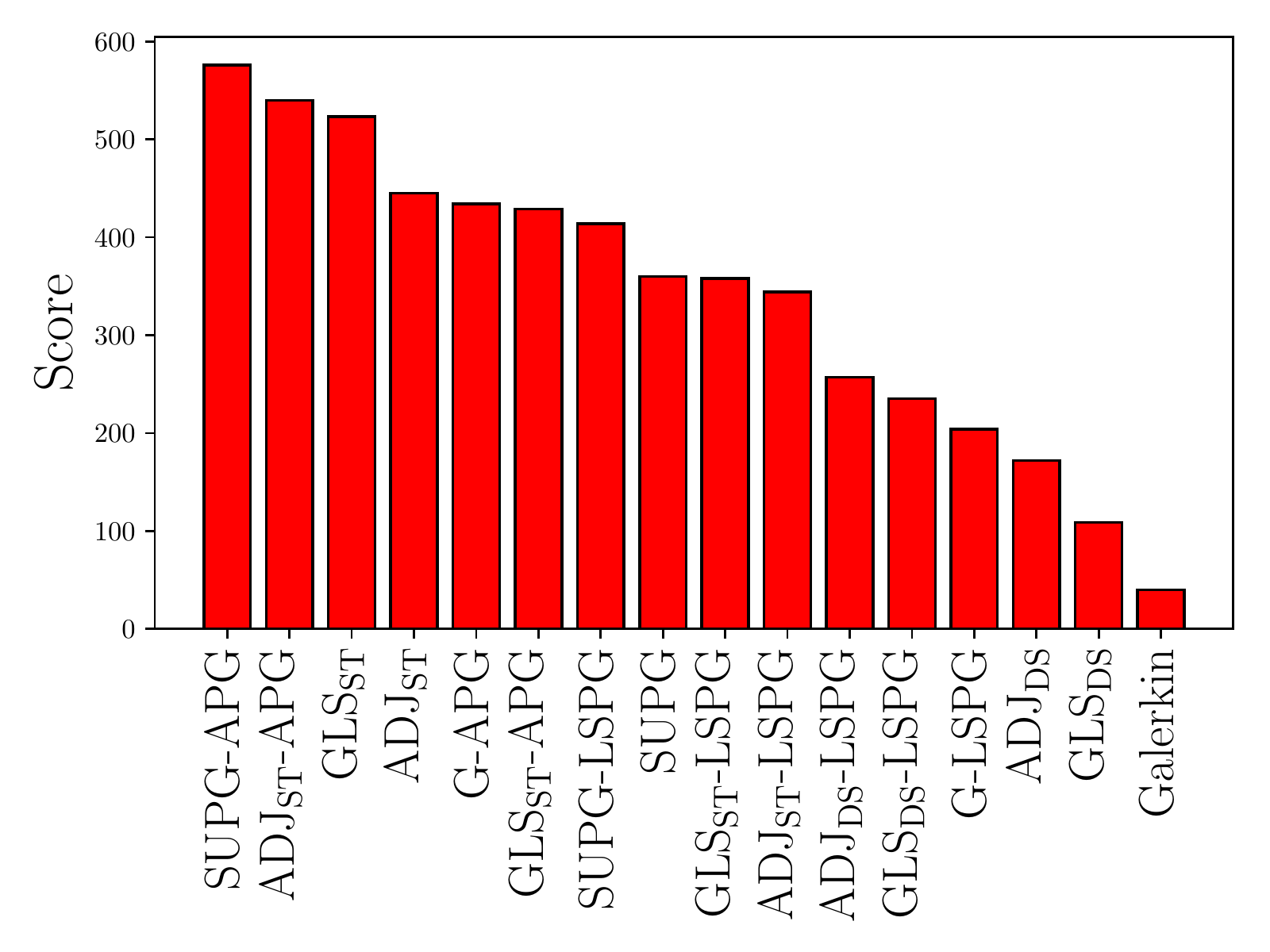}
\caption{Score based on the $\HOne$-norm}
\end{subfigure}
\caption{Summary of numerical experiments. Score as measured by the $\LTwo$-norm (left) and $\HOne$-norm (right). 
For a given basis dimension, we scored a ROM on a scale of $1$ to $N_{\text{ROMS}}$, where $N_{\text{ROMS}}=16$ is the number of ROM formulations considered. The best ROM gets a score of $N_{\text{ROMS}}$, the second best ROM gets a score of $N_{\text{ROMS}}-1$, and so on until the worst ROM gets a score of 1. The total score for each ROM formulation is computed by summing the individual scores across all basis dimensions and both numerical experiments.}
\label{fig:rom_score_summary}
\end{center} 
\end{figure}

Lastly, Figure~\ref{fig:rom_score_summary} attempts to rank the various ROM formulations by scoring their performance. For a given basis dimension, we scored a ROM on a scale of $1-N_{\text{ROMS}}$, where $N_{\text{ROMS}}=16$ is the number of ROM formulations considered. The best ROM gets a score of $N_{\text{ROMS}}$, the second best ROM gets a score of $N_{\text{ROMS}}-1$, and so on until the worst-performing ROM gets a score of 1. The total score for each ROM formulation is computed by summing the individual scores across all basis dimensions and both numerical experiments. By this scoring system, we find that \APGSUPG\ is the best performing ROM in both error measures, followed closely by \APGSTADJ, \STGLS, and \APGSTGLS. The \STGLS\ ROM is the best continuous ROM considered in this work, while the Galerkin ROM performs the worst.  

We end the discussion with a brief note regarding the computational cost of the various ROMs and methods evaluated.  While we did not report wall-clock times for running the various models considered in this work, we remark that, for linear problems, all ROM methods considered herein have similar online computational costs.  The costs will differ for nonlinear problems, and will depend on the choice of hyper-reduction method used to handle the online evaluation of the nonlinearities in the governing problem.  Extension of the analysis presented herein to nonlinear problems will be the subject of a subsequent publication.  

%
%

\section{Conclusions} \label{sec:conclude}
The development of robust ROMs for time-critical and
many-query scenarios remains an active research area. This work
outlined the construction of stabilized ROMs for the transient
convection-diffusion-reaction equation via two differing approaches
that have emerged within the community: discrete and continuous
projection. We outlined the standard Galerkin, \SUPG, \GLS, \ADJ, \STGLS, and \STADJ\
continuous ROMs developed via traditional stabilized finite
elements. We additionally outlined the construction of the discrete
Galerkin, LSPG, and APG ROMs. These discrete ROMs can be constructed from a standard Galerkin FEM model, or
they can also be constructed from a stabilized FEM model. We highlighted the well-established
equivalences between constructing ROM basis vectors at the discrete
and continuous levels. We additionally highlighted the established
equivalence conditions between discrete ROM formulations and continuous
ROM formulations. Lastly, a brief summary of existing numerical analyses was provided, where we discussed consistency, stability, and error bounds of the various methods.

Numerical experiments were conducted for two configurations of the CDR system. These experiments demonstrated that all stabilized ROMs result in superior performance over a standard Galerkin ROM built via continuous projection. APG and \STGLS-based ROMs proved to be the best performing methods, while \GLS, \ADJ, and \LSPGG\ proved to be the worst performing stabilized methods. In particular, we found that equipping a stabilized FEM model with \APG\ or \LSPG\ projection can result in more accurate solutions. This improvement in accuracy comes at the cost of a more complex dependence on the stabilization parameters and time step. In the case of \LSPG, results are optimal at an intermediate time step which is hard to select \textit{a priori}. Further, the accuracy of \LSPG-based methods degrades in the small time step limit. In the case of \APG, we observed a non-trivial relationship between $\tau$ in the stabilized FEM model and $\tauApg$ in the APG projection. Low values of $\tau$ in the stabilized FEM model required high values of $\tauApg$, and vice versa. There exist minimal methods for the \textit{a priori} selection of these optimal stabilization parameters. However, \APG-based ROMs were robust for small time steps. Lastly, we observed that ROMs built via continuous projection from the ``space--time" approach were clearly superior to ROMs built via the ``discretize-then-stabilize" approach: the space–time stabilization formulations had lower errors, were well-behaved in the low time-step limit (the discretize-then-stabilized methods were not), and had a smoother, more intuitive dependence on the time step and stabilization parameter.

In addition, our study highlighted several points (most of which are  well-established in the literature) which we reiterate here:
\begin{itemize} 
\item In the case of discrete projection, the POD basis must be obtained in an inner product that is $\mass$-orthogonal to recover the POD basis obtained via continuous projection. This is well established in the literature.
\item Both the Galerkin and stabilized ROMs obtained via continuous projection can be obtained via the Galerkin ROM obtained via discrete projection. This is additionally well established.
\item \LSPGG\ approximates a continuous minimization principle if the discrete norm is minimized in the $\mass^{-1}$ inner product.
\item All stabilized methods depend on a stabilization parameter, $\tau$, the time step, $\Delta t$, or both. The dependency on these parameters is complex, and more work needs to be done for proper \textit{a priori} selection of the stabilization parameters.
\end{itemize}
Future work should focus on three aspects. First, theory for discretely stabilized ROMs such as LSPG and APG is lacking for linear problems: minimal analyses exist   
studying the stability of these methods as well as their accuracy, and future work should address this. Second, future work should focus on the \textit{a priori} selection of the stabilization parameters and optimal time steps (or, alternatively, the development of methods that eliminate this sensitivity as in Ref~\cite{parish_wls}). Minimal analyses exists for the appropriate \textit{a priori} selection of these parameters within the context of ROMs: SUPG is studied in Refs.~\cite{GIERE2015454,JoMoNo22}, LSPG is studied in~\cite{carlberg_lspg_v_galerkin}, and APG is studied in~\cite{parish_apg}. Lastly, future work should focus on extension to nonlinear and vector-valued systems. Here, hyper-reduction is important and the extension of stabilization techniques to this setting offers numerous interesting and important questions.  

\section*{Acknowledgments}
E. Parish acknowledges funding from the John von Neumann Postdoctoral Fellowship, ASC V\&V 103723/05.30.02, and ASC V\&V 131792/04.09.02. 
I. Tezaur acknowledges funding from her Presidential Early Career Award for Scientists and Engineers
(PECASE), awarded by the U.S. Department of Energy (DOE), as well as support from the U.S. Office of Science, Office of Advanced Scientific Computing Research, Mathematical Multifaceted Integrated Capability Centers (MMICCS) program, under Field Work Proposal 22025291 and the Multifaceted Mathematics for Predictive Digital Twins (M2dt) project.
T. Iliescu acknowledges funding from the NSF grants DMS-2012253 and CDS\&E-MSS-1953113.
This paper describes objective technical results and analysis.
Any
subjective views or opinions that might be expressed in the paper do
not necessarily represent the views of the U.S. Department of Energy
or the United States Government. Sandia National Laboratories is a
multimission laboratory managed and operated by National Technology \&
Engineering Solutions of Sandia, LLC, a wholly owned subsidiary of
Honeywell International Inc., for the U.S. Department of Energy’s
National Nuclear Security Administration under contract DE-NA0003525.

\begin{appendices}
\section{Correspondence of LSPG to a continuous minimization principle}\label{appendix:lspg_proof}
This section outlines the equivalence between LSPG and a continuous minimization principle.  
We define the time-discrete, spatially continuous residual of the PDE~\eqref{eq:cdr} as
$$\cdcrResidTimeDiscreteSpaceContinuous : (\cStateTimeDiscreteDumN ;\cStateTimeDiscreteDumNm) \mapsto
\frac{\cStateTimeDiscreteDumN - \cStateTimeDiscreteDumNm}{\Delta t} - \nu \nabla^2 \cStateTimeDiscreteDumN + \wavespeed
\cdot \nabla \cStateTimeDiscreteDumN + \reaction
 \cStateTimeDiscreteDumN - \forcing.
$$
Under the assumption that the state is sufficiently
regular\footnote{We remark that this assumption does \textit{not}
hold for standard $\mathcal{C}^0(\cDomain)$ FEM discretizations, in
which case the state is not twice differentiable.}, the residual of
the Galerkin FOM O$\Delta$E can be written as
$$
\residGalerkinFEMArg{i} : (\genStateTimeDiscreteSpaceDiscreteDumN;\genStateTimeDiscreteSpaceDiscreteDumNm) \mapsto \cipGen{\cBasis_{i}}{ \cdcrResidTimeDiscreteSpaceContinuous(\cBasisVec \genStateTimeDiscreteSpaceDiscreteDumN ;\cBasisVec \genStateTimeDiscreteSpaceDiscreteDumNm}) , \qquad \forall \cBasis \in \cSpaceFEM.
$$
The \femCoefficients\ of the $L^2(\cDomain)$ orthogonal projection of this residual onto the trial space $\cSpaceFEM$ are given by
$$
\residGalerkinFEMOrtho : (\genStateTimeDiscreteSpaceDiscreteDumN;\genStateTimeDiscreteSpaceDiscreteDumNm) \mapsto
 \mass^{-1} \residGalerkinFEM (\genStateTimeDiscreteSpaceDiscreteDumN;\genStateTimeDiscreteSpaceDiscreteDumNm).
$$
Analogously, at the spatially
continuous level
$$
\cdcrResidTimeDiscreteSpaceContinuousOrtho : (\cStateTimeDiscreteDumN ;\cStateTimeDiscreteDumNm) \mapsto
\cBasisVec \mass^{-1} \cipGen{\cBasisVec}{ \cdcrResidTimeDiscreteSpaceContinuous(\cStateTimeDiscreteDumN ;\cStateTimeDiscreteDumNm)}.
$$
The square residual integrated over the domain is then given by
\begin{equation}\label{eq:min_projected_residual}
\int_{\cDomain} \left(
\cdcrResidTimeDiscreteSpaceContinuousOrtho(\cStateTimeDiscreteArg{n}
;\cStateTimeDiscreteArg{n-1}) \right)^2 dx = \left[ \mass^{-1}
\residGalerkinFEM
(\genStateTimeDiscreteSpaceDiscreteArg{n};\genStateTimeDiscreteSpaceDiscreteArg{n-1})
\right]^T \mass
 \mass^{-1} \residGalerkinFEM (\genStateTimeDiscreteSpaceDiscreteArg{n};\genStateTimeDiscreteSpaceDiscreteArg{n-1}).
\end{equation}
Setting $\resid \leftarrow \residGalerkinFEM$ and $\discreteInnerProductType \leftarrow \mass^{-1}$ in optimization problem~\eqref{eq:LSPG_minproblem}, LSPG corresponds to the continuous minimization principle for $\romStateTimeDiscreteSpaceDiscreteArg{n}$, $n=1,\ldots,\nTimeSteps$
\begin{equation}\label{eq:LSPG_minproblem_continuous}
\romStateTimeDiscreteSpaceDiscreteArg{n} =
\underset{\cStateTimeDiscreteDum \in \cSpaceRom  }{\text{arg\,min}}
\; \int_{\cDomain} \left(
\cdcrResidTimeDiscreteSpaceContinuousOrtho(\cStateTimeDiscreteDum;\romStateTimeDiscreteSpaceDiscreteArg{n-1})
\right)^2 dx.
\end{equation}
LSPG computes the solution $\romStateTimeDiscreteSpaceDiscreteArg{n}$ within the ROM trial space $\cSpaceRom$ that minimizes the $\LTwo$-norm of the time-discrete, spatially continuous residual projected onto the finite element trial space $\cSpaceFEM$.
Analogously, LSPG can be defined on the stabilized FOM O$\Delta$E~\eqref{eq:stab_fom_odeltae}. Setting $\resid \leftarrow \residStabilizedFEM$, the optimality conditions become
$$\dipGen{\frac{\mass \dBasisROMMat}{\Delta t} + \matFull \dBasisROMMat + \stabilizationMat \dBasisROMMat}{\residStabilizedFEM(\dBasisROMMat \ddgenStateDum;\dBasisROMMat \genDiscreteRomStateTimeDiscreteSpaceDiscreteArg{n-1})}{} = \mathbf{0}.$$

\section{Derivation of the Adjoint Petrov--Galerkin method}\label{appendix:derive_apg}
The APG method was derived in Ref.~\cite{parish_apg} in the case where the coarse- and fine-scale bases are orthogonal in a standard $\ell^2$ discrete inner product. In FEM discretizations it is more appropriate to construct the spaces to be orthogonal in an $\mass$ inner product. As such, in this section we derive the APG method for the case where the fine and coarse scales are $\mass$-orthogonal. We consider application to the dynamical system given by
$$\mass \frac{d \mathbf{x}}{dt} + \mathbf{A} \mathbf{x} - \mathbf{f} = \bz,$$
where $\mass \in \mathbb{S}^{\fomdim}$ is the mass matrix, $\mathbf{x}: [0,T] \rightarrow \RR{\fomdim}$ is the state, $\mathbf{A} \in \RR{\fomdim \times \fomdim}$ is the system matrix, and $\mathbf{f}$ is a forcing term.
The derivation begins by decomposing $\mathbb{R}^{\fomdim}$ into a coarse-scale space $\dSpaceROM$ and a fine-scale space $\dSpaceROMFine$ such that $\dSpaceROM \oplus \dSpaceROMFine \equiv \mathbb{R}^{\fomdim}$. The coarse-scale space $\dSpaceROM \subset \mathbb{R}^{\fomdim}$ corresponds to the standard ROM space and is of dimension $\text{dim}(\dSpaceROM) = \romdim$, while the fine-scale space $\dSpaceROMFine$ comprises the $\mass$-orthogonal complement of the coarse-scale space and has dimension $\text{dim}(\dSpaceROMFine) = \fomdim - \romdim$. We equip the coarse- and fine-scale spaces with $\mass$-orthogonal bases $\dBasisROMMat \in \RR{\fomdim \times \romdim}$ and $\dBasisROMMat' \in \RR{\fomdim \times (\fomdim - \romdim)}$. Note that $\dBasisROMMat^T \mass \dBasisROMMat' = \bz$ by definition.

The APG derivation proceeds by expressing the dynamical system in terms of the generalized coordinates associated with the coarse and fine-scale bases. This process results in a coupled system for the coarse and fine scales
\begin{equation*}
\begin{split}
& \dBasisROMMat^T \mass \frac{d}{dt} \dBasisROMMat \genDiscreteRomStateTimeDiscreteSpaceDiscreteArg{} + \dBasisROMMat^T \mathbf{A} \left[ \dBasisROMMat  \genDiscreteRomStateTimeDiscreteSpaceDiscreteArg{} +  \dBasisROMMat^{'} \genDiscreteRomStateTimeDiscreteSpaceDiscreteArg{'}\right] - \dBasisROMMat^T \mathbf{f} = \bz \\ 
&\left[\dBasisROMMat^{'}\right]^T \mass \frac{d}{dt} \dBasisROMMat^{'} \genDiscreteRomStateTimeDiscreteSpaceDiscreteArg{'} + \left[\dBasisROMMat^{'}\right]^T \mathbf{A}  \left[ \dBasisROMMat \genDiscreteRomStateTimeDiscreteSpaceDiscreteArg{} + \dBasisROMMat^{'} \genDiscreteRomStateTimeDiscreteSpaceDiscreteArg{'}\right] - \left[ \dBasisROMMat^{'}\right]^T \mathbf{f} = \bz,
\end{split}
\end{equation*}
where $\genDiscreteRomStateTimeDiscreteSpaceDiscreteArg{} :[0,T] \rightarrow \RR{\romdim}$ are the coarse-scale generalized coordinates and $\genDiscreteRomStateTimeDiscreteSpaceDiscreteArg{'} : [0,T] \rightarrow \RR{\fomdim - \romdim}$ are the fine-scale generalized coordinates. The APG method proceeds to approximate the fine-scales via the Mori--Zwanzig formalism and a perturbation analysis, which here results in the quasi-static aproximation 
$$
 \genDiscreteRomStateTimeDiscreteSpaceDiscreteArg{'}(t) \approx -\tau \left[\dBasisROMMat^{'}\right]^T \mathbf{A}  \dBasisROMMat \genDiscreteRomStateTimeDiscreteSpaceDiscreteArg{}(t),$$
where $\tau \in \RPlus$ is a stabilization constant. 
Injecting in the approximation to the fine-scale state into the coarse-scale equation results in
\begin{equation*}
 \dBasisROMMat^T \mass \frac{d}{dt} \dBasisROMMat \genDiscreteRomStateTimeDiscreteSpaceDiscreteArg{} + \dBasisROMMat^T \mathbf{A} \left[ \dBasisROMMat  \genDiscreteRomStateTimeDiscreteSpaceDiscreteArg{} - \tau \dBasisROMMat^{'}  \left[\dBasisROMMat^{'}\right]^T \mathbf{A}  \dBasisROMMat \genDiscreteRomStateTimeDiscreteSpaceDiscreteArg{}\right] - \dBasisROMMat^T \mathbf{f} = \bz.
\end{equation*}
Next we use the property $\dBasisROMMat \dBasisROMMat^T \mass +  \dBasisROMMat^{'}  \left[\dBasisROMMat^{'}\right]^T \mass = \mathbf{I}$ 
to remove the dependence on the fine-scale basis functions and get
 \begin{equation*}
 \dBasisROMMat^T \mass \frac{d}{dt} \dBasisROMMat \genDiscreteRomStateTimeDiscreteSpaceDiscreteArg{} + \dBasisROMMat^T \mathbf{A} \left[ \dBasisROMMat  \genDiscreteRomStateTimeDiscreteSpaceDiscreteArg{} - \tau \PiFine \mathbf{A}  \dBasisROMMat \genDiscreteRomStateTimeDiscreteSpaceDiscreteArg{}\right] - \dBasisROMMat^T \mathbf{f} = \bz,
\end{equation*}
where $\PiFine = \mass^{-1} - \dBasisROMMat \dBasisROMMat^T$. Next, assuming the forcing to be zero on the fine-scale space such that $\PiFine \forcingVec = \bz$,  we write the above in a Petrov--Galerkin form
\begin{equation*}
\left[ \left( \mathbf{I} - \tau  \PiFine^T  \mathbf{A}^T \right) \dBasisROMMat \right]^T 
\left[\mass \frac{d}{dt} \dBasisROMMat \genDiscreteRomStateTimeDiscreteSpaceDiscreteArg{} + \mathbf{A} \dBasisROMMat \genDiscreteRomStateTimeDiscreteSpaceDiscreteArg{}  - \forcingVec \right] = \bz 
\end{equation*}
where we have leveraged $\PiFine \mass \dBasisROMMat  = \bz$. This Petrov--Galerkin projection is what we refer to as the Adjoint Petrov--Galerkin method.
\section{Proper orthogonal decomposition algorithm}
\label{sec:appendix_podsvd}

Algorithm~\ref{alg:pod} presents the algorithm for computing the trial basis via proper orthogonal decomposition.
\begin{algorithm}
\caption{Algorithm for generating POD Basis.}
\label{alg:pod}
\textbf{Input:} Snapshot matrix, $\snapshotMatrixDiscrete \in \RR{\fomdim \times \nSnap}$ ; cutoff energy tolerance,  $\energyCutoff$ ; symmetric positive definite inner product matrix, $\discreteInnerProductType$ \\
\textbf{Output:} POD Basis $\dBasisROMMat \in \RR{\fomdim \times \romdim}$ \;
\textbf{Steps:}
\begin{enumerate}
\item Compute the ``covariance matrix"
$$\mathbf{K} = \snapshotMatrixDiscrete^T \discreteInnerProductType \snapshotMatrixDiscrete.$$
\item Compute the eigenvalue decomposition 
$$ \mathbf{K} = \dSnapEigVecMat \dSnapEigMat [\dSnapEigVecMat]^{-1}.$$
\item Compute the statistical energy
$$\energy_K =\frac{ \sum_{i=1}^K \dSnapEig_i } { \sum_{i=1}^{\nSnap} \dSnapEig_i}, $$
where $K \le \nSnap$.
\item Determine basis dimension from cutoff criterion
$$ \romdim =  \text{Card}( \{ \energy_i \}_{i=1}^{\nSnap} | \energy_i \le \energyCutoff ).$$

\item Compute the ROM bases as
 $$\dBasisRomMat =  \snapshotMatrixDiscrete  \dSnapEigVecMatTruncate \sqrt{[\dSnapEigMatTruncate]^{-1}},$$
where $\dSnapEigVecMatTruncate$ and $\dSnapEigMatTruncate$ comprise the first $\romdim$ columns of $\dSnapEigVecMat$ and the first $\romdim$ columns and rows of $\dSnapEigMat$, respectively.

\end{enumerate}
\end{algorithm}

\end{appendices}

\bibliographystyle{siam}
\bibliography{refs}

\begin{thebibliography}{100}

\bibitem{l1}
{\sc R.~Abgrall and R.~Crisovan}, {\em Model reduction using {L1}-norm
  minimization as an application to nonlinear hyperbolic problems},
  International Journal for Numerical Methods in Fluids, 87 (2018),
  pp.~628--651.

\bibitem{afkham2017structure}
{\sc B.~M. Afkham and J.~S. Hesthaven}, {\em {Structure preserving model
  reduction of parametric Hamiltonian systems}}, SIAM J. Sci. Comput., 39
  (2017), pp.~A2616--A2644.

\bibitem{ahmed2021closures}
{\sc S.~E. Ahmed, S.~Pawar, O.~San, A.~Rasheed, T.~Iliescu, and B.~R. Noack},
  {\em On closures for reduced order models $-$ a spectrum of first-principle
  to machine-learned avenues}, Phys. Fluids, 33 (2021), p.~091301.

\bibitem{rozza_stabilized}
{\sc S.~Ali, F.~Ballarin, and G.~Rozza}, {\em Stabilized reduced basis methods
  for parametrized steady {Stokes} and {Navier-Stokes} equations}, arXiv
  e-print,  (2020).

\bibitem{AlnaesBlechta2015a}
{\sc M.~S. Aln{\ae}s, J.~Blechta, J.~Hake, A.~Johansson, B.~Kehlet, A.~Logg,
  C.~Richardson, J.~Ring, M.~E. Rognes, and G.~N. Wells}, {\em The fenics
  project version 1.5}, Archive of Numerical Software, 3 (2015).

\bibitem{amsallem_stab}
{\sc D.~Amsallem and C.~Farhat}, {\em Stabilization of projection-based
  reduced-order models}, International Journal for Numerical Methods in
  Engineering, 91 (2012), pp.~358--377.

\bibitem{azaiez2021cure}
{\sc M.~Aza{\"\i}ez, T.~C. Rebollo, and S.~Rubino}, {\em {A cure for
  instabilities due to advection-dominance in POD solution to
  advection-diffusion-reaction equations}}, J. Comput. Phys., 425 (2021),
  p.~109916.

\bibitem{doi:10.1002/fld.3777}
{\sc J.~Baiges, R.~Codina, and S.~Idelsohn}, {\em Explicit reduced-order models
  for the stabilized finite element approximation of the incompressible
  {Navier–Stokes} equations}, International Journal for Numerical Methods in
  Fluids, 72 (2013), pp.~1219--1243.

\bibitem{BAIGES2015173}
{\sc J.~Baiges, R.~Codina, and S.~Idelsohn}, {\em {Reduced-order subscales for
  POD models}}, Computer Methods in Applied Mechanics and Engineering, 291
  (2015), pp.~173 -- 196.

\bibitem{Balajewicz_rom0}
{\sc M.~Balajewicz and E.~H. Dowell}, {\em Stabilization of projection-based
  reduced order models of the {Navier--Stokes}}, Nonlinear Dynamics, 70 (2012),
  pp.~1619--1632.

\bibitem{basis_rotation}
{\sc M.~Balajewicz, I.~Tezaur, and E.~Dowell}, {\em Minimal subspace rotation
  on the {Stiefel} manifold for stabilization and enhancement of
  projection-based reduced order models for the compressible {Navier--Stokes}
  equations}, Journal of Computational Physics, 321 (2016), pp.~224--241.

\bibitem{ballarin2015supremizer}
{\sc F.~Ballarin, A.~Manzoni, A.~Quarteroni, and G.~Rozza}, {\em Supremizer
  stabilization of {POD--G}alerkin approximation of parametrized steady
  incompressible {N}avier--{S}tokes equations}, Int. J. Numer. Meth. Engng.,
  102 (2015), pp.~1136--1161.

\bibitem{BARONE20091932}
{\sc M.~F. Barone, I.~Kalashnikova, D.~J. Segalman, and H.~K. Thornquist}, {\em
  Stable {Galerkin} reduced order models for linearized compressible flow},
  Journal of Computational Physics, 228 (2009), pp.~1932 -- 1946.

\bibitem{benner_st}
{\sc M.~Baumann, P.~Benner, and J.~Heiland}, {\em Space-time {Galerkin POD}
  with application in optimal control of semilinear partial differential
  equations}, SIAM Journal on Scientific Computing, 40 (2018),
  pp.~A1611--A1641.

\bibitem{structurePreserveBeattie}
{\sc C.~Beattie and S.~Gugercin}, {\em Structure-preserving model reduction for
  nonlinear port-{H}amiltonian systems}, in Decision and Control and European
  Control Conference (CDC-ECC), 2011 50th IEEE Conference on, IEEE, 2011,
  pp.~6564--6569.

\bibitem{willcox_benner_rev}
{\sc P.~Benner, S.~Gugercin, and K.~Willcox}, {\em {A survey of
  projection-based model reduction methods for parametric dynamical systems}},
  SIAM Review, 57 (2015), pp.~483--531.

\bibitem{bergmann2009enablers}
{\sc M.~Bergmann, C.~H. Bruneau, and A.~Iollo}, {\em {Enablers for robust POD
  models}}, J. Comput. Phys., 228 (2009), pp.~516--538.

\bibitem{berkooz_turbulence_pod}
{\sc G.~Berkooz, P.~Holmes, and J.~L. Lumley}, {\em The proper orthogonal
  decomposition in the analysis of turbulent flows}, Annu. Rev. Fluid Mech., 25
  (1993), pp.~539--575.

\bibitem{BOCHEV20042301}
{\sc P.~B. Bochev, M.~D. Gunzburger, and J.~N. Shadid}, {\em Stability of the
  {SUPG} finite element method for transient advection–diffusion problems},
  Computer Methods in Applied Mechanics and Engineering, 193 (2004), pp.~2301
  -- 2323.

\bibitem{borggaard2011artificial}
{\sc J.~Borggaard, T.~Iliescu, and Z.~Wang}, {\em Artificial viscosity proper
  orthogonal decomposition}, Math. Comput. Modelling, 53 (2011), pp.~269--279.

\bibitem{brooks_thesis}
{\sc A.~N. Brooks}, {\em A Petrov-Galerkin Finite Element Formulation for
  Convection Dominated Flows}, PhD thesis, California Institute of Technology,
  1981.

\bibitem{BROOKS1982199}
{\sc A.~N. Brooks and T.~J. Hughes}, {\em {Streamline upwind/Petrov-Galerkin
  formulations for convection dominated flows with particular emphasis on the
  incompressible Navier-Stokes equations}}, Computer Methods in Applied
  Mechanics and Engineering, 32 (1982), pp.~199 -- 259.

\bibitem{bui_thesis}
{\sc T.~Bui-Thanh}, {\em {Model-constrained optimization methods for reduction
  of parameterized large-scale systems}}, PhD thesis, Massachusetts Institute
  of Technology, 2007.

\bibitem{bui_resmin_steady}
{\sc T.~Bui-Thanh, K.~Willcox, and O.~Ghattas}, {\em Model reduction for
  large-scale systems with high-dimensional parametric input space}, SIAM
  Journal on Scientific Computing, 30 (2008), pp.~3270--3288.

\bibitem{bui_unsteady}
\leavevmode\vrule height 2pt depth -1.6pt width 23pt, {\em Parametric
  reduced-order models for probabilistic analysis of unsteady aerodynamic
  applications}, AIAA Journal, 46 (2008), pp.~2520--2529.

\bibitem{caiazzo2014numerical}
{\sc A.~Caiazzo, T.~Iliescu, V.~John, and S.~Schyschlowa}, {\em A numerical
  investigation of velocity-pressure reduced order models for incompressible
  flows}, J. Comput. Phys., 259 (2014), pp.~598--616.

\bibitem{carlberg_thesis}
{\sc K.~Carlberg}, {\em Model Reduction of Nonlinear Mechanical Systems via
  Pptimal Projection and Tensor Approximation}, PhD thesis, Stanford
  University, 2011.

\bibitem{carlberg_lspg_v_galerkin}
{\sc K.~Carlberg, M.~Barone, and H.~Antil}, {\em {Galerkin v. least-squares
  Petrov-Galerkin projection in nonlinear model reduction}}, Journal of
  Computational Physics, 330 (2017), pp.~693--734.

\bibitem{carlberg_lspg}
{\sc K.~Carlberg, C.~Bou-Mosleh, and C.~Farhat}, {\em Efficient non-linear
  model reduction via a least-squares petrov-galerkin projection and
  compressive tensor approximations}, Int. J. Numer. Methods Eng., 86 (2011),
  pp.~155--181.

\bibitem{carlberg_conservative_rom}
{\sc K.~Carlberg, Y.~Choi, and S.~Sargsyan}, {\em Conservative model reduction
  for finite-volume models}, Journal of Computational Physics, 371 (2018),
  pp.~280--314.

\bibitem{carlberg_gnat}
{\sc K.~Carlberg, C.~Farhat, J.~Cortial, and D.~Amsallem}, {\em The {GNAT}
  method for nonlinear model reduction: {Effective} implementation and
  application to computational fluid dynamics and turbulent flows}, Journal of
  Computational Physics, 242 (2013), pp.~623--647.

\bibitem{carlberg2012spd}
{\sc K.~Carlberg, R.~Tuminaro, and P.~Boggs}, {\em Preserving {L}agrangian
  structure in nonlinear model reduction with application to structural
  dynamics}, SIAM J. Sci. Comput., 37 (2015), pp.~B153---B184.

\bibitem{carlberg2015preserving}
\leavevmode\vrule height 2pt depth -1.6pt width 23pt, {\em {Preserving
  Lagrangian structure in nonlinear model reduction with application to
  structural dynamics}}, SIAM J. Sci. Comput., 37 (2015), pp.~B153--B184.

\bibitem{Chan:2020}
{\sc J.~Chan}, {\em Entropy stable reduced order modeling of nonlinear
  conservation laws}, Journal of Computational Physics, 423 (2020), p.~109789.

\bibitem{chaturantabut2016structure}
{\sc S.~Chaturantabut, C.~Beattie, and S.~Gugercin}, {\em {Structure-preserving
  model reduction for nonlinear port-Hamiltonian systems}}, SIAM J. Sci.
  Comput., 38 (2016), pp.~B837--B865.

\bibitem{chekroun2019variational}
{\sc M.~D. Chekroun, H.~Liu, and J.~C. McWilliams}, {\em {Variational approach
  to closure of nonlinear dynamical systems: Autonomous case}}, J. Stat. Phys.,
   (2019), pp.~1--88.

\bibitem{choi_stlspg}
{\sc Y.~Choi and K.~Carlberg}, {\em {Space-Time Least-Squares Petrov-Galerkin
  Projection for Nonlinear Model Reduction}}, SIAM J. Sci. Comput.,  (2019).

\bibitem{Chorin_predictwithmem}
{\sc A.~Chorin, O.~Hald, and R.~Kupferman}, {\em {Optimal prediction with
  memory}}, Phys. D,  (2002), pp.~239--257.

\bibitem{CODINA1998185}
{\sc R.~Codina}, {\em Comparison of some finite element methods for solving the
  diffusion-convection-reaction equation}, Computer Methods in Applied
  Mechanics and Engineering, 156 (1998), pp.~185 -- 210.

\bibitem{codina_oss}
\leavevmode\vrule height 2pt depth -1.6pt width 23pt, {\em {Stabilization of
  incompressibility and convection through orthogonal sub-scales in finite
  element methods}}, Computer methods in applied mechanics and engineering, 190
  (2000).

\bibitem{CODINA20072413}
{\sc R.~Codina, J.~Principe, O.~Guasch, and S.~Badia}, {\em Time dependent
  subscales in the stabilized finite element approximation of incompressible
  flow problems}, Computer Methods in Applied Mechanics and Engineering, 196
  (2007), pp.~2413 -- 2430.

\bibitem{constantine_strom}
{\sc P.~G. Constantine and Q.~Wang}, {\em Residual minimizing model
  interpolation for parameterized nonlinear dynamical systems}, SIAM J. Sci.
  Comput.,  (2012).

\bibitem{CSB03}
{\sc M.~Couplet, P.~Sagaut, and C.~Basdevant}, {\em Intermodal energy transfers
  in a proper orthogonal decomposition--{G}alerkin representation of a
  turbulent separated flow}, J. Fluid Mech., 491 (2003), pp.~275--284.

\bibitem{Dahmen_2014_Double_Greedy}
{\sc W.~Dahmen, C.~Plesken, and G.~Welper}, {\em Double greedy algorithms:
  Reduced basis methods for transport dominated problems}, ESAIM: M2AN, 48
  (2014), pp.~623--–663.

\bibitem{decaria2020artificial}
{\sc V.~DeCaria, T.~Iliescu, W.~Layton, M.~McLaughlin, and M.~Schneier}, {\em
  An artificial compression reduced order model}, SIAM J. Numer. Anal.,
  (2020).
\newblock accepted.

\bibitem{eroglu2017modular}
{\sc F.~G. Eroglu, S.~Kaya, and L.~G. Rebholz}, {\em A modular regularized
  variational multiscale proper orthogonal decomposition for incompressible
  flows}, Comput. Meth. Appl. Mech. Eng., 325 (2017), pp.~350--368.

\bibitem{farhat2014dimensional}
{\sc C.~Farhat, P.~Avery, T.~Chapman, and J.~Cortial}, {\em Dimensional
  reduction of nonlinear finite element dynamic models with finite rotations
  and energy-based mesh sampling and weighting for computational efficiency},
  Int. J. Num. Meth. Eng., 98 (2014), pp.~625--662.

\bibitem{FrVa20}
{\sc L.~Franca and F.~Valentin}, {\em On an improved unusual stabilized finite
  element method for the advective--reactive--diffusive equation}, Computer
  Methods in Applied Mechanics and Engineering, 190 (2000), pp.~1785--1800.

\bibitem{FRANCA1995299}
{\sc L.~P. Franca and C.~Farhat}, {\em Bubble functions prompt unusual
  stabilized finite element methods}, Computer Methods in Applied Mechanics and
  Engineering, 123 (1995), pp.~299 -- 308.

\bibitem{FRANCA1992253}
{\sc L.~P. Franca, S.~L. Frey, and T.~J. Hughes}, {\em Stabilized finite
  element methods: I. {Application} to the advective-diffusive model}, Computer
  Methods in Applied Mechanics and Engineering, 95 (1992), pp.~253 -- 276.

\bibitem{funaro}
{\sc D.~Funaro and D.~Gottlieb}, {\em {Convergence results for pseudospectral
  approximations of hyperbolic systems by a penalty-type boundary treatment}},
  Mathematics of Computation, 57 (1991).

\bibitem{GIERE2015454}
{\sc S.~Giere, T.~Iliescu, V.~John, and D.~Wells}, {\em {SUPG} reduced order
  models for convection-dominated convection–diffusion–reaction equations},
  Computer Methods in Applied Mechanics and Engineering, 289 (2015), pp.~454 --
  474.

\bibitem{gong2017structure}
{\sc Y.~Gong, Q.~Wang, and Z.~Wang}, {\em {Structure-preserving Galerkin POD
  reduced-order modeling of Hamiltonian systems}}, Comput. Methods Appl. Mech.
  Engrg., 315 (2017), pp.~780--798.

\bibitem{Grepl_2007_RB_EIM}
{\sc M.~A. Grepl, Y.~Maday, N.~C. Nguyen, and A.~T. Patera}, {\em Efficient
  reduced-basis treatment of nonaffine and nonlinear partial differential
  equations}, ESAIM: M2AN, 41 (2007), pp.~575--605.

\bibitem{Grepl_2005_RB_Parabolic}
{\sc M.~A. Grepl and A.~T. Patera}, {\em A posteriori error bounds for
  reduced-basis approximations of parametrized parabolic partial differential
  equations}, {ESAIM}: Mathematical Modelling and Numerical Analysis, 39
  (2005), pp.~157--181.

\bibitem{GrFaYo20}
{\sc S.~Grimberg, C.~Farhat, and N.~Youkilis}, {\em On the stability of
  projection-based model order reduction for convection-dominated laminar and
  turbulent flows}, Journal of Computational Physics, 419 (2020), p.~109681.

\bibitem{Gruber:2023}
{\sc A.~Gruber, M.~Gunzburger, L.~Ju, and Z.~Wang}, {\em Energetically
  consistent model reduction for metriplectic systems}, Computer Methods in
  Applied Mechanics and Engineering, 404 (2023), p.~115709.

\bibitem{gunzburger2019evolve}
{\sc M.~Gunzburger, T.~Iliescu, M.~Mohebujjaman, and M.~Schneier}, {\em {An
  evolve-filter-relax stabilized reduced order stochastic collocation method
  for the time-dependent Navier-Stokes equations}}, SIAM-ASA J. Uncertain.,
  (2019), pp.~1162--1184.

\bibitem{Haasdonk_2013_POD_Greedy_Convergence}
{\sc B.~Haasdonk}, {\em Convergence rates of the {POD}-{Greedy} method},
  {ESAIM}: Mathematical Modelling and Numerical Analysis, 47 (2013),
  pp.~859--873.

\bibitem{Haasdonk_2008_RB_FV_Evolution}
{\sc B.~Haasdonk and M.~Ohlberger}, {\em Reduced basis method for finite volume
  approximations of parametrized linear evolution equations}, Mathematical
  Modelling and Numerical Analysis, 42 (2008), pp.~277--302.

\bibitem{HaSt07}
{\sc O.~H. Hald and P.~Stinis}, {\em Optimal prediction and the rate of decay
  for solutions of the euler equations in two and three dimensions},
  Proceedings of the National Academy of Sciences, 104 (2007), pp.~6527--6532.

\bibitem{HARARI20041491}
{\sc I.~Harari}, {\em Stability of semidiscrete formulations for parabolic
  problems at small time steps}, Computer Methods in Applied Mechanics and
  Engineering, 193 (2004), pp.~1491 -- 1516.
\newblock Recent Advances in Stabilized and Multiscale Finite Element Methods.

\bibitem{Hesthaven2016}
{\sc J.~S. Hesthaven, G.~Rozza, and B.~Stamm}, {\em Certified Reduced Basis
  Methods for Parametrized Partial Differential Equations}, Springer
  International Publishing, Cham, 2016.

\bibitem{hijazi2019data}
{\sc S.~Hijazi, G.~Stabile, A.~Mola, and G.~Rozza}, {\em {Data-driven
  POD-Galerkin reduced order model for turbulent flows}}, arXiv preprint,
  \url{http://arxiv.org/abs/1907.09909},  (2019).

\bibitem{HLB96}
{\sc P.~Holmes, J.~L. Lumley, and G.~Berkooz}, {\em Turbulence, Coherent
  Structures, Dynamical Systems and Symmetry}, Cambridge, 1996.

\bibitem{HSU2010828}
{\sc M.-C. Hsu, Y.~Bazilevs, V.~Calo, T.~Tezduyar, and T.~Hughes}, {\em
  Improving stability of stabilized and multiscale formulations in flow
  simulations at small time steps}, Computer Methods in Applied Mechanics and
  Engineering, 199 (2010), pp.~828 -- 840.
\newblock Turbulence Modeling for Large Eddy Simulations.

\bibitem{HuWeDu22}
{\sc C.~Huang, C.~R. Wentland, K.~Duraisamy, and C.~Merkle}, {\em Model
  reduction for multi-scale transport problems using model-form preserving
  least-squares projections with variable transformation}, Journal of
  Computational Physics, 448 (2022), p.~110742.

\bibitem{HUGHES1984217}
{\sc T.~Hughes and T.~Tezduyar}, {\em Finite element methods for first-order
  hyperbolic systems with particular emphasis on the compressible euler
  equations}, Computer Methods in Applied Mechanics and Engineering, 45 (1984),
  pp.~217 -- 284.

\bibitem{brooks_supg0}
{\sc T.~J. Hughes and A.~N. Brooks}, {\em {A multidimensional upwind scheme
  with no crosswind diffusion}}, in Finite Element Methods for Convection
  Dominated Flows, ASME, 1979.

\bibitem{hughes1}
{\sc T.~J. Hughes, G.~Feijoo, L.~Mazzei, and J.~Qunicy}, {\em {The variational
  multiscale method - a paradigm for computational mechanics}}, Computer
  methods in applied mechanics and engineering, 166 (1998), pp.~173--189.

\bibitem{HuFrHu89}
{\sc T.~J. Hughes, L.~P. Franca, and G.~M. Hulbert}, {\em {A new finite element
  formulation for computational fluid dynamics: VIII. The
  galerkin/least-squares method for advective-diffusive equations}}, Computer
  Methods in Applied Mechanics and Engineering, 73 (1989), pp.~173 -- 189.

\bibitem{HUGHES1996217}
{\sc T.~J. Hughes and J.~R. Stewart}, {\em A space-time formulation for
  multiscale phenomena}, Journal of Computational and Applied Mathematics, 74
  (1996), pp.~217 -- 229.

\bibitem{iliescu2018regularized}
{\sc T.~Iliescu, H.~Liu, and X.~Xie}, {\em Regularized reduced order models for
  a stochastic {B}urgers equation}, Int. J. Numer. Anal. Mod., 15 (2018),
  pp.~594--607.

\bibitem{iliescu2013variational}
{\sc T.~Iliescu and Z.~Wang}, {\em Variational multiscale proper orthogonal
  decomposition: {C}onvection-dominated convection-diffusion-reaction
  equations}, Math. Comput., 82 (2013), pp.~1357--1378.

\bibitem{iliescu2014variational}
\leavevmode\vrule height 2pt depth -1.6pt width 23pt, {\em Variational
  multiscale proper orthogonal decomposition: {N}avier-{S}tokes equations},
  Num. Meth. P.D.E.s, 30 (2014), pp.~641--663.

\bibitem{iollo2000stability}
{\sc A.~Iollo, S.~Lanteri, and J.~A. D{\'e}sid{\'e}ri}, {\em {Stability
  properties of POD--Galerkin approximations for the compressible
  Navier--Stokes equations}}, Theoret. Comput. Fluid Dyn., 13 (2000),
  pp.~377--396.

\bibitem{JoMoNo22}
{\sc V.~John, B.~Moreau, and J.~Novo}, {\em Error analysis of a supg-stabilized
  pod-rom method for convection-diffusion-reaction equations}, Computers \&
  Mathematics with Applications, 122 (2022), pp.~48--60.

\bibitem{JoNo11}
{\sc V.~John and J.~Novo}, {\em Error analysis of the supg finite element
  discretization of evolutionary convection-diffusion-reaction equations}, SIAM
  Journal on Numerical Analysis, 49 (2011), pp.~1149--1176.

\bibitem{Johnson_1984_FEM_Hyperbolic}
{\sc C.~Johnson, U.~Nävert, and J.~Pitkäranta}, {\em Finite element methods
  for linear hyperbolic problems}, Computer Methods in Applied Mechanics and
  Engineering, 45 (1984), pp.~285--312.

\bibitem{kalashnikova2014reduced}
{\sc I.~Kalashnikova, S.~Arunajatesan, M.~F. Barone, B.~G. van
  Bloemen~Waanders, and J.~A. Fike}, {\em Reduced order modeling for prediction
  and control of large-scale systems}, Sandia National Laboratories Report,
  SAND,  (2014).

\bibitem{kalash_aiaa}
{\sc I.~Kalashnikova and M.~Barone}, {\em {Stable and Efficient Galerkin
  Reduced Order Models for Non-Linear Fluid Flow}}, AIAA-2011-3110, 6th AIAA
  Theoretical Fluid Mechanics Conference, Honolulu, Hawaii, 2011.

\bibitem{kalash_convergence}
{\sc I.~Kalashnikova and M.~F. Barone}, {\em On the stability and convergence
  of a {G}alerkin reduced order model {ROM} of compressible flow with solid
  wall and far-field boundary treatment}, International Journal for Numerical
  Methods in Engineering, 83 (2010), pp.~1345--1375.

\bibitem{KALASHNIKOVA2014569}
{\sc I.~Kalashnikova, M.~F. Barone, S.~Arunajatesan, and B.~G. van
  Bloemen~Waanders}, {\em Construction of energy-stable projection-based
  reduced order models}, Applied Mathematics and Computation, 249 (2014),
  pp.~569--596.

\bibitem{kalash_eig_reassign}
{\sc I.~Kalashnikova, B.~van Bloemen~Waanders, S.~Arunajatesan, and M.~Barone},
  {\em {Stabilization of Projection-Based Reduced Order Models for Linear
  Time-Invariant Systems via Optimization-Based Eigenvalue Reassignment}},
  Comput. Meth. Appl. Mech. Engng., 272 (2014), pp.~251--270.

\bibitem{kaneko2020towards}
{\sc K.~Kaneko, P.-H. Tsai, and P.~Fischer}, {\em Towards model order reduction
  for fluid-thermal analysis}, Nucl. Eng. Des., 370 (2020), p.~110866.

\bibitem{kaptanoglu2020physics}
{\sc A.~A. Kaptanoglu, K.~D. Morgan, C.~J. Hansen, and S.~L. Brunton}, {\em
  {Physics-constrained, low-dimensional models for MHD: First-principles and
  data-driven approaches}}, arXiv preprint arXiv:2004.10389,  (2020).

\bibitem{koc2019commutation}
{\sc B.~Koc, M.~Mohebujjaman, C.~Mou, and T.~Iliescu}, {\em Commutation error
  in reduced order modeling of fluid flows}, Adv. Comput. Math., 45 (2019),
  pp.~2587--2621.

\bibitem{koc2021optimal}
{\sc B.~Koc, S.~Rubino, M.~Schneier, J.~R. Singler, and T.~Iliescu}, {\em On
  optimal pointwise in time error bounds and difference quotients for the
  proper orthogonal decomposition}, SIAM J. Numer. Anal., 59 (2021),
  pp.~2163--2196.

\bibitem{kondrashov2015data}
{\sc D.~Kondrashov, M.~D. Chekroun, and M.~Ghil}, {\em {Data-driven
  non-Markovian closure models}}, Phys. D, 297 (2015), pp.~33--55.

\bibitem{kragel2005streamline}
{\sc B.~Kragel}, {\em Streamline diffusion {POD} models in optimization}, PhD
  thesis, Universit{\"a}t Trier, 2005.

\bibitem{kunisch1}
{\sc K.~Kunisch and S.~Volkwein}, {\em Galerkin proper orthogonal decomposition
  methods for parabolic problems}, Numerische Mathematik, 90 (2001),
  pp.~117--148.

\bibitem{kunisch2}
\leavevmode\vrule height 2pt depth -1.6pt width 23pt, {\em Galerkin proper
  orthogonal decomposition for a general equation in fluid dynamics}, SIAM
  Journal on Numerical Analysis, 40 (2002), pp.~492--515.

\bibitem{LALL2003304}
{\sc S.~Lall, P.~Krysl, and J.~E. Marsden}, {\em Structure-preserving model
  reduction for mechanical systems}, Physica D: Nonlinear Phenomena, 184
  (2003), pp.~304 -- 318.
\newblock Complexity and Nonlinearity in Physical Systems -- A Special Issue to
  Honor Alan Newell.

\bibitem{layton2012approximate}
{\sc W.~J. Layton and L.~G. Rebholz}, {\em Approximate Deconvolution Models of
  Turbulence: Analysis, Phenomenology and Numerical Analysis}, vol.~2042,
  Springer Berlin Heidelberg, 2012.

\bibitem{legresley_1}
{\sc P.~LeGresley and J.~Alonso}, {\em Airfoil design optimization using
  reduced order models based on proper orthogonal decomposition}.

\bibitem{legresley_3}
\leavevmode\vrule height 2pt depth -1.6pt width 23pt, {\em Dynamic Domain
  Decomposition and Error Correction for Reduced Order Models}.

\bibitem{legresley_2}
{\sc P.~LeGresley and J.~J. Alonso}, {\em Investigation of non-linear
  projection for {POD} based reduced order models for aerodynamics}, 2001.

\bibitem{lin2019data}
{\sc K.~K. Lin and F.~Lu}, {\em {Data-driven model reduction, Wiener
  projections, and the Mori-Zwanzig formalism}}, arXiv preprint
  arXiv:1908.07725,  (2019).

\bibitem{LindsayEtAl2022}
{\sc P.~Lindsay, J.~Fike, I.~Tezaur, and K.~Carlberg}, {\em Preconditioned
  least-squares petrov–galerkin reduced order models}, International Journal
  for Numerical Methods in Engineering, 123 (2022), pp.~4809--4843.

\bibitem{LoggMardalEtAl2012a}
{\sc A.~Logg, K.-A. Mardal, G.~N. Wells, et~al.}, {\em Automated Solution of
  Differential Equations by the Finite Element Method}, Springer, 2012.

\bibitem{LoggWells2010a}
{\sc A.~Logg and G.~N. Wells}, {\em Dolfin: Automated finite element
  computing}, ACM Transactions on Mathematical Software, 37 (2010).

\bibitem{LoggWellsEtAl2012a}
{\sc A.~Logg, G.~N. Wells, and J.~Hake}, {\em DOLFIN: a C++/Python Finite
  Element Library}, Springer, 2012, ch.~10.

\bibitem{loiseau2018constrained}
{\sc J.-C. Loiseau and S.~L. Brunton}, {\em {Constrained sparse Galerkin
  regression}}, J. Fluid Mech., 838 (2018), pp.~42--67.

\bibitem{LoCaLuRo16}
{\sc S.~Lorenzi, A.~Cammi, L.~Luzzi, and G.~Rozza}, {\em {POD-Galerkin} method
  for finite volume approximation of {Navier–Stokes} and {RANS} equations},
  Computer Methods in Applied Mechanics and Engineering, 311 (2016), pp.~151 --
  179.

\bibitem{Maday_2002_RB_Noncoercive}
{\sc Y.~Maday, A.~T. Patera, and D.~V. Rovas}, {\em A blackbox reduced-basis
  output bound method for noncoercive linear problems}, in Nonlinear Partial
  Differential Equations and their Applications - Coll{\`{e}}ge de France
  Seminar Volume {XIV}, Elsevier, 2002, pp.~533--569.

\bibitem{majda2018model}
{\sc A.~J. Majda and N.~Chen}, {\em Model error, information barriers, state
  estimation and prediction in complex multiscale systems}, Entropy, 20 (2018),
  p.~644.

\bibitem{mclaughlin2016stabilized}
{\sc B.~McLaughlin, J.~Peterson, and M.~Ye}, {\em Stabilized reduced order
  models for the advection--diffusion--reaction equation using operator
  splitting}, Comput. Math. Appl., 71 (2016), pp.~2407--2420.

\bibitem{mohebujjaman2019physically}
{\sc M.~Mohebujjaman, L.~G. Rebholz, and T.~Iliescu}, {\em
  Physically-constrained data-driven correction for reduced order modeling of
  fluid flows}, Int. J. Num. Meth. Fluids, 89 (2019), pp.~103--122.

\bibitem{balanced_truncation_moore}
{\sc B.~{Moore}}, {\em Principal component analysis in linear systems:
  Controllability, observability, and model reduction}, IEEE Transactions on
  Automatic Control, 26 (1981), pp.~17--32.

\bibitem{balanced_truncation_roberts}
{\sc C.~T. Mullis and R.~A. Roberts}, {\em Synthesis of minimum roundoff noise
  fixed point digital filters}, IEEE Transactions on Circuits and Systems, 23
  (1976), pp.~551--562.

\bibitem{PACCIARINI20141}
{\sc P.~Pacciarini and G.~Rozza}, {\em Stabilized reduced basis method for
  parametrized advection–diffusion {PDEs}}, Computer Methods in Applied
  Mechanics and Engineering, 274 (2014), pp.~1 -- 18.

\bibitem{parish_wls}
{\sc E.~J. Parish and K.~T. Carlberg}, {\em Windowed least-squares model
  reduction for dynamical systems}, Journal of Computational Physics, 426
  (2021), p.~109939.

\bibitem{parish_apg}
{\sc E.~J. Parish, C.~R. Wentland, and K.~Duraisamy}, {\em The {Adjoint
  Petrov–Galerkin} method for non-linear model reduction}, Computer Methods
  in Applied Mechanics and Engineering, 365 (2020), p.~112991.

\bibitem{rozza_monograph}
{\sc A.~T. Patera and G.~Rozza}, {\em Reduced Basis Approximation and A
  Posteriori Error Estimation for Parametrized Partial Differential Equations},
  MIT Pappalardo Graduate Monographs in Mechanical Engineering, Massachusetts
  Institute of Technology, Department of Mechanical Engineering, 2007.

\bibitem{peng2016symplectic}
{\sc L.~Peng and K.~Mohseni}, {\em {Symplectic model reduction of Hamiltonian
  systems}}, SIAM J. Sci. Comput., 38 (2016), pp.~A1--A27.

\bibitem{krylov_rom}
{\sc L.~T. Pillage, X.~Huang, and R.~A. Rohrer}, {\em Asymptotic waveform
  evaluation for timing analysis}, in Proceedings of the 26th ACM/IEEE Design
  Automation Conference, DAC '89, New York, NY, USA, 1989, ACM, pp.~634--637.

\bibitem{rb_1}
{\sc C.~Prud’homme, D.~V. Rovas, K.~Veroy, L.~Machiels, Y.~Maday, A.~T.
  Patera, and G.~Turinici}, {\em {Reliable Real-Time Solution of Parametrized
  Partial Differential Equations: Reduced-Basis Output Bound Methods }},
  Journal of Fluids Engineering, 124 (2001), pp.~70--80.

\bibitem{quarteroni2015reduced}
{\sc A.~Quarteroni, A.~Manzoni, and F.~Negri}, {\em Reduced Basis Methods for
  Partial Differential Equations: An Introduction}, vol.~92, Springer, 2015.

\bibitem{Quarteroni_2008_PDE}
{\sc A.~Quarteroni and A.~Valli}, {\em Numerical approximation of partial
  differential equations}, Springer, New York, 1997.

\bibitem{doi:10.1137/S0036142901389049}
{\sc M.~Rathinam and L.~R. Petzold}, {\em A new look at proper orthogonal
  decomposition}, SIAM Journal on Numerical Analysis, 41 (2003),
  pp.~1893--1925.

\bibitem{rebollo2017certified}
{\sc T.~C. Rebollo, E.~D. {\'A}vila, M.~G. M{\'a}rmol, F.~Ballarin, and
  G.~Rozza}, {\em {On a certified Smagorinsky reduced basis turbulence model}},
  SIAM J. Numer. Anal., 55 (2017), pp.~3047--3067.

\bibitem{REYES2020112844}
{\sc R.~Reyes and R.~Codina}, {\em Projection-based reduced order models for
  flow problems: A variational multiscale approach}, Computer Methods in
  Applied Mechanics and Engineering, 363 (2020), p.~112844.

\bibitem{Rezaian}
{\sc E.~Rezaian and M.~Wei}, {\em {Impact of Symmetrization on the Robustness
  of POD-Galerkin ROMs for Compressible Flows}}, in AIAA Scitech, Orlando,
  Florida, January 2020.

\bibitem{Rezaian:2021}
{\sc E.~Rezaian and M.~Wei}, {\em A global eigenvalue reassignment method for
  the stabilization of nonlinear reduced-order models}, International Journal
  for Numerical Methods in Engineering, 122 (2021), pp.~2393--2416.

\bibitem{roop2013proper}
{\sc J.~P. Roop}, {\em {A proper-orthogonal decomposition variational
  multiscale approximation method for a generalized Oseen problem}}, Adv.
  Numer. Anal., 2013 (2013).

\bibitem{roos2008robust}
{\sc H.~G. Roos, M.~Stynes, and L.~Tobiska}, {\em Robust Numerical Methods for
  Singularly Perturbed Differential Equations: Convection-Diffusion-Reaction
  and Flow Problems.}, vol.~24 of Springer Series in Computational Mathematics,
  Springer, second~ed., 2008.

\bibitem{rovas_thesis}
{\sc D.~V. Rovas}, {\em {Reduced-basis output bound methods for parametrized
  partial differential equations}}, PhD thesis, Massachusetts Institute of
  Technology, 2003.

\bibitem{rowley2004model}
{\sc C.~W. Rowley, T.~Colonius, and R.~M. Murray}, {\em Model reduction for
  compressible flows using {POD} and {G}alerkin projection}, Phys. D, 189
  (2004), pp.~115--129.

\bibitem{Rozza2008}
{\sc G.~Rozza, D.~B.~P. Huynh, and A.~T. Patera}, {\em Reduced basis
  approximation and a posteriori error estimation for affinely parametrized
  elliptic coercive partial differential equations}, Archives of Computational
  Methods in Engineering, 15 (2008), p.~229.

\bibitem{rozza2008reduced}
\leavevmode\vrule height 2pt depth -1.6pt width 23pt, {\em Reduced basis
  approximation and a posteriori error estimation for affinely parametrized
  elliptic coercive partial differential equations}, Arch. Comput. Method. E.,
  15 (2008), pp.~229--275.

\bibitem{Rozza_2007_RB_Stokes}
{\sc G.~Rozza and K.~Veroy}, {\em On the stability of the reduced basis method
  for {Stokes} equations in parametrized domains}, Computer Methods in Applied
  Mechanics and Engineering, 196 (2007), pp.~1244--1260.

\bibitem{rubino2020numerical}
{\sc S.~Rubino}, {\em {Numerical analysis of a projection-based stabilized
  POD-ROM for incompressible flows}}, SIAM J. Numer. Anal., 58 (2020),
  pp.~2019--2058.

\bibitem{sabetghadam2012alpha}
{\sc F.~Sabetghadam and A.~Jafarpour}, {\em $\alpha$ regularization of the
  {POD-G}alerkin dynamical systems of the {K}uramoto--{S}ivashinsky equation},
  Appl. Math. Comput., 218 (2012), pp.~6012--6026.

\bibitem{sagaut2006large}
{\sc P.~Sagaut}, {\em Large Eddy Simulation for Incompressible Flows},
  Scientific Computation, Springer-Verlag, Berlin, third~ed., 2006.

\bibitem{san2015stabilized}
{\sc O.~San and T.~Iliescu}, {\em A stabilized proper orthogonal decomposition
  reduced-order model for large scale quasigeostrophic ocean circulation}, Adv.
  Comput. Math.,  (2015), pp.~1289--1319.

\bibitem{san2018machine}
{\sc O.~San and R.~Maulik}, {\em Machine learning closures for model order
  reduction of thermal fluids}, Appl. Math. Model., 60 (2018), pp.~681--710.

\bibitem{San2018}
{\sc O.~San and R.~Maulik}, {\em Neural network closures for nonlinear model
  order reduction}, Advances in Computational Mathematics, 44 (2018),
  pp.~1717--1750.

\bibitem{serre}
{\sc G.~Serre, P.~Lafon, X.~Gloerfelt, and C.~Bailly}, {\em {Reliable
  reduced-order models for time-dependent linearized Euler equations}}, Journal
  of Computational Physics, 231 (2012).

\bibitem{SHAKIB199135}
{\sc F.~Shakib and T.~J. Hughes}, {\em A new finite element formulation for
  computational fluid dynamics: Ix. {Fourier} analysis of space-time
  {Galerkin}/least-squares algorithms}, Computer Methods in Applied Mechanics
  and Engineering, 87 (1991), pp.~35--58.

\bibitem{singler}
{\sc J.~Singler}, {\em {New POD error expressions, error bounds, and asymptotic
  results for reduced order models of parabolic PDEs}}, SIAM Journal on
  Numerical Analysis, 52 (2014).

\bibitem{Sockwell:2019}
{\sc K.~Sockwell}, {\em Mass Conserving Hamiltonian-Structure-Preserving
  Reduced Order Modeling for the Rotating Shallow Water Equations Discretized
  by a Mimetic Spatial Scheme}, PhD thesis, Florida State University, 2019.

\bibitem{sotomayor_thesis}
{\sc R.~R. Sotomayor}, {\em On approaching real-time simulations for fluid
  flows}, PhD thesis, Universitat Polit\`ecnica de Catalunya, 2020.

\bibitem{stabile_vms_rom}
{\sc G.~Stabile, F.~Ballarin, G.~Zuccarino, and G.~Rozza}, {\em A reduced order
  variational multiscale approach for turbulent flows}, Advances in
  Computational Mathematics, 45 (2019), pp.~2349--2368.

\bibitem{StHiMoLo17}
{\sc G.~Stabile, S.~Hijazi, A.~Mola, S.~Lorenzi, and G.~Rozza}, {\em
  Pod-galerkin reduced order methods for cfd using finite volume
  discretisation: vortex shedding around a circular cylinder}, Communications
  in Applied and Industrial Mathematics, 8 (2017).

\bibitem{StRo18}
{\sc G.~Stabile and G.~Rozza}, {\em Finite volume pod-galerkin stabilised
  reduced order methods for the parametrised incompressible navier–stokes
  equations}, Computers \& Fluids, 173 (2018), pp.~273 -- 284.

\bibitem{strazzullo2022consistency}
{\sc M.~Strazzullo, M.~Girfoglio, F.~Ballarin, T.~Iliescu, and G.~Rozza}, {\em
  Consistency of the full and reduced order models for evolve-filter-relax
  regularization of convection-dominated, marginally-resolved flows}, Int. J.
  Num. Meth. Eng., 123 (2022), pp.~3148--3178.

\bibitem{St05}
{\sc M.~Stynes}, {\em Steady-state convection-diffusion problems}, Acta
  Numerica, 14 (2005), p.~445–508.

\bibitem{TEZDUYAR19911}
{\sc T.~Tezduyar}, {\em Stabilized finite element formulations for
  incompressible flow computations}, vol.~28 of Advances in Applied Mechanics,
  Elsevier, 1991, pp.~1 -- 44.

\bibitem{thomee2006galerkin}
{\sc V.~Thom{\'e}e}, {\em {Galerkin finite element methods for parabolic
  problems}}, Springer Verlag, 2006.

\bibitem{URBAN2012203}
{\sc K.~Urban and A.~T. Patera}, {\em A new error bound for reduced basis
  approximation of parabolic partial differential equations}, Comptes Rendus
  Mathematique, 350 (2012), pp.~203 -- 207.

\bibitem{Urban_2014_Space_Time_RB}
\leavevmode\vrule height 2pt depth -1.6pt width 23pt, {\em An improved error
  bound for reduced basis approximation of linear parabolic problems},
  Mathematics of Computation, 83 (2014), pp.~1599--1615.

\bibitem{rb_3}
{\sc K.~Veroy and A.~T. Patera}, {\em Certified real-time solution of the
  parametrized steady incompressible {Navier–Stokes} equations: rigorous
  reduced-basis a posteriori error bounds}, International Journal for Numerical
  Methods in Fluids, 47 (2005), pp.~773--788.

\bibitem{rb_2}
{\sc K.~Veroy, C.~Prud'homme, D.~Rovas, and A.~Patera}, {\em A Posteriori Error
  Bounds for Reduced-Basis Approximation of Parametrized Noncoercive and
  Nonlinear Elliptic Partial Differential Equations}.

\bibitem{volkwein2013proper}
{\sc S.~Volkwein}, {\em Proper orthogonal decomposition: Theory and
  reduced-order modelling}, Lecture Notes, University of Konstanz,  (2013).
\newblock
  \url{http://www.math.uni-konstanz.de/numerik/personen/volkwein/teaching/POD-Book.pdf}.

\bibitem{Wang:269133}
{\sc Q.~Wang, N.~Ripamonti, and J.~S. Hesthaven}, {\em Recurrent neural network
  closure of parametric {POD-Galerkin} reduced-order models based on the
  {Mori-Zwanzig} formalism}, Journal of Computational Physics,  (2019).

\bibitem{Wang_ROM_thesis}
{\sc Z.~Wang}, {\em {Reduced-Order Modeling of Complex Engineering and
  Geophysical Flows: Analysis and Computations}}, PhD thesis, Virginia
  Polytechnic Institute and State University, 2012.

\bibitem{wang2012proper}
{\sc Z.~Wang, I.~Akhtar, J.~Borggaard, and T.~Iliescu}, {\em Proper orthogonal
  decomposition closure models for turbulent flows: A numerical comparison},
  Comput. Meth. Appl. Mech. Eng., 237-240 (2012), pp.~10--26.

\bibitem{wells2017evolve}
{\sc D.~Wells, Z.~Wang, X.~Xie, and T.~Iliescu}, {\em An evolve-then-filter
  regularized reduced order model for convection-dominated flows}, Int. J. Num.
  Meth. Fluids, 84 (2017), pp.~598--–615.

\bibitem{wentland_apg}
{\sc C.~R. Wentland, C.~Huang, and K.~Duraisamy}, {\em Closure of Reacting Flow
  Reduced-Order Models via the Adjoint Petrov-Galerkin Method}.

\bibitem{xie2018data}
{\sc X.~Xie, M.~Mohebujjaman, L.~G. Rebholz, and T.~Iliescu}, {\em Data-driven
  filtered reduced order modeling of fluid flows}, SIAM J. Sci. Comput., 40
  (2018), pp.~B834--B857.

\bibitem{xie2017approximate}
{\sc X.~Xie, D.~Wells, Z.~Wang, and T.~Iliescu}, {\em Approximate deconvolution
  reduced order modeling}, Comput. Methods Appl. Mech. Engrg., 313 (2017),
  pp.~512--534.

\bibitem{xie2018numerical}
\leavevmode\vrule height 2pt depth -1.6pt width 23pt, {\em Numerical analysis
  of the {L}eray reduced order model}, J. Comput. Appl. Math., 328 (2018),
  pp.~12--29.

\bibitem{Yano2014stBoussinesq}
{\sc M.~Yano}, {\em A space-time {Petrov}-{Galerkin} certified reduced basis
  method: application to the {Boussinesq} equations}, SIAM Journal on
  Scientific Computing, 36 (2014), pp.~A232--A266.

\bibitem{Yano_2014_Space_Time_RB_Boussinesq}
\leavevmode\vrule height 2pt depth -1.6pt width 23pt, {\em A space-time
  {Petrov}-{Galerkin} certified reduced basis method: application to the
  {Boussinesq} equations}, SIAM Journal on Scientific Computing, 36 (2014),
  pp.~A232--A266.

\bibitem{zoccolan2023streamline}
{\sc F.~Zoccolan, M.~Strazzullo, and G.~Rozza}, {\em {A streamline upwind
  Petrov-Galerkin reduced order method for advection-dominated partial
  differential equations under optimal control}}, arXiv preprint,
  \url{http://arxiv.org/abs/arXiv:2301.01973},  (2023).

\bibitem{zoccolan2023stabilized}
\leavevmode\vrule height 2pt depth -1.6pt width 23pt, {\em Stabilized weighted
  reduced order methods for parametrized advection-dominated optimal control
  problems governed by partial differential equations with random inputs},
  arXiv preprint, \url{http://arxiv.org/abs/arXiv:2301.01975},  (2023).

\end{thebibliography}

\end{document}